\documentclass[a4paper,11pt]{amsart}
\usepackage{amscd}
\usepackage{calrsfs}
\usepackage{amssymb}

\def\N{\mathbb{N}}
\def\Z{\mathbb{Z}}
\def\Q{\mathbb{Q}}
\def\R{\mathbb{R}}
\def\C{\mathbb{C}}
\def\F{\mathbb{F}}
\def\G{\mathbb{G}}
\def\P{\mathbb{P}}
\def\A{\mathbb{A}}

\def\N{\mathbb{N}}

\def\Fr{\mathrm{Fr}}
\def\st{\mathrm{st}}
\def\lc{\mathrm{loc}}

\def\associate{\mathrm{ass}}
\def\op{\mathrm{op}}
\def\tr{\mathrm{trace}}
\def\Ind{\mathrm{Ind}}
\def\residue{\mathrm{res}}

\def\rank{\mathrm{rank}}

\def\Res{\mathrm{Res}}
\def\Sel{\mathrm{Sel}}
\def\BC{\mathrm{BC}}
\def\GL {\mathrm{GL}}
\def \modular{\mathrm{mod}}
\def\tw{\mathrm{tw}}
\def\exc{\mathrm{exc}}
\def\Hecke{\mathrm{Hecke}}
\def\parabolic{\mathrm{par}}
\def\old{\mathrm{old}}
\def\Type{\mathrm{Type}}
\def\stack{\mathrm{stack}}
\def\disc{\mathrm{disc}}
\def\JL{\mathrm{JL}}
\def\Norm{\mathrm{N}}
\def\new{\mathrm{new}}
\def\noeth{\mathrm{noeth}}
\def\artin{\mathrm{artin}}
\def\lift{\mathrm{lift}}
\def\flatcoh{\mathrm{fl}}
\def\unr{\mathrm{unr}}
\def\ev{\mathrm{ev}}
\def\mult{\mathrm{mult}}
\def\ob{\mathrm{ob}}
\def\ffgpsch{\mathrm{GSch ^{ff}}}

\def\deform{\mathrm{def}}
\def\PSL{\mathrm{PSL}}

\def\Ann{\mathrm{Ann}}

\def\tame{\mathrm{tame}}
\def\ab{\mathrm{ab}}

\def\sp{\mathop{\mathrm{sp}}\nolimits}
\def\Spec{\mathop{\mathrm{Spec}}\nolimits}

\def\Sym{\mathop{\mathrm{Sym}}\nolimits}
\def\Aut{\mathop{\mathrm{Aut}}\nolimits}
\def\Hom{\mathop{\mathrm{Hom}}\nolimits}

\def\Ext{\mathop{\mathrm{Ext}}\nolimits}
\def\Tor{\mathop{\mathrm{Tor}}\nolimits}

\def\ker{\mathop{\mathrm{ker}}\nolimits}
\def\coker{\mathop{\mathrm{coker}}\nolimits}

\def\id{\mathrm{id}}

\def\ren{\mathrm{ren}}
\def\val{\mathrm{val}}
\def\slope{\mathrm{slope}}

\def\cycle{\mathrm{cycle}}
\def\cond{\mathrm{cond}}
\def\Art{\mathrm{Art}}
\def\ord{\mathrm{ord}}
\def\GL{\mathrm{GL}}
\def\PGL{\mathrm{PGL}}
\def\SL{\mathrm{SL}}
\def\det{\mathrm{det}}
\def\ad{\mathrm{ad}}
\def\End{\mathrm{End}}
\def\Gal{\mathrm{Gal}}

\def\unr{\mathrm{unr}}
\def\ss{\mathrm{ss}}
\def\der{\mathrm{der}}

\def\et{\mathrm{\acute{e}t}}
\def\univ{\mathrm{univ}}

\def\pr{\mathop{\mathrm{pr}}\nolimits}
\def\length{\mathrm{length}}

\def\Lie{\mathrm{Lie}}

\newtheorem{thm}{Theorem}[section]
\newtheorem{prop}[thm]{Proposition}
\newtheorem{lem}[thm]{Lemma}
\newtheorem{cor}[thm]{Corollary}
\newtheorem{rem}[thm]{Remark}
\newtheorem{dfn}[thm]{Definition}

\newtheorem{quest}[thm]{Question}

\newtheorem{sublem}[thm]{Sublemma}

\newtheorem{conj}[thm]{Conjecture}
\newtheorem{hyp}[thm]{Hypothesis}

\def\longhookrightarrow{\lhook\joinrel\longrightarrow}

\def\ad{\mathop{\mathrm{ad}}\nolimits}

\def\sp{\mathop{\mathrm{sp}}\nolimits}

\def\int{\mathop{\mathrm{int}}\nolimits}

\def\Spec{\mathop{\mathrm{Spec}}\nolimits}

\def\Sym{\mathop{\mathrm{Sym}}\nolimits}
\def\Aut{\mathop{\mathrm{Aut}}\nolimits}
\def\Hom{\mathop{\mathrm{Hom}}\nolimits}

\def\Ext{\mathop{\mathrm{Ext}}\nolimits}
\def\Tor{\mathop{\mathrm{Tor}}\nolimits}
\def\ker{\mathop{\mathrm{ker}}\nolimits}
\def\coker{\mathop{\mathrm{coker}}\nolimits}

\def\id{\mathrm{id}}

\def\et{\mathrm{\acute{e}t}}

\def\Sets{\mathbf{Sets}}

\def\coh{\mathrm{coh}}

\def\frame{\mathrm{Frame}}

\def\longhookrightarrow{\lhook\joinrel\longrightarrow}

\font\manual=manfnt 
\def\dbend{{\manual\char237}} 
\def\d@nger{\medbreak\begingroup\clubpenalty=10000
  \def\par{\endgraf\endgroup\medbreak} \noindent\hangindent3zw\hangafter=-2
  \hbox to0pt{\hskip-\hangindent\dbend\hfill}}
\outer\def\danger{\d@nger}
\def\dd@nger{\medbreak\begingroup\clubpenalty=10000
  \def\par{\endgraf\endgroup\medbreak} \noindent\hangindent4zw\hangafter=-2
  \hbox to0pt{\hskip-\hangindent\dbend\kern1pt\dbend\hfill}}
\outer\def\ddanger{\dd@nger}

\title {Deformation rings and Hecke algebras in the totally real case}
\author{Kazuhiro Fujiwara}

\begin{document}
\maketitle
\begin{abstract} In this paper, the $\ell$-adic Hecke algebras of $\GL_2$ over totally real fields are studied. In particular we show that they are identified with deformation rings of mod $\ell$ Galois representations in important cases, assuming that the representation is absolutely irreducible. 
\end{abstract}

\tableofcontents
\setcounter{section}{-1}

\section{Introduction}
One of the basic questions in number theory is to
determine semi-simple $\ell$-adic representations of the absolute Galois group of a number
field $F$. In case of abelian representations, a satisfactory answer is known for compatible
system of $\ell$-adic representations, and these types of abelian representations are
obtained from algebraic Hecke characters. In this paper, we discuss the question for two
dimensional representations over a totally real number field.   

\medskip

For a totally real field
$F$, let
$G_F=
\Gal (
\bar F / F )
$ be the absolute Galois group, and $I_{F, \infty}  =
\{
\iota : F
\hookrightarrow
\R
\}
$ be the set of the infinite places of $F$. We take a pair $(k, w)$ of integral vector
$ k = ( k _{\iota} ) _{
\iota
\in I _{F, \infty} }
\in
\Z ^{I_{F, \infty} }$ and an integer $w \in \Z$ called discrete infinity type (cf. Definition \ref{dfn-shimura21}). For a
cuspidal representation $\pi$ of 
$\GL _2 (\A_F)
$, we assume that the infinity part $\pi _{\infty}
$ is isomorphic to
$\bigotimes _{ \iota \in I_{F, \infty } } D_{k \iota , w} $. Here $D_{k_\iota , w }$ is an
essentially square integrable representation of $\GL_2 (\R)$ (see Conventions for our
normalization).  These types of cuspidal representations are generated by holomorphic
Hilbert modular forms. \par

Fix a prime
$\ell$, and an isomorphism $\C \simeq \bar \Q_\ell$ by the axiom of choice.
 It is known that the finite part $\pi _f 
$ of $\pi$ as above is defined over some $\ell$-adic field
$E_{\lambda}$ with the integer ring $o_{E_ \lambda}$ and the residue field $k_\lambda$, and there is a two dimensional
continuous
$\lambda$-adic Galois representation 
$$
\rho _{\pi, E_\lambda} : G _F \to \GL_2 (o_{E_\lambda} )
$$ associated to
$\pi$, which is pure of weight $w +1$ with respect to geometric Frobenius elements, and
$\det
\rho _{
\pi ,
E_\lambda } \cdot \chi ^ {w + 1}  _{\cycle} $ is a character of finite order (see
\cite{O1}, \cite{Car2}, \cite{W1}, \cite{T1}, \cite{BR} for $F \neq \Q$). Here
$\chi _{\cycle }: G_F
\to
\Z^{\times} _{\ell} $ is the cyclotomic character. We call this class of
$\lambda$-adic representations {\it modular } $\lambda$-adic representations. Conjecturally, modular
$\lambda$-adic representations in our sense should cover most {\it motivic} and {\it totally
odd } Galois representations (global Langlands correspondence).  When the $\ell$-adic representation $\rho$ is obtained from an elliptic curve
$E$ over
$F$, this conjecture is a generalized form of the Taniyama-Shimura conjecture. \par

In \cite{W2}, in case of $\Q$, Wiles has studied the problem of modularity via
the deformation theory of mod $\ell$-Galois representations. We take the same approach in this paper. For a mod $ \ell$-representation $ \bar \rho : G_F \to \GL_2 ( k_\lambda )$ having values in $k_\lambda $, assume that there is a cuspidal representation
$\pi
$ of
$\GL _2 (\A _F )
$ of infinity type $ ( k, w) $, and
$\bar
\rho$ is obtained as the reduction modulo $\lambda$ of the associated Galois representation,
that is, 
$$
\bar 
\rho \simeq
\rho _{ \pi , E_\lambda} \mod \lambda.
$$ 
We call $\bar \rho$ {\it modular} if this condition is satisfied. 

Instead of the global Langlands correspondence itself, we consider the following question,
which asks the stability of the notion of the modularity of $\lambda$-adic representations under perturbation:
\begin{quest} Assume that mod $\ell$-representation $\bar \rho$ is modular, and take a $\lambda$-adic representation
$\rho$ which lifts $\bar \rho$. Is $\rho $ modular?
\end{quest}

Without any restriction, $\rho$ can not be modular, since
$\rho _{\pi, E_\lambda}
$ is expected to be {\it motivic}, which is known in most cases (\cite{Car2}, \cite{BR}). Thus there are strong restrictions
on local monodromies of $\rho$, in particular at the places dividing $\ell$.  \par

On the other hand, $\rho$ is seen as an $\ell$-adic {\it deformation} of $\bar \rho$ from
the viewpoint of Mazur if $\bar \rho$ is irreducible \cite{M}. So one may expect that if the local conditions are imposed appropriately,
all reasonable $\lambda$-adic deformations of $\bar
\rho$ are modular representations. Hida introduced ordinary $\ell$-adic Hecke algebras and constructed Galois representations having values in it when $F= \mathbb Q$ \cite{H0}. Mazur conjectured these Hecke algebras introduced by Hida are actually
the universal deformation rings which control all deformations with suitable local
monodromy conditions.   \par

 In \cite{W2}, together with \cite{TW}, it was shown affirmatively that the
$\ell$-adic Hecke algebra is the universal deformation ring in the case where $F=\mathbb Q$, including some non-ordinary
cases (see \cite{D1} for a generalization). Moreover, Wiles applied this result to the modularity of $\lambda$-adic
representations, and proved the Taniyama-Shimura conjecture over $\mathbb Q$ in the
semistable case.\par

 Our main theorem in this article is the following (Theorem \ref{thm-final1}):

\begin{thm}[$R=T$ theorem]
Let $F$ be a totally real number field of degree $d$, $ \bar \rho: G_F \to \GL_2( k) $ an absolutely irreducible mod $
\ell$-representation. We fix a deformation type $\mathcal  D$, and assume the following conditions.
\begin{enumerate}
\item $\ell \geq 3$, and $\bar \rho
\vert _{{ F(\zeta_\ell) }}$ is absolutely irreducible. When $\ell= 5$, the following case is excluded: the projective image $\bar G$ of $ \bar \rho $ is isomorphic to $\PGL_2 (\mathbb F_5)$, and the mod $\ell$-cyclotomic character $\bar \chi _{\cycle} $ factors through $G_F \to  \bar G ^{\ab} \simeq \mathbb Z/ 2$ (in particular $[F(\zeta_5): F ] = 2$). 
\item For $ v\vert
\ell$, the deformation condition for 
$\bar
\rho
\vert _{G_{F_v} } $ is either nearly ordinary or flat. 
When the condition is nearly ordinary (resp. flat) at $v$, we assume that $\bar \rho \vert _{G_{F_v} } $ is $G_{F_v} $-distinguished (resp. $F_v$ is absolutely unramified).
\item There is a minimal modular lifting $\pi$ of $ \bar \rho $ in Definition \ref{dfn-nearlyordinary32}.

\item Hypothesis \ref{hyp-nearlyordinary21} is satisfied.
\item If $\mathcal  D $ is not minimal, we assume Hypothesis \ref{hyp-coh31} when $q_{\bar \rho } = 1$. 
\end{enumerate}
Then the universal deformation ring $R_{\mathcal  D} $ of $ \bar
\rho$ of type $\mathcal  D$ is a complete intersection, and is isomorphic to the Hecke algebra $T_{\mathcal  D} $.
\end{thm}

As is already remarked, this is proved in \cite{TW}, \cite{W2}, \cite{D1} in the case where $F= \Q$ (our method in this paper gives a substantial simplification, which is
also due to Diamond \cite{D2}). Theorem \ref{thm-final1} is a basic tool to study the modularity questions.  As in \cite{W2}, one deduces the modularity of some $3$ and $5$-adic
representations from the theorem. \par

As another application of Theorem \ref{thm-final1}, we obtain the finiteness of the Selmer group for the
adjoint representation.
\begin{cor}\label{cor-final1}
Under the same assumption as in Theorem \ref{thm-final1}, the Selmer group $\Sel
_{\mathcal  D} ( F, \ \ad \rho ) $ of
$\rho  =
\rho _{\pi, E_ \lambda} $ for $\pi $ appearing in $T_{\mathcal  D} $ is finite. 
\end{cor}

For other applications of the main theorem, see \cite{Fu2}.

Here is the explanation of the conditions of the main theorem.\par
In (1), the exceptional case when $\ell=5$ does not happen in the application to elliptic curves (see Proposition \ref{prop-exceptional}). Even in the exceptional case, Theorem \ref{thm-final1} holds true by a slight generalization of our method, which will be discussed on another occasion. 

In (2),  $F_v$ can be absolutely ramified if the deformation condition is nearly ordinary at $v$. \par

For (3), if one only assumes that $\bar \rho$ is modular, the existence of a minimal lift is
satisfied outside $\ell$ by \cite{Ja1}, \cite{Ja2}, \cite{Fu1}, and
\cite{Raj}. At places dividing $ \ell$, the existence is satisfied if the condition is nearly ordinary. In general, only a partial answer is known.  \par

In (4), Hypothesis \ref{hyp-nearlyordinary21} (local monodromy hypothesis) is satisfied if $ k = (2, \ldots, 2) $, or the deformation condition is
nearly ordinary at all
$v \vert \ell$. So the theorem is applied to Hida's nearly ordinary Hecke
algebras \cite{H2}. \par

As for (5), Hypothesis \ref{hyp-coh31} is what is called as ``Ihara's Lemma for Shimura curves''. Ribet \cite{Ri1} proved it when $F= \mathbb Q$ and modular curves using a
result of Ihara. If one is only interested in the modularity question, Hypothesis \ref{hyp-coh31} can be avoided
by choosing a totally real quadratic extension
$F'$ of $F$ which depends on $\bar \rho$, and using a base change argument. 
We will discuss the local monodromy hypothesis and Ihara's Lemma for Shimura curves on another occasion. \par
\bigskip 
Our proof of Theorem \ref{thm-final1} consists of two basic steps, similarly to the argument in \cite{W2}. One first
proves the theorem when the deformation condition $\mathcal D$ is minimal, and reduces the
general case to the minimal case by a level raising argument. \par
 In the minimal case, the problem of showing that $R _{\mathcal D }$ is equal to $  T_{\mathcal
D} $ ($R=T$ theorem, a generalized form of Mazur's conjecture) is reduced to the construction of a family of rings and modules
$\{ R_Q, M_Q \}_{Q
\in
X}  
$ with a weak compatibility, which we call a {\it Taylor-Wiles system}. Here an element $Q \in X$ of the index set $X$ is a finite set of finite places, and $\emptyset \in X$.  $R_{\emptyset } = R_{\mathcal D} $,
$R_Q$ is a complete noetherian local $o_{E_\lambda}$-algebra with a surjective map
$ R _Q\twoheadrightarrow R_{\mathcal D } $, and 
$M_Q
$ is an
$R_Q$-module. The axioms of Taylor-Wiles systems and the fundamental theorem (Theorem \ref{thm-TW21}) called the {\it complete intersection-freeness
criterion} are discussed in
\S\ref{sec-TW}. $R=T$ theorem in the minimal case is proved by applying this criterion.  A similar observation was also made independently by Diamond \cite{D2}. This is one
of the main innovations which made possible our generalization to totally real fields.
 In \cite{TW}, Taylor-Wiles system $\{ R_Q, M_Q \} _{Q \in X} 
$ was constructed for Hecke algebras in case of $\mathbb Q$. There $ R_Q$ is a
certain Hecke algebra $T_Q$ dominating $T_{\mathcal D}$, and $M_Q$ is known to be free over  $R_Q$ by an application of the $q$-expansion principle (the freeness theorem).  The freeness of $M_Q$ over $T_Q$ (which implies the Gorensteinness of $T_Q$ in the modular curve case) goes back to Mazur's work.  Though the $q$-expansion principle is also available for Hilbert modular forms,
a direct application of the method used for modular curves seems more difficult, and to show the freeness theorem and the 
Gorensteinness of the Hecke algebras in the totally real case, Taylor-Wiles systems are the most powerful tools at present.\par
In our formulation, notice that $M_Q$ is not assumed to be free over
$R_Q$, but is finite free over some {\it group ring} $o _{E_\lambda}[\Delta_Q]$.
The freeness over a group ring plays an essential role in showing that a well-chosen limit of $R_Q$ is a power series ring, and behaves as if $R_Q$ were smooth rings. 
 \par
In our general totally real case, a Taylor-Wiles system is constructed for universal deformation rings. A direct use of the universal deformation rings is suggested by an observation of Faltings (\cite{TW}, appendix). $M_Q$ is constructed from the middle dimensional
cohomology group of some modular variety
$S_Q$ attached to a quaternion algebra $D$ over $F$ of complex dimension $\leq 1 $.
$S_Q$ is a Shimura curve, or a zero-dimensional variety associated to definite quaternion
algebra whose arithmetic importance was found by Hida \cite{H1}. In the latter case,
$S_Q$ is {\it not} a Shimura variety in the sense of Deligne, which we call a {\it Hida variety}. To verify the axiom TW4, 
we make use of an argument based on a property of perfect complexes. The perfect complex
argument is extremely useful in the study of congruences. For example, if one can control the mod
$\ell$-cohomology groups except the middle dimension, one may expect that the argument given in the paper is effective for other modular varieties, 
and Taylor-Wiles systems are obtained from the middle dimensional cohomology groups. \par
After the construction of a Taylor-Wiles system, we deduce $R=T$ theorem and the complete intersection
property of the Hecke algebra in the minimal case at the same time by the
complete intersection-freeness criterion. In particular we prove
that the minimal Hecke algebra is a complete intersection without knowing a priori  the freeness nor the Gorenstein property. Moreover, the freeness of the middle dimensional cohomology also follows {\it as a consequence
}(Theorem \ref{thm-final1}). 
\begin{thm}[Freeness theorem] Under the same assumption as in $R=T$ theorem,
$ M _{\phi} =M_{\mathcal D }  $ is a free $T_{\mathcal D}
$-module. 
\end{thm}
For general deformation type $\mathcal D$, the three properties
\begin{itemize}
\item $R=T$,
\item the local complete
intersection property of
$T_{\mathcal D}$,
\item the freeness theorem
\end{itemize}
are proved inductively from the minimal case {\it
at the same time} by calculating cohomological congruence modules in the sense of Ribet
\cite{Ri1}. Hypothesis \ref{hyp-coh31} is needed in this calculation. The level raising formalism is given in \S\ref{sec-TW} (one may also use results of Diamond \cite{D2}).
Since we do not know the Gorenstein property of
$T_{\mathcal D}$ in advance, to deduce $R= T$, we use Lenstra's criterion \cite{Len} instead of Wiles' criterion
\cite{W2}. The argument of the freeness given in \S\ref{sec-TW} is due to T. Saito.  \par
\bigskip
Throughout this paper, we make a systematic use of the homological algebra on Shimura varieties. Some foundational results are treated as general as possible for the applications in future including absolutely reducible representations. \par
This paper is organized as follows.\par

In \S\ref{sec-TW}, a formulation and a generalization of a beautiful argument in \cite{TW} is given.
The notion of Taylor-Wiles systems is introduced. Our formulation is influenced by an
argument of Faltings (\cite{TW}, appendix). As is remarked already, a similar observation was made independently by Diamond
\cite{D2}. We also give a formulation of the level raising argument. 

In \S\ref{sec-galdef}, we discuss deformations of local and global Galois deformations. We include the case treated by 
Diamond \cite{D1}, though we make some modifications using G\'eradin's result \cite{Ger}. \par

In \S\ref{sec-shimura}, basic results on modular varieties are discussed. In particular  we treat general coefficient sheaves.  The duality formalism is needed in the later sections. \par
In \S\ref{sec-coh}, the universal exactness of cohomological sequences is formulated and studied. The universal exactness (Ihara's Lemma) is proved for definite quaternion algebras except a class of absolutely reducible representations which is called of residual type. \par

In \S\ref{sec-nearlyordinary}, nearly ordinary automorphic representations are discussed. We formulate a property on local monodromies of Galois representations attached to such representations as Local monodromy hypothesis (Hypothesis \ref{hyp-nearlyordinary21}). The hypothesis is known to hold in many cases. \par
In \S\ref{sec-modular}, we construct and study a Hecke algebra $T_{\mathcal D}$ and a Galois representation
$\rho _{\mathcal D} ^{\modular} $ having values in
$\GL_2 (T_{\mathcal D})$ for deformation type $\mathcal D$ of $\bar \rho$ using the irreducibility of $\bar
\rho$. For the construction of $\rho _{\mathcal D} ^{\modular} $, one applies the method of
pseudo-representation of Wiles, and reduces it to the existence of $
\lambda$-adic representation $\rho_{\pi , E_\lambda} $ attached to a cuspidal representation
$\pi$. $\rho _{\mathcal D} ^{\modular} $ is seen as the universal modular deformation of $\bar
\rho$. To have a ring homomorphism $ \pi_{\mathcal D} : R_{\mathcal D} \to T_{\mathcal D} $ from the universal
deformation ring, one needs to check the local property of $ \rho ^{\modular} _{\mathcal D} $ at
all places. The necessary information at
$v \nmid \ell$ follows from the following fact: the global Langlands correspondence $
\pi \mapsto \rho _{\pi, E_\lambda} $ is compatible with the local Langlands correspondence
\cite{K}. This compatibility is shown in
\cite{Car2}, \cite{W2},
\cite{T1}. For places dividing $\ell$, the situation is much subtle. If the deformation
condition is flat, we prove it by an arithmetic geometric analysis of Shimura curves using
\cite{Car1} if the quaternion algebra
$D$ is indefinite. If $D$ is definite, one can not treat it directly, and we use the approximation method of
Taylor
\cite{T1} , \cite{T2}. The claim is reduced to the case of Shimura curves associated with some indefinite and
ramified quaternion algebras.\par

In \S\ref{sec-const}, we construct a Taylor-Wiles system from universal deformation rings and cohomology
of modular varieties of complex dimension $\leq 1$ as is remarked already. We do not use arithmetic geometrical properties of these ``modular varieties''.\par 

In \S\ref{sec-minimal}, we combine the results obtained in the previous sections, and deduce the main
theorem in the minimal case. The deformation rings are controlled by Galois cohomology groups.\par

In \S\ref{sec-congruent}, the calculation of cohomological congruence modules is done as in
\cite{W2}, \S2. 

In \S\ref{sec-final}, we summarize the previous results, and the main theorem is
proved.\par
\bigskip

 {\bf Acknowledgement:} 

The first version of this paper was written during the author's stay at M\"unster university in 1996. There, only the minimal case was treated. The author would like to thank P. Schneider for hospitality. The second version was written during the author's stay at the Institute for Advanced Study for the academic year
1998-1999, which included the non-minimal cases. The present paper is based on the second version. The author would like to thank P. Deligne, F. Diamond, H. Hida, B. Mazur, R. P. Langlands, D. Ramakrishnan, T. Saito, C. Skinner, J. Tilouine, and R. Taylor for helpful comments.  In particular the author is grateful to T. Saito for his encouragement during preparation of this paper. He has given the author detailed comments on earlier versions of this paper. The argument in \S\ref{sec-TW} of showing the freeness is due to him. \par
Finally, the author is grateful to A. Wiles. Without his support, this paper would have never been written.
\section{Conventions}

For a field $F$, $G_F = \Gal (\bar F  /F)$ means the absolute
Galois group.
For a prime $\ell$, we denote the group of $\ell$-power roots of unit by $\mu _{\ell ^{\infty } } (F) $. Let
$\chi_{\cycle}: G_F \to \Z  ^{\times} _\ell  $
be the cyclotomic character of $F$, $\chi_{\cycle} =  \chi _{\cycle} ^\ell  \cdot \chi _{\cycle, \ell } $ the decomposition into prime-to
$\ell$ and pro-$\ell$-part. We denote $\chi _{\cycle} ^\ell $ by $\omega _F$.
\par

For a non-archimedean local field $F$, the ring of the integers of $F$ is denoted by $o_F$, with the maximal ideal $m_F$ and the residue field $k_F$. A uniformizer of $F$ is denoted by $p_F$. $W_F\subset G_F$ the Weil group of $F$, $I_F\subset G_F$ the inertia subgroup of $G_F$. $W'_F$ is the Weil-Deligne group of $F$.\par
\bigskip
For a non-zero ideal $I$ of $o_F$, subgroups $ K_{11} (I) \subset K_1 ( I ) \subset K_0 ( I ) $ of $\GL_2 (o_F)$ are defined by 
$$K _ {11} ( I ) = \{ g \in 
\GL _2 ( o _F ),\ g \equiv
\left(
\begin{array}{cc}
1 &* \\
0&1
\end{array}
\right)
\mod I \} ,
$$
$$
K _ 1 ( I ) = \{ g \in 
\GL _2 ( o _F),\ g \equiv 
\left(
\begin{array}{cc}
1 &* \\
0&*
\end{array} 
\right) 
\mod I  \} ,
$$
$$
K _ 0 ( I ) = \{ g \in 
\GL _2 (  o _F ),\ g \equiv 
\left(
\begin{array}{cc}
* &* \\
0&*
\end{array}
\right) 
\mod I \} .
$$
\par
By the local class field theory, $W^{\ab}_F \simeq F ^{\times} $, where a uniformizer $p_F$ corresponds to a geometric Frobenius element of $W^{\ab}_F$. \par
The local Langlands correspondence for $\GL_{2, F}$ is a bijection between
the isomorphism classes of $F$-semisimple representation $\rho$ of the Weil-Deligne group $W' _F$ and the isomorphism classes of admissible
irreducible representation $\pi$ of $\GL_2 (F)$ \cite{K}. \par
\bigskip
When $F$ is a global field, the set of all places (resp. finite places) is denoted by $\vert F\vert  $ (resp. $\vert F\vert_f$) . For a place $v$ of $F$, $F_v$ means the local field at
$v$. If $v$ is a finite place, the ring integers of $F_v$ is denoted by $o_{F_v}$, with the maximal ideal $m_{F_v}$, and a uniformizer
$p_v$. $k (v ) = o_{F_v} /m_{F_v}$ the residue field, $q_v $ its cardinality.\par
For a finite set of finite places $\Sigma$, $G_{\Sigma } = \pi^{\et} _1 (\Spec o_F \setminus \Sigma) $ is the Galois group of the maximal Galois extension of $F$ which is unramified outside $\Sigma $.\par
\bigskip
$\A_F$ is the ad\'ele ring of $F$, $\A _{F, f}$ is the ring of finite adeles, and
$\A _{F, \infty} $ is the infinite part of $\A_F$.
For a reductive group $G$ over $F$, a compact open subgroup $K$ of the adelic group $ G (\A_{F, f })$ is factorizable if $ K =
\prod_v K_v $, $ K _v \subset G ( F _v)$. In this case, for a finite set of places
$\Sigma $, 
$$ K _\Sigma = \prod _{ v \in  \Sigma } K_v,\ K ^ \Sigma =  \prod _{ v \not \in 
\Sigma } K_v.
$$
$ K = K _\Sigma \cdot K ^\Sigma $. \par
 For a non-zero
ideal
$I$ of
$o_F
$, we define compact open subgroups of
$\GL _2 (
\A _{F , f} )$ by  
$$K(I ) = \{ g \in 
\GL _2 (
\prod _{v\in \vert F\vert _f }  o _{F_v} ),\ g \equiv 1 \mod I  \} ,
$$

$$K _ {11} ( I ) = \{ g \in 
\GL _2 (
\prod _{v \in \vert F \vert _f}  o _{F_v }),\ g \equiv
\left(
\begin{array}{cc}
1 &* \\
0&1
\end{array}
\right)
\mod I \} ,
$$
$$K _ 1 ( I  ) = \{ g \in 
\GL _2 (
\prod _{v \in \vert F \vert _f} o _{F_v} ),\ g \equiv 
\left(
\begin{array}{cc}
1 &* \\
0&*
\end{array} 
\right) 
\mod I  \} ,
$$
$$
K _ 0 ( I ) = \{ g \in 
\GL _2 (
\prod _{v \in \vert F \vert _f}  o _{F_v} ),\ g \equiv 
\left(
\begin{array}{cc}
* &* \\
0&*
\end{array}
\right) 
\mod I \} .
$$

For a non-archimedian local field $F$, let $B(F)$ be the standard Borel subgroup of $ \GL_2 ( F) $ consisting of the upper triangular matrices. For two quasi-characters
$\chi _i :  F ^\times \to \mathbb C^{\times}$ ($i = 1, 2 $)
$\chi: B ( F) \to \mathbb C ^{\times} $ is defined by $\chi  (
\left(
\begin{array}{cc}
a &b \\
0&d
\end{array} \right) 
 ) =\chi _1 (a ) \chi_2 (d) $.

The non-unitary induction $\pi ( \chi _1 , \chi _2 )= \Ind _{B(F)}  ^{G(F)} \chi$ is given by 
$$
\Ind _{B(F)}  ^{G(F)} \chi = \{ f : \GL_2 ( F ) \to
\mathbb C  ,\ 
f (
\left(
\begin{array}{cc}
a &b \\
0&d
\end{array} \right) 
 g ) =\chi _1 (a ) \chi_2 (d) \vert a \vert  f (g)  \} .
$$

For a prime $\ell$, we fix an isomorphism $ \C\simeq \bar \Q_{\ell}$ by the
axiom of choice. \par
Let $F$ be a totally real field. For an infinity type $(k , w)$ (see Definition \ref{dfn-shimura21}), which satisfies $ k _\iota \equiv w
\mod 2 $ for $\iota \in I_{F, \infty} $, $k ' \in \Z ^{I_{F, \infty} } $ is defined by the formula
$$ k + 2 k' = (w + 2 )
\cdot  (1, \ldots, 1 ). $$ 

For a cupspidal representation $\pi$ of $\GL_2 ( \mathbb A_F) $ of discrete infinity type $(k, w) $ such that $ \pi _f$ is defined over an $\ell$-adic field $E_{\lambda} $, $\rho _{\pi , E_{\lambda} } : G_F \to \GL _2 (E_\lambda ) $ is the associated $\lambda$-adic representation. 
For a finite place $v \nmid
\ell$ such that the $v$-component 
$\pi _v $ is spherical, a geometric Frobenius element
$
\Fr _ v
$ satisfies
$$ \tr \rho _{\pi ,  E_\lambda} (\Fr _v ) = \alpha _v + \beta _v
$$ where $ (  \alpha _v ,\ \beta _v ) $ is the Satake parameter of $\pi _v $ seen as a
semi-simple conjugacy class in dual group $ \GL_2 ^\vee ( \bar \Q_\ell)$, and
$\pi _v$ is a constituant of non-unitary induction $\Ind _{B(F_v)}  ^{G(F_v)} \chi _{\alpha _v, \beta _v}
$. Here $\chi _{\alpha _v, \beta _v}: B (F_v) \to
\bar
\Q ^{\times}_\ell
$ is the unramified character defined by $\chi _{\alpha _v, \beta _v} ( \left(
\begin{array}{cc}
a &b \\
0&d
\end{array} \right ) = \alpha _v ^{\ord _v a}
\beta _v ^{\ord _v d } 
$.

At infinite places, the $\GL_2 (\R)$-representation $D_{ k, w} $ corresponds to the unitary
induced representation
$$\Ind _{B(\R)} ^{G(\R)} ( \mu _{k, w} , \nu _{k, w} ) _{\bold u }
$$ for two characters of the split maximal torus
$$\mu _{ k, w} (a) = \vert a \vert ^{\frac  1 2 - k ' } (\operatorname{sgn}a) ^{- w}
$$
$$
\nu_{ k, w} (d) = \vert d \vert ^{\frac  1 2 -w + k'} 
$$ for $k'$ satisfying $ k -2 + 2 k ' = w $.
\medskip This normalization, which is the $\vert  \cdot \vert _v  ^{ \frac  1 2 }
$-twist of the unitary normalization, has the merit that it preserves the field of definition.
The central character of
$\pi$ corresponds to $
\det
\rho _{\pi,
E_\lambda} (1)
$. Our normalization is basically the same as that in \cite{Car2}, except one point. In
\cite{Car2}, an arithmetic Frobenius element corresponds to a uniformizer.

The global correspondence $\pi \mapsto \rho _{\pi , E_ \lambda} $ is compatible with the local
Langlands correspondence for $ v \nmid \ell$ (\cite{Car2} th\'eor\`eme (A), see
\cite{W1}, \cite{T1} theorem 2 for the missing even degree cases), namely, if we take the
$F$-semisimplification of
$
\rho_{\pi , \lambda} \vert _{G_{F_v}}
$, this corresponds to $\pi_v$ by the local Langlands correspondence normalized as above.

\section{Taylor-Wiles systems}\label{sec-TW}

We present an abstract formulation of
the argument of Taylor and Wiles \cite{TW}. Our proof is influenced by an argument of Faltings on their work (\cite{TW}, Appendix).  A similar method was
found by Diamond independently \cite{D2}.

\subsection{Definition of Taylor-Wiles systems} 
\label{subsec-TW1}
For a global field $F$, let $\vert F \vert _f $ be the set of finite places of $F$, and for $v
\in \vert F \vert _f $, $ q _v $ means the cardinality of the residue field $k (v) $ at $v$.

Let $ o_\lambda  $ be a complete noetherian local algebra with the maximal ideal $m_\lambda$. We assume that the residue field $k_\lambda =  o_\lambda /m_\lambda  $ is a finite field of characteristic
$\ell$.
\begin{dfn}\label{dfn-TW11} 
 Let $H $ be a torus over $F$ of dimension
$d$, 
$X $ a set of finite subsets of
$\vert F \vert _f $ which contains $\emptyset$. 
We take a pair $(R, M )$, where $ R $ is a complete noetherian local $ o _\lambda$-algebra with the residue field $k_\lambda$, and $M$ is a finitely generated $R$-module. 

A Taylor-Wiles system $\{ R _ Q , M _Q
\} _{ Q \in X }$ for $ (R,  M )$ 
consists of the following data:

\begin{itemize}
\item [TW1:] For $Q \in X  $ and $v \in Q$, $H$ is split at $v$, and $ q_ v  \equiv 1 \mod \ell $. We denote by $ \Delta _v $ the
$\ell$-Sylow subgroup of $ H (k (v) ) $, and $
\Delta _Q $ is defined as $\prod _{ v\in Q } \Delta _v $ for $Q \in X$.
\item [TW2:] For $Q\in X$,  $R _Q$ is a complete noetherian local $o _\lambda[\Delta _Q]
$-algebra with the residue field $k_\lambda$, and $ M_Q$ is an $R_Q$-module. For $Q = \emptyset $, $( R_{\emptyset} , M _{\emptyset} ) = (R, M )  $. 
\item [TW3:] A surjection 
$$
 R _Q /I_Q R_Q\twoheadrightarrow  R 
$$
of local $o _\lambda $-algebras for each $ Q \in X$. Here $I_Q \subset
o_\lambda [\Delta_Q]$ denotes the augmentation ideal of
$o_\lambda [\Delta_Q]$. For $Q = \emptyset $, it is the identity of $R$.

\item [TW4:]
The homomorphism $R _Q / I_Q R_Q \to \End _{o _\lambda}  M _Q / I_Q   M _Q$
factors through $R$, and $M _Q / I_Q   M _Q $ is isomorphic to $M$ as an $R$-module.

\item [TW5:]
$ M_Q  $ is free of finite rank $\alpha$ as an  $ o _\lambda [ \Delta_Q] $-module for a fixed integer
$ \alpha
\geq 1$. 
\end{itemize}
\end{dfn}

In \cite{TW}, the conditions that $R_Q$ is Gorenstein and
$M_Q$ is a free $R_Q$-module are required. \par 
Unlike Kolyvagin's Euler systems, we do not 
impose a functoriality when the index set grows. 

\begin{dfn}\label{dfn-TW12}
Let $\{ R _ Q , M _Q
\} _{ Q \in X }$ be a Taylor-Wiles system for $ (R,  M )$ with the coeffiient ring $o_\lambda$, $o_\lambda\to o ' _{\lambda '}$ a local homomophism between complete noetherian local rings with finite residue fields. Then the scalar extension  $\{ R _ {Q, o'_{\lambda'}  }  , M _{Q, o'_{\lambda'} } \} _{Q\in X} $ of $\{ R _ Q , M _Q \} _{ Q \in X }$ is defined by $ R_{Q, o'_{\lambda '} } = R_Q\hat \otimes _{o_{\lambda } }  o'_{\lambda '} $, $  M_{Q, o'_{ \lambda '} } = M_Q \otimes _{o _\lambda }  o'_{\lambda '}$. This is a system for $( R_{o'_{\lambda'} } , M _{o' _{\lambda '}} ) = ( R\hat \otimes _{o_{\lambda } } o'_{\lambda '} ,   M \otimes _{o _\lambda }  o'_{\lambda '} )  $. 

\end{dfn}

\subsection {Complete intersection-freeness criterion}
\label{subsec-TW2}
\begin{thm}\label{thm-TW21}[Complete intersection and freeness criterion]  For a Taylor-Wiles system $\{ R _ Q , M _Q
\} _{ Q \in X }$ for $ (R,  M )$ and a torus $H$ of dimension $d$, assume the following conditions.
\begin{enumerate}
\item 
For any integer $ n\in \N$, the subset $X_n  $ of $X$ defined by 
$$
X_n  = \{ Q \in X ;\  v  \in Q
\Rightarrow q_ v 
\equiv 1
\mod \ell ^n \}. 
$$ 
 is an infinite set. 
\item  $r=  \sharp Q   $ is independent of $ Q \in X$ if $Q \neq \emptyset$. 

\item $R _Q $ is generated by at most $d r$ many elements as a complete
local $
o _\lambda$-algebra for all $ Q\in X$.
\end{enumerate}
 Then under (1)-(3),  
\begin{itemize}
\item $R$ is $o _\lambda$-flat
and of relative complete intersection of dimension zero. 
\item$M$ is a free $R$-module.
\end{itemize}
In particular $R$ is isomorphic to the image T in $ \End_{o_\lambda} M $.

\end{thm}

\begin{proof}[Proof of Theorem \ref{thm-TW21}] First we treat the case when
$o_\lambda $ is a regular local ring. We choose the following data for each element $ Q \in  X\setminus \{ \emptyset \}$:
\begin{itemize}
\item An isomorphism 
$$
\alpha _Q : \Delta_Q \stackrel{\sim} {\longrightarrow} \bigoplus _{ v \in Q}  (\mathbb Z/ \ell ^{ n(v) }\mathbb Z )^{\oplus d}
$$ as finite abelian groups. Here $\ell ^{ n(v) }$ is the order of the $\ell$-Sylow subgroup of $k (v) ^\times$.
\item 
An isomorphism 
$$
 \beta _Q :  M _Q \stackrel{\sim} {\longrightarrow}( o _\lambda [
\Delta _Q ] ) ^\alpha
$$
 as $ o _\lambda [ \Delta _Q ] $-modules.
\item
A surjection $\gamma _Q:  o _\lambda [[ T_1,.., T _{dr }]]
\twoheadrightarrow R _Q $ as complete local $o _\lambda$-algebras. 
\end{itemize}
The existence of $\beta _Q$ and $\gamma _Q $ is assured by TW5 and \ref{thm-TW21}, (3).\par
By the assumption (1), $X_n $ is an infinite set. \par
For $\Gamma = \mathbb Z _{\ell} ^{\oplus dr}$, let $S_{\infty}$ be the complete group ring $o _\lambda   [[   \Gamma ]]$, and $I_{\infty} $ the augmentation ideal. \par

For an element $Q \in  X _n$, where $X_n $ is defined in \ref{thm-TW21} (1), 
$\alpha _Q$ induces an isomorphism $\Delta _Q / \ell ^n \Delta _Q \stackrel{\sim}{\rightarrow}  \Gamma /  \ell ^n \Gamma $ of abelian groups, so the $o _\lambda $-algebra isomorphism
$$
\epsilon _{ Q, n } : o _\lambda [ \Delta _Q ] / J_{Q, n}  \stackrel{\sim}{\longrightarrow} S_n 
$$
is induced. Here $J_{Q, n } $ is the ideal of $o _\lambda [
\Delta _Q ]$ generated by $m_\lambda  ^n$ and $
\delta  ^{\ell ^n } -1$ for $ \delta \in \Delta_Q $, $ S _n = (o_{\lambda} / m_\lambda  ^n ) [ \Gamma / \ell ^ n \Gamma ]  $, which is viewed as a quotient of $S_{\infty}$. $I_n = I _{\infty } S_n $ is the augmentation ideal of $S_n $ as the group ring of $\Gamma / \ell ^ n \Gamma $ over $ o_\lambda / m^n _{\lambda}$. \par
$\beta _Q$ and $\epsilon _{ Q, n }$ induce an isomorphism
$$
M _Q / J _{Q, n} M _Q \stackrel{\sim}{\longrightarrow}S_n ^\alpha 
$$
as $S_n$-modules.
$R_{Q, n  }$ is defined as the image of $ R _Q / J _ { Q, n} R _Q$ in $ \End _{ S_n } ( M _Q / J _{Q, n} M _Q) \stackrel{\sim}{\rightarrow}  \End _{S_n} S _n ^\alpha  $, which is viewed as an $S_n$-algebra by $\epsilon_{Q, n} $. \par
\bigskip
Now we apply the method of Taylor and Wiles to construct a projective system which approximates a power series ring.  
For an integer $n \in \N$ and $Q \in X _n$, we consider the couple $ ( (R_{Q,  n },   \iota _{Q, n }) , p _{Q, n}  )$ defined as follows:
\begin{itemize}
\item[($\ast$)]
$ R_{Q,  n } $ is a finite $S_n$-algebra with a faithful action on $S_n ^\alpha $
by $ \iota_{Q, n }: R_{Q, n } \hookrightarrow \End _{S_n} S _n ^\alpha $.  

\item[($\ast \ast$)] $p_{Q, n }: o_\lambda [[ T _1 \ldots, T_{dr} ]]\twoheadrightarrow  R_{Q, n}  $ is the surjective homomorphism 
$ o_\lambda [[ T _1 \ldots, T_{dr} ]] \overset {\gamma _Q } \to R_Q \to R_{Q, n}$ as local $o_\lambda $-algebras.
\end{itemize}
 
\medskip 
The isomorphism classes of
the couples 
$( (R_{ Q, n} ,
 \iota _{Q, n} ) , p_{Q, n} )
$ are finite as $Q$ varies in $ X _n $, since the cardinality of
$ R_{Q, n  } $ is bounded. Thus there is a sequence $\{ Q( n ) \} _{n \in \N } $ which satisfies the following properties:
\begin{itemize}
\item $ Q( n ) $ is an element of $X _n $.
\item For $m \leq n $,  $ ( (R_{ Q ( n ) , m } ,
\iota _{Q(n ) , m}) , p_{ Q(n ) , m } ) $
is isomorphic to $ ( (R_{ Q ( m ) , m } ,
\iota _{Q(m ) , m}) , p_{Q(m ) , m } )$.
\end{itemize}
For the sequence $\{  Q (n ) \} _{ n \in \N }
$ thus obtained, 
$( (R_{ Q(n) , n} ,
 \iota _{Q(n ) , n} ) , p_{ Q(n) , n })_{ n
\in \N}  $ forms a projective system. The transition map
$$
R _{ Q(n+1 ), n+1 }
 \longrightarrow R _{ Q ( n+1 ),n} 
\simeq R _{Q(n) , n } 
$$
is surjective for any $n \in \N$, since
$  R _{Q(n +1) ,  n} $ is a quotient of
$  R _{Q(n+1 ) ,  n+1 } \otimes _{ S_{n +1} }
S_{n}$.\par
By taking the projective limit, we define
$$
P = \varprojlim _n R_{ Q(n) , n  }, 
$$
which has a structure of $S_\infty $-algebra. By ($\ast \ast$), there is a surjection $
o _\lambda [[ T_1,.., T _{dr} ]] \twoheadrightarrow  P $. By the definition, $P$ has a
faithful non-zero module 
$$
L=\varprojlim_n  {S_n } ^{\oplus \alpha } \simeq \varprojlim _n  M _{Q(n) } / J _{Q(n) , n} M _{Q(n)} ,
$$
which is $S_\infty$-free and finitely generated over $S_\infty$.\par
\begin{lem} \label{lem-TW21} For the $o_\lambda $-algebra $P$ and the $P$-module $L$ thus defined, the following holds:
\begin{enumerate}
\item The local  $o_\lambda$-homomorphism $o _\lambda [[ T_1,.., T _{dr} ]] \twoheadrightarrow  P$ is an isomorphism. 
\item $ P $ is $S_\infty $-flat.
\item $ L$ is a non-zero free $P$-module. 
\end{enumerate}
\end{lem}
\begin{proof}[Proof of Lemma \ref{lem-TW21}]
We prove (1). The $o_\lambda$-algebra homomorphism $S_\infty 
\rightarrow P $ is injective since any element in the kernel must annihilate a non-zero free
$S_\infty $-module $ L $. Since
$ P $ is an $S_{\infty}$-subalgebra of $ \End _{S_\infty} L  $, $P$ is a finite
$  S_\infty $-algebra. It follows that $\dim P = dr +\dim o_\lambda$, and the surjection 
$o _\lambda [[ T_1, .., T _{dr} ]]\twoheadrightarrow P $ must be an isomorphism. \par
Since $P$ is regular and finite over $S_\infty$, (2) follows from \cite{GD2}, Chap. IV, proposition 6.1.5. \par 
$L$ is finitely generated and flat as an $S_\infty 
$-module. Then $ L $ is a finitely generated $P$-module, and $P$-free by the following sublemma. 
\begin{sublem}\label{sublem-TW21} Let $A$, $B$ be regular local rings, and $ M$ a finitely generated
$B$-module. If
$B$ is finite flat over $A$, and $M$ is $A$-flat, then $M$ is $B$-flat.
\end{sublem}
\begin{proof}[Proof of Sublemma \ref{sublem-TW21}] By the assumption, $ \operatorname {depth} _B M =  \operatorname {depth} A = 
\operatorname {depth} B$. Since the projective dimension of $M$ is finite, by the
Auslander-Buchsbaum formula \cite{GD1}, Chap.0, proposition 17.3.4, $
\operatorname{proj.dim}_B M = 
\operatorname {depth} B -  \operatorname {depth} _B M = 0 $ and hence $M$ is projective. 
\end{proof}
\end{proof}
By Lemma \ref{lem-TW21} (2), 
$$
\tilde R = P  /  I _{\infty} P 
$$
is finite flat as an
$S_\infty / I_{\infty } =o _\lambda$-module.\par

The ideal $I =  I_{\infty} P  
$ of $P$ is generated by $dr$ many elements. This means that $
\operatorname{ht} I\leq dr
$. Since $P$ is catenary, $ \dim \tilde R = \dim P -  \operatorname{ht} I \geq dr+
\dim o_\lambda  - dr = \dim o_\lambda $. On the other hand, $\dim
\tilde R \leq \dim o_\lambda  $ because $\tilde R$ is finitely generated as an
$o _\lambda$-module. Thus we have the equalities $ \dim \tilde R = \dim
o_\lambda 
$ and 
$\operatorname{ht} I = dr$. This implies that $dr$ is the minimal number of generators of $I$, so
$\tilde R$ is of relative complete intersection over $o_{\lambda }$.\par

By Lemma \ref{lem-TW21} (3), $\tilde R $-module $ L / I_{\infty}  L$ is a non-zero free module. $L_\infty / I _{\infty} L $ is naturally isomorphic to 
$$
 \varprojlim _n M _{ Q(n) }/ ( m_ \lambda ^n, I _{Q(n )} \
) M _{ Q(n) } \simeq \varprojlim _n  M /m_ \lambda ^n =M  
$$ 
by the definition and TW4. Since the pair $(\tilde R, M )$ has the desired complete intersection-freeness property, to conclude the proof of Theorem \ref{thm-TW21} when $o_\lambda $ is regular, it is sufficient to identify $R$ and $\tilde R$.
\par
\bigskip 
\begin{lem}\label{lem-TW22}
\begin{enumerate}
\item The canonical homomorphism $ P \otimes _{S_\infty} S_n   \twoheadrightarrow  R _{  Q (n ) , n } $ is an isomorphism, and  
$ M _{ Q (n ) , n  } / J _{ Q(n ), n} M _{ Q ( n ) ,n }$ is a non-zero free
$ P \otimes _{S_\infty} S_n 
$-module.
\item $\tilde R $ is a quotient of $R$.
\end{enumerate}
\end{lem}
\begin{proof} [Proof of Lemma \ref{lem-TW22}] 
$  L \otimes _{S_\infty} S_n  =  M _{ Q ( n )  } / J _{Q(n ), n} M _{ Q (n )}
$ is a free $ P \otimes _{S_\infty} S_n  $-module by Lemma \ref{lem-TW21}, (3). Since the $P\otimes _{S_\infty} S_n 
$-action on $  L \otimes _{S_\infty} S_n   $ factors through the quotient
$R _{ Q(n ), n} $, (1) is shown. \par
Since $M _{ Q (n) } /  J_{Q(n), n } M _{ Q (n) }$ is $R _{ Q(n ), n}$-free by (1), 
$$
  M _{ Q (n) } /  J_{Q(n), n } M _{ Q (n) } \otimes _{S_n } S_n / I_n =  M /m^n _\lambda  M  ,
$$
is a free $R _{ Q(n ), n}  \otimes _{S_n } S_n / I_n  =(P \otimes _{S_\infty} S_n) \otimes _{S_n } S_n / I_n = \tilde R / m^n _\lambda \tilde R $-module.  Thus $ \tilde R / m^n _\lambda \tilde R $ coincides with the image
of
$ R _{Q(n ), n} $ in $\End _{o_\lambda / m^n _{\lambda }  }  M /m^n _\lambda  M  $. The $R_{Q(n ) }$-action on $ M_{Q(n) } / I _{Q( n) }  M_{Q(n) } = M $ factors through $R$ by TW4, and hence $\tilde R / m^n _\lambda \tilde R $ is regarded as a quotient of $R$.  By passing to the projective limit, (2) is shown. 
\end{proof}
\bigskip

Fix an integer $N \geq 1$. Since the transition maps in the projective system $ \{ R_{ Q(n), n}
\} _{ n \in \N } $ are surjective, there is some integer $n
\geq N $ such that
$$
R_{Q(n ) ,  n} /m ^N _{  R_{Q(n ) ,  n } } \simeq P / m ^N _P 
\simeq o _\lambda [[ T_1,.., T _{dr} ]]/ (m_\lambda, T_1,.., T _{dr})^N .
$$
On the other hand, $R_{Q(n) } /m ^N _{R_{Q(n) } } $ is a quotient of
$ 
o _\lambda [[ T_1,.., T _{dr }]]/ (m_\lambda, T_1,.., T _{dr})^N
$ by the condition \ref{thm-TW21}, (3). As $ R_{Q(n ) ,  n} $ is a quotient of $R_{Q(n) }  $, $ R_{Q(n) } /m ^N _{R_{Q(n) } }  = R_{Q(n ) ,  n} /m ^N _{  R_{Q(n ) ,  n } }$ holds. \par

We define the $o_\lambda $-algebra $ \tilde R _{ Q(n) }$ by $ \tilde R _{ Q(n) } = R_{Q (n ) }  / I_{Q(n) } R_{Q(n)}$. Since $R_{Q(n ) ,  n} /  I_{Q(n) }R_{Q(n ) ,  n}   $ is isomorphic to $\tilde R / m^n _{\lambda } \tilde R$ by Lemma \ref{lem-TW22} (1), it follows that $ 
 \tilde R _{ Q(n) } /m^N _{\tilde R _{ Q(n) }} \stackrel{\sim}{\rightarrow} \tilde R / m ^N _{
\tilde R } $. The isomorphism $\tilde R _{ Q(n) } /m^N _{\tilde R _{ Q(n) }}
\stackrel{\sim}{\rightarrow} \tilde R / m ^N _{ \tilde R } $ is factorized as 
$$
\tilde R _{ Q(n) } /m^N _{\tilde R _{ Q(n)}} 
\twoheadrightarrow R / m ^N _R
 \twoheadrightarrow
\tilde R / m ^N _{ \tilde R } .
$$
This implies that $R / m ^N _R
 \stackrel{\sim}{\rightarrow}
\tilde R / m ^N _{ \tilde R }$, and hence
$R \stackrel{\sim}{\rightarrow} \tilde R$ by passing to the limit with respect to $N$.\par
\bigskip
We prove the general case.  Let $ k= k_\lambda$ be the residue field of $o_\lambda$, and $ (  R _{Q} \otimes _{o _\lambda} k, M_{Q} \otimes _{o _\lambda} k)_{Q\in X}$ the Taylor-Wiles system for $(R  \otimes _{o _\lambda} k, M \otimes _{o _\lambda} k  )$ obtained by the scalar extension. Since $k$ is regular, we have shown that $R\otimes _{o_\lambda} k$ is a finite $k$-algebra of complete intersection, and $M \otimes _{o_\lambda } k$ is $R\otimes _{o_\lambda} k$-free. In particular this implies that $R$ is finite over $o_\lambda$. \par
We take an integer $\beta \geq 1 $ and a surjective  
$R$-homomorphism $g:  R ^{\oplus \beta } \twoheadrightarrow M $ which induces an isomorphism $ (R/m_\lambda R
) ^{\oplus
\beta }
\simeq M /
m_\lambda M $. Let $M ' $ be the kernel of $g$, which is finitely generated as an $o_\lambda$-module. Since $M$ is $o_\lambda$-flat, $ \Tor _1^{o_\lambda } ( M , k_\lambda ) = \{0 \}$, and the sequence
$$
0  \longrightarrow M ' \otimes _{o_\lambda} k  \longrightarrow ( R ^{\oplus \beta }) \otimes _{o_\lambda} k \overset {g \otimes \id _k }\longrightarrow M  \otimes _{o_\lambda} k \longrightarrow  0 
$$
is exact. By Nakayama's lemma, $M '= \{0\} $, and $g$ is an isomorphism, that is, $M$ is $R$-free. \par
By the $o_\lambda$-flatness of $M$, $R$ is $o_\lambda $-flat.
Since $ R$ is $o_{\lambda}$-flat and $R \otimes _{o_\lambda} k  $ is of complete intersection, $R$ is of relative complete intersection over $o _\lambda$.
\end{proof}

\subsection{Hierarchy of complete intersection and freeness}\label{subsec-TW3}

In this subsection, we give a formulation of a level raising argument of \cite{W2}. \par
 
 Let $E_\lambda $ be an $\ell$-adic field, $ o_\lambda = o_{E_\lambda}$ the integer ring,  $m_\lambda$ the maximal ideal.

\begin{dfn}\label{dfn-raising1}
\begin{enumerate}
\item An admissible quintet is a quintet $(R, T, \pi , M ,  \langle \ , \
\rangle) $, where $R$ is a complete noetherian local
${o _\lambda}$-algebra, $T$ is a finite flat ${o _\lambda}$-algebra, $ \pi : R
\to T$ is a
surjective
${o _\lambda}$-algebra homomorphism,  $M$ is a faithful finitely generated
$T$-module which is $o _\lambda $-free, and $ \langle \ , \ \rangle  :M \otimes _{{o
_\lambda}} M \to {o _\lambda}$ is a perfect pairing which induces $ M \simeq \Hom _{o
_\lambda} ( M, {o _\lambda}) $ as a $T$-module. 
\item  An admissible quintet $(R, T, \pi , M ,  \langle \ , \ \rangle)$ is distingushed
if
$R$ is of complete intersection, and $M$ is a non-zero free $R$-module (it follows that $\pi $ is an isomorphism). 
\item  An admissible morphism from $ ( R', T' ,\pi' , M' , \langle \ , \
\rangle ' )
$ to $ ( R, T ,\pi , M , \langle \ , \ \rangle ) $ is a triple $(\alpha,\beta,\xi )$.
Here
$
\alpha : R' \to R  $, $ \beta: T' \to T $ are surjective ${o _\lambda}$-algebra
homomorphisms making the following diagram
$$
\CD 
R ' @> \alpha >> R \\
 @V \pi ' VV @V \pi VV \\ 
T' @> \beta >> T \\
\endCD
$$ 
commutative, and $\xi : M  
\hookrightarrow M '   
$ is an injective $T'$-homomorphism
 onto an
${o _\lambda}$-direct summand.
\end{enumerate}
\end{dfn}
A Taylor-Wiles system gives rise to a distinguished admissible quintet for a suitably chosen pairing on $M$ under the condition of Theorem \ref{thm-TW21}. Note that we do not assume that the restriction of $\langle \
, \
\rangle ' $ to $\xi ( M ) $ is $\langle \ , \ \rangle $ in the definition of admissible morphisms. \par
\bigskip
For an admissible morphism $(\alpha,\beta,\xi )$ from $ ( R', T' ,\pi' , M' , \langle \ , \
\rangle ' )
$ to $ ( R, T ,\pi , M , \langle \ , \ \rangle ) $, there is a criterion for $(R ', T',
\pi ' , M ', 
\langle
\ ,
\ \rangle ')$ to be a distinguished quintet under the condition that $(R , T,
\pi , M , 
\langle
\ ,
\ \rangle ) $ is distinguished.\par
For a finitely generated $o_\lambda $-module $ L$, we denote $\Hom _{ {o _\lambda}
} (L, {o _\lambda} )  $ by $L ^{\vee}$.

By the perfect pairings on $ M$ and  $M'$, we make the identifications $ M   \simeq M  ^{\vee} $ and $ M  ' \simeq (M ' ) ^{\vee} $,
and $\xi ^{\vee}$ is defined as  
$$
 \xi ^{\vee}:M  ' \stackrel{\sim}{\longrightarrow} (M ' ) ^{\vee} \longrightarrow M ^{\vee} \stackrel{\sim} {\longrightarrow} M ,
$$
such that 
$$ 
 \langle \ \xi ( x),\ y  \ \rangle ' =  \langle \ x, \  \xi^{\vee} (y)\
\rangle \quad  ( \forall x \in M ,\ \forall y \in M ')
$$
holds.\par

 For any noetherian local
$o_\lambda$-algebra $R $ with an $o_\lambda$-algebra homomorphism $f : R
\to o_\lambda$, one attaches a numerical invariant to it, following \cite{Len}: \par 
Consider the annihilator $\Ann_R ( \ker \ f) $ of $\ker \ f $ in $R$. Define an ideal
$\eta_f$ of $o_\lambda$ by
$$
\eta_f = \text{ the image
of }
\Ann_R ( \ker \ f) \text{ by }f .
$$
When $R$ is a finite flat Gorenstein $o _\lambda$-algebra, 
$\Ann _R ( \ker f) $ is generated by the image of $1$ of the $ o _\lambda$-module homomorphism 
$$
o_\lambda \overset { f ^{\vee}}
\longrightarrow
\omega_ R \stackrel{\sim}{\longrightarrow}R , 
$$
where $\omega_ R = R ^{\vee} $ is the dualizing module of $R$, and $\eta_f$ is generated by the image of $1$ of
$$
o_\lambda \overset { f ^{\vee}}
\longrightarrow
\omega_ R \stackrel{\sim}{\longrightarrow} R\longrightarrow o_\lambda.
$$

We fix an $o _\lambda$-algebra homomorphism $f_T : T \to o _\lambda$. $ f
_R$ (resp.$f_{T'} $, resp. $f '_{R' }$) is defined as $f _T \circ
\pi
$ (resp. $f_T \circ \beta $, resp. $f_T 
\circ
\beta \circ \pi ' $).

\begin{thm}\label{thm-TW2} (abstract level raising formalism)
For an admissible morphism between admissible quintets $(R ',
T',
\pi ' , M ', 
\langle
\ ,
\ \rangle ') \to (R, T, \pi , M ,  \langle \ , \ \rangle) $, we assume the following conditions:
\begin{enumerate}
\item $(R, T, \pi , M ,  \langle
\ , \ \rangle) $ is distinguished.
\item 
$T$ and $T'$ are reduced,  $ M'
\otimes_{o _\lambda} E_{\lambda} $ is
$ T'\otimes_{o _\lambda} E_\lambda 
$-free, and the rank is equal to $\rank _T M$.
\item The equality
$$ 
 \xi ^{\vee} \circ \xi  (M )  = \Delta \cdot M 
$$
holds for some non-zero divisor $\Delta $ in $T$.
\item The inequality 
$$
\length _{o _\lambda}
\ker f_{R'}  / ( \ker f_{R'}    ) ^2   \leq \length _{o _\lambda} \ker f_R / ( \ker
f_R ) ^2  + \length _{o _\lambda} {o _\lambda} / \ f_T ( \Delta) o_\lambda
$$
holds.
\end{enumerate}
 Then $(R' , T' , \pi ' , M'  ,  \langle \ , \ \rangle ') $ is also distinguished, that is, $\pi ' : R ' \simeq T ' $, $R'$ is of complete intersection, and 
$M' $ is $T' $-free.
\end{thm}
\begin{rem}\label{rem-TW1}
The freeness of $M'
$ follows from theorem 2.4 of \cite{D2}. The argument for the freeness in the following is due to T. Saito.
\end{rem}

We will show a slightly stronger statement than Theorem \ref{thm-TW2}. The deduction of Theorem \ref{thm-TW2} from Theorem \ref{thm-TW3} is left to the reader.
\begin{thm}\label{thm-TW3}(Lifting theorem)
For an admissible morphism between admissible quintets $(R ',
T',
\pi ' , M ', 
\langle
\ ,
\ \rangle ' ) \to (R, T, \pi , M ,  \langle \ , \ \rangle) $, we assume the following conditions:
\begin{enumerate}
\item $(R, T, \pi , M ,  \langle
\ , \ \rangle) $ is distinguished.
\item 
$(M '/ \xi (M )) _{E_\lambda } $ does not have any non-zero subquotient which is a $T_{E_\lambda}$-module. 
\item The equality
$$ 
 \xi ^{\vee} \circ \xi  (M )  = \Delta \cdot M 
$$
holds for some element $\Delta $ in $T$.
\item $ \eta _{f_T} \neq \{0 \} $, and the inequality 
$$
\length _{o _\lambda}
\ker f_{R'}  / ( \ker f_{R'}    ) ^2   \leq \length _{o _\lambda} \ker f_R / ( \ker
f_R ) ^2  + \length _{o _\lambda} {o _\lambda} / \ f_T ( \Delta)o_\lambda 
$$
holds.
\end{enumerate}
 Then $\pi ' : R ' \to T ' $ is an isomorphism, $R$ is of complete intersection, and there is a $T'$-free direct summand $ F ' $ of $M'$ such that the restriction of $\langle \  ,\  \rangle ' $ to $F'$ is perfect, and $\rank _{T'} F' = \rank _T M $ holds. 
\end{thm}

\begin{lem}\label{lem-raising1}
Assume that we are given an admissible morphism $(R ', T', \pi ' , M ', 
\langle \ , \ \rangle ') \to (R, T, \pi , M ,  \langle \ , \ \rangle) $ which satisfies conditions (1) and (2) of \ref{thm-TW3}.
If 
$  \xi^{\vee}  \cdot \xi  (M )  = \Delta \cdot M $ for some element $\Delta $ in
$T$, then
$$
\length _{o _\lambda} {o _\lambda} / \eta _{f_{T'}}\geq \length _{o _\lambda}
{o _\lambda} / f_T (\Delta ) o_\lambda +
\length _{o _\lambda} {o _\lambda} /
\eta_{f_T} $$ holds.
\end{lem}
\begin{proof}[Proof of Lemma \ref{lem-raising1}] First we show that $\Ann _{T'} ( \ker f_{T'} )M ' \subset \xi ( \Ann _M (\ker f_T) )$. Any element $f$ in $ \ker f_{T'} $ acts as zero on $\Ann _{T'} ( \ker f_{T'} )M '$, so that $ \Ann _{T'} ( \ker f_{T'} )M '$ admits a structure of $T ' /  \ker f_{T'} \simeq T / \ker f_{T} $-module. Since $ (M ' / \xi ( M ) )_{E_\lambda } $ does not contain any non-zero $T_{E_\lambda}$-modules by \ref{thm-TW3} (2), $\Ann _{T'} ( \ker f_{T'} )M '\subset \xi ( M ) $ because $\xi (M) $ is an
$o _\lambda$-direct summand of $M'$. If we identify $\Ann _{T'} ( \ker f_{T'} )M ' $ as a submodule $N$ of $ M$, any element in $N$ is annihilated by $ \ker f_{T}  $, and $N $ is contained in $ \Ann _M (\ker f_T)$.\par

By \ref{thm-TW3} (1), 
$ M $ is
$T$-free, so we have the equality
$\Ann _M (\ker  f_T)  = (  \Ann _T \ker f_T)  M$. By applying $\xi ^{\vee}$, we obtain
$
\xi^{\vee}  (\Ann _{T'} (
\ker f_{T'}) M')
\subset  \xi ^{\vee}  \circ \xi    ((  \Ann _T \ker f_T)    M)=( \Delta  \cdot   \Ann _T \ker f_T)  M $.\par

$h : M ' \otimes _{f_{T'} }o_\lambda \to M \otimes _{f_{T'}} o_\lambda = M \otimes _{f_T} o_\lambda $ is surjective. Since $M \otimes _{f_T} o_\lambda  $ is non-zero and $o_\lambda $-free, we find a surjective $o_\lambda$-homomorphism $ g: M \otimes _{f_T} o_\lambda \twoheadrightarrow L $ such that $L \simeq o_\lambda $. 
Because $o_\lambda $ is PID, there is an element $x$ of $M'$ such that $g \circ h ( x \otimes _{f_{T'} } 1)  $ generates $L $.  The image of $ \Ann _{T'} ( \ker f_{T'} ) x $ in $L$ is $\eta _{f_{T'}} L$. \par
$g \circ h ( \Ann _{T'} ( \ker f_{T'} ) x  \otimes _{f_{T'} } 1)
\subset g ((
\Delta \cdot   \Ann _T \ker f_T)  M  \otimes _{f_{T} } {o_\lambda}) =f_T (\Delta ) \cdot \eta _{f_T} g (M \otimes _{f_T} o_\lambda)= f_T (\Delta ) \cdot \eta _{f_T}L $. The claim is shown.

\end{proof}
\begin{proof}[Proof of Theorem \ref{thm-TW3} ]
First we show that $R'$ is of complete
intersection, and $ \pi ' :R '\rightarrow T '$ is an isomorphism. \par
Since $T$ is of complete intersection, the equality 
$$
\length _{o _\lambda}
{o_\lambda} / \eta _{ f_T} = \length _{o _\lambda} \ker
 f_T  / ( \ker f_T  ) ^2 
$$
holds by \cite{Len}.\par
We show $ \xi ^{\vee}_{E_\lambda}
\circ
\xi _{E_\lambda}$ is an isomorphism. By assumption \ref{thm-TW3}, (2), the restriction of $ \xi ^{\vee} : (M') ^{\vee} \to M ^{\vee}$ to $(M ' /\xi ( M )  ) ^{\vee} $ is zero. So the induced map $ \xi ^{\vee} _{E_\lambda } \vert _{\xi (M ) _{E_\lambda } } $ must be surjective. By compairing the dimension, it is an isomorphism. \par
In particular $\Delta $ is invertible in $ T _{E_\lambda }$, and $  \ f_T ( \Delta)  \neq 0$.\par
By assumption \ref{thm-TW3}, (4) and Lemma \ref{lem-raising1}, we have 
$$
\length _{o _\lambda}
\ker f_{R'}  / ( \ker f_{R'}    ) ^2 \leq \length _{o _\lambda} \ker f_R / ( \ker
f_R ) ^2  + \length _{o _\lambda} {o _\lambda} / \ f_T ( \Delta) o_\lambda 
\leqno{(\ast_1)}
$$
$$
 \leq  \length _{o _\lambda} {o _\lambda}/ \eta
_{f_{T '}}.
$$
Then the claim follows from the isomorphism criterion of Lenstra
\cite{Len}.\par
As a consequence, we have the equality 
$$
\length _{o _\lambda} o _\lambda/ \eta
_{f_{T '}}= \length _{o_\lambda}
o _\lambda / \eta _{f_ T}  + \length _{o _\lambda} o _\lambda / \ f_T (
\Delta)o_\lambda , 
\leqno{(\ast_2)}
$$
since the equality holds in $ (\ast_1)$. \par
\bigskip

We construct a $T'$-free direct summand of $M' $. 
There is a canonical $T'$-homomorphism 
$$
 \beta^ \vee :  \omega_T \longrightarrow \omega _{T'}
$$
by the duality. \par
Since $T$ and $T'$ are complete intersections, we choose isomorphisms $\delta_T :  T \stackrel{\sim}{\rightarrow}  \omega _T $ and $\delta _{T'} : T ' \stackrel{\sim}{\rightarrow}   \omega _{T'}$ as $T'$-modules. 
We define $\Delta _0 $ by 
$$
\Delta _0 :  T\overset {\delta _T } {\stackrel{\sim}{\longrightarrow} }  \omega_T  \overset { \beta ^\vee } \longrightarrow \omega _{T'} \overset {\delta ^{-1} _{T'} }  {\stackrel{\sim}{\longrightarrow}   }T ' \overset {\beta} \longrightarrow T .
\leqno{(\ast_3)}
$$

$\Delta _0$ is a multiplication by an element $\delta_0 $ in $ T $ since it is a $T'$-homomorphism.\par

 By assumption \ref{thm-TW3}, (1),
$M$ is
$T$-free. Let
$ ( e_j ) _{j
\in J} $ be a basis of
$M$, which yields an isomorphism
$e  : F=  T ^{\oplus c  }\stackrel{\sim}{\rightarrow} M $. 
For $j \in J$, we choose an element $\tilde e  _j  $ of $ M '$ such that  $\beta (\tilde e  _j ) = e_j$. $ ( \tilde e_j ) _{j
\in J}  $ gives a $T'$-homomorphism $ e' : F ' \to M'$, where $ F ' = {T '} ^{\oplus c } $. 
The diagram 
$$
\CD 
F'  @>\beta ^{\oplus c }  >> F  \\ 
@V e ' VV @V  eVV \\ 
M '  @> \xi >> M 
\endCD
$$
is commutative. \par
Let $B: F \to F'  $ be a $T'$-homomorphism defined by
$$
B:  F\overset  {(\delta _T) ^{\oplus c }  } {\stackrel{\sim}{\longrightarrow}}(\omega_T  ) ^{\oplus c }\overset {( \beta ^{\vee})^{\oplus c}  } \longrightarrow (\omega _{T'} ) ^{\oplus c} \overset {(\delta ^{-1} _{T'}) ^{\oplus c}  }  {\stackrel{\sim}{\longrightarrow}}  F ' . 
$$
By the definition, the composition $F\overset B  \to  F' \overset { \beta ^{\oplus c }} \to F  $ is $(\Delta _0 )^{\oplus c } $, and hence it is the multiplication by $\delta _0 $. \par 

Consider the $T'$-submodule $ ( e ' \circ B ) ( F) $ of $M'$. By assumption \ref{thm-TW3}, (2), $\coker ( \xi _{E_\lambda } ) $ does not have any non-zero $T' _{E_\lambda}$-subquotients which are $ T _{E_\lambda }$-modules, which implies that $  ( e ' \circ B ) ( T) $ is a $T'$-submodule of $ \xi _{E_\lambda } ( M _{E_\lambda} ) $. Since $ \xi ( M)$ is an $o_\lambda$-direct summand of $M'$, $ ( e ' \circ B ) ( T)  \subset  \xi ( M)$, 
and there is a $T'$-homomorphism $ f : F  \to M $ that
makes the following diagram commutative:
$$
\CD 
F  @>B >> F ' @>\beta ^{\oplus c } >> F \\ 
@V f VV @V e '  VV @V e VV\\ 
M @> \xi >> M ' @>
 \xi ^{\vee}>> M .\\ 
\endCD
\leqno{(\ast_4)} 
$$
 Since $\Delta $ is a non-zero divisor in $T$, there is a $T$-automorphism $\mu$ of $M$ such  that $  \xi ^{\vee}
\circ
\xi  = \Delta \cdot \mu 
$. 

We have 
$$
e \circ (\beta ^{\oplus c } ) \circ B =e \circ (( \Delta _0 )^{\oplus c } ) =  \delta _0 \cdot e  =  \xi  ^{\vee} \circ \xi \circ f =  \Delta \cdot \mu \circ f.
\leqno{(\ast_5)}
$$
This implies that $( \delta _0 ) \subset (\Delta )$ as ideals of $T$ by evaluating $(\ast_5)$ at some
$e_i
$, and
$$
\delta _0 =
\alpha
\cdot \Delta 
$$
 for some $ \alpha \in T$.  We show that $\alpha$ is a unit.

Since $T$ and $T'$ are Gorenstein, $\eta_{f_T}$ and $\eta_{f_{T'}} $ are generated by the images of 1 by
$o_\lambda \overset {( f_T) ^{\vee} } \to 
\omega _T\overset{\delta _T ^{-1} } {\stackrel{\sim}{\rightarrow}} T
\overset { f_{T} } \to  
o_\lambda 
$ and $o_\lambda \overset {( f_{T'} ) ^{\vee}} \to 
\omega _T '  \overset{\delta _{T'}  ^{-1} }  {\stackrel{\sim}{\rightarrow}} T'
\overset { f_{T'} } \to  
o_\lambda 
$, respectively.  By the definition of $\delta_0 $ $(\ast_3)$, we have $ \eta
_{f_{T '} }= f_T(\delta_0 ) \cdot \eta_{f_T}$, and hence 

$$
\length _{o_\lambda} {o_\lambda}/ \eta
_{f_{T '}}= \length _{o_\lambda}
{o_\lambda} / \eta _{ f_T}  + \length _{o_\lambda} {o_\lambda} / \ f_T (
\delta _0 ) o_\lambda
\leqno{(\ast_6)}
$$
holds. Comparing $(\ast_6)$ with $(\ast_2)$, $ f_T ( \delta _0 )o _{\lambda}  = f_T( \Delta )o _{\lambda}  $. This implies that $ f _T (\alpha ) $ is a unit in $o_\lambda $, then 
$\alpha $ must be a unit in $T$. Replacing $\Delta$ by $\alpha \cdot  \Delta $, we may assume that $\delta _0
=\Delta$.\par

From $(\ast_5)$, $ f : F \to M  $ is an isomorphism since $\delta _0
=\Delta $. Taking the dual of $(\ast_4)$, we have a commutative
diagram 
$$
\CD 
F  @>B >> F ' @>\beta ^{\oplus c }  >> F \\ 
@V f VV @V e' VV @V e VV\\ 
M @> \xi >> M ' \simeq {M'}^{\vee}@> \xi ^ \vee>> M\simeq  M ^{\vee}  \\
@.  @V (e ' )^{\vee} VV @V f ^{\vee}  VV\\
  @.  (F ' ) ^{\vee}@>B ^{\vee}>> F^{\vee} . \\
\endCD
$$
By the definition of $B$, $B ^{\vee} $ induces a $T$-isomorphism $  F ^{\vee} \simeq (F ' ) ^{\vee}\otimes _{T'}T $, and is identified with $(F ') ^{\vee} \to (F ' ) ^{\vee} \otimes _{T'} T$. 

Since $ f^{\vee} $ and $e$ are isomorphisms and $\beta ^{\oplus c } $ is surjective, $(  (e ' ) ^{\vee}\circ  e ' ) (F ') + I  ( F ' ) ^{\vee}= (F ') ^{\vee}  $, where $I = \ker ( T ' \to T)$. By Nakayama's lemma, $( (e') ^{\vee} \circ  e ' ) (F ') = (F') ^{\vee}$, and hence $ (e' ) ^{\vee} \circ  e ' $ is surjective. Since the surjectivity implies the injectivity for a homomorphism between finite free modules of the same rank over a commutative ring, $  (e' ) ^{\vee}  \circ  e' $ is an isomorphism, and $ e ' ( F' )$ is a $T'$-direct summand of $M'$. The identification $(M' ) ^{\vee} = \Hom_{o_\lambda} ( M' , o_\lambda ) \simeq M '$ is made by the pairing $\langle \ , \ \rangle ' $, and the composition of $F ' \overset {e'} \to  M' \stackrel{\sim} {\rightarrow} (M') ^{\vee} \overset{(e') ^{\vee}} \to (F') ^{\vee} $ is an isomorphism. Thus the restriction of the pairing to $F'$ is perfect. 
 
\end{proof}
 
\section{Galois deformations}\label{sec-galdef}

For a prime $\ell$ and an $\ell$-adic field $E_\lambda$ with the integer ring $o_{E_\lambda}$, the
maximal ideal of $o_{E_\lambda} $ is denoted as $\lambda$, and $k_\lambda = o_{E_\lambda} / \lambda $ is the residue field. 
For a perfect field $k$ of characteristic $\ell$, $W(k)$ denotes the Witt ring of $k$. \par
\bigskip
In this section, we discuss deformation properties of a Galois representation 
$$ 
\bar \rho : G _F \longrightarrow \GL _2 ( k ) , 
$$
where $k$ is a finite field of characteristic $\ell$, in particular when $F$ is a local field.\par
For a local field $F$, we denote the integer ring by $o_F$ and the residue field by $k_F$. Let $p$ be the residue
characteristic, $p_F$ a uniformizer, and $ q =
\sharp  k_F
$. By $G_F$ and $I_F$, we denote the absolute Galois group of $F$ and the inertia subgroup, respectively.\par

\subsection {Representations outside $ \ell$}\label{subsec-galdef2}

In this subsection, assume that $ \ell \neq p$.  The most basic classification is whether
$\bar
\rho
$ is absolutely irreducible or absolutely reducible. 

\medskip In the absolutely irreducible case, there are two subcases:
\begin{itemize}
\item[$0_E$.] $\bar \rho \vert _{I_F}$ is absolutely irreducible, $ q
\equiv -1
\mod \ell $. In this case,  $ \bar \rho$ is isomorphic to $  \Ind ^{G _F } _{  G _{\tilde  F}  } \bar \psi$, where $\tilde  F
$ is the degree two unramified extension of $F$, $\bar  \psi : G _{\tilde F} \to k ^{\times}
$ is a character which does not extend to $G_F$. We call this case {\it exceptional}. 

\item [$0_{NE}$.] $\bar \rho  $ is absolutely irreducible but not of type
$0_E$.
\end{itemize}
\medskip

In the absolutely reducible case, by extending
$k $ if necessary, we may assume that $\bar \rho$ is reducible over $k$. There is a 
character
$\bar
\mu : G_F
\to k ^{\times}
$ such that the twist $\bar \rho \otimes \bar \mu ^{-1}$ has the non-trivial
$I_F$-fixed part. There are three subcases:
\begin{itemize}
\item [$1_{SP}$.] For some character $ \bar
\mu : G_F \to k ^{\times}$, 
 $(\bar \rho \otimes \bar \mu ^{-1}) ^{I_F} $ is one dimensional, and
the semi-simplification
$(\bar
\rho
\otimes
\bar
\mu ^{-1} ) ^{\ss}$ is unramified.

\item[$1_{PR}$.] For some characters $ \bar
\mu, \bar \mu '  $, $\bar \rho \simeq   \bar \mu   \oplus \bar \mu '$, and $  \bar \mu ' /
\bar
\mu
 $ is ramified. 
\item [$2_{PR}$.]  For some character 
$\bar \mu $, $\bar \rho \otimes \bar \mu ^{-1} $ is unramified.
\end{itemize}

\bigskip Note that this classification is stable under a twist by a character $\bar
\rho \mapsto \bar \rho \otimes \bar \chi$.

\begin{dfn}\label{dfn-galdef21} Assume that $\ell \neq p $. Let $\bar \rho: G_F \to \GL_2 (k) $ be a Galois representation where $k$ is a finite field. 
\begin{enumerate}
\item For a reducible $G_F$-representation $\rho : G_F \to \GL_2 ( k ) $, a character $\bar \mu : G_F \to k^{\times}   $ is called a twist character for $\bar \rho$ if $(\bar \rho \otimes \bar \mu ^{-1}) ^{I_F} \neq \{ 0\} $. The restriction $\bar \kappa = \bar \mu \vert _{I_F} $ to $I_F$ is called a twist type of $\bar \rho$. 
\item In the case of $0_E$, assume that $\bar \rho \vert _{G_{\tilde F} }  $ is reducible over $k$. A character $\bar \psi$ which appear as a constituant of $\bar \rho \vert _{G_{\tilde F} }   $ is called an inertia character of $\bar \rho$. The restriction $\bar \psi \vert _{I_{\tilde F} } $ is called an inertia type of $\bar \rho$. 
\end{enumerate}
\end{dfn}
Since $\bar \kappa $ extends to a character of $G_F$, it is identified with a character of the $G_F$-coinvariant $(I^{\ab}_F  ) _{G_F}$, which is canonically isomorphic to $ o^{\times} _F $ by the local class field theory. \par
 In the cases of $1_{SP}$, and $2_{PR}$, a twist type $\bar \kappa$ is unique for a given $\bar \rho$.  In the case of $ 1_{PR}$, there are two choices of twist types. We choose one twist type in this case.

\medskip
\begin{rem}\label{rem-galdef21}
\begin{enumerate}
\item The notation `` case $A_a $'' has the following meaning. The
number $A$ is the maximum of dimensions of the
$I_F$-fixed part of twists of $\bar \rho$. $E$ (resp. $NE$, $PR$, $SP$) means exceptional (resp. non-exceptional, principal series, special). In the cases of $PR$ and $SP$, typical situations arise
from admissible irreducible representations of these types by the local Lanlglands correspondence.

\item For $ F=\Q_p \ (p \neq \ell)$, $0_{NE}$, $1_{SP}$, and $1_{PR}$ cases belong to case C, A,
B of \cite{W2}, respectively. $0_E$-case is treated by Diamond \cite{D1}.
\end{enumerate}
\end{rem}

\subsection{Deformations outside $\ell$}\label{subsec-galdef3}
Let $o_\lambda $ be a complete noetherian local ring with the maximal ideal $m_{o_\lambda }$ and the finite residue field $k_{o_\lambda}$ of characteristic $\ell$.\par
By $ \mathcal C^{\noeth}_{o_{\lambda} }$ (resp. $ \mathcal C^{\artin}_{o_{\lambda} }$), we mean the category of the complete noetherian (resp. artinian) local $o_{\lambda}$-algebras $A$ with the maximal ideal $m_A$ such that the residue field $k_A  $ is $k_{o_\lambda }$, where morphisms are local $o_{\lambda}$-algebra homomorphisms. 
By a deformation of $\bar \rho: G_F \to \GL _n (k_{o_\lambda }) $, we mean a continuous representation
$$ 
\rho : G _F \longrightarrow\GL _n ( A )
$$
such that $ \rho \mod m _A  = \bar \rho$. Two deformations $\rho $ and $\rho '$ of $\bar \rho$ is isomorphic if there is a $G_F$-isomorphism $f:  \rho \simeq \rho '$ such that $f \mod m_A$ is the identity. 
In this paper, the case of $n \leq 2$ is considered.

\begin{dfn}\label{dfn-galdef31}(finite and unrestricted deformations)  
\begin{enumerate}
\item A deformation $\rho$ of $\bar \rho$ without any conditions is called an unrestricted deformation.
\item Assume that $\ell \neq p$.  
Let $\chi : G_F \to  A^{\times}$ be a deformation of $ \bar \chi : G_F \to k^{\times} _{o_\lambda} $. We say $\chi$ is a finite deformation if $\chi / \bar \chi_{\lift} $ is unramified. Here $ \bar \chi_{\lift} $ is the Teichm\"ller lift of $\bar \chi$. 
\item  Assume that $\ell \neq p$.  
 Let 
$\rho : G_F 
\to
\GL_2 (A)$ be a deformation of $\bar \rho:G_F 
\to
\GL_2 (k_{o_\lambda}) $. We say $\rho$ is a finite deformation if
$\det \rho $ is a finite deformation of $ \det \bar \rho 
$, and satisfies one of the following conditions:
\begin{itemize}
\item If $\bar \rho$ is reducible with a twist type $\bar \kappa$, then the $I_F$-fixed part $ (\rho \vert _{I_F} \otimes \kappa ^{-1}   )^{I_F}$ is $A$-free and is an $A$-direct summand, and the rank is equal to $ \dim _{k_\lambda } (\bar \rho \vert _{I_F} \otimes \bar \kappa ^{-1} ) ^{I_F} $. Here $\kappa= \bar \kappa _{\lift} : I_F \to W ( k_{o_\lambda } ) ^{\times } \hookrightarrow o^{\times}_{\lambda} $ is the Teichm\"uller lift of $\bar \kappa$. 

\item  If $ \bar \rho $ is absolutely irreducible and of type $0_{NE}$, we only require the condition on the determinant. 
\item  If $ \bar \rho $ is absolutely irreducible and of type $0_E$ of the form
$
\bar
\rho=
\Ind ^{ G_F }  _{ G _{\tilde F} } \bar  \psi $, then we require that 
$
\rho
\vert _{ I _{\tilde F} } $ is the sum of $ \psi \vert _{I_{\tilde F}} $ and its Frobenius twist. Here $\psi= \bar \psi_{\lift} :  G_{\tilde F } \to W ( k_{o_\lambda } ) ^{\times } \hookrightarrow o^{\times}_{\lambda}$ is the
Teichm\"uller lift of $\bar \psi$. 
\end{itemize}
\end{enumerate}
\end{dfn}

By $ F^{\bold u} _{ \bar \rho } : \mathcal C^{\noeth}_{o_{\lambda} } \to \Sets  $ (resp. $  F^{ \bold f}_{\bar \rho } : \mathcal C^{\noeth}_{o_{\lambda} } \to \Sets $) we denote the
functor defined on $ \mathcal C^{\noeth}_{o_{\lambda} } $ consisting of the isomorphism classes of the unrestricted deformations (resp. finite deformations). For a continuous character $ \chi : G_F \to o ^{\times} _\lambda $, $ F^{\bold u}_{ \bar
\rho ,  \chi} $ (resp. $ F^{ \bold f}_{\bar \rho , \chi } $) is the subfunctor
of $ F^{\bold u}_{  \bar
\rho  } $ (resp. $ F^{ \bold f } _{ \bar \rho  } $) consisting of the deformations whose  determinant is
$\chi$. 
\medskip
\begin{rem}\label{rem-galdef31}
Since we do not assume $ \dim _{k_{o_\lambda} } \Hom_{G_F}  ( \rho , \rho ) = 1$ in
general, these local deformation functors may not be representable, though the restriction to $\mathcal C^{\artin}_{o_{\lambda} }  $ always has
a versal hull in $\mathcal C^{\noeth}_{o_{\lambda} } $ by \cite{Sch}. A versal hull represents the deformation functor if $ \dim _{k_{o_\lambda}} 
\Hom_{G_F}  ( \rho , \rho ) = 1$.  
\end{rem}

\subsection{Tangent spaces outside $\ell$}\label{subsec-galdef4}

 As usual, we set
 
$$
 \ad \bar \rho = Hom _{k_{o_\lambda} }  ( \bar\rho ,\ \bar\rho ) .
$$

The tangent space $ F^{\bold u }_{  \bar \rho } ( k_{o_\lambda}[\epsilon] ) $ of $F^{\bold u }_{  \bar \rho , 
}$ is canonically isomorphic to $H ^1 (F, \ \ad \bar \rho )$ as in \cite{M}, which is of finite dimension over $k_{o_\lambda}$ by the finiteness of local Galois
cohomology groups.  Here $k_{o_\lambda} [\epsilon] $ is the ring of dual numbers. The trace
$0$-part $ \ad ^0 \bar
\rho$ of $\ad \bar \rho$ gives the tangent space of $F^{\bold u} _{ \bar \rho ,  \chi } $ consisting of the deformations with the determinant $\chi$.
\medskip
The finite deformations are controlled by the mod $ \ell $-version of the finite
part of Bloch-Kato \cite{BK}. For any $G_F$-module $M$ with order prime to the residual
characteristic, we define the finite part by

$$
H ^ 1 _f (F,\ M )  = H ^ 1 ( F^{\unr}/ F ,\ M ^{I_F} ) = \ker (H ^ 1  (F, M ) \longrightarrow  H ^1
(F^{\unr}, \ M  ) ). 
$$

\begin{prop}\label{prop-galdef41}

An element $x$ in  $ F^{ \bold u}_{ \bar \rho  } ( k
_{o_\lambda} [\epsilon] )  $ corresponds to a finite deformation if and only if it belongs to $ H ^ 1 _f(
F,\ \ad  \bar \rho ) $.
\end{prop}
\begin{proof}[Proof of Proposition \ref{prop-galdef41}] In the cases of $2_{PR}$, $0_E$, and $ 0_{NE}$, this is clear from the definition. 

 In the case of $ 1_{PR}$, any deformation splits by the vanishing of $ H^1 (F,\ \bar \mu /
\bar
\nu)$, so $\rho\vert _{I_F}$ is a constant deformation.

In the case of $1_{SP}$, the $I_F $-action factors through the maximal pro-$\ell$ quotient $ t_\ell
:I_F
\twoheadrightarrow
\Z_\ell (1)$, so the monodromy group is topologically generated by any element
$\sigma
 \in I_F$ so that $t_\ell(\sigma)$ is a topological generator. We take a basis so that 
$$ 
\rho (\sigma ) =
\begin{pmatrix} 
1 &1\\ 
0 &1 
\end{pmatrix} 
.
$$ 
This shows that $\rho \vert_{I_F}$ is constant. 
\end{proof}
\bigskip

The calculation of the dimension of tangent spaces is done as follows. First assume that $\bar
\rho$ is absolutely irreducible. Note that the deformation functors are representable in this
case.

\begin{lem}\label{lem-galdef1}
Assume that $ \bar \rho $ is absolutely irreducible. Then the finite part
$H ^1 _f ( F,\  \ad ^0 \bar \rho ) $ always vanishes, and 
$$ 
H ^1 ( F , \ad ^0 \bar \rho ) = \{ 0 \}
$$ 
 in the case of $0_{NE}$. In the case of $0_E$, 
$H ^1 ( F, \ad ^0 \bar \rho )$ is one dimensional.
\end{lem}
\begin{proof}[Proof of Lemma \ref{lem-galdef1}]
Since $\bar \rho $ is absolutely irreducible, $ H ^0 ( F ,\ad ^0 \bar \rho )=\{ 0\}$. By the Euler characteristic formula 
$$
 \dim _{k_{o_\lambda}} H ^1 ( F , \ad ^0 \bar \rho )
=  \dim _{k_{o_\lambda}} H ^0 ( F ,\ad ^0 \bar \rho )
+  \dim _{k_{o_\lambda}} H ^2 ( F , \ad ^0 \bar \rho)=  \dim _{k_{o_\lambda}} H ^2 ( F , \ad ^0 \bar \rho),
$$ and $ \dim _{k_{o_\lambda}} H ^2 ( F , \ad ^0 \bar \rho)=  \dim _{k_{o_\lambda}} H ^0 ( F , \ad ^0 \bar \rho (1))$ by the local duality. if $H ^0 ( F , \ad ^0 \bar \rho (1))\neq \{ 0\} $, $\bar \rho \simeq \bar \rho (1) $, and it is easily seen that $\bar \rho $ is of type $0_E$ and $H ^0 ( F , \ad ^0 \bar \rho (1)) $ is one dimensional. 
\end{proof}

\subsection{Local deformation spaces: examples}\label{subsec-galdef5}

Let $\chi : G _F \to o^{\times} _\lambda $ be an unramified continuous character. 
In the case
of $0_E$, $ F^{\bold u } _{ \bar \rho , \chi }$ is representable by remark \ref{rem-galdef31}.

The universal deformation ring $R $ for $ F ^{\bold u} _{ \bar \rho , \chi } $ is determined in
\cite{D1}, p.141. We recall the result.

$
\bar
\rho =
\Ind ^{ G_F} _{ G _{ \tilde F  } } \bar \psi $ for the unramified extension
$\tilde  F
$ of
$F$ and a character $\bar \psi :  G _{ \tilde F }  \to k ^{\times}_{o_\lambda}
$. Take the Teichm\"uller lift $  \psi : G _{ \tilde F } \to W (k _{o_\lambda} ) ^{\times} \hookrightarrow  o^{\times}_{\lambda} 
$ of
$\bar \psi$. If we consider a deformation
$
\rho : G _F
\to \GL _2 (A) $ for a local $o_{E_\lambda} $-algebra $A$ fixing the determinant,
$
\rho
\vert _{G _{ \tilde F } }
$ decomposes into sum of different characters, and for a suitable lift $
\phi$ of
$
\bar \psi$, 
$
\rho =
\Ind ^{ G _F} _ {G {\tilde F} } \phi$. $ \xi =  \phi  / \psi $ is a character of $ G _{
\tilde F } $ of $\ell$-power order. By the local class field theory
$\xi $ is a character of $ o ^{\times}  _{\tilde F}
$ of $\ell$-power order.

By the consideration, 

\begin{prop}[Diamond]\label{claim-galdef1}
$R$ is isomorphic to the group ring $o_{E_\lambda}[ 
\Delta _{\tilde F} ]
$. Here $ \Delta _{\tilde F} $ is the $\ell$-Sylow subgroup of $ k _{\tilde F} ^{ 
\times}.
$
\end{prop}

We need one more case where a versal hull is explicitly determined. 

\begin{prop}[Faltings]\label{claim-galdef2}
 Let $ \bar \rho$ be an unramified representation with different Frobenius eigenvalues, and
$ q
\equiv 1
\mod \ell$.  Then the versal hull of $F^{\bold u }  _{\bar \rho , \chi } $ over $o_{\lambda}$ fixing the determinant is isomorphic to 
$
o_\lambda [[ \widehat {(F^\times)}  _\ell ]] $. Here $ \widehat {(F^\times)}  _\ell$ is the
pro-$\ell$-completion of
$F ^\times$.
\end{prop}

This is found in \cite{TW}, Appendix, p. 569.  In fact, one shows that any
lifting $\rho$ to an artinian local $o_\lambda $-algebra
$A$ is decomposable under the assumptions of \ref{claim-galdef2}. We construct an isomorphism explicitly
for our later needs. 

We choose a decomposition $\bar \rho =\bar \chi _1 \oplus \bar \chi _2 $. $\chi _1 ,\
\chi _2$ denote the Teichm\"uller liftings of $ \bar \chi _1 ,\  \bar \chi _2$. We may assume that $ \rho$ takes the form
$$
\rho = 
\begin{pmatrix}
 \tilde \chi _1   & 0 \\ 
0 &  \tilde \chi _2  \\
\end{pmatrix}
$$
for liftings $\tilde \chi _1, \ \tilde \chi _2$ of $\bar \chi _1 ,\  \bar \chi _2 $.

$\delta = \tilde \chi _2  / \chi _2 : G_F \to A^\times $ is a character of
$\ell$-power order. So it factors through $G_F \overset {\delta^{\univ} }\to
o_\lambda [[
\widehat {(F^\times)}  _\ell ]]
\to A$ by the local class field theory.  
Note that $\widehat {(F^\times)}  _\ell $ is isomorphic to $ \Z _\ell \times \Delta _F$ if we choose a
uniformizer of $F$. Here $\Delta _F $ is the $\ell$-Sylow subgroup of $ k_F ^{\times} $. \par

In these examples, the local deformation space is the group ring of a
commutative group. 

\subsection {Deformations at $\ell$: nearly ordinary and flat deformations}\label{subsec-galdef6}
In this subsection, we assume that the residual characteristic $p$ of $o_F$ is equal to $\ell$. 
\begin{dfn}\label{dfn-galdef61}
Let $A$ be a complete local noetherian ring with the maximal ideal $m_A$ and the finite residue field $k_A$ of characteristic $\ell$. 
\begin{enumerate}
\item A pair $(\rho , \mu ) $ of a continuous $G_F$-representation $\rho : G_F \to \GL_2(A)$ and a continuous character $\mu : G_F \to A ^{\times } $ is called nearly ordinary if $\rho$ is isomorphic to 

$$
\rho   \sim  
\begin{pmatrix} 
\chi _1 & * \\ 
0 & \chi _2
\end{pmatrix} 			, 
$$
where $\chi _1 = \mu$. $\mu$ is called the nearly ordinary character for $(\rho , \mu)$, and the restriction $\kappa = \mu \vert _{I_F} $ is called the nearly ordinary type of $(\rho , \mu) $. 
\item When $A $ is a field, a reducible representation $\rho: G_F \to \GL_2 (A)  $ is called $G_F$-distinguished if one of the following conditions holds.  $\rho$ is indecomposable, or it is semi-simple and has two distinct constituants. 
\end{enumerate}
\end{dfn}
By the local class field theory, $ (I^{\ab} _F) _{G_F} \simeq o^{\times }_F  $, so a nearly ordinary type is regarded as a character of $ o^{\times }_F$. 

When a nearly ordinary character is clearly specified, we say $\rho$ is a nearly ordinary representation with the nearly ordinary character $\mu$ for short.

\begin{dfn}\label{dfn-galdef615}
For any complete local noetherian ring $A$ with the maximal ideal $m_A$ and the finite residue $k_A$ of characteristic $\ell$, 
\begin{enumerate}
\item We define
$$
H ^1 _{\flatcoh} (  o_{F^{\unr} }, \ A(1) ) 
\overset{def} =  \varprojlim _n H^ 1 _{\flatcoh} (  o_{F^{\unr} },\  A/ m ^n _A  (1)
)  = \varprojlim _n  o_{F^{\unr} } ^\times / ( o_{F^{\unr} } ^\times ) ^{\ell ^ n
} \otimes _{\Z_\ell} A ,
$$
$$
H ^1  (  F^{\unr} , \ A(1) ) 
\overset{def} =  \varprojlim _n H^ 1  ( F^{\unr} ,\  A/ m ^n _A  (1)
)  = \varprojlim _n ( F^{\unr} ) ^\times /( ( F^{\unr}  )^\times ) ^{\ell ^ n
} \otimes _{\Z_\ell} A .
$$
\item The exact sequence
$$
1 \longrightarrow o_{F^{\unr} } ^\times \longrightarrow (F^{\unr})  ^\times \longrightarrow  \mathbb Z \longrightarrow 0 
$$ 
induces
$$
 0 \longrightarrow H ^1 _{\flatcoh} ( o_{F^{\unr} } ,
\ A(1) ) \longrightarrow  H ^1 ( F^{\unr} , \ A (1) ) ) \overset {\ev _A } \longrightarrow A \longrightarrow 0 .
$$
We call $\ev _A  $ the evaluation map. 
\end{enumerate}
\end{dfn}
$H ^1 _{\flatcoh} (  o_{F^{\unr} }, \ A(1) )$ (resp. $H ^1(  F^{\unr} , \ A(1) ) $ classifies (pro-systems of) the extensions
$0 \to A \otimes _{\Z_\ell }  \varprojlim _n \mu _{\ell ^ n } 
\to M
\to A
\to 0  
$ over $ o_{F ^{\unr}} $ (resp. $ F ^{\unr} $) with the flat (resp. \'etale) topology.

\begin{dfn}\label{dfn-galdef62} Let $A$ be a complete local noetherian ring with the maximal ideal $m_A$ and the finite residue $k_A$ of characteristic $\ell$,  $\rho : G_F \to \GL_2 (A) $ a nearly ordinary representation of the form 
$$
0 \longrightarrow \chi_1 \longrightarrow \rho \longrightarrow \chi_2 \longrightarrow 0 , 
$$
where $\chi _1 : G_F \to A ^{\times }$ is the nearly ordinary character. 
\begin{enumerate}
\item Assume that $\rho$ satisfies $\chi _1 \vert _{I_F } = \chi_2( 1) \vert_{I_F} $. 
Then $c_\rho $ is  the extension class in 
$ H ^1 ( F ^{\unr} , \ A (1) ) $ defined by $\rho \otimes \chi _2  ^{-1} $. 

\item A nearly ordinary representation $\rho $ is called nearly ordinary finite if  $\chi _1 \vert _{I_F } = \chi_2( 1) \vert_{I_F} $, and if the class $c_\rho$ belongs to
$  H ^1 _{\flatcoh} ( o_{F^{\unr} } ,
\ A(1) ) \subset  H ^1 ( F^{\unr} , \ A (1) ) $. Equivalently, $\ev _A (c_{\rho} ) = 0 $. 
\end{enumerate}
\end{dfn}

\begin{prop}\label{prop-galdef61}
Let $E_\lambda$ be an $\ell$-adic field, $A$ an $o_{E_\lambda}$-finite flat local algebra, $\rho :G_F \to
\GL_2 (A) 
$ a nearly ordinary representation with the nearly ordinary character $\chi$ which satisfies
$\det \rho \vert _{I_F} = (\chi \vert _{I_F} )^2 ( -1) $.  Assume that there is an
$\ell $-divisible group $ M$ over $ o_F$ whose Tate module $T_\ell (M_ F)$ satisfies 
$$
T_\ell (M_ F ) \otimes _{o_{E_\lambda}} E_\lambda \simeq (\rho \otimes \chi  ^{-1} (1) )
\otimes  _{o_{E_\lambda}} E_\lambda  . 
$$
Then $\rho$ is a nearly ordinary finite
representation. 
\end{prop}

\begin{proof}[Proof of Proposition \ref{prop-galdef61}] We may assume that $\chi= 1$ by a twisting. Let $c _{\rho} $ be the extension class in $ H ^1 ( F ^{\unr} , \ A (1) ) $ which corresponds to $\rho\vert _{F ^{\unr} }  $. It suffices to prove $\ev _A ( c _{\rho})= 0   $ for the evaluation map $ \ev _A $ in Definition \ref{dfn-galdef615} (2). 

By the nearly ordinariness and the
assumption on the determinant of $\rho$, the generic fiber $ M_F
$ admits a filtration by an $\ell$-divisible subgroup over $F$
$$
0 \longrightarrow M_{1,  F} \longrightarrow M_F \to M_{2, F}  \longrightarrow 0 , 
$$
where $T_\ell ( M_{i , F}) $ are $ A\otimes _{o_{E_\lambda} } E_\lambda $-modules of rank one, and both
$T_\ell ( M _{1, F } (-1) )$ and $T_\ell (M_{2, F} )$ are unramified as Galois modules. We view
$M_{1, F} $ as the generic fiber of an $\ell$-divisible group of
multiplicative type $M_1$. By \cite{Ta}, Proposition 12, there is a morphism $i:  M _1 \to M $ as
$\ell$-divisible groups which gives the exact sequence over $F$ as above. This implies that the multiplicative part $M ^{\mult} $ of $M$ satisfies $T_{\ell} ( M ^{\mult} ) _F \simeq T_{\ell} ( M _1)_F  $. Thus we may assume that $i $ is
a closed immersion as an $\ell$-divisible group, and the quotient $M_2 $ is \'etale. \par 
The exact sequence of $\ell$-divisible groups over $o_{F ^{\unr} } $
$$
0 \longrightarrow M_1 \longrightarrow M \longrightarrow M_2 \longrightarrow 0 
$$
defines, in the fppf-topology, an element $c_M$ of $ \varprojlim _n H ^ 1 _{\flatcoh} ( o_{F ^{\unr} } , \ \Hom_A (M_2 [\ell ^n ] , M_1 [ \ell ^n ]  ) )=    H ^ 1 _{\flatcoh } ( o_{F ^{\unr} } , \ \Z_\ell (1) ) \otimes _{\mathbb Z_\ell} N
$. Here $N =  \Hom _A (   T_{\ell} (M_2 ) , T_{\ell} (M_1 (-1) ) )$.  On the other hand, since for the rank $2$-lattice $L = A ^{\oplus 2}$, the representation $\rho(1) $ is defined, and the subspace $L_1 $ where $\rho$ acts by $\chi_1 (1)$ satisfies $L _1 \otimes _{o_{E_\lambda} } E_\lambda = M_1  \otimes _{o_{E_\lambda} }  E_\lambda$. By this identification, 
$$    ( H ^ 1  ( F ^{\unr}  , \ \Z_\ell (1) ) \otimes
_{\Z_\ell} N   ) \otimes _{o_{E_\lambda} } E_\lambda =   H ^ 1 ( F ^{\unr}  , A (1) )  \otimes _{o_{E_\lambda} } E_\lambda , 
$$
and
$$    (H ^ 1 _{\flatcoh } ( o_{F ^{\unr} } , \ \Z_\ell (1) ) \otimes
_{\Z_\ell} N   ) \otimes _{o_{E_\lambda} } E_\lambda =   H ^ 1 _{\flatcoh } ( o_{F ^{\unr} } , A (1) )  \otimes _{o_{E_\lambda} } E_\lambda 
$$
hold, and $c_\rho $ spans the same $ E_\lambda$-subspace as $ c_M$ in $ H ^ 1  ( F ^{\unr}  , A (1) )\otimes _{o_{E_\lambda} } E_\lambda $. This implies that $\ev _A ( c_\rho )  $ is zero in $ A\otimes  _{o_{E_\lambda} } E_\lambda$ because $c_M$ belongs to $H ^ 1 _{\flatcoh } ( o_{F ^{\unr} } , \ \Z_\ell (1) ) \otimes _{\mathbb Z_\ell} N$. Since $A$ is $\ell$-torsion free, $\ev _A (c_\rho )$ is zero in $A$, and hence $\rho$ is nearly ordinary finite. 
\end{proof}
Now we introduce deformation functors for nearly ordinary representations.  
\begin{dfn} \label{dfn-galdef63} Let $o_\lambda $ be a complete local noetherian ring with the maximal ideal $m_{o_\lambda} $ and the finite residue field $k_{o_\lambda}$ of characteristic $\ell$. 
We fix continuous characters $\kappa :( I  ^{\ab} _F )_{G_F} \to o^{\times}_{\lambda}  $, and  $\chi : G_F \to o ^{\times}_{\lambda} $. Let $\bar \rho  : G_F \to \GL_2 (k_{o_\lambda}) $ be a $G_F$-distinguished nearly ordinary representation whose nearly ordinary type $ \bar \kappa$ is $ \kappa \mod m_{o_\lambda}$, and the determinant is $\chi \mod m_{o_\lambda} $. 
\begin{enumerate}
\item We define the functor
$$ 
F ^{\bold {n.o.} }_{ \bar \rho ,  \kappa }  :  \mathcal C^{\noeth} _{o_{\lambda} }\to \Sets
$$
as follows: for each $A \in \ob  \mathcal C^{\noeth}_{o_{\lambda} }$, $ F ^{\bold {n.o.} }_{ \bar \rho ,  \kappa } (A)$ is the set consisting of isomorphism classes of $ (\rho  ,  \xi ) $. Here $\rho : G_F \to \GL_2 (A) $ is a nearly ordinary representation with the nearly ordinary type $\kappa$, $\xi : (\rho \mod m_A) \simeq \bar \rho $ is an isomorphism as a $G_F$-representation. \par
$ F ^{\bold {n.o.} }_{ \bar \rho ,  \kappa , \chi } $ denotes the subfunctor of $ F ^{\bold {n.o.} }_{ \bar \rho ,  \kappa } $ consisting of the deformations whose determinant is $\chi$. \par

\item Assume moreover that $\bar \rho \vert_{I_F} $ has a form 
$$
0 \longrightarrow \bar \kappa
\longrightarrow  \bar \rho 
\longrightarrow \bar \kappa (-1) \longrightarrow 0.
$$
Then 
$$ 
F ^{\bold {n.o.f} }_{ \bar \rho ,  \kappa }  :  \mathcal C^{\noeth}_{o_{\lambda} }\longrightarrow \Sets
$$
is the subfunctor of $F ^{\bold {n.o.} }_{ \bar \rho ,  \kappa } $ consisting of the nearly ordinary finite deformations. $ F ^{\bold {n.o. f} }_{ \bar \rho ,  \kappa , \chi } $ denotes the subfunctor of $ F ^{\bold {n.o.} }_{ \bar \rho ,  \kappa } $ consisting of the deformations whose determinant is $\chi$. 
\end{enumerate}
\end{dfn}

\bigskip
Let $F$ be a local field of mixed characteristic $p$. We assume that $F$ is absolutely unramified, and $p\geq 3$.  We recall basic results on finite flat commutative group schemes over $ o_F$. \par
For a noetherian scheme $S$, let $ \ffgpsch _S$ be the category of commutative finite flat group schemes over $S$. \par

By Raynaud's result \cite{Ray}, $ \ffgpsch _{o_F}$ is an abelian category, and the restriction functor 
$$
\Res^{ \ffgpsch} \vert _F:  \ffgpsch _{o_F} \longrightarrow \ffgpsch _F   \label{eq-galdef64}
$$
is fully-faithful. Since the charateristic of $F $ is zero, all finite flat group schemes over $F$ are \'etale, and hence are regarded as a commutative finite group with a $G_F$-action. 
For $ G  \in \ffgpsch _{o_F} $, we define the associated $G_F$-module $DT(G)$ by 
$$
DT(G) = \Hom (G(\bar F ), \mathbb Q / \mathbb Z ) . 
$$
$DT$ is a cofunctor from  $\ffgpsch _{o_F} $ to the category of $G_F$-modules $\mathcal G_F $, which is fully faithful. The essential image of $DT$ is denoted by $\mathcal G_{[0, 1], F }  $. \par

For a commutative ring $A$, an $A$-action on a finite flat group scheme $G$ over $S$ is a ring homomorphism $A \to \End _{S}  G  $. By $\ffgpsch _{A , S}$, we denote the category of commutative finite flat group schemes with an $A$-action over $S$. \par
Let $\mathcal G ^A _{[0, 1], F }  $ be the subcategory of $\mathcal G _{[0, 1], F}  $ consisting of $G_F$-modules with an $A$-aciton. $ \mathcal G ^A _{[0, 1] , F}$ is equivalent to $\ffgpsch _{A , o_F}$.  \par

Assume that $A$ is a commutative ring of characteristic $p$. The category $M^{FL}_{[0, 1] , A}  $ is defined as follows.
$M^{FL}_{[0, 1] , A} $ consists of triples $D=(F^0 (D), F^1(D) , \varphi ) $ such that 
\begin{itemize}
\item $F^0 (D)$ is a finite $o _F \otimes _{\mathbb Z_p} A$-module. 
\item $ F^1(D)$ is an $o _F
\otimes _{\mathbb Z _p} A $-submodule of $F ^ 0 (D) $. 
\item $\varphi$ is a $\sigma$-linear $A$-isomorphism $F^1(D)  \oplus ( F^0 (D)/F^1(D) )
\simeq  F^0 (D)$. Here $\sigma : o_F \simeq o_F $ is the lift of the absolute frobenius of $ k_F$. 
\end{itemize}

By \cite{FL}, there is an equivalence of categories
$$
D: \mathcal G^A _{ [0, 1], F  }   \simeq  M^{FL}_{[0, 1] , A}    .
$$
A quasi-inverse functor of $D$ is denoted by $V$. 

\begin{dfn}\label{dfn-galdef64}
We assume that $F$ is absolutely unramified local field of the residual characteristic $p$, $ p\geq 3$.
Let $A$ be a complete local noetherian $o_\lambda$-algebra with a finite residue field of characteristic $p$. 
\begin{enumerate}
\item A continuous representation $\rho: G_F \to \GL_2 (A)$ is called flat if there exists a continuous character $ \mu: G_F \to A ^{\times}$ satisfying the following conditions. 
\begin{itemize}
\item For any integer $n \geq 1 $, $(\rho \otimes \mu ^{-1}) \mod m^n _A $ belongs to $\mathcal G_{F, [0,1]} $.
\item $\det (\rho \otimes \mu ^{-1})  \vert _{I_F} = \chi^{-1} _{\cycle} \vert _{I_F}   $. 
\end{itemize}
We call $\mu$ a twist character of $\rho$, and $\mu \vert _{I_F}$ a twist type of $\rho$. 
\item Let $\bar \rho : G_F\to \GL_2 ( k_{o_\lambda } ) $ be a $G_F$-representation, $\kappa : (I_F) _{G_F} \to o^{\times}_\lambda $
a continuous character. 
$$
F^{ \bold{fl} } _{ \bar \rho , \kappa }: \mathcal C ^{\noeth} _{o_\lambda} \longrightarrow \Sets
$$
is the functor defined as follows: for each $A \in \ob \mathcal C ^{\noeth} _{o_\lambda}  $, $ F^{ \bold{fl} } _{ \bar \rho , \kappa } (A) $ is the set of  isomorphism classes of flat representations with the twist type $\kappa$, equipped with an isomorphism $\rho \mod m_A \simeq \bar \rho $. 
For a continuous character $\chi: G_F \to o^{\times}_{\lambda} $, $ F^{\bold {fl}}  _{
\bar \rho ,\kappa , \chi } $ is the subfunctor of $ F^{ \bold{fl} } _{ \bar \rho , \kappa }$ consisting of the deformations with the determinant $\chi$.
\end{enumerate}
\end{dfn}

\begin{rem}\label{rem-galdef61}
\begin{enumerate}
\item Assume that $F$ is absolutely unramified with the residual characteristic $p\geq 3 $. A nearly ordinary finite representation with a nearly ordinary type $\kappa$ is a flat representation with a twist type $\kappa$. 
\item For flat deformations, we allow the case when $\bar \rho$ is absolutely reducible ($\bar \rho
$ can be split). This becomes especially important when $F\neq
\Q_p$.  
\end{enumerate}
\end{rem}

\subsection {Tangent spaces at $\ell$}\label{subsec-galdef7}

We discuss the nearly ordinary case first. Since we allow ramifications of $F$, the tangent space calculation is discussed in detail.  \par
We fix characters $\kappa : (I^{\ab} _F) _{G_F} \to o^{\times}_{\lambda} $, and $\chi : G_F \to o^{\times}_{\lambda} $. 
Let $\bar \rho : G_F \to \GL _2 (k_{o_\lambda })$ be a nearly ordinary representation whose nearly ordinary type $\bar \kappa $ and $\det \bar \rho $ are equal to $\kappa \mod m _{o_ \lambda} $ and $\chi \mod m _{o_ \lambda}  $ respectively. 
Let $G = \GL_{2, k_\lambda}  $ as an algebraic group over $k_\lambda$, $Z $ the center of $G$, $B$ the standard Borel subgroup consisting of the upper triangular matrices in $G$, $U$ the uniportent radical of $B$, $T = B /U$. 
By $k_{o_\lambda} [\epsilon]  $ we mean the ring of dual numbers. \par
Let $\rho :G_F \to G (k_{o_\lambda} [\epsilon]  )  $ be a deformation of $\bar \rho$ whose nearly ordinary type and the determinant are fixed to $\kappa $ and $\chi$.  Then $\rho$ defines the Kodaira-Spencer class $c_G (\rho ) $ in $H ^1 ( F, \ad ^0 \bar \rho)$. 
 If we regard $G$-representation $\ad _{G/Z}$ as a $G_F$-representation by $ \rho$, it is identified with $\ad ^0 \bar \rho $ (projective deformations are the same as the determinant fixed deformations if $\ell \neq 2$).  

Since $\bar \rho$ is nearly ordinary, we may assume that $\bar \rho$ factors through $ B(k_{o_\lambda}) \hookrightarrow G (k_{o_\lambda})$ by taking a suitable basis. We consider the filtration as a $B$-representation $ 0 \subset \ad _{B/Z} \subset \ad _{G/Z}$, where we identify $\ad _{B/Z} \subset \ad _{G/Z} $ with $W_1 =  \bar \rho \otimes (\bar \chi _2) ^{-1}  \subset W = \ad ^0 \bar \rho $ as $G_F$-representations. \par

Since $ \rho : G_F \to G ( k_{o_\lambda} [\epsilon]  )$ is a nearly ordinary deformation of $\bar \rho$, by a suitable choice of a basis, $\rho $ takes the values in $B (k_{o_\lambda}[\epsilon] )  $, and hence the cohomology class $c_G(\rho)$ in $H ^1 ( F, \ad ^0 \bar \rho ) $ associated to $\rho $ is the image of the class $c_B (\rho  ) $ in $H ^1 (F, \ad _{B/Z} )  $, defined similarly as $c_G(\rho) $ as the Kodaira-Spencer class. 

$c_B (\rho) $ must satisfy one more condition since we fix the nearly ordinary type.   $\ad _{B/Z}$ admits a filtration $ \ad _U \subset \ad _{B/Z} $ such that we have the exact sequence
$$
0 \longrightarrow  \ad _U \longrightarrow \ad _{B/Z}\longrightarrow  \ad _{T/Z}\longrightarrow  0   
$$
of $B$-representations. It induces the exact sequence
$$
 0 \longrightarrow \bar \chi _1 /\bar \chi _2 \longrightarrow  W_1 \longrightarrow k_{o_\lambda} \to 0   
$$
as $G_F$-representations. 
Since we are considering deformations with a fixed nearly ordinary type, the deformations in the direction of $T/Z $ should be unramified, and hence the image of $c_B (\rho ) $ in $ H ^1 (F, k_{o_\lambda} ) $ belongs to $ H ^1 ( F ^{\unr}/ F, k_{o_\lambda}) $.

\begin{thm}\label{thm-galdef71} Let $\bar \rho$ be a $G_F$-distinguished nearly ordinary representation. 

\begin{enumerate}
\item The dimension of $ F ^{\bold {n.o.}}_{ \bar \rho  , \kappa , \chi } (k_{o_\lambda}
[\epsilon] ) $ is at most $\dim_{k_{o_\lambda} } H ^ 1 ( F,
\  \bar \chi _1 /\bar \chi _2) +   \dim _{k_{o_\lambda} }H ^0 (F, \ \ad ^0 \bar \rho ) -\dim _{k_{o_\lambda} }H ^0 (F,\  \bar \chi _1 /\bar \chi _2)
  $. (This upper bound is equal to $   \dim _{k_{o_\lambda} }H ^0 (F,\  \ad ^0 \bar \rho )+  [F: \mathbb Q_p ]+  \dim _{k_{o_\lambda} } H ^ 0 (F, \ ( \bar \chi _2 /\bar \chi _1) (1)) $ by the Euler characteristic formula of Tate). 
\item 
 If $\bar \chi _2 =\bar \chi _1 (-1)$, and $\bar \rho $ is not nearly ordinary finite, then the dimension of $ F ^{ \bold {n.o.} }_{ \bar \rho ,\kappa , \chi } (k_{o_\lambda}
[\epsilon] ) $
 is at most $\dim _{k_{o_\lambda}} H ^0  (F, \ad ^0 \bar \rho  ) + [F: \mathbb Q_p ] $.

\item If $\bar \chi _2 =\bar \chi _1 (-1)$, and $\bar \rho $ is nearly ordinary finite, then the dimension of $ F^{ \bold {n.o.f.}  } _{ \bar \rho ,\kappa , \chi }  (k_{o_\lambda}
[\epsilon] )$
 is at most $\dim _{k_{o_\lambda}} H ^0  (F, \ad ^0 \bar \rho  ) + [F: \mathbb Q_p ] $.

\end{enumerate}
\end{thm}
Before beginning the proof of \ref{thm-galdef71}, we prepare several lemmas. 
\begin{lem}\label{lem-galdef71} Let $\bar \rho$ be a $G_F$-distinguished nearly ordinary representation. 
\begin{enumerate}
\item  If $\bar \rho $ is semi-simple, then $\dim _{k_{o_\lambda} }H ^0 (F, \ad ^0 \bar \rho ) = \dim _{k_{o_\lambda}}H ^0 (F, W_1  )  = 1+ \dim _{k_{o_\lambda} }H^0 (F,   \bar \chi _1 /\bar \chi _2)= 1$, and $H ^1 ( F, \bar \chi _1 /\bar \chi _2) \to H ^1 ( F, W_1 )   $ is injective. 
\item If $\bar \rho $ is not semi-simple, then $\dim _{k_{o_\lambda} }H ^0 (F, \ad ^0 \bar \rho ) = \dim _{k_{o_\lambda} }H ^0 (F, W_1  )  =  \dim _{k_{o_\lambda} }H^0 (F,   \bar \chi _1 /\bar \chi _2)$, and the kernel $L_{\bar \rho}$ of $H ^1 ( F, \bar \chi _1 /\bar \chi _2) \to H ^1 ( F, W_1 )   $ is one dimensional. \end{enumerate}
\end{lem}
The verification of this lemma is easy, and left to the reader. \par
\begin{lem}\label{lem-galdef72}
Assume that $\ell \geq 3$, $\bar \chi _2 =\bar \chi _1 (-1)$, and $\bar \rho $ is not nearly ordinary finite. For any nearly ordinary deformation $\rho $ over $k _{o_\lambda} [ \epsilon] $ with a fixed nearly ordinary character and determinant of the form
$$
0  \longrightarrow \tilde \chi _1 \longrightarrow \rho \longrightarrow \tilde  \chi_2 \longrightarrow 0,
$$
$\tilde \chi _2 = \tilde \chi_1 (-1 )  $ holds. 
\end{lem}
\begin{proof}[Proof of Lemma \ref{lem-galdef72}]
  When we view $\bar \chi_1 $ and $\tilde \chi_2$ as characters of $G_F$ with the values in $ k_{o_\lambda} [\epsilon ]$ by $ k_{o_\lambda} \hookrightarrow k_{o_\lambda}[\epsilon ] $, we denote them by $\chi_1$ and $\chi_2 $. \par

$ \mu = \tilde \chi _1 / \chi _1  : G _F \to k_{o_\lambda}[ \epsilon ] ^\times$ is an unramified
character lifting the trivial character. $\tilde \chi_2 =  \chi_2  \cdot (\mu
)^{-1} 
$ since we fix the determinant. $ \rho \otimes (\tilde \chi_2 )  ^{-1}$ defines an element
$c (\rho)  \in \Ext ^ 1 ( F,\  \tilde \chi _ 1 / \tilde \chi _2 ^{-1})= H ^ 1  (F,\ \mu ^2 \chi_1 / \chi_2 )= H ^1 ( F, \mu ^2 (1) )
$ which lifts the class $c (\bar \rho) $ in $H ^ 1  (F,\ \bar \chi _1 /\bar  \chi _2 )=  H ^1 ( F, k_{o_\lambda} (1) )
$ defined by
$\bar \rho \otimes  \bar \chi _2 ^{-1}$ by the canonical homomorphism 
$$
H ^ 1  (F,\  \mu ^2 (1) ) \longrightarrow H ^ 1 ( F,\  k_{o_\lambda} (1) )
$$
deduced from
$$
 \mu ^2 (1) \longrightarrow k_{o_\lambda} (1)  \longrightarrow 0 .
$$
$c(\bar \rho ) \vert _{F^{\unr}} = c_{\bar \rho }$ in $H ^ 1 ( F^{\unr},\  k_{o_\lambda} (1) ) $.\par
Assume that $\mu  \neq 1 $.  
Consider the following homomorphism
$$ 
H ^1 ( F ^{\unr} , \mu ^ 2 (1)  ) \simeq (F ^{\unr}) ^{\times}  /((F ^{\unr}) ^{\times}  ) ^{\ell} \otimes _{\mathbb F_\ell } \mu ^2 \overset {\tilde \beta }  \to \mu ^ 2 , 
$$
where $\tilde \beta $ is deduced from the valuation $(F ^{\unr}) ^{\times} \to \mathbb Z$ (evaluation map, see Definition \ref{dfn-galdef615}).\par
$\tilde \beta (c(\rho ) \vert _{F ^{\unr} } ) $ is fixed by $\Gal ( F ^{\unr} / F)$. 
Since $\ell \neq 2$, $\mu ^2 \neq 1 $, and $ (\mu ^2 ) ^ {\Gal ( F ^{\unr} / F)} $ is the subspace $ k_{o_\lambda} \epsilon$ spanned by $\epsilon$. Then $\tilde \beta (c(\rho ) \vert _{F ^{\unr} } )$ is mapped to zero by $ \mu ^2 \to \mu ^2 \otimes _{k_{o_\lambda}[\epsilon] }k_{o_\lambda} = k_{o_\lambda}$, which is equal to $ \beta (c_{\bar \rho })$ where $\beta = \ev _{ k_{o_\lambda}} : H ^1 ( F ^{\unr} , k_{o_\lambda} (1)  ) \simeq (F ^{\unr}) ^{\times}  /((F ^{\unr}) ^{\times}  ) ^\ell \otimes _{\mathbb F_\ell } k_{o_\lambda} \overset { \ev_{ k_{o_\lambda}} }  \to k_{o_\lambda}$ is the evaluation map for $k_{o_\lambda}$ (Definition \ref{dfn-galdef615}).
This implies that $ \bar \rho$ is nearly ordinary finite.  
\end{proof}

\bigskip
\begin{proof}[Proof of Theorem \ref{thm-galdef71}]
For (1), assume that $ \bar \rho $ is semi-simple. Then $H ^1 ( F, \bar \chi _1 /\bar \chi _2) \to H ^1 ( F, W_1 )   $ is injective by Lemma \ref{lem-galdef71}, (1). So the possible maximal dimension is $\dim _{k_{o_\lambda}} H ^1 ( F, \bar \chi _1 /\bar \chi _2) +1 $. If $\bar \rho$ is not semi-simple, the kernel $L_{\bar \rho}$ of $H ^1 ( F, \bar \chi _1 /\bar \chi _2) \to H ^1 ( F, W_1 )   $ is one dimensional (this is the subspace spanned by the extension determined by $\bar \rho$). So the dimension is at most $\dim _{k_{o_\lambda}} H ^1 ( F, \bar \chi _1 /\bar \chi _2)  $.\par

For (2), by Lemma \ref {lem-galdef72}, there are no deformations in the direction of $ \ad _{T/Z}$, and the possible maximal dimension is $ \dim _{k_{o_\lambda}} H ^1 ( F, \bar \chi _1 /\bar \chi _2)/L_{\bar \rho} =   \dim _{k_{o_\lambda}} H ^1 ( F, \bar \chi _1 /\bar \chi _2)-1$, which is equal to $ [ F: \mathbb Q_p ] + \dim _{k_{o_\lambda}} H ^ 0 (F,  \bar \chi _1 /\bar \chi _2) =  [ F: \mathbb Q_p ] + \dim _{k_{o_\lambda} }H ^0 (F, \ad ^0 \bar \rho ) $ by Lemma \ref{lem-galdef71} and the Euler characteristic formula. \par
For (3), first assume that $\bar \rho$ is semi-simple.  For a class $c  $ which belongs to the tangent space, the set of the classes in the tangent space with the same projection to $H ^1 ( F, k_{o_\lambda} ) $ is a torsor under $ H ^1 _{\flatcoh}  (o_F , k_{o_\lambda} (1)  ) \simeq o_F  ^{\times} / (  o_F  ^{\times})^{\ell} \otimes _{\mathbb F_\ell}  k_{o_\lambda}$. So the dimension is estimated from above by $1 + [F: \mathbb Q_p ] + \dim _{\mathbb F_\ell } H ^0 ( F, \mu _\ell  )  $. Since we assume that $\bar \rho  $ is $G_F$-distinguished, $ H ^0 ( F, \mu _\ell  ) = \{ 1 \}$, and the result follows. \par
If $\bar \rho$ is not semi-simple, the kernel $L_{\bar \rho}$ of $H ^1 ( F, \bar \chi _1 /\bar \chi _2) \to H ^1 ( F, W_1 )  $ belongs to $H ^1 _{\flatcoh} (o_F,  \mu _\ell \otimes _{\mathbb F_\ell}  k_{o_\lambda} )  $ since $\bar \rho$ is nearly ordinary finite. 

So the dimension of the tangent space is estimated from above by $ \dim _{k_{o_\lambda}} H ^1 _{\flatcoh} (o_F,  \mu _\ell \otimes _{\mathbb F_\ell }{k_{o_\lambda} }) /L_{\bar \rho} + 1 = (\dim _{\mathbb F_\ell} H ^1 _{\flatcoh} (o_F,  \mu _\ell) -1 ) + 1  =  [F: \mathbb Q_p ] + \dim _{\mathbb F_\ell } H ^0 ( F, \mu _\ell  ) = [F: \mathbb Q_p ] + \dim _{k_{o_\lambda}}H ^0 (F, \ad ^0 \bar \rho  ) $ by Lemma \ref{lem-galdef71}, (2).
\end{proof}
\bigskip

\begin{rem}\label{rem-galdef71}
In the terminology of \cite{W1}, in the case of \ref{thm-galdef71} (2), Selmer
deformations are strict by a theorem of Diamond, and the dimension estimate follows from proposition 1.9, (iv).
\end{rem}
\bigskip 
Next we consider flat deformations. The case was treated by Ramakrishna \cite{Ram} when $ F =
\Q_p
$ and $\bar \rho $ is absolutely irreducible. Conrad discusses flat deformations in detail for $F=\Q_p$ (\cite{C}, theorem 5.1), though he excludes the case when $\bar \rho $ is split. 

Here we only discuss the tangent space of flat deformation functors by a use of the theory of Fontaine-Laffaille (cf. \cite{Ram},
\cite{TW}).  

\begin{thm} \label{thm-galdef72} Assume that $F$ is absolutely unramified, with the residual characteristic $p\geq 3$. 
 Let $\bar \rho$ be a flat representation with the twist character $\bar \kappa$, $\kappa:( I_F)_{G_F}   \to o^{\times}_{\lambda} $ (resp. $\chi : (I _F )_{G_F} \to o^{\times}_{\lambda}$) a continuous character lifting $\bar \kappa$ (resp. $\det \bar \rho\vert _{I_F}$). 

\begin{enumerate}
\item The dimension of $ F ^{\bold {fl}}_{ \bar \rho  , \kappa  } (k_{o_\lambda}
[\epsilon] ) $ is equal to $   \dim _{k_{o_\lambda} }H ^0 (F,\  \ad  \bar \rho )+  [F: \mathbb Q_p ]$. 
\item The dimension of $ F ^{\bold {fl}}_{ \bar \rho  , \kappa , \chi   } (k_{o_\lambda}
[\epsilon] ) $ is equal to $   \dim _{k_{o_\lambda} }H ^0 (F,\  \ad ^0 \bar \rho )+  [F: \mathbb Q_p ]$.
\end{enumerate}
\end{thm}

\begin{lem}\label{lem-galdef73} Let $\rho: G_F \to \GL_2( A)$ be a flat deformation of $\bar \rho: G_F \to \GL_2 (k) $ with the trivial twist type, where $A$ is an object of $ \mathcal C ^{\artin} _{o_\lambda }$ of characteristic $p$. Then 
$F^0( D(\rho  ) ) $ is free of rank $2$ over $ o_F \otimes _{\Z_p}A$, and $ F^1(D (\rho ) )$ is free of rank $1$ over $ o_F \otimes _{\Z_p} A$. 
\end{lem}
\begin{proof}[Proof of Lemma \ref{lem-galdef73}]
First we treat the case of $A= k$. 
Assume that  $\bar \rho$ is absolutely reducible. By extending $k$ if necessary, we may assume that $\bar \rho  $ is reducible of the form 
$$
0 \longrightarrow \bar \chi _1 \longrightarrow \bar  \rho \longrightarrow \bar  \chi_2 \longrightarrow 0 , 
$$
where $\bar \chi _1 $ and $\bar \chi_2$ are characters, and there is an embedding $k _F \hookrightarrow k$. For $i = 1, 2 $, 
$\bar \chi _{i}  \vert _{I_F}$ is expressed as 
$$
\bar \chi _{i } \vert _{I_F}= \prod _{\sigma \in \Gal ( F/ \mathbb Q_p )  }(\chi ^{\sigma  } _{F } ) ^{-\epsilon _{i, \sigma} } ,
$$
where $\chi _F: I_F \to (I^{\ab}_F ) _{G_F} \to k^\times _F  \hookrightarrow k ^{\times}$ is the fundamental character of $F$, and $\chi ^{\sigma}_F$ denotes the $\sigma $-twist of $\chi _F$. Since $\bar \chi _i $, $i = 1, 2 $, is associated with a finite flat group scheme, $ \epsilon _{i, \sigma } \in \{ 0, 1 \} $ for any $\sigma \in \Gal ( F/ \mathbb Q_p ) $. \par
By our assumption,  $\det \bar \rho \vert _{I_F} = \bar \chi ^{-1} _{\cycle} \vert _{I_F}$, and hence we have an equality
$$
\epsilon _{1, \sigma } + \epsilon _{2, \sigma} = 1.
\leqno{(\ast)}
$$
for $ \sigma \in \Gal ( F/ \mathbb Q_p )$. \par
 The exact sequence
$$
0 \longrightarrow D ( \bar  \chi _1   ) \longrightarrow D (\bar \rho ) \longrightarrow D (\bar \chi _2  ) \longrightarrow 0 
$$
is strictly compatible with the filtrations. $F^0 (D (\bar \chi_i ))$ for $i= 1, 2 $ are free of rank 1 over $k_F \otimes k$, and hence $F ^0(D(\bar \rho  ) )$ is free of rank 2. For $i = 1,2 $, $F ^1 ( D ( \bar \chi _i ) )$ is the part of $F^0 ( D ( \bar \chi _i ) ) $ where $ k^\times _F$ acts as $\prod _{\sigma \in \Gal ( F/ \mathbb Q_p )  }(\iota  _{\sigma } ) ^{\epsilon _{i, \sigma} } $. Here $\iota _{\sigma} :k^{\times }_F \overset {\sigma }\to k^{\times}_F \hookrightarrow k ^{\times}$. By $(\ast)$, $F^1 (D (\bar \rho ))$ is free of rank 1 over $ k_F \otimes k$. \par
In the absolutely irreducible case, let $\tilde F $ be the unramified extension of $F$ of degree 2.  To show \ref{lem-galdef73} for $\bar \rho$, it suffices to prove it for $\bar \rho \vert _{\tilde F} $. Since $\bar \rho  \vert _{\tilde F} $ is absolutely reducible, we are reduced to the absolutely reducible case, which is already treated. \par

In the general case, we proceed by an induction on $n=\length A$.  For $n=1$, the claim is already shown. 
For $n \geq 2$, take a quotient
 $ A ' = A / I $ of $A$ with $\length I = 1$ and $I ^2 = \{ 0\} $. The kernel of $ \rho
\twoheadrightarrow  \rho ' = 
\rho \otimes _A A' $ is isomorphic to $\bar \rho$. Since $D$ is an exact
functor,
$$
0 \longrightarrow D( \bar \rho  ) \longrightarrow D (\rho ) \longrightarrow D(\rho'    ) \longrightarrow 0
\leqno{(\dagger)}
$$
is exact. Moreover, $(\dagger)$ is strictly compatible with the filtrations. This gives us an equality of the length

$$
\length  F ^{j } ( D(  \rho  )  )
= \length F ^{j } ( D( \rho' ))  + \length F ^{j}( D (\bar \rho  ) )
$$
for $j = 0, 1  $.
By the assumption on the induction, $F ^{j } ( D ( \rho'  ))$ is free of rank $2-j$ over $ o_F \otimes _{\Z_p} A'$ for $j = 0, 1 $, so $F ^j (D(\rho) )$ is a quotient of $(o_F \otimes _{\Z_p} A) ^ {\oplus (2-j )}$ by Nakayama's Lemma. Since $\length  F ^j( D (\rho)  ) = \length (o_F \otimes _{\Z_p} A) ^ {\oplus (2-j )}$, we have the freeness of $ F ^j( D (\rho ) )$. 

\end{proof}

\begin{dfn}\label{dfn-galdef71} Let $A$ be an object of $ \mathcal C ^{\artin} _{o_\lambda }$ of characteristic $p$.
\begin{enumerate}
\item $\mathcal E _A $ is the full subcategory of $ M^{FL}_{[0, 1] , A} $ consisting of objects $D=(F^0(D), F^1(D) , \varphi) $ which satisfy the following conditions:
\begin{itemize}
\item $F^0 (D) $ is free of rank $2$ over $o_F \otimes _{\Z_p} A $.
\item $ F^1 (D )$ is free of rank $1$ over $ o_F \otimes _{\Z_p} A$. 
\end{itemize}
\item For $ L _A = (o_F \otimes _{\Z_p} A )^{\oplus  2}  $, we define a filtration by $F ^ 0 ( L_A ) = L_A$, $F^1 (L_A ) = $ the $A
\otimes _{\Z_p} o_F$-submodule of $  ( o_F \otimes _{\Z_p} A )^{\oplus  2}$ generated by $(1, 0)$. \par
For an object $D=(F^0(D), F^1(D) , \varphi) $ of $\mathcal E_A$, a frame $f $ of $D $ is an isomorphism $f: F^ 0 (D)  \simeq L_A$ as an $o_F \otimes _{\Z_p} A $-module which preserves the fltrations.  $(D , f)  $ is called a framed object of $ \mathcal E_A$. 
The set of all frames for $D$ is denoted by $ \frame (D)$.
\item A framed isomorphism $( D_1,\ f_1 ) \to (D_2,\ f_2) $ is an isomorphism $\alpha  : D_1 \to D_2 $ in $M^{FL}_{[0, 1] , A}  $ such that $f _2 \circ \alpha =  f_1 $. 
\end{enumerate}

\end{dfn}
\begin{rem}\label{rem-galdef72} By Definition \ref{dfn-galdef71} (3), a framed automorphism of a framed object of $ \mathcal E_A$ reduces to the identity, and hence any framed isomorphisms between two framed objects are unique. 
\end{rem}
For a flat deformation $\rho$ of $\bar \rho$ with the trivial twist type, by Lemma \ref{lem-galdef73}, a frame always exists for $D (\rho)  $. \par
Let $G$ be the Weil restriction $ \Res_{o_F /
\mathbb Z_p }\GL_{2, o_F}$ of $\GL_{2, o_F} $, $B$ the Weil restriction of the standard Borel subgroup of $\GL_{2, o_F}$ consisting of the upper triangular matrices. These groups are smooth over $\bold Z_p$ since $F$ is absolutely unramified. \par
For $D \in \ob \mathcal E_A$, $\frame (D)$ forms a trivial $B(A)$-torsor. 
By Definition \ref{dfn-galdef71} (3),  two frames $ f_1 $ and $f_2$ of $D$ are isomorphic if and only if there is an automorphism $\alpha : D\simeq  D
$ in
$M^{FL}_{[0, 1] , A} 
$ such that $f_2 = f_1 \circ \alpha $.  

\begin{lem}\label{lem-galdef74}
Let $\rho: G_F \to \GL_2( A)$ be a flat deformation of $\bar \rho: G_F \to \GL_2 (k) $ with the trivial twist type.
\begin{enumerate}
\item Assume that $ pA=0 $. Let $X(A)$ be the set of framed isomorphism classes of framed objects of $ \mathcal E_A $.
If $X(A)$ is non-empty, it has a structure of a trivial torsor under $ G (A)$.
 \item  We fix a frame of $ D (\bar \rho )$.  For $A= k [\epsilon]$, the set $X _{\bar \rho}( A) $ of the framed isomorphism classes of framed objects $D$ of $ \mathcal E_A $ with a framed isomorphism $D \otimes _A k \simeq  D (\bar \rho ) $ has a structure of a $\Lie G _k $-torsor. 
\item We fix a frame of $D (\bar \rho )$.  For $A= k  [\epsilon]$, the set of the isomorphism classes of the frames of $D ( \rho)  $ lifting the frame of $ D ( \bar \rho ) $ is a torsor under $ \Lie B _k / K$, where $K$ is a subgroup of $ \Lie B _k$ which is canonically isomorphic to $ H^0 ( F, \ad \bar \rho )$. 
\end{enumerate}
\end{lem}
\begin{proof}[Proof of Lemma \ref{lem-galdef74}] For (1), by the definition of $ M^{FL}_{[0, 1] , A} $, $X(A)$ is identified with 
$$
X(A) = \{ h : (o_F \otimes _{\Z_p} A )^ 2 \rightarrow  (o_F \otimes _{\Z_p} A )^ 2,\ h \text{ is a } \sigma \text{-linear } A\text{-isomorphism} .\}
$$
(1) is clear from this description with the $G(A)$-action given by $h \mapsto h \circ g^{-1} $ for $g \in G(A)$ .\par
 For (2), note that $X_{\bar \rho } = X_{\bar \rho} (A)$ is non-empty since it contains the constant deformation $\bar \rho \otimes _k A $. As in the case of (1), $X_{\bar \rho} $ has a structure of a trivial torsor under $\hat G(A)$, where $\hat G(A) $ is the subgroup of $G(A)$ consisting of the elements whose mod $m_A$-reduction is $1$. $\hat G(A) $ is canonically isomorphic to $\Lie G _k $. \par
For (3), by Lemma \ref{lem-galdef73}, $D_{\rho} = D ( \rho  )  $ admits a frame $f_0$. By changing it by an element of $B(A)$, we may assume that $f_0$ lifts the frame $\bar f_0  $ of $D ( \bar \rho ) $. Then the set $\frame (D_\rho ) _{\bar f_0} $ of all frames of $D_{\rho}  $ lifting $\bar f_0$ is a torsor under $\hat B (A) $, where $\hat B (A) = \hat G (A) \cap B (A) $. $\hat B (A) $ is canonically isomorphic to $\Lie B _k  $. To determine the isomorphism classes of the frames of $D_{\rho} $, it suffices to determine the automorphism group $K$ of $ D_{\rho} $ whose mod $m_A$-reduction is the identity of $D ( \bar \rho )  $, since the isomorphism classes of the frames in $\frame (D_\rho ) _{\bar f_0}$ is identified with $\frame (D_{\rho}  )_{\bar f_0 }  / K \simeq \Lie B_k / K $. By the categorical equivalence, $K$ is isomorphic to the group $K'$ of the automorphisms of $\rho$ which induce the identity on the mod $m_A$-reduction $\bar \rho$. It is well-known that $K'$ is canonically isomorphic to $H^0 (F,  \ad \bar \rho)$.
\end{proof}

\begin{proof}[Proof of Theorem \ref{thm-galdef72}] We may assume that the twist character $\kappa$ is trivial. $k = k_{o_\lambda}
$, $A= k [\epsilon] $.  We fix a frame of $D ( \bar \rho) $. There is a canonical map
$$
\pi : X _{\bar \rho } \longrightarrow F ^{\bold {fl}}_{ \bar \rho  , \kappa  } (k [\epsilon] ) 
$$
which associates the Galois module $V ( D) _F $ to a framed object $(D, f) \in  X _{\bar \rho }$. 
By Lemma \ref{lem-galdef74} (3), $\pi$ has a structure of $ \Lie B_k / K$-torsor. By Lemma \ref{lem-galdef74} (2), $X _{\bar \rho }$ is a $\Lie G_k $-torsor. Thus we have
$$\dim _k  F ^{\bold {fl}}_{ \bar \rho  , \kappa  } (k [\epsilon] ) = \dim _k \Lie G_k- \dim_k \Lie B_k  + \dim _k  K = [F: \mathbb Q _p ] + \dim _k H ^0 (F, \ad \bar \rho ).
$$
 For \ref{thm-galdef72}, (2), we just note that unramified twists $ \rho
\mapsto
\rho \otimes \mu$ give a one dimensional family of deformations, so one should diminish 1 in the case of $ F ^{\bold {fl}}_{ \bar \rho  , \kappa , \chi } (k [\epsilon] )$ because the determinant is fixed. 
\end{proof}
\subsection {$K$-types} \label{subsec-galdef8}

The discussion in this paragraph localizes a global
argument in
\cite{W2}, proposition 2.15.
$0_E$-case is analyzed as in \cite{DT2}, \cite{D1}, though we make a modification using a
result of
\cite{Ger} when the relative conductor is even. \par
\bigskip

First we recall a result of G\'eradin \cite{Ger}, which gives the local Langlands correspondence explicitly in some supercuspidal cases (cf. \cite{DT2}, \S2, and \cite{D1}, \S4 for other applications to Hecke algebras). 

 \par
Let $\tilde F$ be the unramified quadratic extension of $F$, $ \sigma$ the
non-trivial element in
$\Gal (\tilde F / F) $. Take a character $\psi: G_{\tilde F} \to E_\lambda ^{\times}$ which
does not extend to $G_F$. $\rho= \Ind _{G_{\tilde F }} ^{G_ F} \psi $. The conductor $c$ of $ \psi / \psi ^\sigma$ depends only on $\rho$, and it is called the relative conductor of $\rho$. Here $ \psi ^\sigma $ is the $\sigma
$-twist of $\psi$, and $c \geq 1$ holds.\par
For an integer $m$, let $(m)_2 \in \{0, 1  \} $ be the integer which satisfies $m \equiv  (m)_2 \mod 2 $. $(m)_2$ is called the parity of $m$.  
For the relative conductor $c$, define a quaternion algebra $D_c$ central over $F$ by 
$$
D_c = \tilde F + \tilde F \Pi ,\quad \Pi ^2 = (p_F )^ {(c) _2 } ,\  \Pi x  =\sigma ( x )\Pi  \quad (\forall x \in \tilde F) ,
$$
and $G_c = D ^{\times }_c $ is the multiplicative group of $D_c$. \par
$D_c$ depends only on the parity of $c$, with the invariant $ \frac {c} {2}\mod \mathbb Z $, and $o _{D_c} = o_{\tilde F} + o_{\tilde F}\Pi $ is a maximal order of $D_c$. 
We define an open subgroup $K_c$ of $G_c(F) $ and a character $\mu _{\psi } $ of $K_c$ by
$$
K _c = \tilde F ^\times ( 1 + m_{\tilde F} ^d \Pi ) , 
$$
$$
\mu_\psi  \vert _{\tilde F ^{\times}}
 = \psi \cdot \chi _{\tilde F} ,\ \mu_\psi \vert _{  1 + m_{\tilde F }^d \Pi  } 
= \chi \circ
\Norm _{D/F} .
$$
Here $d = \frac {c- (c) _2} {2}$,  $\chi _{\tilde F}$ is the unramified character corresponding to $ \tilde F/ F$, and $\Norm _{D/F} $ is the reduced norm (when $d= 0 $, we regard $m_{\tilde F }^0$ as $o_{\tilde F} $). $\chi$ is a character of $G_F$ such that 
$\psi \cdot{ \chi  \vert _{G_{\tilde F} }}^{-1} $ and $ \psi / \psi ^{\sigma}$ has the same
conductor (see \cite{Ger}, 3.2 for the existence of $\chi$).  $\psi $ and $\chi$ are regarded as characters of $\tilde F ^\times $ and $F^\times$ respectively by the local class field theory. $\mu _{\psi} $ does not depend on a choice of $ \chi$. \par

The induced representation $\Ind _{K _{c} }^{G_c  (F) } \mu_{\psi} $ is irreducible and supercuspidal (\cite{Ger}, 3.4 and 5.1). 
Moreover, the representation $\pi $ which corresponds to $\rho$ by the local
Langlands correspondence (resp. the local Langlands correspondence composed with the Jacquet-Langlands correspondence) is $ \Ind _{K _{c} }^{G_c  (F) } \mu_{\psi}$ when $c$ is even (resp. odd) by \cite{Ger}, 3.8 and 5.4.\par
\bigskip
Assume that $\ell \neq p $, and we fix an $\ell$-adic field $E_\lambda$. 
For $\bar \rho : G_F \to \GL_2(k _\lambda)$, and $ * = \bold f , \bold u$, we attach an inner twist
$G_{\bar \rho}$ of $\GL_{2, F}$, compact open subgroups $K _{*} (\bar \rho) $ in $G_{\bar \rho}(F)$, and characters
$\nu_{*}( \bar \rho)
$ of $K_{* } (\bar \rho)$ having values in $o^{\times}_{E_\lambda }$. \par
\bigskip

In the case of $0_E$, $\bar \rho$ is expressed as $ \Ind _{G_{\tilde F}} ^ {G_F}\bar \psi  $. By enlarging $k_\lambda$ if neccesary, we may assume that $\bar \psi$ is defined over $k_\lambda$. Let $\psi=\bar \psi _{\lift}  : G_{\tilde F} \to o_{E_\lambda} ^\times$ be the Teichm\"uller lift of $\bar \psi$. We also view
$\psi $ as a character of
$\tilde F ^\times $ by the local class field theory. \par

Let $c$ be the relative conductor of $\bar \rho$. Then $ G _{\bar \rho } = G_c$, 
$$
K _{\bold f } (\bar \rho)  
= o ^{\times}_{\tilde F}   \cdot ( 1 + m_{\tilde F} ^ d 
 \Pi  ),
$$
$$
K_{\bold u } (\bar \rho)  
= Z ^{\ell}_{\tilde F}  \cdot ( 1+  m_{\tilde F} ^ d 
 \Pi) .
$$
Here $d= \frac {c-(c) _2 } {2}$, and $Z ^{\ell}_{\tilde F}  $ is the minimal subgroup of $o^{\times}_{\tilde F} $ containing $1+ m_{\tilde F} $ of $\ell$-power index. 
$\nu_{*}( \bar \rho): K_{* } (\bar \rho) \to o_{E_\lambda} ^{\times}  $ is the restriction of $\mu _{\psi } $ to $ K_{* } (\bar \rho) $ for $ *= \bold f, \bold u $. 
\par
\bigskip 
Next we treat the cases other than $0_{E}$. When $\bar \rho $ is absolutely reducible, by enlarging $k_\lambda $ if necessary, we may assume that $\bar \rho$ is reducible over $k_\lambda$, and $\bar \kappa$ be the twist type of $\bar \rho$.  In the case of $0_{NE}$, we define the twist type $\bar \kappa$ of $\bar \rho$ as the trivial character. 
 
We define the integers $c(\bar \rho )  $ and $d(\bar \rho )$ by
$$
c (\bar \rho ) = \Art  (\bar \rho\vert _{I_F}  \otimes \bar \kappa^{-1} ) ,
$$
$$
d (\bar \rho ) = \dim _{k_\lambda } (\bar \rho\vert _{I_F}  \otimes \bar \kappa^{-1}) ^{I_F}.
$$
Here $  \Art $ denotes the Artin conductor. 
$ G_{\bar \rho}$ is $\GL_{2, F} $, and $K _{\bold f}(\bar \rho)$ and $K _{\bold u}(\bar \rho) $ are defined by
$$
K _{\bold f}(\bar \rho)= \mu_{\ell ^{\infty}} ( F) \cdot K _1 ( m^{c(\bar \rho )} _F ), 
$$
$$
K _{\bold u}(\bar \rho)= K _1 (m^{c(\bar \rho )+ d (\bar \rho ) } _F ) . 
$$
For $ * =\bold f,  \bold u $, the $K$-character $\nu_{*}( \bar \rho): K_* ( \bar \rho) \to o_{E_\lambda} ^\times$ of $\bar \rho$ is the composition of 
$$
K_* ( \bar \rho)  \longrightarrow \GL _2 (o_F) \overset {\det}
 \longrightarrow  o_F^\times \overset {\kappa} \longrightarrow o_{E_\lambda} ^\times .
$$
Here $ \kappa $ is the Teichm\"uller lift of the twist type $\bar \kappa $ of $\bar \rho$, which is regarded as a character of $o_F^\times$ by the local class
field theory. \par

\begin{dfn} \label{dfn-galdef81} Assume that $\ell \neq p$, and let $\bar \rho : G_F \to \GL_2 (k_\lambda) $ be a $G_F$-representation.  
The data $( G_{\bar \rho} , (K _* (\bar \rho ) ,\nu_{*}( \bar \rho))_{ * \in \{ \bold f, \bold u \} } )$ defined for $\bar \rho$ is called the type of $\bar \rho$.  $ \nu_{*} (\bar \rho)$ is called the $K$-character of $\bar \rho$.
\end{dfn}

Take a deformation $ \rho $ of $\bar \rho$ over $o_{E_\lambda} $. 
 Let $ \rho^{WD}_{E_\lambda}  : W ' _F  \to \GL_2 (E_\lambda) $ be the representation of the Weil-Deligne group $W'_F$ of $F$ associated to $\rho _{E_\lambda} $, and $\pi $ the irreducible admissible representation of $\GL _2 (F) $ defined over $ \bar E_\lambda$ associated to the $F$-semi-simplification of $ \rho^{WD}_{E_\lambda}  $ by the local Langlands correspondence.
If $  G_{\bar \rho } $ is isomorphic to $\GL_{2, F} $ (resp. the multiplicative group of a division quaternion algebra), we denote $\pi$ (resp. the Jacquet-Langlands correspondent of $\pi$) by $\pi_{\rho } $. 

\begin{dfn}\label{dfn-galdef82}
Assume that $\ell\neq p$, and let $\bar \rho : G_F \to \GL_2 (k_\lambda) $ be a $G_F$-representation of type $( G_{\bar \rho} , (K _* (\bar \rho ) ,\nu_{*}( \bar \rho))_{ * \in \{ \bold f, \bold u \} } )$, and $\rho$ a deformation of $\bar \rho $ over $o_{E_\lambda} $.  
\begin{enumerate}
\item An irreducible admissible representation $\pi $ of $G_{\bar \rho } (F) $ defined over $\bar E _\lambda $ is associated to $\rho$ if it is isomorphic to $\pi_{\rho} $. 
\item 
For an admissible irreducible representation $\pi : G _{\bar \rho} (F) \to \Aut _{\bar E_\lambda }V $ associated to $\rho$, and for $ * = \bold f , \bold u $, 
$$
 I _* (\bar \rho, \pi )  = \Hom _{K _* (\bar \rho )  }( \nu_{*}( \bar \rho) , \ V ) .
$$  
\end{enumerate}
\end{dfn}

\begin{prop} \label{prop-galdef81}
Assume that $ \ell\neq p  $, and let $\bar \rho : G_F \to \GL_2 (k_\lambda) $ be a $G_F$-representation of type $( G_{\bar \rho} , (K _* (\bar \rho ) ,\nu_{*}( \bar \rho))_{ * \in \{ \bold f, \bold u \} } )$, $\rho$ a deformation of $\bar \rho $ over $o_{E_\lambda} $, and $\pi $ an admissible irreducible representation associated to $\rho$. \par
Then $  I _{\bold f}  (\bar \rho , \pi )$ is at most one dimensional.  It is non-zero if and only if $\rho$ is a finite
deformation of
$\bar
\rho$. 
\end{prop}
\begin{proof}[Proof of Proposition \ref{prop-galdef81}] Note that $ K_{\bold f} (\bar \rho)$ contains $\mu _{ \ell ^{\infty} }(F) $, so $\det \rho _{\pi} $ is a finite deformation of $\det \bar \rho$. \par
First we treat the case of $0_E$. Assume that $  I _{\bold f}  (\bar \rho , \pi )\neq \{ 0 \}$. $\bar \rho =\Ind_{G_{\tilde F} } ^{G_F }\bar \psi $,
$\psi $ the Teichm\"uller lift of $\bar \psi$. The relative conductor of $\psi$ is $ c$. Let $V$ be the representation space of $\pi$. $  I _{\bold f}  (\bar \rho , \pi)$ is identified with the subspace of $V^{\ker \nu _{\bold f}( \bar \rho  ) } $ where $K_{\bold f } (\bar \rho )  $ acts as $ \nu _{\bold f} (\bar \rho  ) $. By twisting by an unramified character, we may assume that the central element $ p_F$ acts on $V$ by $\mu _{\psi} ( p_F ) $. So we may assume that the action of $ p^{\mathbb Z} _F \cdot K_{\bold f } (\bar \rho ) = K_c  $ on $ I _{\bold f}  (\bar \rho , \pi) \subset V^{\ker \nu _{\bold f}( \bar \rho  ) } $ is $\mu _{\psi}$. 

 Since $K_{\bold f } (\bar \rho )$ is a compact open subgroup and $ \pi$ is admissible, the $K_{\bold f } (\bar \rho )$-action on $V$ is semi-simple, and hence $ I _{\bold f} (\bar \rho , \pi  ) $ appears as a quotient representation of 
$\pi \vert _{K_{\bold f } (\bar \rho )} $. By the Frobenius reciprocity, 
$$
\dim _{ \bar E_\lambda } \Hom _{K_c  }(\pi \vert _{K_c}, \mu _{\psi}  ) = \dim _{ \bar E_\lambda } I _{\bold f}  (\bar \rho , \pi) ^{\vee}\cdot   \dim _{ \bar E_\lambda } \Hom _{G_{\bar \rho} (F)  }(\pi,  \Ind _{K_c} ^{G_{\bar \rho}(F) } \mu _{\psi}  ) .
$$ Since 
$\pi$ and $\Ind _{K_c} ^{G_{\bar \rho} (F) } \mu _{\psi} $ are both irreducible, they are isomorphic, and $\dim _{ \bar E_\lambda } I _{\bold f} (\bar \rho , \pi)= 1 $ by Schur's lemma. $\rho_{E\lambda } $ is isomorphic to $\Ind_{G_{\tilde F} } ^{G_F }\psi  $, and hence is a finite deformation of $\bar \rho$. The argument also shows that $\dim _{ \bar E_\lambda } I _{\bold f} (\bar \rho , \pi)= 1  $ in the case of finite deformations. \par

In the other cases, let $\nu$ be a character of $ F ^\times $ such that $ \nu \vert _{I_F} $ is the Teichm\"uller lift $ \bar \kappa _{\lift} $ of twist character $\kappa$. By replacing $\pi $ by $\pi
\otimes
\nu^{-1}$, we may assume that $ \nu_{\bold f}  ( \bar \rho )$ is trivial. Let $c (\pi)$ be the conductor of $\pi$.  Assume that $  I _{\bold f}  (\bar \rho , \pi )\neq \{ 0 \}$. Since $K_{\bold f }(\bar \rho ) $ contains $ K_1 ( m^{c(\bar \rho)}  _F )$, 
$$
c (\bar \rho ) \geq c (\pi )  
$$
since $\pi $ has a non-zero $K_1 (m ^ {c (\bar \rho )} )$-fixed vector. $c(\pi ) $ is equal to $ \Art \rho_{E_\lambda} $ since $\pi$ is associated to $\rho$ (the local Langlands correspondence preserves the Galois and automorphic conductors), and hence 
$$
 \Art \bar \rho   \geq \Art \rho _{E_\lambda }.
$$  
On the other hand, for the Artin conductors
$$
\Art \bar \rho   - \Art \rho _{E_\lambda }   =  \dim _{E_\lambda }(\rho _{E_{\lambda} } \vert _{I_F} )^{I_F}  -  \dim _{k_\lambda} (\bar \rho \vert _{I_F} )^{I_F}
\leq 0 $$
holds. Thus we have $ \Art \rho _{E_\lambda }= \Art \bar \rho $, and hence the equality $ \dim _{k_\lambda} ( \rho _{E_\lambda } \vert _{I_F} )^{I_F}
= \dim _{E_\lambda }(\bar \rho  \vert _{I_F} )^{I_F}$.  $\rho ^{I_F} $ is $o_{E_\lambda}$-free and is an $ o_{E_\lambda}$-direct summand, and hence the last equality implies that $ \rho$ is a finite deformation of $\bar \rho$. \par
If $\rho$ is a finite deformation, by an argument as above $c(\pi ) = \Art \bar \rho $ holds. $\pi ^{K_{\bold f}  (\bar \rho) }=  \pi ^{K_1 ( m^{c(\pi )}  _F) }$ is the space of new vectors, and hence $I _{\bold f} (\bar \rho , \pi)$ is one dimensional. \end{proof}
\bigskip

Let $\rho $ be a deformation of $\bar \rho$ over $o_{E_\lambda}$, and $\pi$ the admissible representation associated to $\rho$. Assume that $\bar \rho$ is not of type $0_E$. In the unrestricted case, we need to choose a one dimensional subspace from $  I _{\bold u} (\bar \rho , \pi ) $. We use the $U (p_F)$-operator for this purpose (\cite{W2}) .\par
In general, for a reductive group $G$ and a compact open subgroup $K$ of $G(F)$, let $ H ( G (F), K ) _{\bar E_\lambda } $ be the convolution algebra formed by the $\bar E _\lambda $-valued compactly supported $K$-biinvariant functions on $G(F)$. For any irreducible admissible representation $\pi$ of $G(F) $ defined over $\bar E_\lambda $, $H ( G (F), K ) _{\bar E_\lambda }  $ acts on $ \pi ^K$. 

\begin{dfn}\label{dfn-galdef83} Let $\pi$ be an irreducible admissible representation of $\GL _2 (F)$ defined over $\bar E_\lambda $. For a uniformizer $p_F$ of $F$,  
$$
U(p_F) ,\ U (p_F, p_F) :\pi ^K \longrightarrow  \pi ^K 
$$
are defined by the characteristic functions of the double cosets
$K
\begin{pmatrix} 1 & 0
\\ 0 &  p _F
\end{pmatrix} 
K  $ and $K
\begin{pmatrix} p_F  & 0
\\ 0 &  p _F
\end{pmatrix} 
K  $, repectively.
\end{dfn}
Note that
$U(p_F)$ and $U(p_F, p_F)$-operators thus defined may depend on a choice of a uniformizer.
\par
In our situation, we choose $ K = \ker \nu _{\bold u } (\bar \rho) $. $ U(p_F)$ acts on  $  I _{\bold u} (\bar \rho, \pi) $ by identifying $ I _{\bold u} (\bar \rho, \pi) $ as a subspace of $ \pi ^K$.  
\begin{prop}\label{prop-galdef82} Assume that $ \ell\neq p  $, and let $\bar \rho : G_F \to \GL_2 (k_\lambda) $ be a $G_F$-representation of type $( G_{\bar \rho} , (K _* (\bar \rho ) ,\nu_{*}( \bar \rho))_{ * \in \{ \bold f, \bold u \} } )$, $\rho$ a deformation of $\bar \rho $ over $o_{E_\lambda} $, and $\pi $ an admissible irreducible representation associated to $\rho$. 
$\kappa= \bar \kappa _{\lift}$ is the Teichm\"uller lift of the twist character of $\bar \rho $.\par
\begin{enumerate}

\item  If $\bar \rho$ is of type $0_E$, $ I _{\bold u}  (\bar \rho , \pi)$ is one dimensional. 

\item If $\bar \rho$ is not of type $0_E$, then $  I _{\bold u}  (\bar \rho , \pi) $ is isomorphic to
$$
\mathcal L _{\bar  \rho, p_F } =\bar E_\lambda [ U] /(U \cdot L (U,\
\pi \otimes ( \nu _{\bar \rho,  p_F } \circ \det ) ^{-1} ) )
$$ 
as an
$\bar E_\lambda[ U  ]$-module. Here $U $ acts as $U( p_F) $ on $ I _{\bold u}  (\bar \rho , \pi)$,  $\nu _{\bar \rho , p_F }$ is the character of $F^{\times} $ such that $\nu _{\bar \rho , p_F} ( p_F ) = 1$ and $ \nu _{ \bar \rho , p_F }\vert _{o^{\times} _F }= \kappa$, and $L(T,\   \pi \otimes   ( \nu _{\bar \rho,  p_F } \circ \det ) ^{-1} )$ is the standard
$L$-function of
$ \pi \otimes   ( \nu _{\bar \rho,  p_F } \circ \det ) ^{-1} $. The localization of $\mathcal L _{\bar  \rho, p_F } $ at $ U= 0$ is one dimensional.
\end{enumerate}
\end{prop}
\begin{proof}[Proof of Proposition \ref{prop-galdef82}]
For (1),  first show that $ I _{\bold u}  (\bar \rho , \pi)\neq \{ 0\} $. Since $\rho$ is a deformation of $\bar \rho$, $\rho _{\bar E_\lambda } $ is written as $\Ind^{G_F}  _{G_{\tilde F} } \phi $, where $\phi$ is a character of $ G_{\tilde F}$. We may assume 
that $\phi$ is a lift of $\bar \psi$. The relative conductor $c$ of $\phi$ is the same as that of $\bar \psi$. 
By the result of G\'eradin, $\pi$ is isomorphic to $\Ind^{G_{\bar \rho}(F)}  _{K_c} \mu _{\phi} $, and $\pi \vert _{K_c } $ contains a non-zero subspace $W$ where $K_c$ acts by $\mu _{\phi}$. It is clear that $ W \subset I _{\bold u } (\bar \rho , \pi ) $. \par
$I _{\bold u } (\bar \rho , \pi  )$ is identified with the subspace of $ V^{\ker \nu _{\bold u } ( \bar \rho ) } $ where $K_{\bold u }(\bar \rho ) $ acts as $  \nu _{\bold u } (\bar \rho) $. We consider the action of $\mu_{\ell ^{\infty} }(\tilde F ) $ on $I _{\bold u } (\bar \rho ) $, and let $\alpha $ be a character of $\mu_{\ell ^{\infty} }(\tilde F ) $ which appears as a subrepresentation of $I _{\bold u } (\bar \rho , \pi  ) $. We identify $\alpha$ with a character of $k^{\times} _{\tilde F } $ of $\ell$-power order. 
Let $\pi _{\psi} $ be the $G(F)$-representation corresponding to $ \Ind_{G_{\tilde F} } ^{G_F }\psi$, $ \beta $ be the
character of $o_D ^\times $ defined by $  o_D ^\times\to
 o_D ^\times/( 1 + \Pi o_D )= o^{\times}_{\tilde F}/ (1+ m_{\tilde F} )= k^{\times} _{\tilde F}  \overset {\alpha} \to o_{E_\lambda} ^\times
$. By twisting by an unramified character, we may assume that the central character of $\pi$ takes the same value as the central character of $\pi _{\psi} $ at $p_F  $. The action of $\tilde F ^{\times } $ is thus determined, and it follows that there is a subspace of $I _{\bold u } (\bar \rho , \pi ) $ where $K_c  $ acts by $ \mu _{\psi } \cdot \beta \vert _{K_c}  $. \par

By the definition, $\mu _{\psi \cdot \alpha } = \mu _{\psi } \cdot \beta \vert _{K_c} $, and hence $\pi$ corresponds to $\Ind_{G_{\tilde F} } ^{G_F }(\psi \cdot \alpha)$ by the local Langlands correspondence. This shows the uniqueness of $\alpha$, 
and hence $ K_c $ acts on $I _{\bold u} (\bar \rho , \pi) $ by $\mu _{\psi \cdot \alpha }  $. The calculation of the dimension is the same as in the case of $I _{\bold f} (\bar \rho , \pi) $ by using the Frobenius reciprocity.\par

For (2), we may assume that $\bar \kappa $ is trivial by twisting by a character of order prime to $\ell$.
$$
\Art \rho _{E_\lambda } + \dim _{E_\lambda} \rho _{E_\lambda } ^{I_F}= 2+ \operatorname{sw}
\rho_{E_\lambda }   =2+
\operatorname{sw} \bar\rho  =  \Art \bar \rho +\dim _{k_\lambda} \bar
\rho ^{I_F} = c (\bar \rho ) + d (\bar \rho )
$$ 
by the formula for Artin conductors. Here
$
\operatorname{sw}
$ means the swan conductor, which remains unchanged by a mod $\lambda$-reduction. $\dim_{E_\lambda} \rho _{E_\lambda  } ^{I_F}= \deg L (T, \pi ) 
$ and $\Art \rho _{E_\lambda} = \cond \pi$,
 since the local Langlands correspondence preserves the $L$ and $\epsilon$-factors. Thus \ref{prop-galdef82}, (2)
is reduced to the following lemma. 
\begin{lem}\label{lem-galdef81}
For an admissible irreducible representation $\pi$ of $\GL _2 (F) $ defined over $\bar E_{\lambda} $,  let $ c$ be the conductor of $ \pi $, $d$ the degree of $ L (T, \ \pi)$, where $ L (T, \ \pi)$ is the standard $L$-function of $\pi$. 
We regard $\pi ^ {K_1 (m
^{ c+ d} )  }$ as an $ \bar E_\lambda [ U] $-module, where the $U$-action is given by
the $U(p_F)$-operator. Then it is isomorphic to $\bar E_\lambda [ U] /(U \cdot L (U,\
\pi ) )
$. 

\end{lem} 
\begin{proof}[Proof of Lemma \ref{lem-galdef81}]
When $\pi$ is supercuspidal, this is well-known. 
Let $ v$ be a non-zero vector in $\pi  ^ {K_1(m ^c )} $ (new vector). Then $\pi 
^ {K_1(m ^{c + d})} $ has a basis by
$
\begin{pmatrix}  
p_F & 0 \\ 
0 & 1
\end{pmatrix}
 ^ i v ,\ 0 \leq i \leq d
$. By writing down the $U(p_F)$-action explicitly, the lemma follows. We omit the details
(for $F=\Q$ and in the global setting, this is found in \cite{W2}). 
\end{proof}
\end{proof}
 
\medskip

\subsection{Global deformations}\label{subsec-galdef9}
 In this subsection, $F$ is a global field. 
For a continuous representation $ \rho  $ of the global Galois group $G_F= \Gal (\bar F /F)$ and a finite place $v$ of $F$, we denote the
restriction to the decomposition group $\rho \vert _{G_{F_v}} $ by $\rho\vert_ {F_v}$. For a finite set of finite places $\Sigma$, $G_{\Sigma } = \pi^{\et} _1 (\Spec o_F \setminus \Sigma) $ be the Galois group of the maximal Galois extension of $F$ which is unramified outside $\Sigma $.\par
Let $k$ be a finite field, $ \bar \rho : G _F \to
\GL _ 2 ( k ) $ be an absolutely irreducible representation of $G_F$. Assume the following conditions:
\begin{itemize}
\item[{\bf DC1}] For $v \vert \ell$, $\bar
\rho\vert _{F_v} 
$ is flat or nearly ordinary. If $\bar
\rho\vert _{F_v} 
$ is flat (resp. nearly ordinary) at $v$, a flat twist type (resp. nearly ordinary type) $ \bar \kappa _v $ is specified, and is defined over $k$.
\item[{\bf DC2}] For $v \nmid \ell$ and $\bar \rho \vert _{F_v}$ is absolutely reducible, it is reducible over $k$, and a twist type $\bar \kappa_v $ (Definition \ref{dfn-galdef21}) is specified, and is defined over $k$.
\item[{\bf DC3}] For $v \nmid \ell$ and $\bar \rho \vert _{F_v}$ is of type $0_E$, an inertia character $\bar \psi $ (Definition \ref{dfn-galdef21}) is specified, and is defined over $k$.
\end{itemize}

\begin{dfn}\label{dfn-galdef91}
Let $ \bar \rho : G _F \to
\GL _ 2 ( k) $ be an irreducible mod $ \ell$-representation which satisfies {\bf DC1}-{\bf 3}.
\begin{enumerate}
\item A deformation function $ d $ of $\bar \rho$ is a map $ \vert
F \vert _f  \to  \{ {\bf f} ,{\bf u }, {\bf n.o.f.}, {\bf n.o.}, {\bf fl}\} $ which satisfies the following properties:

\begin{itemize}
\item If $  v \nmid \ell $, then  $
   d  (v) \in \{
\text{\bf f},
\text{\bf u}
\}
$.
\item If $v \vert \ell $ and $\bar \rho \vert _{F_v}$ is nearly ordinary, 
$  d ( v) \in 
\{\text{\bf n.o.f.} ,   \text{\bf n.o.} \}$.  
\item If $v \vert \ell $ and $\bar \rho \vert _{F_v}$ is flat, 
$  d ( v)  = \bf fl $.
\item  For almost all places $v$ of $F$, $ d (v) = \bold f $.
\end{itemize}

\item A deformation type $ \mathcal D $ of $\bar \rho$ is a quartet
$(d 
 ,\  o_{\lambda} , \{ \kappa _v\}_{v \in \bold {NO} (\bar \rho)} , \{ \kappa _v\}_{v \in \bold {FL} (\bar \rho)}  ) 
$ such that 
\begin{itemize}
\item $d$ is a deformation function of $\bar \rho$.  $d$ is called the deformation function of $\mathcal D$, and denoted by $\deform _{\mathcal D} $.
\item $
o_{\lambda}
 $ is a complete noetherian local ring with the maximal ideal $m _{o_\lambda} $ and the residue field $k_{o_\lambda} $. $k_{o_\lambda}$ is isomorphic to $k$.  $o_{\lambda} $ is called the coefficient ring of $\mathcal D$, and denoted by $o_{\mathcal D} $. 
\item $\bold {NO} (\bar \rho)$ (resp. $ \bold {FL} (\bar \rho)$) is the set of finite places where $\rho \vert _{F_v} $ is nearly ordinary with the nearly ordinary type $\bar \kappa _v$ (resp. flat with the twist character $\bar \kappa_v$). $\kappa _v : G_{F_v} \to o^{\times} _{\lambda}$ for $v \in \bold {NO} (\bar \rho)$ (resp. $ v\in  \bold {FL} (\bar \rho)$) is a continuous character such that $\kappa _v \mod m_{o_\lambda } = \bar \kappa_v $. For $v \in \bold {NO} (\bar \rho)$ (resp. $ v\in  \bold {FL} (\bar \rho)$),  $\kappa _v$ is called the nearly ordinary type of $\mathcal D$ (resp. flat twist type of $\mathcal D $), and is denoted by $\kappa _{\mathcal D, v} $.
\end{itemize}

\item The ramification set $\Sigma _{\mathcal D} $ of a deformation type $\mathcal D $ of $\bar \rho$ is 
$$
\Sigma _{\mathcal D}=\{ v \vert \ell \} \cup  \{ v :\bar \rho
\text{ is ramified at } v \} \cup  \{ v :\bar \rho
\text{ is unramified at } v,\  \deform _{\mathcal D} (v) = \bold u 
\}.$$

\end{enumerate}
\end{dfn}

\begin{dfn}\label{dfn-galdef92}
\begin{enumerate}
\item For a deformation type $\mathcal D = (d 
 ,\  o_{\lambda} , \{ \kappa _v\}_{v \in \bold {NO} (\bar \rho)} , \{ \kappa _v\}_{v \in \bold {FL} (\bar \rho)}) 
$ of $\bar \rho$ and a local homomorphism $o_{\lambda}\to o'_{\lambda '} $ of complete noetherian local rings, the scalar extension $ \mathcal D_{o' _{\lambda ' }}$ is defined by $ (d 
 ,\  o'_{\lambda'} , \{ \kappa ' _v\}_{v \in \bold {NO} (\bar \rho)},   \{ \kappa ' _v\}_{v \in \bold {FL} (\bar \rho)}) $. Here $  \kappa '_v $ is the composition of $G_{F_v} \overset {\kappa _v }\to o^{\times } _{\lambda}  \to ( o '_{\lambda '} )^{\times}$ for $v \in  \bold {NO} (\bar \rho)$ (resp. $v \in \bold {FL} (\bar \rho )$). 
\item We define a partial order $\leq $ on the set $ \{ {\bf f} ,{\bf u }, {\bf n.o.f.}, {\bf n.o.}, {\bf fl}\}  $ by $ \bold f
\leq
\bold u $, $\bold {n.o.f} \leq \bold {n.o.} $, $\bold {fl}\leq \bold {fl} $. 
A deformation type $\mathcal D$ of $\bar \rho$ is minimal if $ \deform _{\mathcal D}  ( v ) $ takes the minimal possible value at any $ v \in \vert F \vert _f $ for the partial order.
\item A morphism $\mathcal D \to \mathcal D '$ between deformation types of $\bar \rho$ is a local ring homomorphism $ f: o_{\mathcal D ' } \to o _{\mathcal D } $ and a condition on deformation functions which satisfy the following properties:
\begin{itemize}
\item $ G _{F_v}\overset {\kappa _{\mathcal D ' , v} }  \to o^{\times}_{\mathcal D ' } \overset {f} \to o^{\times}_{\mathcal D }$ is $ \kappa _{\mathcal D , v} $ for $v \in \bold {NO} (\bar \rho)\cup \bold {FL} (\bar \rho)$.
\item $\deform _{\mathcal D } ( v) \leq \deform_{\mathcal D'} (v) $ for any $v \in \vert F \vert _f $. 
\end{itemize}
The category of deformation types of $\bar \rho$ is denoted by $\Type (\bar \rho)$. 
\end{enumerate}
\end{dfn}

We now define the global deformation functor associated with a deformation type $\mathcal D  $ of $\bar \rho$. 

\begin{dfn}\label{dfn-galdef93} Let $\bar \rho: G_F\to \GL_2 ( k ) $ be an absolutely irreducible representation over a finite field of characteristic $\ell$, $\mathcal D$ a deformation type of $\bar \rho$.\par
For the coefficient ring $o_{\mathcal D }$ of $\mathcal D$ and $A \in \ob \mathcal C ^{\noeth} _{o_{\mathcal D }}$, 
$ \rho : G _{\Sigma_{\mathcal D} }  \to \GL _2 ( A) $ is a deformation of $\bar \rho$ of type $\mathcal D$ if the following conditions are satisfied:
\begin{itemize}
\item $\rho \mod m_A \simeq \bar \rho$. 
\item If $\deform _{\mathcal D} (v) =\bold f $ (resp. $\bold u$) at $v \in \vert F\vert_f $, $\rho \vert _{F_v } $ is a finite (resp. unrestricted) deformation of $\bar \rho \vert _{F_v} $. 
\item If $\deform _{\mathcal D} (v) =\bold {n.o.f.}$ (resp. $\bold {n.o.}$) at $v \in \vert F\vert_f $, $\rho \vert _{F_v} $ is a nearly ordinary finite (resp. nearly ordinary) deformation of $\bar \rho \vert _{F_v} $ of nearly ordinary type $ \kappa _{\mathcal D, v}$. 
\item If $\deform _{\mathcal D} (v) =\bold {fl} $ at $v \in \vert F\vert_f $, $\rho \vert _{F_v } $ is a flat deformation of $\bar \rho \vert _{F_v} $ of twist type $ \kappa _{\mathcal D, v}$. 
\end{itemize}
By $ F_{\mathcal D} (A) $ we denote the set of isomorphism classes of deformations of $\bar \rho$ of type $\mathcal D$ over $A$, and the functor
$$
 F_{\mathcal D} :  \mathcal C ^{\noeth} _{o_{\mathcal D }} \longrightarrow \Sets
$$
is the deformation functor of $\bar \rho$ of type $\mathcal D$. \par
For a continuous character $ \chi: G _{\Sigma _{\mathcal D} } \to o^{\times}_{\mathcal D} $, $F_{\mathcal D, \chi }$ is the subfunctor of $ F_{\mathcal D}$ consisting of the deformations with the determinant $\chi$.
\end{dfn}
\begin{thm} \label{thm-galdef91} Let $\bar \rho: G_F \to \GL_2 ( k ) $ be an absolutely irreducible representation over a finite field of characteristic $\ell$, $\mathcal D$ a deformation type of $\bar \rho$,  $ \chi: G _{\Sigma _{\mathcal D} } \to o^{\times}_{\mathcal D} $ a continuous character. Then the deformation functor $ F_{\mathcal D}$ (resp. $ F_{\mathcal D, \chi }$) of $\bar \rho$ is representable. 
\end{thm}
\medskip
The representability follows as in
\cite{M}, by the use of the Grothendieck-Schlessinger criterion, using $ H^0 ( G_{\Sigma}, \ad
^0  \bar \rho )= \{ 0 \}$ to assure the universality. The proof is so standard, and the details are omitted.   \par
\bigskip
By $ R_{\mathcal D }
$ (resp. $R_{\mathcal D, \chi }$), we denote the universal deformation ring representing
$ F_{\mathcal D} $ (resp. $F_{\mathcal D, \chi}$). 
$\rho  ^{\univ}_{\mathcal D}  : G _{\Sigma _{\mathcal D}  }
\to \GL _2 ( R _{ \mathcal D } )$ (resp. $\rho  ^{\univ}_{\mathcal D , \chi }  : G _{\Sigma _{\mathcal D}  }
\to \GL _2 ( R _{ \mathcal D , \chi })$ is the universal representation. \par
 As in \cite{W2}, the universal deformation
ring behaves nicely under the change of the coefficient ring by
$o_{\lambda} 
\to
o' _{\lambda'} $: $ R _{\mathcal D _{o' _{\lambda'}}} \stackrel{\sim} {\rightarrow} R _{\mathcal D }\otimes _{o_{\lambda} } o' _{\lambda'}$. \par

By the discussions in \S\ref{subsec-galdef4} and \S\ref{subsec-galdef7}, we have 
\begin{prop}\label{prop-galdef92}
\begin{enumerate}
\item The tangent space $ F_{\mathcal D} (k_{o _{\mathcal D}}[\epsilon ] ) $, which is identified with the Zariski tangent space $\Hom_{k_{o _{\mathcal D}}}  ( m_{R_{\mathcal D}} 
/( m _{ R _{\mathcal D} } ^2, m _{o_{\mathcal D} } )  ,k_{o _{\mathcal D}} )  
$ of
$ R_{\mathcal D}
$ over $ o _{\mathcal D}  $, is canonically isomorphic to 
$$
H ^ 1 _{\mathcal D} ( F, \ \ad \bar \rho ) = \ker ( H ^ 1 ( F,\
\ad
\bar  \rho ) \longrightarrow \bigoplus _{v \in \vert F\vert _f } H ^ 1 ( F_v ,\ \ad   \bar \rho \vert _{F_v}   ) /L_v ).
$$
Here $L_v$ is the local tangent space for $  \bar \rho \vert _{F_v} $ defined by $F ^{\deform _{\mathcal D} (v)} _{\bar \rho \vert _{F_v }  , \kappa _{\mathcal D, v} } ( k_{o _{\mathcal D}} [ \epsilon ] )$ if $\deform_{\mathcal D} (v ) \in \{ \bold {n.o.f} ,  \bold {n.o.} , \bold {fl} \}  $, $ L_v = F ^{\deform _{\mathcal D} (v)} _{\bar \rho \vert _{F_vs}   } ( k_{o _{\mathcal D}} [ \epsilon ] )$ in the other cases. 
\item For a continuous character $\chi: G_F \to o^{\times} _{ \mathcal D } $, the statement corresponding to (1) holds for $ F _{\mathcal D, \chi }$ and $R_{\mathcal D, \chi} $ using $\ad ^0 \bar \rho  $ instead of $\ad
\bar \rho $.
\end{enumerate}
\end{prop}

\subsection{Selmer groups}\label{subsec-galdef10}
To apply the level raising formalism in \S\ref{sec-TW}, we need
to use a variant of Selmer group associated to $\ad \rho $ as in \cite{W2}.
Since it is an easy translation of \S1 of \cite{W2}, we briefly discuss it here for our later use. \par
\bigskip
For an absolutely irreducible representation $\bar \rho$ with a deformation type $\mathcal D$, let $R_{\mathcal D} $ be the
universal deformation ring.  We assume that the coefficient ring $o _{\mathcal D } $ is an $\ell$-adic integer ring with the fraction field $E_{\mathcal D} $, and that there is an $o_{\mathcal D} $-homomorphism
$ f: R_{\mathcal D}
\to
o _{\mathcal D} $. By the universality of $R_{\mathcal D } $, $f$ corresponds to a deformation $\rho: G_{\Sigma_{\mathcal D}} \to \GL_2 (o_{\mathcal D} )  $
of $\bar \rho$ of type $\mathcal D$.\par

Let  $L\simeq o^{\oplus 2} _{\mathcal D} $ be the representation space of $\rho$. We define $M $ and $M ^ 0$ by 
$$ 
M = \ad \rho= L ^\vee \otimes_{o_{\mathcal D } }L,
$$
$$
M ^0 =
\ad ^ 0 \rho . 
$$
For each integer $ n \geq 1$, 
$$
M_n =  M
\otimes _{o_{\mathcal D} } m _{\mathcal D} ^{-n}   / o  _{\mathcal D}. 
$$
Note that the natural inclusion $ M_n \hookrightarrow  M\otimes _{o_{\mathcal D} } E_{\mathcal D}  /
o_{\mathcal D}  $ induces an injection 
$$
 H^ 1 ( F, \ M_ n ) \longhookrightarrow H ^ 1( F,\ M\otimes _{o _{\mathcal D}  } E_{\mathcal D } / o _{\mathcal D} )
$$
since $ H ^ 0 (F, \ \ad^0  \bar \rho  ) =  \{ 0 \} $.\par

For a finite place $v$ and an integer $n$, we define a local
subgroup
$H ^ 1 _{\deform _{\mathcal D} (v) } ( F_v, \ M_n )
\subset H ^ 1 (F_v ,\ M_n  )
$ (\cite{W2}, \cite{D1}, p.142) which reduces to the local tangent space if $n= 1$.\par
\bigskip
First assume that $ \deform _{\mathcal D}(v) = \bold {n.o.} $. Let $ \kappa _{\mathcal D, v}$ be the nearly ordinary type of $\bar \rho $. $\rho \vert _{F_v}$ has the form 
$$ 
0 \longrightarrow \chi _{1, v}  \longrightarrow   \rho \vert _{F_v }\longrightarrow  \chi_{2, v} \longrightarrow  0 .
$$
Here $\chi_{1, v} $ is the nearly ordinary character of $\rho \vert _{F_v} $. 

As in \S\ref{subsec-galdef6}, we define $W_1 M^0  _n \subset M^0  _n$ by 
$$
W_1 M^0  _n = \rho \vert _{F_v }\otimes \chi ^{-1} _{2, v} . 
$$
There is an exact sequence 
$$
0 \longrightarrow   W_2  M ^ 0 _n \longrightarrow   W_1 M ^ 0 _n \longrightarrow   \mathcal O _n \longrightarrow  0 
$$
for $W_1 M ^ 0 _n $. Here $ \mathcal O _n  = o _{\mathcal D } / m_{\mathcal D}  ^ n$, and $W_2  M ^ 0 _n = \chi _{1, v}\cdot \chi^{-1} _{2, v} \mod m ^n _{\mathcal D} $. \par

We define a subgroup $N_{\bold {n.o.} , n} $ of $H ^ 1 ( F^{\unr} _v,  M ^ 0_n )  $ as the image of $ H ^1 ( F^{\unr} _v , W_2 M ^0 _n )$ induced by $ W_2 M ^0 _n \to M ^0 _n $. Then 
$$
H ^1 _{\bold {n.o.} } ( F_v,\ M_n^0 )  =\ker  ( H ^ 1 (
F_v, \ M_n ^ 0 )
\to  H ^ 1 ( F^{\unr} _v,
\
 M ^ 0_n )  /N_{\bold {n.o.} , n}  ) , 
$$
$$
H ^1 _{\bold {n.o.} } ( F_v,\ M_n )= H ^1 _{\bold {n.o.} } ( F_v,\ M_n^0 )\oplus H ^1 
( F_v ^{\unr} / F_v ,  \mathcal O _n )  .
$$

Let $\mathcal O_n [\epsilon]$ be the ring of dual numbers over $\mathcal O_n $. Using that $\bar \rho $ is $G_{F_v} $-distinguished, as in \cite{W2}, proposition 1.1,
we may regard
$H ^1 _{\bold {n.o.} } ( F_v,\ M_n ^0)
$ as the group of extensions of the form 
$$
0 \longrightarrow \tilde  \chi  _{1, v}  \longrightarrow  \tilde  \rho   \longrightarrow  \tilde  \chi _{2, v}  \longrightarrow  0 
$$
such that $\tilde  \rho   \mod \epsilon \simeq  \rho \vert _{F_v}  \mod m_{\mathcal D}  ^n$, and $\det  \tilde \rho = \det \rho \vert _{F_v}$. Here
$\tilde \chi _{i, v}  : G_{F_v}
\to
\mathcal O_n [\epsilon]^\times $ is a character which lifts $\chi _{i , v} \mod m_{\mathcal D} ^ n$ for $i = 1, 2 $, $ \chi _{1, v} 
 \vert _{I_{F_v}} $ is equal to the nearly ordinary type $ \kappa_{\mathcal D, v}$. \par
\bigskip
Next assume that $ \deform _{\mathcal D}(v) = \bold {n.o.f.} $.  $\chi _1\vert _{I_{F_v}}  = \chi_2 (1)
\vert_{I_{F_v}}$ holds by the definition, and hence $\rho\vert_{F_v} $ defines an element $ c _{\rho\vert_{F_v} } \in H ^ 1 (
F_v ^{\unr} ,\  o_{\mathcal D}  (1) )$.

The subgroup $ H ^1 _{\bold
{n.o.f.}}( F^{\unr} _v,
\ W_2  M ^ 0 _n ) $
of
$ H ^ 1 ( F^{\unr} _v,
\ W_2 M ^ 0_n  )$ is defined as follows. Any element $ \rho '$ of $H ^ 1 ( F^{\unr} _v,
\ W_2 M ^ 0_n  ) $ is regarded as an extension
$$
 0 \longrightarrow  \kappa _{\mathcal D, v} 
\longrightarrow  \rho ' \longrightarrow  \kappa _{\mathcal D, v}  (-1) \longrightarrow 0
$$
over $\mathcal O_n [\epsilon]$, where the the extension class $ c_{\rho '} $ in $
 H ^ 1 ( F_v ^{\unr},\ \mathcal O_n [\epsilon]  (1) ) $ lifts $c_{\rho \vert _{F_v}}
\mod m_{\mathcal D}  ^n$. Then $\rho ' $ belongs to $ H ^1 _{\bold
{n.o.f.}}( F^{\unr} _v,
\ W_2  M ^ 0 _n ) $ if and only if $ \ev_{ \mathcal O_n [\epsilon]}  (c_{\rho '} ) = 0 $ for the evaluation map (Definition \ref{dfn-galdef615})\par
Let $ N _{\bold {n.o.f. } , n}  $ be the image of $ H ^1 _{\bold {n.o.f.}}( F^{\unr} _v,
\ W_2 M ^ 0 _n )$ in $H ^ 1 ( F^{\unr} _v,  M ^ 0_n )  $.
$$
H ^1 _{\bold {n.o.f.} } ( F_v,\ M_n^0 )=\ker  ( H ^ 1 ( F_v, \ M_n ^ 0 )
\to  H ^ 1 ( F^{\unr} _v,
\
 M ^ 0_n )  / N _{\bold {n.o.f. } , n} ) ,
$$
$$
H ^1 _{\bold {n.o.f.} } ( F_v,\ M_n )= H ^1 _{\bold {n.o.f.} } ( F_v,\ M_n^0 )\oplus H ^1
( F_v ^{\unr} / F_v ,\ \mathcal O _n )  .
$$
\begin{lem}\label{lem-galdef101}
Assume that $\deform _{\mathcal D } (v) = \bold {n.o.f.}$. 
$$
\length_{o_\mathcal D}  H ^1 _{\bold {n.o.} }
( F_v,\ M_n )  /H^1 _{\bold {n.o.f.} } ( F_v,\ M_n )\leq \length_{o_\mathcal D}   o _{\mathcal D }/ ( \chi
_{1, v}  (\sigma  ) ^ 2 - \det \rho \vert _{F_v}  (-1) (\sigma )).
$$
Here $\sigma $ is an element of $ G_{F_v} $ which lifts the geometric Frobenius element of $G_{k(v) } $.
\end{lem}
\begin{proof}[Proof of Lemma \ref{lem-galdef101}] By the definition, 
$ H ^1 _{\bold {n.o.} } ( F_v,\ M_n )/ H ^1 _{\bold {n.o.f.} } ( F_v,\ M_n )$ is seen as
a subgroup of $ H ^ 0( F ^{\unr} _v , \ \chi _{1, v}  / \chi_{2, v} 
(-1) ) $ by
$$ 
 H ^1 _{\bold {n.o.} } ( F_v,\ M_n )/ H ^1 _{\bold {n.o.f.} } ( F_v,\ M_n )
\longhookrightarrow  N _{\bold {n.o. } , n}   /N _{\bold {n.o.f. } , n} \simeq   H ^1 ( F^{\unr} _v , W_2 M ^0 _n )/ H ^1_{\bold{n.o.f.}} ( F^{\unr} _v , W_2 M ^0 _n )
$$
$$
 \overset{\ev}  \longhookrightarrow H ^ 0( F ^{\unr} _v , \ \chi _{1, v}  / \chi_{2, v} 
(-1) ) .
$$ 
Since the image belongs to the $G_{F_v}$-invariant subspace, the claim follows.
\end{proof}
Assume that $\deform_{\mathcal D}  (v) =\bold {fl} $. Since $M$ remains unchanged by a twist by a character, we may assume that the twist type $\kappa _v$ is trivial.\par
Then we define $ H ^ 1 _{\bold {fl}}( F_v,\ M_n ) $ as the
group of extensions 
$$
0 \longrightarrow  L_n  \longrightarrow  E \longrightarrow  L_n \longrightarrow  0 
$$
in $\mathcal G^{\mathcal O _n } _{F, [0,1]}$. \par
\bigskip
If $\deform _{\mathcal D} (v) = \bold f $, 
$$ 
H^ 1_{\bold f} ( F_v,\ M _n)=H^ 1_f
( F_v ,\ M_n ) = H ^1 ( F^{\unr}_v / F_v , M _n ) 
$$
 is the finite part. If $\deform _{\mathcal D} (v) = \bold u $, then
$$
H^ 1_{\bold u} ( F_v,\ M _n)=H^
1( F_v ,\ M_n ) .
$$
By the definition of $H^ 1_{\bold f} $ and $H^ 1_{\bold u}  $, the following lemma follows easily.
\begin{lem}\label{lem-galdef102} Assume that $\deform _{\mathcal D} (v) = \bold f $. Then 
$$
\length_{o_{\mathcal D} }  H^ 1_{\bold u} ( F_v,\ M _n)/H^ 1_{\bold f} ( F_v,\ M _n)\leq \length_{o_{\mathcal D} }   H^
1( F_v ,\ M_n ) /  H^1_f ( F_v ,\ M_n )
$$
$$ \leq \length_{o_{\mathcal D} }   H^0 ( F_v , \ M _n (1) ).
$$
\end{lem}

With these local groups defined for $v \in \vert F\vert _f$, we define the Selmer group for $ \ad \rho$. 
\begin{dfn}\label{dfn-galdef1 01} Let $\rho $ be a deformation of type $\mathcal D $ of $\bar \rho$. 
\begin{enumerate}
\item For any integer $n \geq 1$,  
$$ 
H ^ 1 _ {\mathcal D } ( F,\ M_n ) = \ker ( H ^ 1 ( F ,\ M_n ) \longrightarrow
\bigoplus_{ v \in \vert F\vert _f  } H ^ 1 (F_v, \ M _n ) / H ^ 1_{\deform _{\mathcal D} (v)} (F_v, \
M_n )  ).
$$
\item The Selmer group is defined by 
$$
\Sel _{\mathcal D} (F, \ M ) =\cup _n \operatorname{Image} (H ^ 1 _ {\mathcal D } ( F,\
M_n )  \to   H ^ 1 ( F,\ M\otimes _{o_{\mathcal D} } E_{\mathcal D} / o_{\mathcal D}  
)).
$$
\end{enumerate}
\end{dfn}
As in \cite{W2}, proposition 1.1, we have
\begin{prop}\label{prop-galdef101} Let $\rho $ be a deformation of type $\mathcal D $ of $\bar \rho$ which corresponds to $ f : R_{\mathcal D } \to o _{\mathcal D }$. Then
$$
\Hom _{o_{\mathcal D}  } ( \ker f / ( \ker f ) ^ 2 , E_{\mathcal D}  / o_{\mathcal D} )
\simeq \Sel _{\mathcal D} (F, \ M ) .
$$
\end{prop}
\bigskip
\begin{rem}\label{rem-galdef101}
Our Selmer group is related to the Selmer group in the sense of Bloch-Kato. As is discussed in \cite{D1}, p.142, under the canonical map $ H ^ 1 ( F_v
,\ M_n )
\to H ^1 (F_v , \ M_{n+1} ) $ induced from the inclusion $M_n \hookrightarrow M_{n+1} $, the
inverse image of $ H ^ 1 _{\bold{fl}} (F_v ,\ M _ {n+1}  )$ is $ H ^ 1  _{\bold{fl}} (F_v ,\
M _ n  )$. We define 
$ H ^ 1 _{\bold{fl}}( F_v , \ M )  $ as the union of images of $ H ^ 1  _{\bold{fl}} (F_v, \
M _n ) $ in $ H ^ 1 (F_v , \ M  )$. This subspace is equal to the finite part $ H ^1_f (F_v,
M ) $ in the sense of Bloch-Kato: by Wiles' argument in \cite{W2}, proposition 1.3, 1), $ H ^
1 _f ( F_v,
\ M )\subset H ^
1 _{\bold{fl}} ( F_v,
\ M )
$. Then the equality follows since $H ^
1 _f ( F_v,
\ M )$ is divisible, and $\length_{o_{\mathcal D}}  H ^
1 _f ( F_v,
\ M_1 ) \geq \length_{o_{\mathcal D}}   H ^
1 _{\bold{fl}} ( F_v,
\ M_1 )$. The last inequality follows from \cite{BK}, proposition 1.9, and Theorem \ref{thm-galdef72}.
\end{rem}

\section {Modular varieties and automorphic representations associated to quaternion algebras}\label{sec-shimura}

 In this section we briefly review modular varieties and automorphic representations
associated to quaternion algebras. See
\cite{W2}, \cite{H1}, \cite{H2} for further studies of Hecke algebras. 

\subsection{Modular varieties associated to quaternion algebras}
\label{subsec-shimura1}

Let $F$ be a totally real number field of degree $ [ F : \mathbb Q] = d$. By
$I _{F, \infty}$ we denote the set of all field embeddings $\iota :F \hookrightarrow \mathbb R$. $I_{F, \infty}$ is identified with the set $\vert F\vert _{\infty} $ of infinite places of $F$. \par
We take a quaternion algebra $ D$ which is central over $F$. By $D ^{\times}$ we mean the multiplicative group of $D$, which is  regarded as an algebraic group over $F$, and $G _D= \operatorname {Res} _{ F/ \mathbb Q } D ^{\times}$ is the Weil restriction to $\mathbb Q$. Let $Z$ be the center of $G_D$. \par

Let $I_D\subset I_{F, \infty} $ be the set of infinite
places of $F$ where $ D$ is split. We fix identifications
$$
D \otimes _{ \iota} \R \simeq M _2 ( \R ) \
\text{ for }
\iota
\in I _D,\quad
D \otimes _{ \iota} \R \simeq \mathbb H \ \text{ for } \iota \in  I _{F, \infty} \setminus I _D.
$$
Here $ \mathbb H $ is the Hamilton quaternion algebra. $G_D (\mathbb R )$ is isomorphic to $  \GL _2 ( \R )  ^{ I _D} \times  (\mathbb H ^{\times} )  ^{ I _{F, \infty} 
\setminus I _D  }$.\par
\bigskip
Let $X_D$ be a $G_D(\mathbb R )$-homogeneous space defined by
$$
 X _D = (\mathcal  H ^{\pm} ) ^{ I _D}  .
$$
Here $\mathcal  H ^{\pm} =  \P ^1 (\C) \setminus \P ^1 (\R)  $ is the double half plane, and the $G_D(\mathbb R)$-action is given by $G_D ( \mathbb R) \to  \GL_2 ( \mathbb R ) ^{I_D} \to \PGL_2 ( \mathbb R )^{I_D}  \simeq  ( \Aut  \mathcal  H ^{\pm} ) ^{I_D} $ by taking the projection to $\GL_2 ( \mathbb R ) ^{I_D}$. The complex dimension of $X_D$ is equal to $ \sharp I_D$, which we denote by $q_D$. 
We denote the stabilizer of $G _D (\mathbb R ) $ at $ ( \sqrt {-1}, \ldots,  \sqrt {-1} )\in X_D$ by $K _{\infty} $. \par
\bigskip

For a compact open subgroup $K\subset G_D (\A _{\mathbb Q, f})$, the associated modular variety $S _K ( G_D, X _D) (\mathbb C ) $ is defined by 
$$
 S _K ( G_D, X _D) (\mathbb C ) 
= G_D (\mathbb Q ) \backslash G _D( \A_ f ) \times X _D / K =  G_D (\mathbb Q )
\backslash G _D( \A  )/ K \times K_\infty.
$$
$S _K ( G_D, X _D) (\mathbb C )  $ is viewed as the set of the $\mathbb C$-valued points of a $\mathbb C$-scheme $ S _K ( G_D, X _D)_{\mathbb C} $ of finite type. \par
For an inclusion of subgroups $K \to K '  $, a natural projection 
$$
\pi_{K ', K }:  S _K ( G_D, X _D)_{\mathbb C} \longrightarrow  S _{K'} ( G_D, X _D)_{\mathbb C}  
$$
is induced, and $\{S _K ( G_D, X _D)_{\mathbb C} \} _{K \subset G_D (\mathbb A_{\mathbb Q, f } ) }$ forms a projective system. $G_D (\mathbb A_{\mathbb Q , f }) $ acts on 
$$
S ( G_D, X_D ) _{\mathbb C} = \varprojlim _{K \subset G_D (\mathbb A_{\mathbb Q, f } )  } S_K ( G_D, X_D ) _{\mathbb C} 
$$
by the right translation. The action of $ g \in  G_D (\mathbb A_{\mathbb Q, f})$ is denoted by $ R(g)$. \par
In this paper, we mainly consider the following two cases. 

\medskip
\begin{itemize}
\item
$q_D = \sharp I_D = 1 $: $ D$ is ramified at any infinite places except $\iota _0\in I_{F, \infty} $.
\end{itemize}

In this case, $S _K ( G_D, X _D) _{\mathbb C} $ is a Shimura curve, which is proper unless $F = \mathbb Q$ and $ D= M_2 (\mathbb Q)$. By
the theory of canonical models of Shimura \cite{Sh}, there is a model
$S _K(G_D , X_D )_F $ canonically defined over $F$, $S _K(G_D , X_D )_F  \times _{F,   \iota_0 } \Spec \mathbb C = S _K ( G_D, X _D)_{\mathbb C} $. Since these models are canonical, 
$$ 
S ( G_D, X_D ) _F=
\varprojlim _K 
S_K ( G_D, X_D )_F
$$ 
is defined over $F$, and hence the $G _D (\A _{\mathbb Q, f}  )$-action is also defined over $F$.

\medskip
\begin{itemize}
\item 
$q_D = \sharp I_D = 0 $: $ D$ is a quaternion algebra over $F$ which is ramified at all infinite places. 
\end{itemize}
\medskip

In this case, $ S _K ( G_D, X _D) _{\C }$ is a zero-dimensional scheme over $\C$, which we call the {\it Hida variety} associated to $ ( G_D, X _D)$. A Hida variety is not a (non-connected) Shimura variety in the sense of Deligne (since it has a compact factor defined over
$\mathbb Q$). The variety was first considered by Hida in his study of Hecke algebras for $\GL _2$ \cite{H1}.

\subsection {Equivariant sheaves on modular varieties}
\label{subsec-shimura2}
In this subsection, $D$ is a division algebra, and we denote $S _K ( G_D, X _D) $ by $S_K $ for short. 
For a prime $\ell$ and an $\ell$-adic field $E_\lambda $, there is an $ E_{\lambda}$-smooth sheaf $\bar{\mathcal F} ^K _{(k,w) , E_\lambda} $ on $S_K$ (see
\cite{Car2} in the Shimura curve case, \cite{H1} for the Hida variety case).  We discuss it here with
a $\Z_{\ell}$-structure.\par
For a finite place $v\vert \ell$, $I _{F, v } $ is the set of field embeddings $F_v \hookrightarrow \bar E_\lambda $ over $\mathbb Q_{\ell}$,  $ I _{F, \ell}= \coprod _{v \vert \ell } I _{F, v} $. By an isomorphism $\bar E_{\lambda } \simeq \mathbb C$, $ I_{F, \ell } $ is identified with $I _{F, \infty}$.
\begin{dfn} \label{dfn-shimura21} Let  $ ( k , w) $ be a pair of an element $ k$ in $ \mathbb Z ^{I_{F,\infty } } $ and an integer $w$. 
\begin{enumerate}
\item $(k, w ) $ is called an infinity type if $k = ( k_\iota ) _{\iota \in I_{F, \infty} }  $ and $w$ satisfy $k_\iota \geq 1  $ and $k_\iota \equiv w \mod 2 $ for any $\iota \in I _{F, \infty}$. 
\item An infinity type $(k, w)  $ is called a discrete type if $ k _{\iota } \geq 2 $ for any $\iota \in I_{F, \infty}$. 
\item For $ \iota \in I_{F, \infty} $, $k ' _{\iota } =  \frac {w- k_{\iota }}{2} +1$.
\end{enumerate}
\end{dfn}
By the above identification, an infinity type $(k,w )$ is also viewed as a pair $k \in \mathbb Z ^{I_{E, \ell} }  $ and $w\in \mathbb Z$. \par
We assume that $D$ is split at all $v \vert \ell$, and choose an isomorphism 
$ D \otimes _{\mathbb Q} \mathbb Q_{\ell} \simeq \prod _{v \vert \ell} M_2 (F_v)$. By this isomorphism we regard $\prod _{v \vert \ell }\GL_2 (o_{F_v} )$ as a compact open subgroup of $G_D(\mathbb Q_\ell)$. \par

Let $ V_{F_v} = F_v ^{\oplus 2} $ be the standard representation of $ \GL _{2, F_v} $ over $F_v$, $ V_{o_{F_v}}= o^{\oplus 2}_{F_v }$ the standard lattice in $V_{F_v} $ fixed by $\GL_2 ( o_{F_v})$. 
We take $ E_\lambda$ so that any $\iota: F _v \hookrightarrow \bar E _{\lambda}$ in $I _{F, \ell }$ factors through $ E_{\lambda }$. For an infinity type $ (k, w ) $ of discrete type, the representation 
$$
V_{(k, w) , E_\lambda } = \bigotimes _{\iota : F_v \hookrightarrow E_\lambda, v \vert \ell  } (\iota \det) ^{-k' _{\iota}}  \cdot \Sym ^{ k _{\iota} - 2 }(V_{F_v} \otimes _{\iota}  E_\lambda )^{\vee} 
$$
of $ G_D (\mathbb  Q _\ell)$ is defined over
$E_\lambda$, and has an
$o_{E_\lambda}$-lattice 
$$
V_{ ( k,w ) , o_{E_\lambda}} = 
\bigotimes _{\iota  : F_v \hookrightarrow E_\lambda , v \vert \ell } (\iota \det) ^{-k' _{\iota}} \cdot \Sym   ^{ k
_{\iota}-2}
(V _{o_{F_v}} \otimes _{\iota} o  _{E_\lambda} ) ^{\vee} ,
$$
which is stable under $\prod _{v \vert \ell }\GL_2 (o_{F_v} ) $.\par

\begin{dfn}\label{dfn-small1} Let $K$ be a compact open subgroup of $G_D (\mathbb A _{\mathbb Q, f}  ) $. $K$ is called small if for any compact open normal subgroup $ K ' \subset K $,  the action of  $(K\cap \overline  {F^\times }) \cdot K '  \backslash K $ on $S_{K '} $ induced by the right action of $K$ is free. 
\end{dfn}
It is easily seen that if $K$ is a small subgroup, then any compact open subgroup $K'$ of $K$ is again small. We discuss the smallness of a compact open subgroup in subsection \ref{subsec-shimura5}.\par
For a small compact open subgroup $ K= K _\ell \cdot K ^{\ell}$ of $ G _D ( \mathbb A_{\mathbb Q, f} )  $, let 
$$
\pi _{\ell} : \tilde S _{\ell} \longrightarrow S _{K} 
$$
be the Galois covering corresponding to $  (K \cap \overline { F ^{\times} })_{\ell} \backslash K _{\ell}$, where $(K \cap \overline { F ^{\times} } )_{\ell} $ is the image of $ K \cap \overline { F ^{\times}} $ by the projection $K \to K_{\ell}$.\par

If the action of $ K\cap
\overline{ F^\times }$ on the representation  $V_{(k, w), o_{E_\lambda} } $ is trivial, an $o _{E_\lambda}$-smooth sheaf $ \bar {\mathcal F} ^K  _{(k, w) , o_{E_\lambda}}$ on $S_K$ is obtained from the covering
$\pi _{\ell}$ and the representation $V_{ ( k,w ) , o_{E_\lambda}}$. 
Let $g $ be an element of $G_D (\mathbb A _{\mathbb Q , f } )  $. By the isomorphism $S_{g  Kg ^{-1}  } \overset { R(g)  }  {\stackrel{\sim}{\rightarrow}} S_K$ induced by the right $G_D (\mathbb A_{\mathbb Q, f} )$-action, 
$\bar {\mathcal F}  ^{g  Kg ^{-1} }  _{( k, w) , o_{E_\lambda}}=R(g) ^* \bar {\mathcal F} ^K _{( k, w) , o_{E_\lambda} } $ is defined by the lattice $g _{\ell}  V _{(k, w ) , o_{E_\lambda} }$, where $g_{\ell}$ is the $\ell$-component of $g$. \par

 For the right $G_D (\A _{\mathbb Q, f }) $-action on $S = \varprojlim _K S_K $, the $E_\lambda
$-sheaf 
$$ 
\bar { \mathcal F }^K  _{(k, w),\ E_\lambda }=  \bar {\mathcal F}^K  _{(k, w) , o _{E_\lambda} }\otimes _{o_{E_\lambda} }
E_\lambda 
$$
forms a $ G_D (\A _{\mathbb Q, f }) $-equivariant system of smooth sheaves as $K$ varies by the construction.
The $o_{E_\lambda}$-lattice structure is preserved by $G_D (\A _{\mathbb Q, f } ^{\ell }  )  $-action. \par
When $q_D $ is one, the sheaf $ \bar {\mathcal F}^K   _{(k, w)}$ is canonically defined over
$F$ by the theory of canonical models, which we denote by $ \mathcal F^K _{(k, w)}$. This canonical
$F$-structure gives a continuous $G_F$-action on the \'etale cohomology groups over $\bar F$. \par
The Betti-version of $\bar {\mathcal F} _{(k, w)} $ is discussed in \cite{Car2}, p.418--419. By the comparison
theorem in
\'etale cohomology, those two cohomology theories are canonically isomorphic, so we do not make
any distinction unless otherwise stated. \par
\begin{rem} In the case of Shimura curves, $ \mathcal F ^K_{(k, w) , E_\lambda}$ is pure of weight $w$. 
\end{rem}

\subsection{Hecke correspondences and the duality formalism} 
\label{subsec-shimura3}
For a commutative ring $R$ and a compact open subgroup $K  $ of $G _D (\mathbb A _{\mathbb Q, f })  $, let $H _{K  , R}= H ( G _D (\mathbb A _{\mathbb Q, f }) ,   \ K)_ R$ be the convolution algebra formed by the compactly supported $R$-valued
$K$-biinvariant functions on $G _D (\mathbb A _{\mathbb Q, f }) $. 
The right $G_D (\mathbb A _{\mathbb Q, f } ) $-action on $S= \varprojlim _{ K } S _{ K } $ yields a left $ G_D (\mathbb A _{\mathbb Q, f } ) $-module structure on 
$$ 
H ^* ( S, \bar {\mathcal F} _{(k, w), E_\lambda }  )\overset {def} =\varinjlim H ^* ( S_K, \bar {\mathcal F}^K  _{(k, w), E_\lambda }  )  . 
$$
At each finite level $K$, $H ^* ( S_K, \bar {\mathcal F}^K  _{(k, w), E_\lambda }  )  $ has a left action of the convolution algebra $H_{K, E_\lambda} $. In this subsection, we discuss this action in detail.\par

\begin{dfn} Let $K$ be a compact open subgroup $K$ of $ G_D (\mathbb A _{\mathbb Q, f } ) $.
\begin{enumerate}
\item  $K$ is $\mathbb Q$-factorizable if $K = \prod _{q}K_q $, $K _q \subset G_D ( \mathbb Q_q ) $. For a finite set $\Sigma $ of primes , $K _\Sigma = \prod _{q \in \Sigma }K_q $, $K ^\Sigma = \prod _{q \not \in \Sigma }K_q $.
\item $K$ is $F$-factorizable if $K = \prod _{v \in \vert F\vert _f }K_v$, $K _v \subset D ^{\times } (F_v ) $ by the identification $ G_D (\mathbb A _{\mathbb Q, f } )= D ^{\times } ( \mathbb A_{F, f} ) $.
For a finite set $\Sigma $ of finite places of $F $, $K _\Sigma = \prod _{v\in \Sigma }K_v $, $K ^\Sigma = \prod _{v \not \in \Sigma }K_v $.
\item 
For an $F$-factorizable compact open subgroup $K = \prod _v  K_v $ of $G_D (\mathbb A _{\mathbb Q, f } ) $, $\Sigma _K $ is the set of finite places which satisfies the following property:\par
$v\not \in \Sigma_K $ $\Leftrightarrow$ $D$ is split at $v$, and $K_v$ is a maximal hyperspecial subgroup of $D ^\times
(F_v)$. 
\end{enumerate}
\end{dfn}
When $ v \not \in \Sigma _K$, $K_v$ is isomorphic to $\GL_2 (o_{F_v} ) $ by some group scheme isomorphism
$ D ^\times_v \simeq \GL_{2. F_v} $ suitably taken.\par

For two compact open subgroups $K ,\ K '   \subset G_D (\A _{\mathbb Q , f} ) $, $g \in G_D (\A _{\mathbb Q , f} ) $ define an algebraic correspondence
$$
[K   gK'   ] : \quad  
\CD 
S_{K}  @< \pr  _2<< S_{K\cap g  K ' g^{-1} }\overset {R(g^{-1}) } {\stackrel{\sim }{\longrightarrow}} S _{ g^{-1} K g   \cap K '} 
@> \pr _1 >> 
 S_{K'},
\endCD 
$$
where $ \pr _1 = \pi _{ g^{-1} K g   \cap K ', K '}$ and $ \pr  _2= \pi _{K\cap g  K ' g^{-1} , K } $ correspond to the inclusion of groups, and the direction of the correspondence is from the second projection to the first projection, that is, from $S_K $ to $S _{K '}$.

The algebraic correspondence induced by $K 'g K $ from $S_K $ to $S_{K'} $ is dual to 
$ K g ^{-1} K' $ from $S_{K'} $ to  $S_K $. \par

Since $\mathcal F_{(k, w) , E_\lambda} $ is $G_D (\A _{\mathbb Q , f} ) $-equivariant by the
construction, 
$$
[K gK' ] ^*: R\Gamma ( S_{ K  ' } , \bar {\mathcal F}^{K'}  _{(k, w),  E_\lambda} ) \longrightarrow  R\Gamma ( S_{K } ,
\bar {\mathcal F} ^K _{(k, w)  ,  E_\lambda} ) 
$$
is induced by 
$$
R\Gamma ( S_{ K  ' } , \bar {\mathcal F}^{K'}  _{(k, w),  E_\lambda} ) \overset { \pr_1 ^* }  {\stackrel{\sim}{\longrightarrow}} R\Gamma ( S_{ g ^{-1} Kg \cap K ' } , \bar {\mathcal F}^{ g ^{-1} Kg \cap K '  }  _{(k, w),  E_\lambda} ) 
$$
$$
\overset{R(g^{-1}) ^* } \longrightarrow R \Gamma (S_{K \cap g K' g^{-1}} , \bar {\mathcal F}^{K \cap g K' g ^{-1}}  _{(k, w),  E_\lambda} )   \longrightarrow  R\Gamma ( S_{K } ,
\bar {\mathcal F} ^K _{(k, w)  ,  E_\lambda} ) .
$$
We call $ [K gK'  ] ^*$ the standard action of $[K gK '] $. When $K = K '$, it is defined with the $o_{E_\lambda} $-lattice structure if $ g^{-1} _\ell  ( V_{(k, w ) , o_{E_\lambda }} ) \subset  V_{( k, w ) , o_{E_\lambda} }$. \par
The action of the characteristic function $\chi_{KgK} $ of $H_{K, E_\lambda }$ on $ H ^*  (S_K, \bar {\mathcal F} ^K _{E_\lambda}  ) $ is the standard action $ [ K g K ]^* $. 
Moreover, we have the action of
$H ( G_D (\A ^{\ell}_{\mathbb Q, f } ) ,\ K ^{\ell})_{o_{E_\lambda}} $ on $H ^* ( S_K,  \bar {\mathcal F}^{K}  _{(k, w),  o_{E_\lambda}} ) $ for $ K = K _\ell \cdot K ^\ell$: For $g \in  G_D (\A ^{\ell}_{\mathbb Q, f } )  $, $ K ^{\ell} gK^ \ell$ acts by the standard action of 
$ K \tilde g K $, $\tilde g = ( 1 _{D(F_v)} ) _{ v \vert \ell } \cdot g \in G_D (\A ^{\ell}_{\mathbb Q, f })$. \par

It is possible to lift the action of $H (G_D (\A ^{\ell} _{\mathbb Q, f} ) ,\ K ^{\ell}) _{o_{E_\lambda}}$ on the cohomology groups to $R\Gamma ( S_K,  \bar {\mathcal F}^{K}  _{(k, w),  o_{E_\lambda}} ) $ as the following proposition shows:
\begin{prop}\label{prop-shimura31}
There is a complex of left $H (G_D (\A ^{\ell} _{\mathbb Q, f} ) ,\ K ^{\ell})_{o_{E_\lambda}} $-modules bounded below which represents $ R\Gamma ( S_{ K } , \bar
{\mathcal F} ^K _{(k, w ), o_{E_\lambda}} )$. 
\end{prop}
\begin{lem}\label{lem-shimura31} Let $Y \overset { f_1, f_2} \rightarrow X$ be finite \'etale morphisms between schemes of finite type over $\mathbb C$, $\mathcal F $ an abelian sheaf on $ X (\mathbb C)$. We regard $  X\overset{ f_1}  \leftarrow Y\overset {f_2} \rightarrow X $ as a correspondence, and assume that a cohomological correspondence $ c: f^* _2 \mathcal F\stackrel{\sim} {\rightarrow} f^*_1 \mathcal F  $ is given.  Let $L^{\cdot} $ be the Godement's canonical resolution of $\mathcal F$. Then $c$ induces a morphism of complexes $c ^* : \Gamma ( X(\mathbb C) , L ^\cdot ) \to \Gamma ( X(\mathbb C) , L ^\cdot )$ which gives the endmorphism of $R \Gamma ( X(\mathbb C) , \mathcal F )$ in the derived category $ D ^+(X(\mathbb C)) $ induced by $c$.
\end{lem}
\begin{proof}[Proof of Lemma \ref{lem-shimura31}] 
Since $f_1 $ and $f_2 $ are \'etale, they are local isomorphisms in the analytic category, and $f ^*_ i  L ^*  $ is the canonical resolution of $ f^*_ i \mathcal F  $ for $i = 1,2 $. Thus $c$ induces an isomorphism of complexes $c^* : f^*_ 1 L ^{\cdot} \stackrel{\sim} { \rightarrow}  f^*_2 
L ^{\cdot}  $. Since $f_2 $ is finite  \'etale, the trace map is defined for sheaves, thus $ \operatorname{tr} : (f _{2})_*  f _ 2 ^*  L ^{\cdot} \to L ^{\cdot}$ is defined. The composite of 
$$
\Gamma (X(\mathbb C ) , L ^{\cdot}) \longrightarrow  \Gamma (Y(\mathbb C ) , f ^* _1 L ^{\cdot})
$$
$$\overset {\Gamma ( c ^* ) } 
 \longrightarrow  \Gamma (Y(\mathbb C ) , f ^* _2 L ^{\cdot}) = \Gamma ( X(\mathbb C ) ,   (f _{2})_* f _ 2 ^*  L ^{\cdot}  ) \overset{\Gamma (\operatorname{tr} ) } \longrightarrow \Gamma ( X(\mathbb C ) ,   L ^{\cdot}  ) 
$$
satisfies the desired property. 
\end{proof}
\begin{proof}[Proof of Proposition \ref{prop-shimura31}]
For simplicity, we work with the Betti
realization described in \cite{Car2}, p.418--419. This is sufficient because the Betti and \'etale cohomologies are
canonically isomorphic by the comparison theorem in \'etale cohomology.\par
We assume that the $\ell$-component of $g$ is $1$.  We have a morphism 
$$
 \pr _1 ^* \bar {\mathcal F}^{K}  _{(k, w), o_{ E_\lambda} } = \bar {\mathcal F}^{ g ^{-1} Kg \cap K   }  _{(k, w),  o_{E_\lambda}} \overset{R(g^{-1}) ^* } {\stackrel{\sim}{\longrightarrow}}\bar {\mathcal F}^{K \cap g K g ^{-1}}  _{(k, w),  o_{E_\lambda}} 
=  \pr ^*_2 \bar {\mathcal F} ^K _{(k, w)  , o_{ E_\lambda} } . 
$$
Let $ L^\cdot $ be Godement's canonical resolution of $\bar { \mathcal F} ^K _{ (k, w ), o_{E_\lambda}  } $ on $S_K$. 
Since $\pr _1 $ and $\pr _2$ are finite \'etale, Lemma \ref{lem-shimura31} is applied, 
and the action of $\chi _{KgK } $ is actually defined on
$\Gamma ( S_K , L^\cdot ) $.  
One checks that this action of $\chi _{KgK } $ extends to a left action of $H (G_D (\A ^{\ell} _{\mathbb Q , f} ) ,\ K^{\ell}  )_{o_{E_\lambda}}$, then 
$\Gamma (L^\cdot ) $ defines an object of the derived category of $ H (G_D (\A ^{\ell} _{\mathbb Q , f} ) ,\ K^{\ell}  )_{o_{E_\lambda}}
$-modules bounded below which lifts $ R \Gamma (S_K , \ \bar  { \mathcal F} ^K _{(k, w), o_{E_\lambda} } ) $. 
\end{proof}

\bigskip

For the relation between $H ( G _D (\mathbb A ^{\ell}_{\mathbb Q , f } )  ,\ K^{\ell} ) _{o_{E_\lambda}}$-action and the Verdier duality, we have the following. By the definition of $\bar {\mathcal F}^K _{(k,w), o_{E_\lambda}  }$,
$$
(\bar {\mathcal F}^K  _{(k, w), o_{E_\lambda}} ) ^{\vee}\stackrel{\sim}{\longrightarrow} \bar {\mathcal F} ^K _{(k, -w), o_{E_\lambda} } , 
\leqno{(\ast)}
$$
since $V_{(k,w) , o_{E_\lambda}} $ satisfies $V_{(k,w) , o_{E_\lambda}} ^\vee \stackrel{\sim}{\rightarrow}  V_{ (k, -w) , o_{E_\lambda} } $. 

This gives a perfect pairing in the derived category of $o_{E_\lambda}$-modules
$$
R
\Gamma  (S_{K } ,\ \bar {\mathcal F }^K _{(k, w), o_{E_\lambda}} )\otimes^{\mathbb L} _{o_{E_\lambda}} R
\Gamma  (S_{K }  ,\ \bar { \mathcal F }^K  _{(k, -w), o_{E_\lambda} }) \longrightarrow o _{E_\lambda} ( -q_D )[-2q_D] 
$$
by Poincar\'e duality. We need to know how Poincar\'e duality exchanges $H (G_D(\A ^{\ell} _{\mathbb Q, f } ), K^{\ell} )_{o_{E_\lambda}}$-actions. By the isomorphism ($\ast$), $(R (g ) ^* ) ^{\vee}$ is identified with $ R(g ^ {-1}) ^*$.
This implies that the standard action of $K gK$ on $
 R \Gamma  (S_{K } ,\
\bar {\mathcal F }  _{(k, w) , o_{E_\lambda}} )$ corresponds to the standard action of
$Kg ^{-1}K $ on $ R\Gamma  (S_{K }  ,\ \bar {\mathcal F}   _{(k, -w), o_{E_\lambda}}) $.

\begin{prop}
\label{prop-shimura32}
The standard action 
$$ 
R \Gamma  (S_{K '} ,\
\bar {\mathcal F} ^{K '}_{(k, w), o_{E_\lambda}})\overset {[K  g K ' ] ^*}  \longrightarrow R \Gamma  (S_{K} ,\
\bar {\mathcal F} ^K _{(k, w), o_{E_\lambda}} )
$$
 induced by $ [K g K ']$ is dual to the standard action
$$
 R \Gamma  (S_{K} ,\
\bar {\mathcal F}^K  _{(k, -w), o_{E_\lambda}})(q_D) [2q_D] \overset {[K' g^{-1} K] ^* (q_D ) [2q_D]} \longrightarrow R \Gamma  (S_{K' } ,\
\bar {\mathcal F}^{K '} _{(k, -w), o_{E_\lambda}})(q_D)  [2q_D]
$$ by $ [K'  g ^{-1} K ]$.
\end{prop}

We have two geometric actions of the convolution algebra, which we call the
standard action and the dual action. The standard action of $[KgK] $ is $[K g K ]^* $ we have already introduced. By the dual action of $ [ K g K ]$ on $  R
\Gamma  (S_K ,
\bar {\mathcal F} _{(k, w), o_{E_\lambda}}) $, we mean the standard action of $ [ K g^{-1} K ]$. Proposition \ref{prop-shimura32} implies that
Poincar\'e duality exchanges the standard action to the dual action.  The standard action (resp. dual action) is a left (resp. right) action.

For an $F$-factorizable compact open subgroup 
$K $ and a finite place $v$ of $F$ such that $D$ is split at $v$, choose a uniformizer $ p_v $ of $F_v$. Define $a(p_v ) $ and $b(p_v)
\in G_D (\A _{\mathbb Q, f}  )
$ as the elements having 
$\begin{pmatrix} 1 & 0
\\ 0 & p _v 
\end{pmatrix} ,\
 \begin{pmatrix}
 p_v & 0
\\ 0 & p_v
\end{pmatrix} $ as the $v$-component and the other components are $1$, respectively.\par

As in \S\ref{subsec-galdef8}, $U(p_v) $ and $U ( p_v , p_v )$ are defined by $\chi_{ K a(p_v)K }$ and $\chi_{ K b(p_v)K }$.
These operators are called the standard Hecke operators. They are independent of the choice of a uniformizer if $ K_v $ contains $K_0 ( m_v) $. When $K _v = \GL _2 (o_{F_v} )  $, $U ( p_v ) $ and $U (p_v , p_v)$ are denoted by $T_v$ and $T_{v, v} $. \par

\subsection{The reciprocity law for $S(G_D , X_D) $}\label{subsec-shimura4}

For a quaternion algebra $D$ over $F$ and an infinity type $(k, w)$ of discrete type, let $\mathcal A ^{\disc}_{(k, w)} (G_D(\A _{\mathbb Q}) )$ be the set of isomorphism classes of irreducible essentially square integrable representations of $G_D ( \A _{\mathbb Q}) $ of infinity type
$(k, w)
$: for $\pi \in \mathcal A ^{\disc} _{(k, w)} (G_D ( \A _{\mathbb Q})  )$, $\pi _\infty$ takes the following form  
$$
\pi _\infty =(\bigotimes _{\iota \in I_D } D_{ k _{\iota }, w})\otimes 
(\bigotimes _{ \iota \in I_{F, \infty}\setminus I_D } \bar D_{ (k _\iota, w) }) , 
$$ 
where $\bar D_{( k_\iota , w) } = \Norm_{ \mathbb H / \R } ^{-k' _\iota } \cdot   \Sym^{ (k _\iota -2 ) }V_ \st ^{\vee}   $ is the irreducible representation of $ \mathbb H ^{\times} $ which corresponds to $
D_{(k_{\iota } , w ) } $ by the Jacquet-Langlands correspondence \cite{JL}. 
$\Norm_{ \mathbb H / \R } $ is the reduced norm, $V_{\st} $ is the standard representation of $\GL_2(\C) $, and we view 
$\GL_2(\C)$-representation $ \Sym ^{(k _\iota -2 ) } V_{\st} ^{\vee}  $ as a representation of
$\mathbb H ^\times$. For a Hecke character
$\chi : F ^\times \backslash \A^\times  _F \to \C ^\times $, $\mathcal A ^{\disc} _{(k, w), \chi  }( G_D ( \A _{\mathbb Q}) )
$ is the subset of $\mathcal A^{\disc} _{(k, w)}(G_D ( \A _{\mathbb Q}) ) $ consisting of the representations with the central character $\chi$. \par

By the Jacquet-Langlands \cite{JL}, Shimizu \cite{S} correspondence, we have an injection
$$
\JL: \mathcal A ^{\disc} _{(k, w)} (G_D ( \A _{\mathbb Q})  ) \longhookrightarrow \mathcal A ^{\disc} _{(k, w)} (\operatorname {Res} 
_{F/{\mathbb Q}}
\GL_{2, F} ( \A _{\mathbb Q}) )  , 
$$
and
$\mathcal A _{(k, w)} (G_D  ( \A _{\mathbb Q}) )
$ (resp. $ \mathcal A _{(k, w), \chi } (G_D ( \A _{\mathbb Q})  )$) is defined as the subset of $\mathcal A ^{\disc} _{(k, w)} (G_ D  ( \A _{\mathbb Q}) )$ (resp. $\mathcal A ^{\disc} _{(k, w), \chi } (G_ D )$) which correspond to cuspidal
representations of $\GL_2 (\mathbb A_F) $. The image of $\mathcal A _{(k, w)} (G_ D ( \A _{\mathbb Q})  )
$ by $\JL$ consists of a cuspidal representation $\pi$ which has an essentially square integrable component $\pi _v$ at
$v$ where $D$ is ramified. When $D$ is a split quaternion algebra over $F$, we denote $\mathcal A ^{\disc} _{(k, w)}( G_D ( \A _{\mathbb Q}) )$ (resp. $\mathcal A ^{\disc} _{(k, w), \chi  }( G_D  ( \A _{\mathbb Q}) )$, resp. $\mathcal A _{(k, w) } (G_D ( \A _{\mathbb Q}) )$, resp. $\mathcal A _{(k, w), \chi  } (G_D ( \A _{\mathbb Q}) )$) by $\mathcal A ^{\disc} _{F, (k, w)}$ (resp. $\mathcal A ^{\disc} _{F, (k, w), \chi  }$, resp. $\mathcal A _{F, (k, w) } $, resp. $\mathcal A _{F, (k, w), \chi  }$). \par
\bigskip
 
First consider the Shimura curve case. For an infinity type $(k,w )$ of discrete type, the decomposition of the \'etale cohomology groups as $ G_F \times H (G_D (\mathbb A_{\mathbb Q} ) , K  ) $-bimodules is given by
$$
H ^1 _{\et}( S_{ K, \bar F} ,\ \bar {\mathcal F}^K _{ (k, w),  \bar E_\lambda }) \simeq
\bigoplus _{\pi \in
\mathcal A _{(k, w)} (G_D ( \A _{\mathbb Q})  ) }
\rho _{
\pi,
\bar E_\lambda }\otimes _{ \bar E_\lambda  }  \pi _f   ^{K },
$$
and 
$$
H ^0 _{\et}( S_{ K, \bar F} ,\ \bar {\mathcal F}^K _{ (k, w),  \bar E_\lambda }) \oplus H ^2 _{\et}( S_{ K, \bar F} ,\ \bar {\mathcal F}^K _{ (k, w),  \bar E_\lambda })\simeq 
\oplus _{\pi \in \mathcal A  ^c_{(k, w)} (G_ D  ( \A _{\mathbb Q}) ) }
\rho _{\pi, \bar E_\lambda }\otimes _{ \bar E_\lambda  }  \pi _f   ^{K }.
$$
Here $\mathcal A  ^c_{(k, w)} (G_ D  ( \A _{\mathbb Q}) )=\mathcal  A  ^{\disc}_{(k, w)} (G_ D  ( \A _{\mathbb Q}) )\setminus \mathcal A  _{(k, w)} (G_ D  ( \A _{\mathbb Q}) ) $.
For $\pi  \in
\mathcal A _{(k, w)} (G_D ) $, $\rho _{\pi , \bar E_\lambda} : G_F \to \GL _2 (\bar E _\lambda ) $ is the two dimensional irreducible $\ell$-adic representation associated to $\pi$ (\cite{O1},\cite{Car2}), and we view the finite part $ \pi _f$ of $\pi $ is defined over $\bar E _\lambda $ by the identification $
\bar E _\lambda  \simeq \C$. 
Note that $ \mathcal A  ^c_{(k, w)} (G_ D ( \A _{\mathbb Q}) )$ is
non-empty if and only if $k = (2,
\ldots, 2)$, and $\mathcal A  ^{c} _{(k, w)} (G_D ( \A _{\mathbb Q}) ) $ consists of one dimensional representations which factor through the
reduced norm $D ^\times  ( \mathbb A_F ) \overset {\Norm_{D/F}} \to \mathbb  A_F ^\times \overset {\chi } \to \mathbb C^{\times} $. 
In this case, $\chi$ is an algebraic Hecke character of weight $w$, and  $\rho _{\pi, \bar E_\lambda }$ is $\rho_{\chi ,  \bar E_\lambda} \oplus \rho _ {\chi,  \bar E_\lambda} ( -1)  $, where $\rho_{\chi,  \bar E_\lambda }:G_F \to \bar E^{\times }_{\lambda} $ is the $G_F$-representation attached to the algebraic Hecke character $\chi$.
\par
In the case of Hida varieties, there are no natural Galois actions, still there is a decomposition 

$$
 H ^0_{\et} ( S_{K, \C} ,\  \bar {\mathcal F} ^K _{(k, w) ,  \bar E_{\lambda} } )\simeq 
 I_{(k, w), E_\lambda } ^K  \oplus  (\bigoplus _{ \pi
\in
\mathcal A  _{(k, w)} (G_D  ( \A _{\mathbb Q})) }  \pi  _f ^K  )
$$
as a $H (G_D (\mathbb A_{\mathbb Q} ) , K  ) $-module. Here $I_{(k, w ), \bar E_\lambda } ^K$ is the subspace of $H ^0_{\et} ( S_{K, \C} ,\  \bar {\mathcal F} ^K _{(k, w) ,  \bar E_{\lambda} } )$ consisting of section $f$ such that the lift of $f $ to $D ^\times (\A _{F, f}  ) $ factors
through $D ^\times (\A _{F, f}  ) \overset {\N_{D/F} } \to  (\A_{F , f } )
^\times $.
$ I_{(k, w), \bar E_\lambda } ^K$ is isomorphic to $\bigoplus _{ \pi
\in
\mathcal A  ^{c} _{(k, w)} (G_D ( \A _{\mathbb Q}) )}  {\pi_f 
}^K $.\par
In the both cases (especially when $q_D=0$), the Galois representation $ \rho _{\pi, \bar E_\lambda }$ attached to $\pi \in \mathcal A _{(k, w)} (G_D ( \mathbb A  _{\mathbb Q} )  )$ exists, and is isomorphic to $ \rho _{\JL (\pi), \bar E_\lambda }$ attached to the cuspidal representation $\JL (\pi ) $ of $\GL _2 ( \mathbb A _F) $.\par

\begin{rem}\label{rem-shimura41}
There is an action of the field automorphisms $\Aut \bar E_{\lambda} $ on $\mathcal A _{ (k, w)}(G_D)$: for a representation $\pi \in \mathcal A _{ (k, w)}(G_D)  $ and any element $\tau \in \Aut \bar E_{\lambda} $, let $\pi ^\tau_f $ be the twist of $\pi _f $ by $\tau$: $\pi ^\tau_f : G_D ( \mathbb A  _{\mathbb Q , f} )\overset {\pi _f} \to \Aut _{\bar E_\lambda } V \simeq \Aut _{\bar E_\lambda } V \otimes _{\bar E_\lambda, \tau } \bar E_\lambda$. Then there is a representation $\pi ^\tau \in \mathcal A _{(k, w)}(G_D) $
such that $( \pi ^\tau ) _f  \simeq \pi ^{\tau}_f $. $ \pi ^\tau $ is unique up to isomorphisms. \par
 Since the $\Aut \bar E_{\lambda} $-orbit of $\pi$ is a finite set, for the stabilizer $H_{\pi} $ of the isomorphism class of $\pi$,
$$
E _{\pi } = \bar E_{\lambda}  ^ {H_{\pi}} 
$$
is a number field of finite degree over $\mathbb Q$, which we call the field of
definition of $\pi$. When a subfield $ E$ of $\bar E_{\lambda}  $ contains $ E_{\pi}$, we say that $\pi$ is defined over $E$.
By the strong multiplicity one theorem for $\GL_2 $, $ \pi ^\tau \simeq \pi $ if and only if $ \pi _v ^{\tau} \simeq \pi_v $ for almost all $v$. In particular $\pi _v$ is defined over $E$ if and only if $\pi_v $ is defined over $E$ for almost all $v$. \par
\end{rem}

\subsection{Smallness of compact open subgroups}\label{subsec-shimura5}
\begin{prop} \label{prop-shimura51} Let $K$ be an $F$-facrtorizable compact open subgroup of $ G_D ( \mathbb A _{\mathbb Q, f} ) $, $y$ a finite place of $F$, $S$ a finite set of finite places of $F$. Assume the following conditions: 
\begin{enumerate}
\item $y \not \in S$, $ \Sigma _K \cap ( \{ y\} \cup S) = \emptyset $, and any element in $S$ and $y$ do not divide $2$. 
\item The map $ o^{\times } _F / ( o^{\times } _F ) ^2 \to \prod _{ v \in S} k (v) / (k(v) ) ^2$ is injective. 
\item $G ^{\der} _D ( \mathbb Q )  \cap g^{-1} (K_{11} (m_y) \cdot K ^y  )g  $ is torsion free for any $ g \in G_D  ( \mathbb A _{\mathbb Q, f})$.
\end{enumerate}
For $u \in S$, define a compact open subgroup $U_u $ of $ D^{\times } (F_u ) $ by
$$
U _u = \{ \alpha \in  K _u, \Norm _{ D_{F_u } / F_u }  ( \alpha ) \in (o ^{\times}_{F_u }) ^2  \} .
$$
Then $K ( y, S ) = ( K _{11} ( m_y )\cap K_y) \cdot \prod _{u \in S} U_u  \cdot  K ^{ \{ y\} \cup S }  $ is a small subgroup of $G_D (\mathbb A_{\mathbb Q, f} )  $.
 \end{prop}
To show the proposition, it suffices to show the following lemma. 
\begin{lem}\label{lem-shimura50} For an open subgroup $K'$ of $K(y, S )$ and an element $k $ of $K ( y, S)  $, assume that $k $ normalizes $K'$. If the right action of $k$ on $ S_{K'} $ admits a fixed point, there is a unit $ \delta \in {o_F}  ^{\times}  $ such that $\delta ^{-1} k \in K ' $. 
\end{lem}
\begin{proof}[Proof of Lemma \ref{lem-shimura50}]
We may assume that any local component $K_v$ of $K$ at a finite place $v$ is a maximal compact open subgroup of $ D^ {\times}(F_v)$, and $K _y = K_{11} (m_y )$.\par 
Assume that $k$ fixes the double coset $G _D (\mathbb Q  ) x (K ' \cdot K_{\infty}) $ defined by $x \in G _D (\mathbb A _{\mathbb Q}  ) $. 
There are elements $ \gamma \in G _D (\mathbb Q  ) $ and $k ' \in K '  $ such that 
$$
x\cdot   k = \gamma \cdot x \cdot k '
\leqno{(\ast)}
$$
holds. 
By taking the reduced norm of ($\ast$),  $\epsilon = \Norm _{D/F} ( k ')  \cdot \Norm _{D/F}( k )^{-1} $ for $\epsilon = \Norm_{D/F}( \gamma )$. $\epsilon $ is a unit of $F $ since it belongs to $\widehat {o^{\times}_F} $. At $ u \in S$, $ \epsilon $ belongs to $(o ^{\times}_{F_u }) ^2 $ by the definition of $K ( y, S ) $. By Proposition \ref{prop-shimura51} (2), this implies that $ \epsilon \in  ( o^{\times } _F ) ^2$. We take $ \delta \in  o^{\times } _F $ such that $ \delta ^2 = \epsilon$, and $ \delta \mod m_{F_y } = 1$. $\delta \in \tilde K =K ( y, S )  $, and $\tilde \gamma = \gamma / \delta $ belongs to $ G_D ^{\der} (\mathbb Q  )  $. \par
Equation ($\ast$) reduces to 
$$
x\cdot   \tilde k = \tilde \gamma \cdot x ,
\leqno{(\dagger)}
$$
where $ \tilde k = k / (\delta \cdot k') \in \tilde K  $. \par
When $ q_D = 0$, ($\dagger$) implies that $\tilde \gamma $ fixes a $\mathbb Z$-lattice in some faithful $\mathbb Q$-representation of $G_D ^{\der}  $ since it is contained in $x \tilde K  x ^{-1} $. Thus $\tilde \gamma $ 
 is an element of finite order, because it is contained in the compact group $ G_D ^{\der} (\mathbb R ) $. By our assumption  \ref{prop-shimura51} (3), $\tilde \gamma = 1 $, and $ \tilde k = 1 $ follows from ($\dagger$). This implies the freeness of the action on $S_{K '}$. \par
When $q_D = 1$,  ($\dagger$) implies that the finite part $x_f $ of $x$ satisfies $ \tilde \gamma \in x_f \tilde K x _f^{-1} $, and the class of the infinite part $x_\infty K _{\infty} \in X_D $ is fixed by $\tilde \gamma$ from the left. By Proposition \ref{prop-shimura51} (3), the action of $G ^{\der} _D ( \mathbb Q )  \cap x _f  (K_{11} (m_y) \cdot K ^y )x^{-1}_f   $ on $ X_D$ is free, which implies that $\tilde \gamma = 1 $. It follows that $ \tilde k= 1 $, and the claim is shown. 
\end{proof}
As for the existence of a nice pair $(y, S)$ as in Proposition \ref{prop-shimura51}, we have the following lemma. 
\begin{lem} \label{lem-shimura51} Let $K$ be an $F$-factorizable compact open subgroup of $ G_D ( \mathbb A _{\mathbb Q, f} ) $.
\begin{enumerate}
\item  For any finite set $P_0 $ of finite places of $F$ which contains $ \Sigma _K \cup \{ u; u \vert 2 \}  $, there is a finite set $S$ of finite places that is disjoint from $ P_0 $, and the condition (2) of Proposition \ref{prop-shimura51} is satisfied for $S$. 
\item There is an integer $a _F \geq 1 $, which depends only on $F$, with the following property: for any finite place $y\not \in  \Sigma _K  $ and $q_y \geq  a_F$, the condition (3) of \ref{prop-shimura51} is satisfied for $y$. 
\end{enumerate}
\end{lem}
\begin{proof}[Proof of Lemma \ref{lem-shimura51}]
To construct $S$, let $ F' $ be the Galois extension of $F$ defined by $F' = F ( \sqrt \epsilon, \ \epsilon \in o^{\times}_F  )$. The Galois group $G' = \Gal (F' /F)$ is an abelian group of type $(2, \ldots, 2)$, which we view as an $\mathbb F_2 $-vector space. Let $\{ \sigma _j \}_{j \in J } $ be a basis of $G'$ over $\mathbb F_2$. 
We take $S$ so that the following two conditions are satisfied:
\begin{itemize}
\item $S \cap  P_0 =  \emptyset  $, 
\item For any $j\in J$, there is an element $s_j $ in $S$ such that $F' $ is unramified at $s_j$, and the geometric Frobenius element $\Fr _{s_j} $ at $s_j$ is mapped to $\sigma _j$. 
\end{itemize}
The existence of $S$ is guaranteed by the Chebotarev density theorem, and Proposition \ref{prop-shimura51} (2) is satisfied for this choice of $S$. \par
For (2), the existence of $a_F$ is proved in \cite{H1}, lemma 7.1 (we take $a_F$ so that $a_F > 2 ^d$, which excludes the possibility that $y $ divides $2$). 
\end{proof} 

\begin{rem}\label{rem-shimura51}The notations are as in Lemma \ref{lem-shimura51}. Let $\mathcal E $ be the kernel of $\Norm_{F/\mathbb Q} : o ^ \times_F \to \mathbb Z ^{\times} $. If we add a finite place $s$, which is split in $ F( \sqrt {\epsilon } , \epsilon \in \mathcal E )$ but does not split in $ F(\sqrt{-1} ) $, to $S$, $\mathcal F ^{K'} _{(k, w ) , E_\lambda } $ is defined on $S_{K'} $ for any compact open subgroup $K '$ of $K (y, S)$. 
\end{rem}

For an $F$-factorizable compact open subgroup $K$ of $G_D ( \mathbb A _{\mathbb Q, f} ) $, assume that $y$ satisfies the condition (1) of Proposition \ref{prop-shimura51}, and the $y$-component $K_y$ of $ K$ is contained in $K_{11} ( m_y )  $. By Lemma \ref{lem-shimura51}, we take an auxiliary set $S$ which is disjoint from $\Sigma _K \cup \{ u : u \vert 2\}$ and satisfies the condition (2) of Proposition \ref{prop-shimura51} (cf. Remark \ref{rem-shimura51}).\par
When $\ell\geq 3$, for any integer $ q $ and a finite $o _{\lambda}$-algebra $R$, the cohomology group $ H ^{q} _{\stack}( S_K \ \bar { \mathcal F} ^K _{(k, w) , R}) $ of degree $q $ is defined as the $K$-invariant part $ H ^{q} ( S_{K(y, S ) }  \  \bar {\mathcal F } ^{K(y, S )} _{(k, w) , R}) ^{K }$ of $ H ^{q} ( S_{K(y, S ) }  \  \bar {\mathcal F } ^{K(y, S )} _{(k, w) , R})$. 

Since $ K(y, S )\backslash K  $ is an abelian group of type $(2, \ldots, 2)$ and $\ell \geq 3$, $ H ^{q} _{\stack}( S_K \  \bar {\mathcal F} ^K _{(k, w) , R}) $ is an $R$-direct summand of $ H ^{q} ( S_{K(y, S ) }  \  \bar {\mathcal F} ^{K(y, S )} _{(k, w) , R}) $. This property is usually suffcient to make an analysis on $S_{K(y, S ) }$. In particular exact sequences are preserved.  

\begin{rem}\label{rem-shimura52} As the notation suggests, $H ^{q} _{\stack} (S_K)  $ is canonically isomorphic to the (\'etale) cohomology group of $S_K$ when $S_K$ is regarded as a Deligne-Mumford stack, and is independent of the choice of an auxiliary set $S$.  
\end{rem}

\section{Universal exactness of cohomology sequences} 
\label{sec-coh}
We study the exactness of a homomorphism defined by degeneracy maps on cohomology groups. In particular we show the sequences in consideration are universally exact, that is, the sequences remain exact under any extensions of scalars. This is a basic tool in the study of cohomological congruences, in particular in the theory of congruence modules.  \par

In the elliptic modular case and the subgroup is $ K_0 (v) $, Ribet calls the universal injectivity ``Ihara's Lemma''.  \par
Throughout this section, $D$ is a division algebra unless otherwise stated. 

\subsection {Modules of residual type}\label{subsec-coh1}
Let $K $ be an $F$-factorizable compact open subgroup of $ G_D ( \mathbb A_{\mathbb Q, f} ) $. For a finite set of finite places $\Sigma $ which contains $ \Sigma _K $ and the places dividing $\ell$, let $T_\Sigma =  H (  D ^\times ( \A ^{\Sigma } _f ),\
K^\Sigma )_{o_{E_\lambda} }$ be the convolution algebra consting of $o_{E_\lambda }$-valued functions. By the assumption on $\Sigma$, $T_{\Sigma } $ is commutative. \par
\begin{dfn}\label{dfn-coh11} Consider the category $\mathcal C _{T_\Sigma}$ of
$T_\Sigma$-modules which are finitely generated as $o_{E_\lambda}$-modules. 
\begin{enumerate}
\item  We call an object $N $ in $\mathcal C _{T_\Sigma}$ residual type
if any constituant $N'$ of $ N/ \lambda N $ satisfies the relation
$$ 
[T_v ]^2 = [T_{v,v} ] (1 + q_v)^2
$$ 
on $N '$ for any $v\not \in \Sigma$.

\item A maximal ideal $m$ of $T_\Sigma$ is of residual type if $T_\Sigma /m$ is of residual type.
\end{enumerate}

 By
$\mathcal C _{\residue}
$ we denote the subcategory of $ \mathcal C_{T_\Sigma} $ consisting of the $T_\Sigma $-modules
of residual type.
\end{dfn}
$\mathcal C _{\residue}$ is a Serre subcategory of $\mathcal C _{T_\Sigma} $, and is
stable under the dual action of $T_\Sigma$. By
$
\bar {\mathcal C} _{T _\Sigma }$, we mean the quotient category of $\mathcal C _{T_\Sigma}$ by $\mathcal C _{\residue}
$.\par
\medskip
A typical example of modules of residual type is obtained by a one dimensional
representation $\pi : D ^\times ( \A ^{\Sigma  }_f ) \to E ^\times _{\lambda}  $ which factors through the reduced norm $\Norm _{D/F} $.
The induced $ T_\Sigma$-action gives a module of residual type.  \par

For a Galois representation $\bar \rho :G_{\Sigma} \to \GL_2 (k_{\lambda} )  $,  
the maximal ideal $m _{\bar \rho } $ of $T_\Sigma$ associated to $\bar \rho $ is the kernel of $o_{E_\lambda} $-algebra homomorphism $ f_{\bar \rho} : T _{\Sigma } \to k _{\lambda }$ such that $  f_{\bar \rho} (T _v ) = \tr \bar \rho (\Fr _v) $, $ f_{\bar \rho}(T_{v, v} ) =  q_v ^{-1} \det \bar \rho ( \Fr _v) $. The following proposition shows that the maximal ideals of $T_\Sigma$ of residual type coming from Galois representations correspond to very special reducible representations. 

\begin{prop}\label{prop-coh11}
For a continuous representation $\bar  \rho : G_\Sigma \to \GL_2 (\bar k_\lambda
)$, the maximal ideal $m_{\bar \rho}$ corresponding to $\bar \rho$ is of residual type if and only if $\bar \rho ^{\ss}$ satisfies 
$$
\bar \rho ^{\ss} \simeq \bar \chi \oplus \bar
\chi (-1)
$$ 
for some one dimensional character $\bar \chi: G_\Sigma \to \bar k_\lambda ^\times$ over $\bar k _{\lambda}$. 
\end{prop}
A proof is found in \cite{Fu1}, Proposition 3.6.

\medskip
\begin{rem} \label{rem-coh1} 
\begin{enumerate}
\item The notion of modules of residual type is stronger than the
notion of Eisenstein modules in
\cite{DT1}, and is introduced in \cite{Fu1} (there it is called of $\omega$-type).  
\item Let $m$ be a maximal ideal of $T_\Sigma $ which is not of residual type. If a sequence 
$$
0 \longrightarrow M_1 \longrightarrow M_2 \longrightarrow M_3 \longrightarrow 0  
$$ of $T_\Sigma$-modules is exact in $\bar  {\mathcal C}  _{T_\Sigma}$, then the sequence 
$$
0 \longrightarrow ( M_1)_m  \longrightarrow (M_2 )_m \longrightarrow ( M_3 )_m\longrightarrow 0 
 $$
localized at $m$ is exact. 
\end{enumerate}
\end{rem}
\medskip
\subsection{Universality under scalar extensions}\label{subsec-coh2}
For a finite $o_{E_\lambda}$-algebra $R$, we denote $\bar {\mathcal F } ^K _{(k,w) , o_{E_\lambda } }\otimes _{o_{E_\lambda } } R  $ by $ \bar {\mathcal F } ^K _{(k,w) ,R }$.
\begin{lem}\label{lem-coh21} Assume that $D$ is a division algebra which defines a Shimura curve, $K$ is an $F$-factorizable compact open subgroup of $ D ^{\times } (\mathbb A _{F, f} ) $, and $\Sigma $ is a finite set of finite places of $F$ that contains $\Sigma _K \cup \{ v \vert \ell \} $. \par
For any integer $\alpha \geq 1 $ and any $D ^{\times}( \mathbb A ^{\Sigma}_{F, f } )$-equivariant $k_\lambda$ smooth subquotient $\mathcal F$ of $(\bar {\mathcal F } ^K _{(k,w) ,k_\lambda } ) ^{\oplus \alpha }$, $H ^q (S_K , \mathcal F) $ is of residual type for $q= 0, 2$. 
\end{lem}
\begin{proof}[Proof of Lemma \ref{lem-coh21}] First we prove the claim for $q=0 $. For an inclusion of two groups $K ' \hookrightarrow K  $, $ H ^0 (S_K , \mathcal F  ) $ is a subspace of $H ^0 ( S_{K'} , \mathcal F_{K'}) $ since $\pi _{K ' , K }$ is surjective. Here $  \mathcal F_{K'}$ is the pullback of $\mathcal F $ to $ S_{K '}$. So we may assume that $K _u \subset K( m_{F_u}) $ for $u \vert \ell$. By the definition of $\bar {\mathcal F} _{(k, w), o_{E_\lambda}}  $, $\bar {\mathcal F} _{(k, w), k_\lambda }
$ is trivialized in
$D ^{\times}( \mathbb A ^{\Sigma}_{F, f } )$-equivariant way on $S_{K }$. This implies that any subquotient of $ \mathcal F $ is isomorphic to $k_\lambda $ as a $D ^{\times}( \mathbb A ^{\Sigma}_{F, f } )$-equivariant sheaf. If $F= k_\lambda$, the claim follows from the isomorphism $ \pi _0 ( S_K ) \simeq \pi _0 ( F^ \times \backslash \mathbb A^{\times}_{F, f} / \Norm_{D/F} K )$, by taking a finite set of finite places $ \Sigma$ which contains $\Sigma _K$, and all places dividing $\ell$. By the induction on the $k_\lambda$-rank of $F$, the claim is shown in this case. \par
For $q=2$, this follows from the case of $q=0$ by Poincar\'e duality.
\end{proof}

\begin{prop} \label{prop-coh21} Assume that $D$ is a division algebra with $q_D \leq 1$. For any finite $o_{E_\lambda }$-algebra $R$, $H ^ {q_D} ( S_K,  \bar {\mathcal F } ^ K _{(k, w ) , o_{E_\lambda} } )$ is $\lambda $-torsion free in $ \bar {\mathcal C }_{T_\Sigma} $, and 
$$
H ^ {q_D}  ( S_K,  \bar {\mathcal F } ^ K _{(k, w ) , R} ) = H ^ {q_D}  ( S_K,  \bar {\mathcal F } ^ K _{(k, w ) , o_{E_\lambda} } )  \otimes _{o_{E_\lambda} } R
$$
holds in $\bar {\mathcal C }_{T_\Sigma} $. 
\end{prop}
\begin{proof} [Proof of Proposition \ref{prop-coh21}] The claim is clear for definite quaternion algebras. So we may assume that $q_D=1$, and $D$ defines a Shimura curve. In the following exact sequence
$$ 
H ^ 0 ( S_K,\  \bar {\mathcal F } ^ K _{(k, w ) , k_{\lambda} } ) \longrightarrow
H ^ 1 ( S_K,\  \bar {\mathcal F } ^ K _{(k, w ) , o_{E_\lambda} } ) \overset {\lambda}
\longrightarrow  H ^ 1( S_K,\  \bar {\mathcal F } ^ K _{(k, w ) , o_{E_\lambda} } )
$$
$$ 
\longrightarrow H ^ 1 ( S_K,\  \bar {\mathcal F } ^ K _{(k, w ) , k_{\lambda} } )\longrightarrow  H ^ 2( S_K,\  \bar {\mathcal F } ^ K _{(k, w ) , k_{\lambda} } ), 
$$ 
$ H ^ q ( S_K,\  \bar {\mathcal F } ^ K _{(k, w ) , k_{\lambda} } )$ for $q= 0 , 2$ are of residual type by Lemma \ref{lem-coh21}.
Thus $H ^ 1 ( S_K,\  \bar {\mathcal F } ^ K _{(k, w ) , o_{E_\lambda} } ) $ is $\lambda$-torsion free, and $H ^ 1 ( S_K,\  \bar {\mathcal F } ^ K _{(k, w ) , o_{E_\lambda} } )\otimes _{o_{E_\lambda }} k _\lambda =   H ^ 1 ( S_K,\  \bar {\mathcal F } ^ K _{(k, w ) , k_{\lambda} } )$ in $\bar {\mathcal C }_{T_\Sigma} $. For general coefficient $R$, we reduce to the case when $R$ has a finite length, and in that case it follows from an induction on $\length R$. 
\end{proof}

\begin{prop} \label{prop-coh22}
For a finite complex $ L = (L ^{\cdot}, d ^{\cdot} ) $ of finite free $o_{E_\lambda}$-modules, define $L_R$ by $L_R = ( L ^{\cdot} \otimes _{o_{E_\lambda}}R, d ^{\cdot}   \otimes \id_R) $ for any finite $o_{E_\lambda}$-algebra $R$.  
Then the following properties are equivalent:
\begin{enumerate}
\item For any finite $o_{E_\lambda}$-algebra $R$, $L_R $ is exact. 
\item $ L _{k_\lambda } $ is exact. 
\item $L  $ is exact, and the image of $ d ^i : L ^i \to L^{i+1}$ in $L ^{i+1}$ is an $o_{E_\lambda } $-direct summand for any $i$. 
\end{enumerate}
\end{prop} 
The verification is left to the reader. 
\begin{dfn}  \label{dfn-coh21} Let $L$ be a finite complex of finite free $o_{E_\lambda}$-modules. If $L$ satisfies the equivalent conditions in Proposition \ref{prop-coh22}, we say $L$ is universally exact. For an $o_{E_\lambda}$-homomorphism  $f : L _0 \to L _1$ between finite free $o_{E_\lambda }$-modules, $f$ is called universally injective if $[ 0 \to L_0 \overset {f} \to L_1  ]$ is universally exact. 
\end{dfn}
\subsection {Cohomological universal injectivity}
\label{subsec-coh3}  

\begin{prop}\label{prop-coh31}
Let $D$ be a division quaternion algebra with $q_D \leq 1$. For an $F$-factorizable small compact open subgroup
$K$ of $G_D ( \mathbb A_{\mathbb Q, f} ) $ and a compact open subgroup $K '$ of $K$, 
$$
H ^{q_D} ( S_K ,\bar { \mathcal F } ^K _{(k , w ) , o_\lambda }) \overset{\pi^* _{K ' , K }} \longrightarrow H ^{q_D} ( S_{K'} ,\bar { \mathcal F } ^{K'} _{(k , w ) , o_\lambda }) 
$$
is universally injective up to modules of residual type.
\end{prop}
\begin{proof}[Proof of Proposition \ref{prop-coh31}] By Proposition \ref{prop-coh21} and \ref{prop-coh22}, it is sufficient to show the kernel of 
$$
H ^{q_D} ( S_K ,\bar { \mathcal F } ^K _{(k , w ) , k_\lambda }) \overset{\pi^* _{K ' , K }} \longrightarrow H ^{q_D} ( S_{K'} ,\bar { \mathcal F } ^{K'} _{(k , w ) , k_\lambda })
$$
is a module of residual type. If $q_D=0$, this is clear. So we may assume that $q_D = 1$. Moreover, by replacing $K'$ by an $F$-factorizable open subgroup, it is enough to consider the case when $K '$ is a normal subgroup of $K$. $ \pi^* _{K ' , K }$ is an \'etale torsor under $G = (\overline {F^ \times} \cap K ) \cdot  K ' \backslash K$, so we have the following exact sequence
$$
H ^1 ( G,  H ^0 (S_{K'} ,\bar { \mathcal F } ^{K'} _{(k , w ) , k_\lambda } ) ) \longrightarrow H ^1 ( S_K ,\bar { \mathcal F } ^K _{(k , w ) , k_\lambda }) \overset{\pi^* _{K ' , K }} \longrightarrow H ^1 ( S_{K'} ,\bar { \mathcal F } ^{K'} _{(k , w ), k_\lambda} ) .
$$ 
For $\Sigma = \Sigma _{K ' } \cup \{ v \vert \ell \}$, the action of $T_\Sigma  $ on $ N= H ^0 (S_{K'} ,\bar { \mathcal F } ^{K'} _{(k , w ) , k_\lambda } ) $ commutes with the $G$-action on $N$, and is of residual type by Lemma \ref{lem-coh21}. So $H ^1 ( G, N ) $ is also of residual type.  
\end{proof}

Let $D$ be a division quaternion algebra with $q_D \leq 1$. For an $F$-factorizable compact open subgroup $K$ of $G_D ( \mathbb A _{\mathbb Q, f} ) $, assume that $D$ is split at $v$, and the $v$-component  of $K$ is $\GL _2 (o_{F_v}) '= 
\ker ( \GL _2 ( o_{F_v})\overset { \det} \to o_{F_v} ^\times \to k(v ) ^\times)
$. \par
Let
$
\pr _{ i, v } : S_{ K _0 (v)
\cap K }\to S_K 
\ (i = 1, 2)
$ be two degeneracy maps defined as follows. $\pr_2= \pi _{K_0 (v) \cap K , K } $ is the canonical projection corresponding to 
the inclusion
$ K \cap K _0(  v ) \subset  K $, $ \pr _1= \pi_{a( p_v) ^{-1} (K_0 (v) \cap K) a( p_v)  , K } \circ R (a ( p_v ) ^{-1} ) $ is the projection
twisted by the conjugation by 
$ 
\begin{pmatrix}
 1 & 0\\ 
0 & p_v
\end{pmatrix}
 $
at $v$, that is, 
$$ 
K_v \cap K_0 ( m _v ) \longrightarrow K_v \cap \begin{pmatrix}
 1 & 0\\ 
0 & p_v
\end{pmatrix}
^{-1}
K_0 ( m_{F_v}   )
\begin{pmatrix}

1 & 0
\\ 0 & p_v 
\end{pmatrix}
 = K _v  \cap K ^ {\op} ( v )   \subset K_v. 
$$
Here $ K ^ {\op} ( v )= K ^{\op } (m_{F_v}  ) \cdot K ^v $, and $ K ^ {\op} ( m_{F_v} ) = \{  ^t g ; \ g \in K_0 (m_{F_v}) \} $.\par
By the definition, $ S_ K\overset {\pr _1} \leftarrow S_{K\cap K_0 (v) }\overset {\pr_ 2} \rightarrow S_K $ is equal to $[K a(p_v ) K] $ as a correspondence.\par

Consider the map 
$$
H ^{q_D} ( S _{K} ,\ \bar {\mathcal F} _{(k, w), o_{E_\lambda}}) ^{\oplus 2}   \overset { \pr ^* _1 +\pr^ * _2} \longrightarrow
H ^{q_D} ( S _{K_0 (v)
\cap  K } ,\
\bar {\mathcal F} _{(k, w), o_{E_\lambda}}). 
$$
Later in \S\ref{sec-congruent}, we need a universal injectivity (Ihara's Lemma) in the calculation of congruence modules. We state it as a hypothesis here. 
\begin{hyp}[Cohomological universal injectivity]\label{hyp-coh31}
Let $D$ be
a quaternion algebra over $F$ with $q_D \leq 1$ which is split at
a finite place $v$. Assume that $ K$ is an $F$-factorizable compact open subgroup of $G_D ( \mathbb A _{\mathbb Q , f} )  $, with the $v$-component $ K_v = \GL_2 (o_{F_v}) '$. For a discrete infinite type $ (k, w)$, if $ v \vert \ell$, we further assume that $k_\iota  = 2 $ and $w=0$ for any $\iota \in I_{F, v } $. \par
Then
$$
H ^{q_D} ( S_K,\  \bar {\mathcal F} _{(k, w), o_{E_\lambda}}   ) ^{\oplus 2} \overset { \pr ^* _1 +\pr^ * _2} \longrightarrow
H ^{q_D} (
S_{K
\cap K _0 (v) },\
\bar {\mathcal F} _{(k, w), o_{E_\lambda}} ) 
$$
is universally injective. 
\end{hyp}
As for the validity of Hypothesis \ref{hyp-coh31}, we prove the following theorem. 
\begin{thm} \label{thm-coh31} Hypothesis \ref{hyp-coh31} is true when $q_D = 0 $. 

\end{thm}

Before beginning the proof of Theorem \ref{thm-coh31}, let us recall the following general fact. 
\begin{lem}
\label{lem-coh31}
Assume that $ f_1: X  \to
 X _1 $, $f_2 : X \to X_2$ are two surjective maps between finite sets. 
\begin{enumerate}
\item 
$$ 
0 \to \ker (  f_1 ^* + f_2 ^* ) \longrightarrow H ^0 ( X_1 ,\  o_{E_\lambda} )\oplus H ^0 ( X_2 ,\  o_{E_\lambda} )
\overset { f_1 ^* + f_2 ^*  } \longrightarrow  H ^0 (X,\   o_{E_\lambda}  )
 $$ is universally exact. 
\item $\ker (  f_1 ^* + f_2 ^* )$ is identified with $H ^0 ( X_1 \amalg_X X_2 ,\  o_{E_\lambda} ) $. Here $X_1 \amalg _X X_2$ is the coproduct with respect to $f_1$ and $f_2$. 
\end{enumerate}
\end{lem}
A proof of (1) is found in
\cite{T1}, lemma 4, case 1. We give a proof based on the theorem of van
Kampen. \par
\begin{proof}[Proof of Lemma \ref{lem-coh31}]
For any
$o_{E_\lambda} $-algebra $R$, define a complex $L_R$ of $R$-modules by
$$
L_R = [H ^0 (X_1 ,\ R )\oplus H ^0 (X_2 ,\ R ) \overset { f_1 ^* +f_2 ^*  } \longrightarrow
H ^0 (X,\ R) ] .
$$
Here the components are placed at degree $0$ and $1$. The formation $ R \mapsto L_R$ commutes with any extension of scalars $R\to R'$: $L_R \otimes _R R' = L_{R'} $ holds. \par
It is easily seen that
$H ^0 (L_R)$ is identified with the cohomology group of the coproduct
$ H ^0 (X_1 \amalg _X X_2,\ R )
$. The identification is 
$$ 
 H ^0 (X_1 \amalg _X X_2 ,\ R ) \overset { (g ^*_1,  - g ^* _2)  }\longrightarrow  H ^0 ( X ,\
R)^{\oplus 2} 
$$ 
for the maps defined by the following commutative diagram.
$$
\CD 
X @>{f_1}>> X_1\\ 
@V{f_2}VV @V{g_1}VV\\ 
X_2 @>{g_2}>> X_1\amalg_X X_2\\
\endCD
$$

It follows that the formation $R \mapsto  H ^ 0 (L_R )$ also commutes with any extension of scalars.
So $H ^1 (L_R)$ must satisfy the same property since $L_{o_{E_\lambda}}$ is a perfect
complex of $o_{E_\lambda} $-modules. This implies that $H ^1 (L_{o_{E_\lambda}})$ is locally free.
\end{proof}

\begin{proof}[Proof of Theorem \ref{thm-coh31}] Define $\Sigma ' $ by $ \Sigma '  = \Sigma _K \cup \{ v\} \cup \{ u : u \vert \ell \}$. 
First we show \ref{thm-coh31} when the infinity type $(k, w)$ is $( (2, \ldots, 2), 0 )$. \par 
Let $T   $ be the $\mathbb Q$-torus $ \Res _{F/ \mathbb Q} \mathbb G _{m, F} $, $U _{\infty } = T ( \mathbb R ) $. For a compact open subgroup $ U$ of $ T ( \mathbb A _{\mathbb Q, f} )$, define $ S_U (T)$ by 
$$
S_U  ( T ) = T (\mathbb Q  ) \backslash T ( \mathbb A _{\mathbb Q}) / U \times U _{\infty } .
$$
For any compact open subgroup $\tilde  K $ of $G_D ( \mathbb A _{\mathbb Q , f} )  $, the reduced norm $ \Norm _{ D/F} $ induces a surjective map 
$$
\alpha _{\tilde K} : S _{\tilde K } \overset { \Norm _{ D/F}} \longrightarrow S _{\Norm _{ D/F} (\tilde  K ) } (T),
$$
which we call the augmentation map. \par
Since the augmentation maps are functorial with respect to $\tilde K$, and $ \Norm _{ D/F}( K ) =  \Norm _{ D/F} ( K \cap K_0 (v) ) $, $\alpha _{K \cap K_0 ( v) } \cdot  \pr _i =  \alpha _K$ for $ i = 1, 2$, which induces a surjective map 
$$
\beta : S_{K} \amalg
_{S_{K  \cap  K_0 (v)  }}  S_{K } \longrightarrow  S _{\Norm _{ D/F} ( K ) } (T).
\leqno{(\ast)}
$$

We show that $\beta$ is an isomorphism. In the decomposition
$$
H ^0 ( S_K, \mathbb C )\simeq I ^K_{(k, w ) , \mathbb C}  \oplus  (\bigoplus _{ \pi
\in
\mathcal A  _{(k, w)} (G_D  ( \A _{\mathbb Q} ) ) }  \pi  _f ^K  ), 
$$
in \S \ref{subsec-shimura4}, the map 
$$
H ^0 ( S_K,\ \C ) ^{\oplus 2}  \overset { \pr ^* _1 +\pr^ * _2} \longrightarrow H ^0 ( S_{K \cap K _0 (v) },\ \C  )
$$
is injective on the part parametrized by $ \mathcal A    _{(k, w)} (G_D  ( \A _{\mathbb Q}))  $. 
By definition, $  I ^K_{(k, w ) , \mathbb C}$ is identified with $ H ^0 ( S _{\Norm _{ D/F} ( K ) } (T) , \mathbb C )$, and is regarded as the kernel $V$ of ($\ast$) via the embedding 
$$  
I ^K_{(k, w ) , \mathbb C}\overset{( i_K , - i _K)} \longrightarrow ( I ^K_{(k, w ) , \mathbb C} ) ^{\oplus 2} \hookrightarrow H ^0 ( S_K,\ \C )  ^{\oplus 2}, 
$$
where $i _K : I ^K_{(k, w ) , \mathbb C}\hookrightarrow H ^0 ( S_K,\ \C ) $ is the inclusion. 
On the other hand, by Lemma \ref{lem-coh31} (2), $V$ is canonically isomorphic to
$ H ^ 0 ( S_{K \cap  K_0 (v) } \coprod
_{S_{K  }}  S_{K
\cap  K_0 (v) } ,\ \C) $. So the source and the target of ($\ast$) have the same cardinalities, and
$\beta$ must be an isomorphism.\par
 
 We apply Lemma \ref{lem-coh31} to $\pr _1, \ \pr _2 : S_{K \cap K_0 (v) }  \to S_K $. 
It follows that the kernel of 
 $$
H ^0 ( S_K,\  k_\lambda  ) ^{\oplus 2} \overset {
\pr ^* _1 +\pr^ * _2} \longrightarrow
H ^0 (
S_{K
\cap K _0 (v) },\
 k_\lambda ) 
$$
is canonically isomorphic to $H ^ 0 (S _{\Norm _{ D/F} ( K ) } (T) ,\ k_\lambda ) $. In particular the action of $  T_{\Sigma '}  $ is of residual type, since the constituants as $D ^{\times } (\A ^{\Sigma '} _f  ) 
$-representations are all one dimensional over $ \bar k_{\lambda}$. \par
We prove the general case. By the definition of $\bar {\mathcal F} _{(k, w), o_{E_\lambda}}  $, $\bar {\mathcal F} _{(k, w), k_\lambda }
$ is trivialized in
$G_D (\A ^{ \ell } _{\mathbb Q , f} ) $-equivariant way on $S_{K '}$. Here $K ' = K '_{\ell}\cdot  K^\ell $, $K' _\ell= \prod _{u\vert \ell} K(m_{F_u})\cap K _u$ if $v \nmid \ell$, $ K ' _ \ell = K_v \cdot \prod _{u \vert \ell, u \neq v}K (m_{F_u} )\cap K_u$ if $v \vert \ell$. 
In the following commutative diagram
$$
\CD
H ^0 ( S_K,\ \bar {\mathcal F}  ^K_{(k, w), k_\lambda }   ) ^{\oplus 2} @>{ \pr ^* _1 +\pr^ * _2} >>
 H ^0 (S_{K \cap K _0 (v) },\
\bar {\mathcal F} ^{K\cap K _0 (v)}_{(k, w), k_\lambda}  )\\
@VVV @VVV\\
H ^0 ( S_{K'},\ \bar {\mathcal F} ^{K'} _{(k, w), k_\lambda } )^{\oplus 2} @>{ \pr ^* _1 +\pr^ * _2} >>
 H ^0 (
S_{K'\cap K _0 (v) },\
\bar {\mathcal F} ^{K'\cap K _0 (v) }_{(k, w), k_\lambda } ), \\
\endCD
$$
the restriction maps in the vertical arrows are injective. Since the kernel is of residual type in the case of constant coefficient $k_\lambda$ and $ \bar {\mathcal F} ^{K'} _{(k, w), k_{\lambda}} $ is trivialized, the general case follows.
\end{proof}
\begin{rem}\label{rem-coh31}
In \cite{DT1}, the universal injectivity (Ihara's Lemma) is proved for division quaternion algebra $D$ over $\mathbb Q$ in the following cases: $D$ is definite, with
``Eisenstein'' instead of  ``residual type''. When $D$ is indefinite, there are some restrictions on $\ell$ and $ ( k, w)$ to apply the $p$-adic Hodge theory. \par
The universal injectivity when $q _D= 1 $ will be treated, by assuming that there is a finite place $u\vert \ell $ where $D$ is split. 
\end{rem}

\subsection{Cohomological universal exactness}
\label{subsec-coh4}  

In this paragraph we show a cohomological universal exactness from $ K\cap K_{11}( v ^
{n-1} ) 
$ to $ K\cap K_{11 }( v ^
{n} ) $ for $n \geq 1$. This case turns out to be easier than the $K_0 (v)$-case. \par
\bigskip
Assume that $D$ is
a quaternion algebra over $F$ which defines either a Shimura curve or a Hida variety, and is split at
a finite place $v$. $ K$ is an $F$-factorizable compact open subgroup of $G_D ( \mathbb A _{\mathbb Q , f} )  $, with the $v$-component $ K_v = \GL_2 (o_{F_v}) '$.

Let
$
\pr _{ i, v } : S_{ K _{11} (v ^n )
\cap K }\to S_{K \cap K_{11} ( v^{n-1} ) } 
\ (i = 1, 2)
$ be two degeneracy maps defined as follows. $\pr_2 $ is the canonical projection corresponding to 
the inclusion
$ K \cap K _{11} (  v^n  ) \subset  K\cap K_{11} ( v ^{n-1} ) $, $ \pr _1$ is the projection
twisted by the conjugation by 
$ 
\begin{pmatrix}
 1 & 0\\ 
0 & p_v
\end{pmatrix}
 $
at $v$, that is, $  K_v \cap 
K_{11} ( m ^n _v ) \to K_v \cap \begin{pmatrix}
 1 & 0\\ 
0 & p_v
\end{pmatrix}
^{-1}
 K_{11} ( m^n _v )  
\begin{pmatrix}

1 & 0
\\ 0 & p_v 
\end{pmatrix}
  \subset K_v\cap K_{11} ( m_{F_v} ^{n-1} )$.\par

$\pr '_1 ,\ \pr '_2: S _{ K \cap K _{11}( v ^n ) \cap K_0 (v ^{ n +1} ) } \to S _{K
\cap  K _{11}(  v ^n) } $ are both defined in a similar way. We use the same  notations in
the case of
$K_{1}$, too.

\begin{prop}
\label{prop-coh41}
Assume that $D$ is a division quaternion algebra over $F$ which satisfies $q_D= \sharp I_D \leq 1 $, and is split at
a finite place $v$. $ K$ is an $F$-factorizable compact open small subgroup of $G_D ( \mathbb A _{\mathbb Q , f} )  $, with the $v$-component $ K_v = \GL_2 (o_{F_v}) '$. For a discrete infinite type $ (k, w)$, assume that $k_\iota  = 2 $ and $w=0$ for any $\iota \in I_{F, v } $ if $ v \vert \ell$.  For $n \geq 1 $, 
$$
 0 \to H ^{q_D} ( S _{K \cap  K _{11}( v ^{n-1}) } ,\  \bar {\mathcal F} _{(k, w), R} )\overset {(\pr_2 ^* , - \pr_1 ^*)}\longrightarrow H ^{q_D} (  S _{K
\cap  K _{11}(  v ^n) } ,\  \bar {\mathcal F} _{(k, w),R} ) ^{\oplus 2}
$$ 
$$
\overset { {\pr'_1} ^* + {\pr'_2 }^*}\longrightarrow H ^{q_D} (  S _{ K \cap K _{11}( v ^n ) \cap K_0
(v ^{ n +1} ) } ,\  \bar {\mathcal F} _{(k, w), R} )
$$
is exact for any finite $o_{E_\lambda} $-algebra $R$ up to residual type modules.
The same is true for $K_{1 }$.
\end{prop}

Proposition \ref{prop-coh41} is proved in
\cite{W2}, lemma 2.5 in the case of
$D=M_2 (\mathbb Q)$ and
$\mathcal F_{(k, w), o_{E_\lambda} } =
o_{E_\lambda} $, assuming
$v\nmid \ell$. 

\bigskip
\begin{proof}[Proof of Proposition \ref{prop-coh41}] We prove the proposition for $K_{11}$. $K_1$-case is proved in the same way. \par
 By Proposition \ref{prop-coh11}, the cohomology groups commutes with the extension of scalars up to
residual type modules. So we may assume that $ R = k_\lambda$. 
$$
\begin{pmatrix} 
1 & 0 \\
 0 &p_v
\end{pmatrix} ^{-1}
K_{11} ( m_{F_v} ^{n-1}) 
 \begin{pmatrix}
1  & 0
\\ 0 & p_v
\end{pmatrix} 
= \{g \in \GL _2 ( o_{F_v} ); g =
\begin{pmatrix} 
a & b \\ 
c & d
 \end{pmatrix}
,\ a , \ d \equiv 1 \mod m ^{n-1}_{F_v} ,\ b \in
m_{F_v},\ c \in m_{F_v} ^ {n -1}
\} 
$$
$$
=K _{11} ( m_{F_v} ^{n-1} ) \cap K ^{\op} (m_{F_v}) 
$$
and
$$ 
\begin{pmatrix}
 1 & 0
\\ 0 &p_v
\end{pmatrix}^{-1}
 ( K_{11} ( m_{F_v} ^ n ) \cap K_0(m^{n+1}_v) )
\begin{pmatrix}
1 & 0\\ 
0 &p_v
\end{pmatrix}
=K_{11} ( m_{F_v} ^ n ) \cap K(m_{F_v})  
$$ 
hold. So the commutative diagram 
$$
\CD 
S _{ K \cap K _{11}( v ^n ) \cap K_0 (v ^{ n +1} ) }@> {\pr' _1} >> S _{K
\cap  K _{11}(  v ^n) }\\ @V \pr'_2 VV @V \pr_1 VV\\ S _{K
\cap  K _{11}(  v ^n) } @> \pr_2>> S _{K \cap  K _{11}( v ^{n-1}) }
\endCD
\leqno{(\ast)}
$$
is factorized into
$$
\CD 
S _{ K \cap K _{11}( v ^n ) \cap K_0 (v ^{ n +1} ) }@> {R (a(p_v) ^{-1} ) } >> S _{K\cap  K _{11}(  v ^n)\cap K (v)  }\\ 
@V \pr '_2 VV @VVV\\ 
S _{K\cap  K _{11}(  v ^n) } @>  {R (a(p_v) ^{-1} ) }   >> S _{K \cap  K _{11}( v ^{n-1})\cap K ^{\op} (v)  }\\
@V \pr_2 VV @VVV\\ 
S_{K \cap K_{11} ( v ^{n-1} ) }  @>>> S_{K \cap K_{11} ( v ^{n-1} )} ,  \\
\endCD
$$
and ($\ast$) is isomorphic to 
$$
\CD 
Z= S_ {K_{11} ( v ^n ) \cap K(v )  \cap K } @>\alpha ' >>Y =S_ {(K_v \cap K_{11} (
m_{F_v} ^{n-1}  ) \cap K ^0 (m_{F_v}) ) \cdot K^v }\\ @V\beta ' VV @V\beta VV \\ Y' = S_{K_{11} (
v ^n ) \cap  K } @>\alpha >> X= S_{K_{11} ( v ^{n-1} )\cap K}.
\endCD
$$

Here the arrows in the second diagram is obtained by inclusions of compact open subgroups.
We write $ f = \alpha \cdot \beta' =  \beta  \cdot \alpha ' $.\par

\begin{lem}
\label{lem-coh41} In the notations as above, 
\begin{enumerate} 
\item The coproduct $Y \amalg _Z Y '  $ is isomorphic to $X$. 
\item For any $o_{E_\lambda} $-sheaf $\mathcal F $ on $X$, 
 $$ 
0 \longrightarrow \mathcal F \longrightarrow \alpha _* \alpha ^*  \mathcal F  \oplus \beta _* \beta ^*  \mathcal F 
 \longrightarrow f_* f^* \mathcal F 
$$ is exact.
\end{enumerate}
\end{lem}
\begin{proof}[Proof of Lemma \ref{lem-coh41}] We check the claim fiberwise. Take a geometric point $x$ of $ X$. Then the fibers $
\alpha ^{-1} (x) $, $\beta ^{-1}(x) $, and $f ^{-1} ( x) $ over $x$ are all
$K_{11} ( m _{F_v} ^{n-1} )\cap K_v$-homogeneous spaces, and are identified with the following left $ K_{11} ( m _{F_v} ^{n-1} )\cap K_v$-spaces, respectively: When
$n = 1$,
$ \SL_2 ( k(v) ) /A \cdot N  $ ($A$ is the group of the scalar matrices represented by a unit of
$F$, $N$ the standard unipotent subgroup consisting of the upper triangular matrices),
$
\SL_2 ( k(v) ) /A \cdot B ^{\op}
$ ($B ^{\op}$ is the opposite Borel subgroup consisting of the lower triangular matrices), $ \SL_2 (k(v) ) / A
$. Since $N$ and $B^{\op}$ generates $\SL_2 (k(v) )$, the claim follows from the
following remark: for a group $G$ and subgroups $H_1$ and $H_2$ which generate $G$, the
coproduct $ G/ H_1 \amalg _G G/H_2 $ is a point (the pullback to $G$ of a $\mathbb C$-valued function on the coproduct is right $H_1$ and $H_2$-invariant, and hence is a constant function). The claim for $ n \geq 2$ is proved in the same way. 
\par
(2) follows from (1) by applying Lemma \ref{lem-coh31} fiberwise. 

\end{proof}

\medskip Consider the complex 
$$
L = [  \alpha _* \alpha ^* \bar {\mathcal F}_{(k, w) , o_{E_\lambda}}  \oplus
\beta _* \beta ^* 
\bar {\mathcal F} _ {(k, w) , o_{E_\lambda}} \longrightarrow f _* f ^* 
\bar {\mathcal F}_{(k, w),  o_{E_\lambda}} ] 
$$
where the components are placed in degree $-1$ and 0. This is a perfect complex of $o_{E_\lambda}$-modules on
$X$ and the cohomology sheaves are locally constant since all components are locally constant sheaves. By Lemma \ref{lem-coh41}, for any
$o_{E_\lambda}
$-sheaf
$\mathcal G$ on $X$, $ H ^ {-1} ( L \otimes^{\mathbb L} _{o_{E_\lambda}} \mathcal G ) $ is isomorphic to $\bar {\mathcal
F }_{(k, w),  o_{E_\lambda} } \otimes  _{o_{E_\lambda}} \mathcal G $, and hence exact in $\mathcal G$. Since
non-trivial cohomology sheaves of
$ L$ are located in degree $-1$ and $0$,
$ \mathcal H = H ^ 0 (L)  $ is a smooth $o_{E_\lambda}$-sheaf.
Let $ \mathcal H ' $ be the image of  $  \alpha _* \alpha ^* \bar {\mathcal F}_{(k, w),  o_{E_\lambda} }  \oplus
\beta _* \beta ^* 
\bar {\mathcal F }_ {(k, w) ,  o_{E_\lambda}} $ in $f _* f ^* 
\bar {\mathcal F}_{(k, w),  o_{E_\lambda}} $.  \par
In summary, we have exact sequences
$$
0 \longrightarrow  \bar {\mathcal
F }_{(k, w),  o_{E_\lambda} } \longrightarrow  \alpha _* \alpha ^* \bar {\mathcal F}_{(k, w),  o_{E_\lambda} }  \oplus
\beta _* \beta ^* 
\bar {\mathcal F }_ {(k, w) ,  o_{E_\lambda}} \longrightarrow  \mathcal H ' \longrightarrow 0 ,
$$
and 
$$
0 \longrightarrow \mathcal H ' \longrightarrow f _* f ^* 
\bar {\mathcal F}_{(k, w),  o_{E_\lambda}} \longrightarrow \mathcal H \longrightarrow 0.
$$
Here, $ \mathcal H $ and $\mathcal H '$ are $o_{E_\lambda}$-smooth sheaves, and hence the exactness of the above sequences is preserved by any scalar extensions. In the case of $q_D= 0 $, the claim follows immediately.\par
 In the case of $q_D=1$, we make use of the action of $H_{o_{E_\lambda } } =  H ( D ^\times ( \A ^{\Sigma  }_{F, f} ), \ K
^\Sigma)_{o_{E_\lambda } } $, and work in $\bar {\mathcal C} _{T_\Sigma}  $.  
As in the proof of Theorem \ref{thm-coh31}, we take a finite set of finite places $ \Sigma$ which contains $\Sigma _K$, $v$, and all places dividing $\ell$. \par
In the exact sequence
$$
H ^ 0 (X,   \mathcal H ' _{k_\lambda } )  \longrightarrow  H ^1 (X,  \bar {\mathcal F }_{(k, w),  k_\lambda  } ) \longrightarrow   H ^1 (X,   \alpha _* \alpha ^* \bar {\mathcal F}_{(k, w),  k_\lambda  }  \oplus
\beta _* \beta ^* 
\bar {\mathcal F }_ {(k, w) ,  k_{\lambda}} )  $$
$$\longrightarrow  H ^1 (X,    \mathcal H '_{k_\lambda } ) \longrightarrow H ^2 (X,  \bar {\mathcal
F }_{(k, w),  k_\lambda  } ) 
$$
$ H ^ 0 (X, \ \mathcal H ' _{k_\lambda } )$ and $H ^ 2 (X, \ \bar {\mathcal
F }_{(k, w),  k_\lambda  })$ are modules of residual type by Lemma \ref{lem-coh41}, and hence vanishes in $ \bar {\mathcal C} _{T_\Sigma} $. \par
To show the injectivity of $H ^1 (X,  \mathcal H ' _{k_\lambda}) \to H ^1 (X, f _* f ^* \bar {\mathcal F}_{(k, w),  k_\lambda}) $, 
it suffices to prove that $H ^0 ( X, \mathcal H _{k_\lambda} )  $ is of residual type. \par
$ f ^* \mathcal H _{k_\lambda} $ is an $D ^{\times} ( \mathbb A ^{\Sigma} _{F, f} ) $-equivariant quotient of $f ^*  f_* f ^* \bar {\mathcal
F }_{(k, w),  k_\lambda  } $, and $f ^*  f_* f ^* \bar {\mathcal
F }_{(k, w),  k_\lambda  } = f ^*  f_* \bar {\mathcal
F } ^{ K \cap K _{11}( v ^n ) \cap K_0 (v ^{ n +1} )}_{(k, w),  k_\lambda  }$ is isomorphic to a direct sum of copies of $ \bar {\mathcal
F } ^{ K \cap K _{11}( v ^n ) \cap K_0 (v ^{ n +1} )}_{(k, w),  k_\lambda  }$ in a $ D ^{\times } (\mathbb A ^{\Sigma} _{F, f} )$-equivariant way. 
Thus $ H ^0 ( Z, f ^* \mathcal H _{k_\lambda} )$ is of residual type by Lemma \ref{lem-coh41}. $H ^0 ( X, \mathcal H _{k_\lambda} )  $ is a submodule of $ H ^0 ( Z,  f ^* \mathcal H _{k_\lambda} )$ and hence vanishes in $ \bar {\mathcal C} _{T_\Sigma}$.

\end{proof}

\section{Nearly ordinary automorphic representations}
\label{sec-nearlyordinary}

\subsection{Renormalization of cohomological Hecke correspondences}\label{subsec-nearlyordinary1}
We fix a place $v$ which divides $\ell$.  For $ K = \GL_2 (o_{F_v} )$, 
we have the Cartan decomposition
$$
\GL_2(F_v) = \cup_{(a, b)  \in  \mathbb Z ^ 2, a\leq b } K  \begin{pmatrix}
p^a_v  &0 \\ 0 & p ^b _v
\end{pmatrix} K.
$$
\begin{dfn} \label{dfn-nearlyordinary11} 
\begin{enumerate}
\item  For an infinity type $ ( k , w ) $ and a place $v $ dividing $\ell$, the $v$-type is defined as a pair $ (k _v , w ) $, where $ k _v = ( k_{\iota } ) _{\iota \in I _{F, v} } $. Here we regard $ I_{F, v} $ (cf. \S \ref{subsec-shimura3}) as a subset of $ I _{F, \infty} $ by $\bar E_{\lambda } \simeq \mathbb C $. 
\item For $h  \in \GL_2(F_v)$, a pair of integers $( a_v (h) ,  b_v (h)  )$ is defined by $ Kh K =  K \begin{pmatrix}
p^{a_v (h)}  &0 \\ 0 & p ^{b_v (h) } _v 
\end{pmatrix}  K $ for $K =\GL_2 (o_{F_v} )$.  $( a_v (h) ,  b_v (h)  )$ depends only on $h$ under the condition that $a_v (h) \leq  b_v (h)$. 
\item For an element $ g \in G_D ( \mathbb A _{\mathbb Q , f} ) $, 
$$
m _v(g)=   \prod _{\iota \in I_{F, v} }  (\iota ( p_v) ) ^{w \cdot a _v ( g_v  ) + k ' _\iota \cdot( b_v ( g_v ) -a_v ( g_v ) )  }    
$$
is the multiplier of $g$ at $v$ with respect to the $v$-type $( k_v , w )$ and $p_v $. Here $ k ' _\iota $ is defined as \ref{dfn-shimura21}, and $k'_v = ( k ' _\iota )_{\iota \in I_{F, v} } $. 
\end{enumerate}
\end{dfn}

Then $( g_v ) ^{-1} ( V_{(k, w ) , o _{E_\lambda} }  )  \subset m _v (g)  \cdot  V_{(k, w) , o _{E_\lambda} }$, and $(R (g ^{-1}  ) ) ^* \bar {\mathcal F} _{(k, w) , o_{E_\lambda} } \subset  \prod _{ v \vert \ell } m _v (g_v) \bar {\mathcal F} _{(k, w) , o_{E_\lambda} }$.
For an $F$-factorizable compact open subgroup $K$ of $ G_D ( \mathbb A_{\mathbb Q, f} )$, one defines a cohomological correspondence of $ \bar {\mathcal F} _{(k, w) , o_{E_\lambda} } $ by using $ (\prod _{ v \vert \ell } m _v (g_v) )^{-1} R(g^{-1})^* $ at $v\vert \ell$ when $K _v $ is a subgroup of $\GL_2 (o_{F_v} ) $. 

\begin{dfn}\label{dfn-nearlyordinary12}
\begin{enumerate}
\item For an $F$-factorizable compact open subgroup $K$ of $G_D ( \mathbb A_{\mathbb Q, f} ) $, the renormalized cohomological correspondence $[Kg K ] ^{\ren} $ is the correspondence given by  $(\prod _{ v \vert \ell } m _v (g_v) )^{-1} R(g^{-1}) ^* $.
\item For a uniformizer $p_v $ of $F_v$, $\tilde U ( p_v ) $ and $\tilde U (p_v , p_v ) $ is the renormalized $U(p_v ) $ and $U ( p_v, p_v ) $-operators defined by $[ K a(p_v ) K]^{\ren}$ and $[K b(p_v ) K] ^{\ren}$, respectively. 
\end{enumerate}
\end{dfn}
Note that this renormalization depends on the choice of infinity type $(k, w) $, and the choice of an isomorphism $\bar {\mathbb Q}_{\ell} \simeq \mathbb C  $.\par

Since the renormalized correspondence preserves the $o_{E_\lambda}$-lattices, we have
\begin{prop}\label{prop-nearlyordinary11} Let $\pi $ be a representation in $\mathcal A _{(k, w) } (G_D  )  $ with a discrete infinity type. Assume that $D$ is split at the places dividing $\ell$.  For any $F$-factorizable compact open subgroup $K$, and any eigenvalue $\alpha $ of $\chi _{K gK }: \pi ^ K \to \pi ^K  $, the inequality
$$
 \val (\alpha  ) \geq \sum _{ v \vert \ell } \val ( m _v  (g) ) 
$$
$$= w \sum _{v \vert \ell }[F_v : \mathbb Q_\ell]  \cdot a_v (g_v) \val (p_v)   +  \sum _{ v \vert \ell} (\sum _{ \iota \in I _{F, v} }  k' _{\iota } \cdot (b_v (g_v ) - a_v (g_v )  ) ) \val (p_v)
$$
holds. 
Here $\val : \bar {E} ^{\times} _\lambda \to \mathbb Q $ is the additive
valuation. 
\end{prop}

At $v$, we normalize $\val $ by $\val ( \iota ( p_v )) =1 $ for $\iota \in I_{F, v}$. This definition does not depend on the choice of a uniformizer.\par
For a quasi-character $\mu : F^{\times} _v \to \bar E _\lambda  ^\times $, the slope of $\mu$ is defined by 
$$
 \slope(\mu)   = \val (\mu ( p _v)   ) , 
$$
where $p_v$ is a uniformizer at $v$.\par
\begin{cor}[Hida, \cite{H1}, theorem 4.11]
\label{cor-nearlyordinary11} Let $\pi $ be a representation in $\mathcal A _{(k, w) } (G_D  )  $ with a discrete infinity type, $\rho _{F_v} $ the $F$-semi-simple representation
of the Weil-Deligne group $W' _F $ over
$\bar {E} _\lambda $ corresponding to $ \pi _v$ by the local Langlands correspondence
\cite{K}. Assume that
$\chi$ is a subrepresentation of $\rho _{F_v}$. Then the inequality 
$$ 
\slope (\chi)  \geq \sum _{\iota \in I_{F, v} } k ' _{\iota}
$$
holds. 
\end{cor}
Hida proves the theorem including limit discrete series by using Hilbert modular varieties, that is, in the case when $
k _\iota \geq 1
$ for all $\iota \in I_{F, \infty} $.

\subsection{Nearly ordinary automorphic representations}\label{subsec-nearlyordinary2}
We consider the local restrictions of the Galois representation attached to automorphic representations when the residual characteristic $p$ is $\ell$, and
define the notion of minimality. In this article, nearly ordinary representations are mainly considered for $v \vert \ell$.

Let $\pi $ be a cuspidal representation $\GL_2 (\A _F)$ of discrete infinity type $ (k, w )
$. The local restriction of the associated
$\ell$-adic representations
$
\rho _{\pi, \lambda}
\vert_{G_{F_v}}$ for
$ v
\vert
\ell$ is a potentially stable representation in the sense of Fontaine unless $ k =
(2, \ldots 2)$ (see 
\cite{Sa2}. It is enough to show this over some finite extension of
$F_v$, and the claim follows from Carayol
\cite{Car2}, Blasius-Rogawski \cite{BR}, and de Jong's method of alteration
\cite{dJ} and Tsuji's theorem \cite{Tsu}). If $ k = (2, \ldots, 2)$, we assume that $d$ is odd, or $\pi$
admits an essentially square integrable local component. For potentially stable
$\ell$-adic representations, T. Saito
\cite{Sa1} has constructed a representation of the Weil-Deligne group $W_{F_v} ' $ over
$\bar E_\lambda $.

Here we make a temporary construction in some one dimensional cases. 
We assume that all continuous field embeddings $ F_v \hookrightarrow \bar E_{\lambda}$ factors through $E_\lambda$. \par

For a quasi-character $\mu: F_v ^{\times}\to E_{\lambda}^{\times}$ and an integral vector $k' _v = ( k ' _\iota ) _{ \iota \in I_{F, v} } \in \Z  ^{I_{F, v}}$ which satisfies
$$
 \slope(\mu)= \sum _{\iota \in I_{F, v}} k ' _\iota, 
$$
we define the associated
$\ell$-adic character
$ L_{k'_v} (\mu )
$. \par

Let $ \mathcal L ( p _ v) $ be the Lubin-Tate formal group over $o_{F_v} $ attached to $p_v$, $ \chi _{p _ v }: 
G_{F_v} \to o  ^\times  _{F_v}\hookrightarrow o^{\times}_{E_\lambda }   $ the Galois representation attached to $p_v$-divisible group $ \mathcal L ( p _ v) [ p^{\infty}_v ] $. It is explicitly given as follows. 
The splitting $  \mathbb Z \to F^{\times }_v $ of 
$$
1 \longrightarrow o^{\times}_{F_v} \longrightarrow  F^{\times }_v  \overset {\val }\longrightarrow \mathbb Z \longrightarrow 0 
$$
given by $1 \mapsto  p_v$ induces the projection $\alpha :\widehat{F _v ^{\times} }\twoheadrightarrow o^{\times} _{F_v}$. $\alpha (x) = x \cdot p^{-\val (x) }_v  $ for $ x\in F _v ^{\times} $.\par

By our normalization of the Galois module associated to $p$-divisible groups, $ \chi _{p_v} $ is the composition of $G^{\ab}_{F_v} \simeq \widehat{F _v ^{\times} }\overset {\alpha  } \to o^{\times} _{F_v}\overset{ (-)^{-1}}  \to o^{\times} _{F_v}\hookrightarrow o^{\times}_{E_\lambda} $. \par
 We define a quasi-character $ \mu _{p_v} : F_v ^ \times \to E_\lambda ^{\times}$ by 
$$
 \mu_{p_v} (x)  =\mu (x) /
\prod _{ \iota \in I_{F, v}} \iota (p_v ^{\val(x)}) ^ { k '_\iota }.
$$
Then $\mu_{ p_v} $ takes the values in $o^{\times}_{E_\lambda} $, extends to a continuous character of $ \widehat {F^{\times} _v}$, and is regarded as a character of $G _{F_v}$ by the local class field theory. 
$$
L _{k'_v}(\mu) = \mu _{ p_v}   \cdot \prod _{\iota \in I_{F, v} } (\iota \circ \chi _{p_v} ) ^ {k '_\iota  }
: G_{F_v} \longrightarrow o_{E_\lambda}  ^\times
$$
is independent of any choice of $ p_v$.

\begin{dfn}[Hida]\label{dfn-nearlyordinary21} Let $\pi $ be a cuspidal representation of $\GL_2 ( \mathbb A _F) $ of infinity type $(k, w)$ defined over $\bar E_{\lambda } $, and $v$ a finite place of $F$ which divides $\ell$.
\begin{enumerate}
\item
Let $\rho _v
$ be the $F$-semi-simple representation of the Weil-Deligne group $W'_{F_v} $ over
$\bar {E} _\lambda $ which corresponds to $ \pi _v$ by the local Langlands correspondence. 
We say $ \pi $ is nearly ordinary at $v$ if $ \rho _v$ contains a character 
$\mu _v $ as a subrepresentation, and the slope of
$\mu _v $ satisfies 

$$
\slope (\mu _v ) = \sum _ {\iota \in I _{F, v}} k ' _\iota .
$$ 
The character $\mu _v $ is called a nearly ordinary character of $\pi $ at $v$.

\item Assume that $\pi$ is nearly ordinary at $v$ with nearly ordinary character
$\mu _v $. A non-zero vector
$z$ in the representation space
$V$ of
$\pi_v$ is called a nearly ordinary vector if 
$$ 
\pi _v(
\begin{pmatrix}
a &*\\ 0&d 
\end{pmatrix}
)
 \cdot z
= \chi _{1 , v} ( a) \vert a \vert _v  \chi _{2 , v} (d) z  ,\  \quad  a, \ d
\in F_v ^\times, \ 
$$
holds for some quasi-character $\chi_{1, v}$, and $\chi_{2, v} =\mu_v$.  $\mu_v \vert _{o^\times_{F_v} }$ is a nearly ordinary type of $\pi $ at $v$. 
\item For a quaternion algebra $D$ central over $F$, assume that $D$ is split at the places dividing $\ell$. Then for a place $v$ dividing $\ell$ and a representation $\pi \in \mathcal A _{(k , w)} (G _D  ) $, $\pi$ is nearly ordinary at $v$ if the Jacquet-Langlands correspondent $\JL ( \pi ) $ is nearly ordinary at $v$. 
\end{enumerate}
\end{dfn} 

In other words, $\pi$ is nearly ordinary at $v$ if and only if
$\tilde U(p_v)$-operator has an $\ell$-adic unit eigenvalue. In the case when $ k_\iota \geq
2$ for any $\iota \in I_{F, v}$, a nearly ordinary character $\mu_v$ is unique, since $\slope (\chi _{1, v} ) = \slope ( \chi _{2, v} ) + \sum _{\iota \in I_{F, v} } ( k _{\iota } -1 )  $ holds.

\begin{conj}\label{conj-nearlyordinary21}
For a cuspidal representation $\pi$ of $\GL_2 ( \mathbb A _F) $ of infinity type $(k, w)$, assume that $\pi $ is nearly ordinary at $v$ with nearly ordinary character
$\mu _v $. Let $\rho = \rho _{ \pi , E_{\lambda}}$ be the $\ell$-adic representation associated to
$\pi$. Then the local representation
$\rho\vert _{F_v}$ is reducible, and contains $ L _{k '_v}(\mu_v ) $ as a subrepresentation.
\end{conj}

\begin{hyp}[Local monodromy hypothesis]\label{hyp-nearlyordinary21} Conjecture \ref{conj-nearlyordinary21} is true for $\pi$ with discrete infinity type $ (k, w) $ if one of the following conditions hold:
\begin{itemize}
 \item The degree
$d$ is odd. 
\item $d$ is even, and $\pi $ has an essentially square integrable component $\pi _u $ for
some finite place $u$.
\end{itemize}
\end{hyp}
Hypothesis \ref{hyp-nearlyordinary21} is known to hold in many cases. 

\begin{thm}
\label{thm-nearlyordinary21}Hypothesis \ref{hyp-nearlyordinary21} is true in the following cases:
\begin{enumerate}
\item $\pi $ is nearly ordinary at all $ v \vert \ell$.
\item The infinity type of $\pi$ is $ ((2,
\ldots, 2), 0 ) $, and the degree $d$ is odd, or $d$ is even and $\pi $ has an essentially square integrable
component at some finite place $u $ (\cite{W2}, lemma 2.1.5). 
\end{enumerate}
\end{thm}
 The case of (1) follows from the case of (2) by Hida's control theorem of nearly
ordinary Hecke algebras \cite{H2}, and the existence of Galois representations. The method of
\cite{W2} is used there. Since
$\pi$ is obtained as a specialization of Hida family, it suffices to check Hypothesis \ref{hyp-nearlyordinary21} for the
Galois representation attached to the nearly ordinary Hecke algebra. Since the algebraic points
corresponding to cuspidal representations with infinity type $ ( (2, \ldots, 2), 0 )$ are dense, the claim follows
from (2). For (2),  see
\cite{W2}, theorem 2.1.4 when $\pi_v$ belongs to principal series. When
$\pi _v$ is a (twisted) special representation, (2) follows from
\cite{Gr2}. \par

\begin{rem}
\label{rem-nearlyordinary21}
A proof of Hypothesis \ref{hyp-nearlyordinary21} will
be given if the residual representation is not of residual type in a forthcoming article. For $v \vert \ell$, we construct an $\ell$-adic family of cuspidal representations which are nearly ordinary at $v$, and reduce to the case of (2).
\end{rem}

\subsection{Modularity and minimality}
\label{subsec-nearlyordinary3}

As in the introduction, an absolutely
irreducible mod $\ell $ Galois representation
$
\bar \rho : G _\Sigma \to \GL _2(k_\lambda )$ is modular if there is an
$\ell$-adic field $E_\lambda $ and a cuspidal
representation
$\pi
$ on
$\GL _2 (\A _F)
$ of discrete infinity type $(k , w) $ such that the finite part $\pi _f $ is defined over $E_\lambda $, and $
\bar
\rho \simeq \rho_{\pi,
E_{\lambda} } 
\mod
\lambda $ over $\bar k_\lambda$.
 In \S\ref{sec-modular}, we need to choose $\pi$ with a distingushed property. We make the following temporaly
definition for such a choice. \par
\bigskip
For $s \in \mathbb C$, $\omega _s : F ^{\times} \backslash \mathbb A ^\times _F \to \mathbb C $ is the Hecke character defined by 
$$
\omega _s ( x) = \vert x \vert   ^s .
$$ If $s$ is an integer, this is an algebraic Hecke character, and the $\ell$-adic representation associated to $\omega _s$ is $ \chi ^s  _{\cycle} $ by our normalization of the reciprocity map. \par
\begin{dfn}\label{dfn-nearlyordinary31} Let $\alpha $ be an integer. 
\begin{enumerate}
\item For a character $\bar \mu : G_F \to k^{\times} _\lambda  $, let $(\bar \mu) _{\lift} $ be the Teichm\"uller lift of $\bar \mu $. 
We view $\nu = (\bar \mu )_{\lift}\cdot (\chi^{\ell} _{\cycle} ) ^{w+1} $ as a Hecke character of finite order by the class field theory. Then the algebraic Hecke character associated to $\bar \mu  $ of weight $2 \alpha $ is defined as
$$
\chi ^{\Hecke} _{\bar \mu , \alpha   }  =  \nu \cdot \omega _{-\alpha } .
$$ 
\item For $\bar \rho  : G _F \to \GL_2 (k_\lambda)$, and an integer $\alpha $,  $\chi ^{\Hecke} _{\det \bar \rho , \alpha} $ is denoted by $ \chi ^{\Hecke} _{ \bar \rho , \alpha} $.

\end{enumerate}
\end{dfn}
By the definition, the $\lambda$-adic representation $\chi_{\bar \mu , \alpha } : G_F\to o ^\times _{E_\lambda} $ of weight $2\alpha $ attached to $\chi ^{\Hecke} _{\bar \mu, \alpha }$ is given by 
$$
\chi_{\bar \rho , \alpha } = (\det \bar \rho )_{\lift} \cdot  (\chi _{\cycle , \ell}   ) ^ { -\alpha  } .
$$

\begin{dfn}
\label{dfn-nearlyordinary32}
Let $\bar \rho$ be a modular mod $\ell$-representation with a minimal deformation type $\mathcal D_{\min}$. A cuspidal representation 
$\pi _{\min}$ of $ \GL_2 (\mathbb A_F )$ with a discrete infinity type $(k, w )$ defined over $E_{\lambda}$ is called a minimal lift of $\bar \rho$ if the following properties hold:
\begin{enumerate}
\item $\rho_{\pi , E_\lambda} $ is a deformation of $\bar \rho $.
\item  $\rho _{\pi, E_ \lambda } \vert _{G_{F_v}} $ for $v \nmid \ell$ is a finite
deformation of $\bar \rho \vert _{F_v} $.
\item The central character of $\pi_{\min}$ is $\chi ^{\Hecke}_{\bar \rho , w}  $.
\item  If $  \deform _{\mathcal D _{\min} } ( v) =  \bold {fl}$, the $v$-type $ (k_v, w)$ is $( (2, \ldots, 2) , w )$, and $\pi _v \otimes \mu ^{-1}_v $ is spherical for some quasi-character $\mu _v $. $\mu _v \cdot \vert \cdot \vert _v ^{\frac w 2}  $ takes the values in $o^{\times} _{E_\lambda} $, $ \mu _v\vert _{o^{\times}_{F_v}} $ has order prime to $\ell$. $\mu _v$ is called a flat twist character, and $\mu_v \vert _{o^{\times}_{F_v}} $ the flat twist type of $\pi$ at $v$. 
\item If $v \vert \ell $ and $ \deform _{\mathcal D _{\min} } ( v) = \bold {n.o.f.} , \bold {n.o. }$, $\pi _{\min}$ is nearly ordinary at $v$. For the nearly ordinary character $\mu _v $, the nearly ordinary type $\mu _v \vert _{o^{\times } _{F_v } }$ has the order prime to $\ell$. $ L _{k'_v }(\mu _v ) $ is a lift of the nearly ordinary type $\bar \kappa _v$ of $\bar \rho \vert _{F_v} $. 
\item If $  \deform _{\mathcal D _{\min} } ( v) =  \bold {n.o.f.}$, the $v$-type $ (k_v, w)$ is $( (2, \ldots, 2) , w )$, and $\pi_v \otimes \mu_v ^{-1}$ is spherical for the nearly ordinary character $\mu_v$.
\end{enumerate}
\end{dfn}
\begin{rem}
\label{rem-nearlyordinary31}
\begin{enumerate}
\item Assume the following condition is satisfied at $v \nmid
\ell$:
$\bar \rho \vert _{G_{F_v}} $ is of type $0_{NE}$, or
$
\bar
\rho
\vert _{G_{F_v}} $ is absolutely reducible and ramified at $v$ with $ \dim _{k_\lambda}
\bar
\rho ^{I_{F_v}} = 1$.\par

Then the condition (2) for $v \nmid  \ell$ is equivalent to the equality of conductors
$$ 
\cond \pi _v = \Art \bar  \rho \vert _{I_{F_v}}.
$$ 
\item If the conditions (4)-(6) at $v\vert \ell$ is satisfied for a modular lift $\pi$, by the results of \cite{Fu1}, \cite{Ja1}, \cite{Ja2}, \cite{Raj}, one finds $\pi _{\min}$ so that (2) and (3) are satisfied. 
\end{enumerate}
\end{rem}
\bigskip
We define the $K$-type $(G _{\bar \rho \vert _{F_v} }, K _* ( \bar \rho \vert _{F_v} ), \nu _* ( \bar \rho \vert _{F_v} )  ) $ for $* = \bold {fl} , \bold {n.o.f.}, \bold {n.o.}$ from $\bar \rho : G_F \to \GL _2 (k_\lambda) $ and the $v$-type $(k_v, w) $. 

\begin{dfn}\label{dfn-nearlyordinary33}
Let $\bar \rho$ be a mod $\ell$-representation with a deformation type $\mathcal D$ and a $v$-type $(k_v , w )$. 
\begin{enumerate}
\item When $ \deform_{\mathcal D} (v) = \bold {fl}$, $k_v = (2, \ldots, 2)$, and for the flat twist character $\bar \kappa _v$ of $\bar \rho \vert_{F_v} $, a quasi-character
$\chi$ is called an automorphic flat twist character if $ \chi \cdot \vert \cdot \vert^{\frac w 2 }  _v $ takes the values in $o_{E_\lambda} $ and the reduction mod $\lambda $ is $\bar \kappa _v \cdot \bar \chi ^{\frac w 2} _{\cycle}  $. $\nu _v = \chi
\vert _{o_{F_v} ^\times} 
$ is called the automorphic flat twist type at $v$.
\item When $ \deform_{\mathcal D} (v) = \bold {n.o.f.}, \bold {n.o.}$, for the nearly ordinary character 
$\bar \kappa _v $ of $\bar \rho \vert_{F_v} $, a quasi-character
$\chi$ of slope $ \sum _{ \iota \in I_{F, v} } k' _\iota $ is called an automorphic nearly ordinary character if $ \chi \vert _{o^\times_{F_v} }$ has order
prime to $\ell$, and $L_{k '_v  } (\chi) \mod \lambda = \bar \kappa _v $. $\nu _v = \chi
\vert _{o_{F_v} ^\times} 
$ is called the automorphic nearly ordinary type at $v$.
\end{enumerate}
An automorphic twist type is unquely determined from $\bar \rho $ and $(k_v, w)$. 
\end{dfn}
\begin{dfn}\label{dfn-nearlyordinary34} Let $\bar \rho$ be a mod $\ell$-representation with a deformation type $\mathcal D$ and a discrete infinity type $(k , w )$. 
\begin{enumerate}
\item $G _{\bar \rho \vert_{F_v} } = \GL _2 (F_v ) $ for $v \vert \ell$. \par
\item If $\deform _{\mathcal D}(v) = \bold {fl} $, let $\nu _v$ be the automorphic flat twist type of $\bar \rho$. Then
$$
K_{\bold {fl}} (\bar \rho \vert _{F_v}  )= \GL_2 (o_{F_v} ) 
$$
with the $K$-character $\nu _{\bold {fl} } (\bar \rho \vert _{F_v} ) : K_{\bold {fl}} (\bar \rho \vert _{F_v} )\overset {\det} \to o ^{\times}_{F_v }\overset{ \nu_v} \to o ^{\times}_{\mathcal D} $.
\item If $\deform _{\mathcal D}(v) = \bold {n.o.f.}$ let $\nu _v$ be the automorphic nearly ordinary type of $\bar \rho $. Then 
$$
K_{ \bold {n.o.f.}} (\bar \rho \vert _{F_v} )= \GL_2 (o_{F_v} )
$$
with the $K$-character $\nu _{\bold {n.o.f.} } (\bar \rho \vert _{F_v} ) : K_{\bold {n.o.f.}} (\bar \rho\vert _{F_v} )\overset {\det} \to o ^{\times}_{F_v }\overset{ \nu_v} \to o ^{\times}_{\mathcal D} $.
\item If $\deform _{\mathcal D}(v) = \bold {n.o.}$ let $\nu _v$ be the automorphic nearly ordinary type of $\bar \rho $. Then 
$$
K_{ \bold {n.o.}} (\bar \rho \vert _{F_v}) = \mu _{\ell ^{\infty} } (F)  \cdot K _1 ( m _{F_v} )
$$
with the $K$-character $\nu _{\bold {n.o.} } (\bar \rho \vert _{F_v} ) : K_{\bold {n.o.}} (\bar \rho \vert_{F_v} )\overset {\det} \to o ^{\times}_{F_v }\overset{ \nu_v} \to o ^{\times}_{\mathcal D} $.
\end{enumerate}
Similarly as in \S\ref{subsec-galdef8}, the intertwining space is defined as 
$$
I_{\deform _{\mathcal D} (v)  } (\bar \rho \vert _{F_v} , \pi _v  ) = \Hom _{K _{\deform _{\mathcal D } (v) }}(\nu_{\deform _{\mathcal D } (v) } (\bar \rho \vert _{F_v} ),    \pi _v ) .
$$
\end{dfn}

As an analogue of proposition \ref{prop-galdef82}, we have 
\begin{prop}
\label{prop-nearlyordinary31} 
For a cuspidal representation $\pi$ of infinity type
$(k, w) $ giving $\bar \rho$, 
 assume that $\pi $ is nearly ordinary at $v$ and we are in case
\ref{dfn-nearlyordinary32}, (4). Then the part $\tilde U(p_v)$-operator has an $\ell$-adic unit eigenvalue in
$I_{\deform _{\mathcal D} (v)  } (\bar \rho \vert _{F_v} , \pi _v  )$ is at most one dimensional.
\end{prop}
\begin{proof}[Proof of Proposition \ref{prop-nearlyordinary31}] By twisting, we may assume that the twist type at $v$ is trivial. The only
non-trivial statement in \ref{prop-nearlyordinary31} is the case when
$\pi_v$ is spherical and $\deform _{\mathcal D } (v) = \bold {n.o.}$, since $I_{\bold {n.o.}  } (\bar \rho \vert _{F_v} , \pi _v  ) $ is two dimensional. For the two eigenvalues $\alpha _v$, $\beta_v$, slopes are different as in \S\ref{subsec-nearlyordinary2}.
So the only one eigenvalue of $ \tilde U(p_v)$-operator is an $\ell$-adic unit.  
\end{proof}

\begin{rem}
If $\deform _{\mathcal D } ( v) = \bold {n.o.}, \bold{n.o.f.}$ and $\bar \rho \vert _{F_v} $ is semi-simple,  we must specify the nearly ordinary type as in \S\ref{subsec-galdef9}, since the possibility is not unique. So the notion of minimal lift depends on this choice. 
\end{rem}

\section{Universal modular deformations}\label{sec-modular}

As is recalled in the introduction, the existence of the
$\lambda$-adic representation $G_F
\to \GL _2 ( T   \otimes _{o_{E_\lambda}} E_\lambda  ) $ is known, where $T $ is a certain
$\ell$-adic Hecke algebra attached to automorphic representations. It is quite important to know the existence of a Galois
representation having the values in $\GL _2 (T  ) $, which can be seen as a deformation of a mod
$\ell$-representation. In this section, we construct it when the residual
representation is absolutely irreducible, and study the local deformation properties. The local behavior is controled by the compatibility of the local and the global Langlands correspondences at the places outside
$\ell$.  At a place dividing $\ell$, we need much finer information beyond the general theory still. \par
The discovery of Galois
representations having values in $\ell$-adic Hecke algebras is due to Hida (in the case of $\mathbb Q$, and for ordinary
Hecke algebras \cite{H0}). We use the method of
pseudo-representations of Wiles \cite{W2} to construct the Galois representation. \par

\subsection{Hecke algebras}\label{subsec-modular1}

Let $ \mathcal S _{(k, w)} (G_D) $ be a $G _D (\A _{\mathbb Q, f} )$-representation over $\bar E_{\lambda} $ defined by
$$
\mathcal S _{(k, w)} (G_D) = \bigoplus _{\pi \in \mathcal A _{(k, w)} (G_D )}
\pi _f. 
$$
Here $ \mathcal S _{(k, w)} (G_D) $ is defined up to isomorphisms, and for a compact open subgroup $K$ of $ G _D (\A _{\mathbb Q, f} )$, the $K$-fixed part $\mathcal S _{(k, w)} (G_D)^K $ is finite dimensional over $\bar E_\lambda $.\par
For an algebraic Hecke character $\chi : F ^{\times } \backslash \mathbb A ^{\times}_F \to \mathbb C^{\times}$, $\mathcal S _{(k, w), \chi} (G_D )$ is defined as the subspace parametrized by $\mathcal A _{(k,
w), \chi } (G_D )$.\par

For an $F$-factorizable compact open subgroup $K$ of $G _D (\A _{\mathbb Q, f} ) $ and a finite set of finite places $\Sigma $ which contains $\Sigma _K $, the reduced Hecke algebra generated
by the standard Hecke operators is defined by
$$
T ^{\ss } _{K , \Sigma }= \Z  [ T_v ,\ T_{v, v} ,\  v \not \in \Sigma]\subset
\End_{\Z } \mathcal  S _{(k, w)} (G_D) ^K, 
$$
and for a commutative ring $R$, 
$$
T ^{\ss } _{K , \Sigma ,  R }= T ^{\ss } _{K , \Sigma }\otimes _\Z R.
$$

$T ^{\ss } _{K ,  \Sigma } $ is commutative by the definition, and $ T^{\ss }_{K, \Sigma , 
\bar E_\lambda }$ acts faithfully on $ \mathcal  S _{(k, w)} (G_D) ^K $. It is known that $T ^{\ss } _{K ,  \Sigma }$ is finite free over $\Z $ since it fixes a
$\Z $-lattice in $\mathcal S _{(k, w)} (G_D) ^K $.  \par
For an $\ell$-adic integer ring $o   _{\lambda } $ and $\Sigma $ which contains all places dividing $\ell$, an $o  _{\lambda }$-algebra homomorphism $ f: T ^{\ss}_{K, \Sigma,  o  _{\lambda }}
\to o _{\lambda  }$ corresponds to an admissible irreducible $ D ^{\times } ( \mathbb A ^{\Sigma }_{F , f} ) $-representation of the form $\pi ^{\Sigma} =  \bigotimes_ { v \not \in \Sigma } \pi _v  $ for some cuspidal representation
$\pi \in \mathcal A _{(k,w) }(G_D) $ defined over
$E  _{\lambda  }$. By the strong multiplicity one theorem for $\GL_2$ and their inner twists,
$\pi $ is uniquely determined from $f$, and $f$ is denoted by $f _{\pi }$. In this way we may identify such an $o_{\lambda}$-homomorphism $f$ with a representation $\pi \in \mathcal A _{(k,
w)}(G_D) $ which is defined over $E_{\lambda }$ and satisfies
$(\pi _f ) ^K \neq \{  0 \} $.\par
In terms of the Hecke algebra homomorphism $f_{\pi }$ attached to $\pi \in \mathcal A _{(k, w)} (G_D)$, the Galois representation $\rho _{\pi , \bar E_\lambda } $ is described as follows. For an $F$-factorizable compact open subgroup $K$ such that $\pi ^K _f \neq \{ 0 \} $, and a finite set of finite places $\Sigma $ which contains $\Sigma _K \cup \{ v: v \vert \ell  \} $, 
$$
\begin{cases}
\tr \rho _{\pi , \bar E_\lambda }  (\Fr _v) =   f _{\pi }(T _v )\\
	    \det  \rho _{\pi , \bar E_\lambda } (\Fr _v ) = q _v \cdot f _{\pi } (T_{v , v} )

	\quad  \text{for any $v \not \in \Sigma  $}
	\end{cases}
$$
holds. The field of definition $E_{\pi}$ is generated by $ f_{\pi } (T_v)  $ and $f_{\pi} (T_{v, v} ) $ for $v \not \in \Sigma $, and the embedding $E_{\pi} \hookrightarrow \bar E_{\lambda}$ determines a valuation $\tilde \lambda $ of $E_{\pi }$. The Galois representation $\rho_{\pi , \bar E_{\lambda}} $ is defined over the completion $ E_{\pi , \lambda '}$ with respect to $\lambda'$.

\subsection{$\ell$-adic Hecke algebras attached to deformation types}\label{subsec-modular2}
Let $F$ be a totally real field, $E_\lambda $ an $\ell$-adic field with $\ell \geq 3$, $\bar \rho : G_F \to \GL_2 ( k_\lambda) $ an absolutely irreducible representation. 
Let $\mathcal D$ be a deformation type of $\bar \rho$. In the following sections, we fix a discrete infinity type $(k, w)$ which satisfies
\begin{itemize}
\item[{\bf IT}] If $\deform _{\mathcal D} (v) = \bold{fl} $ (resp. $\bold {n.o.f.}$) for $v \vert \ell$, the $v$-type $ (k_v , w ) =( ( 2, \ldots, 2) , w ) $. 
\end{itemize}
We restrict ourselves to deformation types which satisfy the following conditions:
\begin{itemize}
\item[{\bf D1}] $o _{\mathcal D} $ is an $\ell$-adic integer ring.
\item[{\bf D2}] If $\deform_{\mathcal D} (v) = \bold {fl}$ for $v \vert \ell$, the flat twist type $\kappa_{\mathcal D , v} $ is $ L _{k' _v } ( \chi _v )\vert _{I_{F_v}}$. Here $\chi_v $ is an automorphic flat twist character (see Definition \ref{dfn-nearlyordinary33}).
\item[{\bf D3}] If $\deform_{\mathcal D} (v) = \bold {n.o.f}, \bold{n.o.}$ for $v \vert \ell$, the nearly ordinary type $\kappa_{\mathcal D , v} $ is $ L _{k' _v } ( \chi _v )\vert _{I_{F_v}}$. Here $\chi_v $ is an automorphic nearly ordinary character (see Definition \ref{dfn-nearlyordinary33}).
\end{itemize}
We view the quotient field $E_\mathcal D$ as a subfield of $\bar E_{\lambda} $.\par
Moreover, we assume the following conditions on $\bar \rho$:
\begin{itemize}
\item[{\bf MM1}] There is a minimal deformation type $\mathcal D _{\min} $ such that $o _{\mathcal D_{\min} } = o _{E_\lambda} $.
\item[{\bf MM2}] There is a cuspidal representation $\pi _{\min}$ of $ \GL_2 (\mathbb A_F )$ with a discrete infinity type $(k, w )$ such that $\pi _{\min , f }$ is defined over $E_\lambda $, and $\pi_{\min}$ is a minimal lift of $\bar \rho$ (see Definition \ref{dfn-nearlyordinary32}). 
\end{itemize}
Note that there is a morphism $\mathcal D _{\min} \to \mathcal D $ in $ \Type (\bar \rho ) $ (Definition \ref{dfn-galdef92}) for some minimal deformation type $\mathcal D _{\min} $.\par
\bigskip
We attach a triple $(G _{\bar \rho} , K _{\mathcal D}, \nu _{\mathcal D}  )$ to $\mathcal D$, where $G _{\bar \rho}  $ is an inner twist of $ \Res _{F/\mathbb Q} \GL _{2, F} $ which depends only on $\bar \rho$, $ K _{\mathcal D}$ is a compact open subgroup of $G _{\bar \rho}  (\mathbb A _{\mathbb Q, f} )  $, and $\nu _{\mathcal D}: K _{\mathcal D } \to o^{\times} _{\mathcal D}  $ is a continuous character. \par
 We choose a division quaternion algebra $D $ which satisfies the following conditions:
\begin{itemize}
\item For a finite place $v$,
$D _{F_v}$ is split unless $\bar \rho \vert _{F_v} $ is of type $0_E$.
\item If $ \bar \rho \vert _{F_v} $ is of type $0_E$ at a finite place $v$, 
$ D_{F_v} $ is a quaternion algebra over $F_v$ with the invariant $\frac {c} {2} \mod \mathbb Z$, where $c$ is the relative conductor of $\bar \rho \vert _{F_v}$. 
\item $ q _{D  } = \sharp I _{D }  \leq 1$. 
\end{itemize}

A quaternion algebra $D$ which satisifies the above conditions exists. We make a choice from such quaternion algebras, and denote it by $D(\bar \rho)$.
$$ 
q_{D(\bar \rho)}  \equiv \sharp \{v: \ \bar \rho \vert_ {F_v }\text{ is of type } 0_E,\ \text{and the relative conductor
is odd}
\} + d \mod 2
$$
holds for $d= [F: \mathbb Q]$, and $D(\bar \rho )$ is split at all places dividing $\ell$. \par
$ G _{\bar \rho}$ is defined as $\Res _{F/ \mathbb Q} D (\bar \rho ) ^{\times}$. 
$K _{\mathcal D}$ is defined as the product $\prod _{v \in \vert F \vert _f } K_ {\deform _{\mathcal D} (v)  } (\bar \rho \vert _{F_v} ) $, and $\nu _{\mathcal D} :  K _{\mathcal D} \to o^{\times} _{\mathcal D} $ is the product of $ \nu _{\deform _{\mathcal D} ( v) }(\bar \rho \vert _{F_v} ):  K_ {\deform _{\mathcal D} (v) } (\bar \rho \vert _{F_v} ) \to o^{\times} _{\mathcal D}$. Here $ (K_ {* } (\bar \rho \vert _{F_v} ),  \nu _{*}(\bar \rho \vert _{F_v} ) )$ is a pair of the compact open subgroup of $D (\bar \rho) ^{\times} (F_v)  $ and the $K$-character defined in \S\ref{subsec-galdef8} for $ * = \bold f , \bold u $, and in \S\ref{subsec-nearlyordinary3} for $ * = \bold {n.o.f.}, \bold{n.o.}, \bold{fl}$.
\par
\begin{dfn}\label{dfn-modular21} For an absolutely irredicible representation $\bar \rho$ which satisfies {\bf MM1} and {\bf MM2}, let $(G _{\bar \rho} , K _{\mathcal D}, \nu _{\mathcal D}  )$ be a triple defined as above.
\begin{enumerate}
\item $I _{ (k, w ) , \mathcal D }$ is the intertwining space defined by
$$
I _{ (k, w ) , \mathcal D } = \Hom _{K_{\mathcal D} } ( \nu _{\mathcal D} , \mathcal S _{(k, w) } (G_{\bar \rho}  ) ) .
$$
\item For $ \pi \in \mathcal A _{(k, w) } (G_{\bar \rho} ) $, $ I _{ (k, w ) , \mathcal D } (\pi ) =  \Hom _{K_{\mathcal D} } ( \nu _{\mathcal D} , \pi _f)$.
\item For an algebraic Hecke character $\chi : F ^{\times } \backslash \mathbb A ^\times _F \to \bar E _\lambda  ^\times$,  $I _{ (k, w ) , \mathcal D , \chi }$ is the subspace which corresponds to $\mathcal S _{(k, w), \chi  } (G_{\bar \rho}  ) $.
\end{enumerate}
\end{dfn}
$ I _{ (k, w ) , \mathcal D } $ is isomorphic to $\oplus _{\pi \in \mathcal A _{(k, w ) } ( G _{\bar \rho} ) }  I _{ (k, w ) , \mathcal D } (\pi )$, and each $ I _{ (k, w ) , \mathcal D } (\pi ) $ is isomorphic to the tensor product of the local intertwining spaces $ \otimes _{v \in \vert F\vert _f } I _{\deform _{\mathcal D}(v) }(\bar \rho \vert _{F_v} , \pi _v ) $. 
\begin{dfn}\label{dfn-modular22} Let $\bar \rho$ be an absolutely irreducible representation which satisfies {\bf MM1} and {\bf MM2}. 
\begin{enumerate}
\item We define a decomposition of $ \Sigma_\mathcal D=  P^{\bold {fl} } \cup P^{\bold {n.o.} } \cup P ^{\bold f }_{\mathcal D}  \cup P_{\mathcal D} ^{\bold u } $ by $ P^{\bold {fl} } =  \deform_{\mathcal D} ^{-1}( 
\bold {fl})$, $ P^{\bold {n.o.} } =   \deform_{\mathcal D} ^{-1}( 
\{ \bold {n.o.}, \bold {n.o.f.} \} )$, $P ^{\bold f }_{\mathcal D}  =   \deform_{\mathcal D} ^{-1}( 
\bold {f} )$, and 
$P _{\mathcal D} ^{\bold u}  =  \deform  ^{-1} _{\mathcal D} (\bold {u}) $. \par
 $ P ^{\exc}= \{ \bar \rho \vert _{F_v} \text{ is of type } 0_E, \text{ and the relative conductor of } \bar \rho \text{ is odd.} \}  $ is a subset of $ P ^{\bold f }_{\mathcal D}  \cup P_{\mathcal D} ^{\bold u } $.

\item The Hecke algebra $\tilde T  _{\mathcal D }$ is the $o_{\mathcal D} $-algebra generated by the following elements in $\End_{o_{\mathcal D}  }  I _{(k, w) , \mathcal D }   $: $ T _v $ and $T_{v, v} $ for $v \not \in \Sigma _{\mathcal D} $, $U(p_v)  $ and $U (p_v, p_v ) $ for $v \in P ^{\bold u }_{\mathcal D }\setminus P ^{\exc}$, and $\tilde U ( p_v )$ and $ \tilde U (p_v, p_v )$ for $v \in P^{\bold {n.o.} } $.

The Hecke algebra $\tilde T_{\mathcal D, \chi }$ with a fixed central character $\chi$ is the image of $\tilde  T  _{\mathcal D   }$ in $\End_{o_{\mathcal D}  }  I _{(k, w) ,  \mathcal D , \chi } $. 

\item The reduced Hecke algebra $ \tilde T ^{\ss}  _{\mathcal D }$ is the $o_{\mathcal D} $-subalgebra of $\tilde T  _{\mathcal D  }$ generated by $ T _v $ and $T_{v, v} $ for $v \not \in \Sigma _{\mathcal D} $. $\tilde T^{\ss}_{\mathcal D, \chi }$ with a fixed central character $\chi$ is the image of $\tilde  T  _{\mathcal D  }$ in $\End_{o_{\mathcal D}  }  I _{(k, w) ,  \mathcal D , \chi } $

\item For any $o_{\mathcal D }$-algebra $R$, $\tilde T  _{\mathcal D  , R } =\tilde  T  _{\mathcal D  } \otimes_{o_{\mathcal D} } R$, $\tilde T  _{\mathcal D  , \chi  , R } = \tilde T  _{\mathcal D , \chi } \otimes_{o_{\mathcal D} } R$. We use the similar notation for the reduced Hecke algebras.
\end{enumerate}
\end{dfn}
By the definition, $ \tilde T ^{\ss}  _{\mathcal D  } $ (resp. $ \tilde T ^{\ss}  _{\mathcal D,  \chi  } $) is a subring of $  \tilde T  _{\mathcal D  } $ (resp. $ \tilde T  _{\mathcal D,  \chi  }  $), and $ \tilde T  _{\mathcal D,  \chi   }  $ (resp. $ \tilde T ^{\ss} _{\mathcal D,  \chi } $) is a quotient ring of $ \tilde T  _{\mathcal D  } $ (resp. $ \tilde T ^{\ss} _{\mathcal D  } $ ).\par

We define a maximal ideal $m_{\mathcal D}$ of $ \tilde T  _{\mathcal D }$ corresponding to
$\bar \rho$ as follows. 

For a minimal lift $\pi _{\min}$, an 
$o_{E_\lambda } $-algebra homomorphism
$f _{\pi _{\min} } : \tilde T^{\ss} _{ \mathcal D_{\min}  , \chi _{\bar \rho} } \to o_{E_\lambda}$ is defined. Here $  \chi _{\bar \rho}= \chi  ^{\Hecke}_{\bar \rho , w  } $ is as in Definition \ref{dfn-nearlyordinary31}. For $ v \in P ^{\bold {n.o.}}$, let $ \mu_{\pi _{\min} , v } $ be the nearly ordinary character of $\pi _{\min } $ at $v$ (Definition \ref{dfn-nearlyordinary21}), $\kappa_{\pi _{\min} , v } =\mu_{\pi _{\min} , v } \vert _{o^{\times}_{F_v}  } $ the nearly ordinary type. \par

By $\tilde m  _{\mathcal D}  $ we denote the maximal ideal of $\tilde  T _{ \mathcal D }$ 
generated by $ m_{o_\mathcal D} $ and the following operators:
\begin{itemize}
\item  $T_v  - f_{\pi _{\min}}  ( T_v ) ,\ T_{v, v}  - f _{\pi _{\min} } (T_{v, v} ) $
for $v \not \in \Sigma _{\mathcal D} $.
\item $ U (p_v),\ U (p_v, p_v ) -  \chi _{\bar \rho } ( p_v) $ for $ v \in P^{\bold u }_{\mathcal D}\setminus P ^{\exc}
$.  
\item $  \tilde  U ( p_v ) -\alpha_{\pi _{\min}, v } $,
$  \tilde  U ( p_v, p_v  ) -\gamma _{\pi _{\min}, v } $ for
$v 
\in P ^{\bold {n.o.}}
$. $ \alpha_{\pi _{\min}, v } = (\prod _{ \iota \in I _{F, v} } \iota (p_v) ^{-k ' _\iota } )\cdot \mu_{\pi _{\min} , v }  (p_v) $, $\gamma_{\pi _{\min}, v } =q^{-w}_v \chi _{\bar \rho } ( p_v) $.
\par
\end{itemize}

Here, $\tilde  m _{\mathcal D  , \chi}$ is the maximal ideal similarly defined for $\tilde   T_{ \mathcal D , \chi } $. 

\begin{lem}\label{lem-modular22}
 $\tilde m_{\mathcal D} $ and $ \tilde m_{\mathcal D , \chi_{\bar \rho}} $ are proper ideals of $ \tilde T_{\mathcal D}$ and $ \tilde T_{\mathcal D, \chi_{\bar \rho} }$, respectively. 
\end{lem}
\begin{proof}[Proof of Lemma \ref{lem-modular22}]
It suffices to prove it for $ \tilde T_{\mathcal D, \chi_{\bar \rho} }$.  For $\pi = \pi _{\min}$, we look at the intertwining space $ I _{ (k, w ) , \mathcal D } (\pi ) =\otimes _{v \in \vert F\vert _f } I _{\deform _{\mathcal D} (v) }(\bar \rho \vert _{F_v} , \pi _v ) $, and calculate the action of generators of $\tilde m_{\mathcal D , \chi _{\bar \rho}} $.   \par
The central character of $\pi _{\min}$ is $\chi _{\bar \rho}$. $T_{v, v} $, $U(p_v, p_v )$, $\tilde U (p_v , p_v) $ are evidently compatible with this central character. So we only consider the actions of $T_v $, $U (p_v) $, and $\tilde U (p_v) $. \par
For $v \not \in \Sigma _{\mathcal D}$, the action of $T_v $ on $I (\bar \rho \vert _{F_v} , \pi _v )  $ is $ f_{\pi _{\min}}  ( T_v )$.  \par
For $v \in P ^{\bold {n.o.}}$, $\alpha_{\pi _{\min}, v } $ is the $U (p_v)$-eigenvalue for a non-zero nearly ordinary vector, So the part in $  I _{\deform _{\mathcal D} (v) }(\bar \rho \vert _{F_v} , \pi _v ) $ where $U ( p_v)$ acts as $\alpha_{\pi _{\min}, v } $ is non-zero.\par
For $v \in  P^{\bold u }_{\mathcal D}\setminus P ^{\exc}$, by Proposition \ref{prop-galdef82} (2), $ I _{\bold u }(\bar \rho \vert _{F_v} , \pi _v ) $ has a non-zero subspace where $U (p_v)$ acts by zero map. So the claim follows. 
\end{proof}

\begin{dfn} \label{dfn-modular23} Let $\mathcal D $ be a deformation type of $\bar \rho$. 
\begin{enumerate}
\item The $\ell$-adic Hecke algebra $T_{\mathcal D} $ of $\bar
\rho$ with deformation type $\mathcal D$ is
$$
T_{\mathcal D} = (\tilde T _{ \mathcal D } )_{ \tilde m _\mathcal D  } , 
$$
and the $\ell$-adic Hecke algebra with a fixed central character $\chi$ is 
$$
T_{\mathcal D , \chi}  = (\tilde T _{ \mathcal D , \chi } )_{ \tilde  m _{\mathcal D , \chi} } . 
$$
\item The $\ell$-adic reduced Hecke algebra $ T^{\ss}  _{\mathcal D }$ (resp. $T^{\ss}  _{\mathcal D  , \chi } $) is the image of $\tilde T ^{\ss} _{\mathcal D } $ (resp. $\tilde T ^{\ss} _{\mathcal D , \chi } $) in $T_{\mathcal D}  $ (resp. $T_{\mathcal D, \chi }  $).
\end{enumerate}
\end{dfn}
By the definition, $ T^{\ss}  _{\mathcal D } $ is the $o_{\mathcal D}$-subalgebra of $ T _{\mathcal D } $ generated by $T _{v } , T_{v, v}  $ for $ v \not \in \Sigma _{\mathcal D} $, and the maximal ideal is generated by $ m_{ o _{\mathcal D} } $ and $T_v  - f_{\pi _{\min}}  ( T_v ) ,\ T_{v, v}  - f _{\pi _{\min} } (T_{v, v} ) $ for $ v \not \in \Sigma _{\mathcal D} $.\par
Let $ e_{\bar \rho}$ (resp. $ e_{\bar \rho , \chi }$)  be the idempotent in $\tilde T_{\mathcal D}
$ (resp. $\tilde T_{\mathcal D , \chi }$) which defines the direct factor $ T_\mathcal D $ (resp. $  T_{\mathcal D , \chi}  $). \par
\begin{dfn}\label{dfn-modular235}
The intertwining spaces which corresponds to $\bar \rho$ are defined by
$$
I _{ (k, w ) , \mathcal D } (\bar \rho )   =  e_{\bar \rho}  I _{ (k, w ) , \mathcal D }  ,
$$
and
$$  
I _{ (k, w ) , \mathcal D , \chi  } (\bar \rho )   =  e_{\bar \rho}  I _{ (k, w ) , \mathcal D  , \chi}  .
$$
\end{dfn}

\begin{lem}\label{lem-modular23}
\begin{enumerate}
\item $I _{ (k, w ) , \mathcal D } (\bar \rho )  $ is free of rank one as a $T _{\mathcal D,  \bar E_{\mathcal D} }=  T _\mathcal D 
\otimes_{o_{\mathcal D } } \bar E_{\mathcal D}$-module. The same is true for $ T _{\mathcal D, \chi _{\bar \rho},  \bar E_{\mathcal D} }=  T _{\mathcal D , \chi _{\bar \rho}}
\otimes_{o_{\mathcal D } } \bar E_{\mathcal D}$.
\item $T^{\ss} _{\mathcal D , E _{\mathcal D} }=T _{\mathcal D , E _{\mathcal D} }   $ and $T^{\ss} _{\mathcal D , \chi _{\bar \rho} , E _{\mathcal D} }=T _{\mathcal D , \chi_{ \bar \rho} ,  E _{\mathcal D} }   $ hold. 
\end{enumerate}
\end{lem}

\begin{proof}[Proof of Lemma \ref{lem-modular23}] We prove the lemma for $T _{\mathcal D}$. $T_{\mathcal D , \chi _{\bar \rho} }$ is treated similarly. \par
First we prove (1). 
By Lemma \ref{lem-modular22}, $I _{ (k, w ) , \mathcal D } (\bar \rho )  $ is non-zero. For an element $\pi$ of $ \mathcal A _{(k, w) }( G_D) $ which appears in the decomposition of $I _{ (k, w ) , \mathcal D } (\bar \rho )  $, we show that $e _{\bar \rho} I _{(k, w) , \mathcal D } (\pi ) $ is one dimensional over $\bar E_{\mathcal D} $. \par
$I _{(k, w) , \mathcal D } (\pi ) =\bigotimes _{v \in \vert F\vert _f } I _{\deform _{\mathcal D} (v)  }(\bar \rho \vert _{F_v} , \pi _v ) $, where $\pi_v$ the $v$-component of the restricted tensor product $ \pi _f = \bigotimes _{v \in \vert F\vert _f } \pi _v$. \par
Take an $\ell$-adic field $E' _{\lambda '}$ so that it contains $E_{\mathcal D} $ and $\pi _f $ is defined over $E' _{\lambda '}$. The associated $E' _{\lambda '}$-representation $ \rho ' = \rho _{\pi ,  E' _{\lambda '}} : G_F \to \GL_2 (  o_{E' _{\lambda '}})$ is a deformation of $\bar \rho$: 
the mod $\lambda '$-reduction $\bar \rho ' = \rho ' \mod \lambda '   $ has the same semi-simplification as $\bar \rho$ since the characteristic polynomial of $\Fr _v $ for $v \not \in \Sigma_{\mathcal D} $ are the same, and hence isomorphic over $k _{\lambda '}$. \par
 For $v \not \in \Sigma _{\mathcal D}   $, $ I _{\bold f} ( \bar \rho \vert _{F_v } , \pi _v )$ is one dimensional since $\pi _v$ is spherical.\par
 For $ v \in P ^{\bold f}_{ \mathcal D} $, by Proposition \ref{prop-galdef81}, it is one dimensional, since $ \pi _v$ corresponds to $\rho '  \vert _{F_v} $ which is a deformation of $\bar \rho$ by the compatibility of the local and the global Langlands correspondence. \par
\item For $v \in P^{\exc}\cap  P ^{\bold u} _{ \mathcal D}  $, $  I _{\bold u} ( \bar \rho \vert _{F_v } , \pi _v )$ is one dimensional by Proposition \ref{prop-galdef82}, (1).\par
 For $v \in P ^{\bold u}_ { \mathcal D}\setminus P ^{\exc} $, we need to consider $U(p_v)$-action on $  I _{\bold u} ( \bar \rho \vert _{F_v } , \pi _v )$. By Proposition \ref{prop-galdef82}, (2), it is isomorphic to $\bar E _{\lambda '} [U  ]  / ( U \cdot L ( U, \ \pi _v \otimes \nu_v ^{-1}) ) $, where $U$ acts as $U (p_v)$, and $\nu _v$ is a character of $F^\times _v $ defined in Proposition \ref{prop-galdef82}. By our definition of $\tilde m_\mathcal D$, $U(p_v )$ acts as zero mod $\tilde m_\mathcal D$. Since the all roots of the local $L$-factor $L ( U, \ \pi _v \otimes \nu_v ^{-1})$ have $\ell$-adic units, the localization at $ (\lambda ' , U(p_v)) $ of the $o _{E' _{\lambda '}} [ U (p_v) ]$-module $  I _{\bold u}( \bar \rho \vert _{F_v } , \pi _v )$ is one dimensional, and this is the part which contributes to $e _{\bar \rho} I _{(k, w) , \mathcal D } (\pi ) $. \par
 For $v \in P ^{\bold {fl}} $, $I _{\bold {fl}} ( \bar \rho \vert _{F_v } , \pi _v )$ is one dimensional since $\pi _v$ is spherical, and $K (\bar \rho \vert _{F_v} ) = \GL _2 (o_{F_v} ) $.\par
 For $v \in P ^{\bold {n.o.}}$, the argument is similar to the case of $ P ^{\bold u}_ { \mathcal D}$, using Proposition \ref{prop-nearlyordinary31} instead,
by specifying the one dimensional subspace where $\tilde U (p_v)$-operator acts by a scalar multiplication by an $\ell$-adic unit. \par

Since the action of $ T _{\mathcal D, \bar E_{\mathcal D} }$ on $I _{(k ,w ) ,\mathcal D  } (\bar \rho )  $ is faithful, the one-dimensionality of $e _{\bar \rho} I _{(k, w) , \mathcal D } (\pi ) $ shows that $ T _{\mathcal D }$ is reduced. It is easy to see that any $\bar E _{\mathcal D} $-homomorphism $f : T _{\mathcal D, \bar E_{\mathcal D} } \to \bar E_{\mathcal D}  $ is obtained in this way, and (1) is shown. (2) follows from (1) since $T ^{\ss}_{\mathcal D} $-action on each non-zero component $e _{\bar \rho} I _{(k, w) , \mathcal D } (\pi ) $ is non-trivial. 
\end{proof}
\begin{cor}\label{cor-modular21} For any element $\pi \in \mathcal A _{(k ,w ) } (G_{\bar \rho})$, 
 $e _{\bar \rho} I _{(k, w) , \mathcal D } (\pi ) $ is isomorphic to $ \bigotimes _{v \in \vert F\vert _f } I' _{\deform _{\mathcal D} (v)  }(\bar \rho \vert _{F_v} , \pi _v ) $. Here, $ I '  _{\deform _{\mathcal D} (v)  }(\bar \rho \vert _{F_v} , \pi _v )$ is the subspace of  $I_{\deform _{\mathcal D} (v)  }(\bar \rho \vert _{F_v} , \pi _v )$ defined as follows.
\begin{itemize}
 \item For $ v\not \in P ^{\bold {n.o.}}\cup P ^{\bold u } _{\mathcal D}  $, $I' _{\deform _{\mathcal D} (v)  }(\bar \rho \vert _{F_v} , \pi _v ) =I _{\deform _{\mathcal D} (v)  }(\bar \rho \vert _{F_v} , \pi _v )  $. 
\item For $v \in  P ^{\bold {u}}_{\mathcal D}$, $I ' _{\bold u } ( \bar \rho \vert _{F_v} , \pi _v ) $ is the subspace of $I _{\bold u  } ( \bar \rho \vert _{F_v} , \pi _v ) $ where $U (p_v) $ acts by zero map.
\item For $v \in  P ^{\bold {n.o.}}$, $I ' _{ \deform _{\mathcal D} (v ) } ( \bar \rho \vert _{F_v} , \pi _v ) $ is the $U (p_v) $-eigensubspace of $I _{\bold u } ( \bar \rho \vert _{F_v} , \pi _v ) $ of the unique $\ell$-adic unit eigenvalue.
\end{itemize}
 \end{cor}
This immediately follows from the proof of Lemma \ref{lem-modular23}.

\subsection {Galois deformations and $\ell$-adic Hecke algebras}\label{subsec-modular3}

We have defined $\ell$-adic Hecke algebras $ T ^{\ss}_\mathcal D \subset T_\mathcal D$. By the existence of the
$\ell$-adic representation attached to the elements of $\mathcal A _{(k,w) }(G_{\bar \rho} ) $ (cf. \S\ref{subsec-shimura4}), there exists a continuous representation
$$
\rho  ^{\modular} _{\mathcal D, E_\lambda} :  G_{\Sigma  _{\mathcal D } } \longrightarrow  \GL _2 ( T^{\ss} _{\mathcal D , E_ \lambda}) 
$$
such that 
$$	\begin{cases} \tr \rho  ^{\modular} _{\mathcal D, E_\lambda}\ (\Fr _v) \ = \ T _v 	\\
	       \det  \rho  ^{\modular} _{\mathcal D, E_\lambda} \  (\Fr _v ) \ = \ q_ v \cdot  T_{v, v}
	\quad  \text{for each $v \not \in \Sigma  _{\mathcal D }$}.
	\end{cases}
$$

\begin{lem}\label{lem-modular31} For $T^{\ss} _{\mathcal D }  $ attached to a deformation type $\mathcal D $, the following holds:
\begin{enumerate}
\item $ T^{\ss} _{\mathcal D } $ contains $\tr_{T ^{\ss} _{\mathcal D , E_\mathcal D}}
\rho  ^{\modular} _{\mathcal D, E_\lambda}(
\sigma )$ 
for $\sigma \in G _{\Sigma  _{\mathcal D } } $.
\item For a finite set of finite places $\Sigma $ which contains $\Sigma _{\mathcal D}$, $ T^{\ss} _{\mathcal D }$ is generated by $ T _v $, $v \not \in \Sigma $ as an $o_\mathcal D $-algebra. 
\end{enumerate}
\end{lem}
\begin{proof}[Proof of Lemma \ref{lem-modular31}] (Cf. \cite{W2}.) 
Let $ \tilde T $ be the normalization of
$T^{\ss} _{\mathcal D }$, $T _{\Sigma } $ the $o_{\mathcal D}$-subalgebra of $T^{\ss} _{\mathcal D } $ generated by $ T _v $, $v \not \in \Sigma $, $ T ' $ the $o_{\mathcal D}$-subalgebra of $\tilde T$ generated by $ \tr
\rho  ^{\modular} _{\mathcal D, E_\lambda}(
\sigma ),
\sigma \in G _{\Sigma  _{\mathcal D }} $. It suffices to see the images of $
T_{\Sigma}$ and $T '  $ in $ \tilde T / \lambda ^ n \tilde T$ is the same for any $ n
\geq 1 $.  $\rho  ^{\modular} _{\mathcal D, E_\lambda} $ is defined over $ \tilde T$, which we denote by $\tilde \rho$. By the Chebotarev density theorem, for any
$\sigma \in G_{\Sigma _{\mathcal D} }$, 
there is a finite place
$v \not \in \Sigma $ such that $\tilde \rho (\Fr _v)$ and $\tilde \rho (\sigma) $ is conjugate in
$\tilde T /
\lambda ^ n \tilde T$. This means that $ \tr_{\tilde T } \tilde \rho (\sigma ) \equiv
\tr_{\tilde T } \tilde \rho (\Fr _v)  \mod \lambda ^ n  $, and $T_{\Sigma} = T ' $. The equality $2q_v T _{v, v} = (\tr \rho (\Fr_v ) ) ^2 - \tr \rho (\Fr^2_v ) $ for $v \not \in \Sigma _{\mathcal D}$ implies that $  T '  = T^{\ss} _{\mathcal D }$. 
\end{proof}
\begin{prop}\label{prop-modular31}(Functoriality) For two deformation types $\mathcal D$ and $\mathcal D '$ of $\bar \rho$, assume that there is a morphism $\mathcal D \to \mathcal D ' $ in $\Type (\bar \rho ) $ (Definition \ref{dfn-galdef92}). Then there is a canonical surjective homomorphism $T _{\mathcal D ' } \to T _{\mathcal D } \otimes _{o_{\mathcal D}} o _{\mathcal D ' }  $. 
\end{prop}
This statement is non-trivial: $T _v $ and $T_{v, v} $-operators are missing in $T _{\mathcal D ' }  $ when $ \deform _{\mathcal D } (v ) = \bold f $, and $\deform _{\mathcal D' } (v) = \bold u $. By Lemma \ref{lem-modular31}, these missing operators are recovered from $ \rho  ^{\modular} _{\mathcal D, E_\lambda} $.\par
\begin{prop} \label{prop-modular32}
Assume that the residual representation $\bar
\rho$ of
$\rho  ^{\modular} _{\mathcal D, E_{\mathcal D} }:  G_{\Sigma  _{\mathcal D }}
\to
\GL _2 ( T ^{\ss} _{\mathcal D , E_{\mathcal D} }) $ is irreducible. Then there exists a representation
$$
\rho ^{\modular} _{\mathcal D } : G _{\Sigma   _{\mathcal D } } \to \GL _2 ( T^{\ss} _{\mathcal D } ) 
$$ 
having the same trace and determinant as $\rho  ^{\modular} _{\mathcal D , E_\lambda} $. 
\end{prop}
\begin{proof}[Proof of Proposition \ref{prop-modular32}] 
 Take a basis of $ ( T ^{\ss}  _{\mathcal D , E_{\mathcal D}}  ) ^{\oplus 2}
$ such that 
$$ 
\rho  ^{\modular} _{\mathcal D, E_{\mathcal D}}  \ ( c) = 
\begin{pmatrix} 
1& 0 \\ 
 0& -1
 \end{pmatrix} .
$$
Here,
$c $ is a complex conjugation. Then the entries of 
$$
\rho  ^{\modular} _{\mathcal D, E_{\mathcal D}}  \ ( \sigma ) = 
\begin{pmatrix} a(\sigma) & b(\sigma )\\ 
c(\sigma) & d(\sigma)
\end{pmatrix}
$$
have the property that $ b(\sigma)  \cdot c(\sigma ) $ is contained in $ T^{\ss} _{\mathcal D } $ and independent of a choice of basis. Using the irreducibility of
$\bar \rho$, one shows the existence of some $
\tau$ such that $ b (
\tau ) c (\tau ) $ is a unit. Change the basis again so that
$ b (
\tau ) = 1
$. Then
$\rho  ^{\modular} _{\mathcal D, E _{\mathcal D }} $ is defined over $  T ^{\ss} _{\mathcal D } $ using
this basis, since 
$$	
\begin{cases}
	a ( \sigma ) =  1/2 ( \tr \rho ^{\modular} _{\mathcal D, E _{\mathcal D }}   ( \sigma ) + \tr
\rho^{\modular} _{\mathcal D, E _{\mathcal D }}  (
\sigma \cdot c )) 	\\ d ( \sigma ) = 1/2 ( \tr \rho ^{\modular} _{\mathcal D, E _{\mathcal D }}  (
\sigma ) -\tr
\rho^{\modular}_{\mathcal D, E _{\mathcal D }}   ( \sigma \cdot c ) )	\\ b ( \sigma )= x ( \sigma,
\tau) / x (
\tau,
\tau) \\ c( \sigma ) = x( \tau, \sigma ) .
	\end{cases}
$$
Here $x( \sigma, \tau ) = a ( \sigma \tau ) -a(\sigma ) a (\tau )= b (
\sigma ) c (\tau) $. So the claim follows.
\end{proof}

The following lemma is proved in a similar way, so the proof is omitted.
\begin{lem}\label{lem-modular32}
Let $G$ be a group with $ c\in G ,\ c ^2 = 1 $, $R$ a local ring such that $2$ is a non-zero divisor,
$\rho, \rho '  : G \to \GL_2 ( R )$ two representations with the irreducible residual representations such that $ \det \rho ( c) = \det \rho '( c) = -1 $.  $\tr \rho ( g ) = \tr \rho ' ( g) $ for any $g \in G$, then $\rho $ and $\rho ' $ are
isomorphic.
\end{lem}
By Lemma \ref{lem-modular31}, $ \rho^{\modular} _{\mathcal D}$ is unique up to isomorphisms. 
\subsection{Cohomology groups of modular varieties as Hecke modules}\label{subsec-modular4}

We construct $T_{\mathcal D}$-module $M_{\mathcal D}$ from a cohomology group of the modular variety
associated to the deformation type $\mathcal D$. 
First we show how to go around the technical difficulty we encounter when a compact open subgroup is not
small by the method of auxiliary places.  

For a given absolutely irreducible $\bar \rho$, a deformation type $\mathcal D$, there are infinitely many finite places $y$ which satisfies the following properties:
\begin{itemize}
\item $ y \not \in \Sigma _\mathcal D $, and $q_y \geq a_F$, where $a_F$ is the integer given in Lemma \ref{lem-shimura51}.
\item $ q_y \not \equiv 1\mod \ell $, and the eigenvalues $\bar \alpha _y,\ \bar \beta _y  $ of $
\bar \rho ( \Fr _y ) $ counted with the multiplicities satisfy
$\bar \alpha _y \neq q_y  ^{\pm 1 } \bar \beta _y$. 
\end{itemize}
Note that we do not exclude the case when $\bar  \alpha _y$ equals $\bar \beta _y $.\par
The existence of finite places with these two properties follows from \cite{DT2}, lemma 11 (in the reference, the linear independence of $F$ and $\mathbb Q(\zeta_\ell) $ is assumed, which can be removed by a slight modification. See also \cite{Fu1} for an argument which is independent of the classification of subgroups of $\GL_2 (\F_q )$). \par

For an auxiliary place $y$ that satisfies the conditions as above, define the deformation function $\deform _{\mathcal D_y} $ by 

$$
  \deform _{\mathcal D _y} ( v ) =
\deform _{\mathcal D} (v ) \text { if } v \neq y ,\  \deform _{\mathcal D _y } ( y )  = \bold u , 
$$
and let $\mathcal D_y$ be the deformation type which has the defomation function $\deform _{\mathcal D_y} $ as above, and the same data as $\mathcal D$ except for the deformation function.

There is a natural surjective homomorphism $T_{\mathcal D _y }  \twoheadrightarrow  T_{\mathcal D} $
by Proposition \ref{prop-modular31}. By the conditions imposed on $y$, this is in fact an isomorphism:
\begin{prop}\label{prop-modular41} There is a natural $o_{\mathcal D } $-algebra isomorphism
$$  
T_{\mathcal D _y  }  \simeq T_{\mathcal D}. 
$$
The same is true for $ T_{\mathcal D _y , \chi _{\bar \rho }  } $ and $T_{\mathcal D , \chi _{\bar \rho} }$.
\end{prop}
\begin{proof}[Proof of Proposition \ref{prop-modular41}] We prove the claim for $T_{\mathcal D _y  }$. The case of $ T_{\mathcal D _y , \chi _{\bar \rho }  } $ follows from it.  Take an element $\pi $ of $\mathcal A _{(k,w)} ( G_D ) $ which contributes to $I _{(k, w), \mathcal D _y } ( \bar \rho )  $ non-trivially. We show that $\pi  $ appears in $I _{(k, w), \mathcal D  } ( \bar \rho ) $. By enlarging $E_{\mathcal D} =E_{\mathcal D_y} $ if necessary, we may assume that $\pi  _f $ is defined
over $E_{\mathcal D}$. \par
 The Galois representation $ \rho = \rho _{\pi , E_{\mathcal D } }: G_F \to \GL _2 ( o _{ E_{\mathcal D} } )   $ attached to $\pi$ is a deformation of $\bar \rho$, and the local representation $\rho \vert _{F_y}  $ corresponds to the $y$-component $\pi _y$ of $\pi$ by the compatibility of the local and the global Langlands correspondences. \par

If $\pi _y $ is supercuspidal, $ \rho _y = \rho \vert _{F_y}  $ must be induced from a tame quadratic extension 
$\tilde F_y$ of $F_y$: $ \rho _y \vert _{G_{\tilde F _y} } \simeq  \chi \oplus \chi '$,
$\chi ' $ is the twist of $\chi$ by the non-trivial element in $ \Gal (\tilde F_y /
F_y)$. \par
$\chi / \chi ' $ must have an $\ell$-power order, otherwise any mod $ m_{o_{\mathcal D}}$ reduction of $\rho _y$ with respect to an $o_{\mathcal D}$-lattice is absolutely irreducible. Thus we have 
$$
(\bar \rho _y) ^{\ss} =  \bar \chi\otimes \Ind ^{G_{F_y}} _{G_{\tilde F_y}} 1 \simeq  \bar \chi \oplus \bar \chi \cdot \chi _{\tilde F_y/F_y}. 
\leqno{(\ast)}
$$
Here $ \bar \chi$ is the mod $ m_{o_{\mathcal D}}$ reduction of $\chi$, and $ \chi _{\tilde F_y/F_y}$ is the quadratic character of $G_{F_y}$ which corresponds to $\tilde F_y/F_y$. 
Since $\bar \rho _y$ is unramified, $ \tilde F_y$ is the unramified quadratic extension
of $F_y $, which implies that $q_y\equiv -1 \mod \ell $. ($\ast$) implies that the two eigenvalues of $\bar \rho _y $ are of the form $\bar \alpha ,
\ - \bar \alpha
$, and this contradicts to our choice of $y$ since $ q_y \equiv -1 \mod \ell$. \par
Since the Swan conductor remains the same under a mod $m_{o_{\mathcal D}}$-reduction, $ \rho _y$ is a tame representation, and hence $ \pi_y$ is a special representation twisted by a tame character, or a tame principal series
representation.  In the former case, the two eigenvalues of any frobenius lift is of the form $
\alpha, \ q_y
\alpha
$, which is a contradiction. Since $q_y \not \equiv 1 \mod \ell$, and $\bar \rho _y $ is unramified, $\pi_y$ is spherical. The eigenvalues of $T_y$ and $T_{y, y} $-operators on $\pi ^{K_{\mathcal D , y} } _y  $ are congruent to $\tr \bar \rho (\Fr _y)  $ and $q^{-1}_y \det \bar \rho (\Fr _y) $ modulo $ m_{o_{\mathcal D}}$, thus $\pi $ must contribute to $ I _{(k, w), \mathcal D  } ( \bar \rho ) $. 
\end{proof}
For an element $\pi \in \mathcal A _{(k, w ) } (G_{\bar \rho})$, we modify the intertwining space $I _{(k, w) , \mathcal D_y } $ and $I _{(k, w) , \mathcal D_y } (\pi) $ to 
$$
I^y _{(k, w) , \mathcal D_y  }= \Hom _{ K _{\mathcal D }(y) } ( \nu _{\mathcal D_ y } \vert _{ K _{\mathcal D }(y) } , \mathcal S _{(k, w)} ( G _{\bar \rho} ) ), 
$$
and 
$$
I^y _{(k, w) , \mathcal D_y  } (\pi ) = \Hom _{ K _{\mathcal D }(y) } ( \nu _{\mathcal D_ y } \vert _{ K _{\mathcal D }(y) } ,  \pi _f ).
$$
$ I^y _{(k, w) , \mathcal D_y  }$ (resp. $I^y _{(k, w) , \mathcal D_y  } (\pi ) $) contains $ I _{(k, w) , \mathcal D_y } $ (resp. $ I _{(k, w) , \mathcal D_y  } (\pi )$) as a subspace. The following lemma shows that the part corresponding to $\bar \rho$ is the same for $I^y _{(k, w) , \mathcal D_y  } $ and $I _{(k, w) , \mathcal D_y  }$.
\begin{lem}\label{lem-modular41} For an element $\pi \in \mathcal A _{(k, w ) } (G_{\bar \rho})$ which is defined over $ E_{\lambda }$, if $I^y _{(k, w) , \mathcal D_y  } (\pi ) \neq \{ 0 \}$ and $\rho _{\pi , E_{\lambda} }$ is a deformation of $\bar \rho$, then $I^y _{(k, w) , \mathcal D_y  } (\pi )= I_{(k, w) , \mathcal D_y  } (\pi )$. 
\end{lem}
\begin{proof}[Proof of Lemma \ref{lem-modular41}] As in the argument of the proof of Proposition \ref{prop-modular41}, the $y$-component $\pi _y$ of $\pi$ is spherical, so it suffices to see $ \pi ^{K_{11} ( m^2_y ) } _y = \pi ^{K_1 ( m^2_y)} _y$. Since the central character of $\pi_y$ is unramified, the action of central elements $ o^{\times}_{F_y }  \cdot 1_{D(\bar \rho ) _{F_y} } $ is trivial. The claim follows. 
\end{proof}
By Lemma \ref{lem-modular41}, we regard $ T_{\mathcal D _y }$ as the Hecke algebra acting on $ I^y _{(k, w) , \mathcal D_y  } $.\par
\bigskip
We define a $T_{\mathcal D}$-module $M_{\mathcal D}$ for deformation type $\mathcal D$. \par
 We modify the $y$-component of $K_{\mathcal D _y }$, and use $\tilde K_{\mathcal D _y } = K _{11} (m^2 _y) \cdot K ^y_{\mathcal D}$. Fix a finite set of finite places $P_0$. We choose a finite set $S$ which satisfies the condition of Lemma \ref{lem-shimura51} for $K =\tilde  K _{\mathcal D_y}$ and $P_0 $. By Proposition \ref{prop-shimura51}, $\tilde K_{\mathcal D_y} ( y, S) = K_{11} ( m^2_y )\cdot \prod  _{u \in S} U_u  \cdot K ^{\{y \} \cup S } _{\mathcal D _y } $ is small, and $H ^{q} _{\stack}( S_K \  \bar {\mathcal F} ^K _{(k, w) , R}) $  is defined for any compact open subgroup $K$ of $ \tilde K_{\mathcal D _y }$ and a finite $o_{\mathcal D}$-algebra $R$ (cf. \S\ref{subsec-shimura5}). \par
For $\Sigma = \Sigma _{\mathcal D_y } \cup S $ and $ K = \ker \nu _{\mathcal D_y } \vert _{\tilde K _{\mathcal D _y } }$, the convolution algebra $H _{K^{\Sigma}} =  H (G_D ( \mathbb A_{\mathbb Q, f}) , K ^{\Sigma} ) _{o_{\mathcal D} } $ acts on $H ^{q} _{\stack}( S_K \  \bar {\mathcal F} ^K _{(k, w) , R}) $. Let $\tilde m _{\Sigma , \bar \rho }  $ be the maximal ideal of $H _{K^{\Sigma}}$ which corresponds to $\bar \rho$ (cf. \S\ref{subsec-coh1}). By Proposition \ref{prop-coh21}, the localization $ H ^{q} _{\stack}( S_K \  \bar {\mathcal F} ^K _{(k, w) })_{\tilde m _{\Sigma , \bar \rho } }$ at $\tilde m _{\Sigma , \bar \rho } $ vanishes unless $q= q_{\bar \rho} $, where $q _{\bar \rho} =q_{D(\bar \rho )}$.  $H ^{q_{\bar \rho}} _{\stack}( S_K \  \bar {\mathcal F} ^K _{(k, w) })_{\tilde m _{\Sigma , \bar \rho } }$ is $o_{\mathcal D_y}$-free, and the formation commutes with scalar extensions. \par
\begin{dfn}\label{dfn-modular405} We define $\tilde M^y _{\mathcal D_y }$ by
$$
\tilde M^y _{\mathcal D_y } = \Hom _{\tilde K_{\mathcal D_y} } ( \nu _{\mathcal D_ y } \vert _{\tilde K _{\mathcal D_y} } , H ^{q _{\bar \rho }}_{\stack} ( S _{ K } ,  \bar {\mathcal F} ^ { K  } _{(k, w) }   ) _{\tilde m _{\Sigma , \bar \rho } }).
$$
\end{dfn}
Since the order of $ \nu _{\mathcal D_ y }$ is prime to $\ell$, $\tilde M^y _{\mathcal D_y } $ is regarded as an $o _{\mathcal D_y}$-direct summand of $ H ^{q _{\bar \rho }}_{\stack} ( S _{ K } , 
 \bar {\mathcal F} ^ { K  } _{(k, w) }   )$ where $\tilde K _{\mathcal D_y }  $ acts by $ \nu _{\mathcal D_ y }$.\par
For $ \tilde M^y _{\mathcal D_y }$ defined as above, the following proposition holds by the decomposition in Proposition \ref{subsec-shimura4}.
\begin{prop}\label{prop-modular415} We have
$$
 \tilde M^y _{\mathcal D_y }\otimes _{o_{\mathcal D_y } } \bar E _{\mathcal D_y } = \bigoplus _{\pi \in \mathcal A _{\bar \rho} } I ^y _{(k, w ) , \mathcal D_y } (\pi ) ^{\oplus 2 ^{q_{\bar \rho }}},  
$$
where $\mathcal A _{\bar \rho}$ is the subset of $ \mathcal A _{(k, w)} ( G_D )  $ consisting of the representations $\pi $ whose associated Galois representation $\rho _{\pi , \bar E _{\mathcal D_y } } $ is a deformation of $\bar \rho$.
\end{prop}
 As $\bar \rho$ is not of residual type, the contribution from $\mathcal A ^c_{(k, w) } ( G_D) $ disappears. 
The standard Hecke operators $T_v , T_{v, v} $ for $v \not \in \Sigma $ and $\tilde U (p_v) , \tilde U (p_v , p_v )$-operators for $v \in \Sigma _{\mathcal D_y}\setminus P ^{\exc}$ acts on $ \tilde M^y _{\mathcal D_y }$. 
By Lemma \ref{lem-modular31} (2), the image of $H _{K^{\Sigma}} $ in $  \End _{o _{\mathcal D_y } }\tilde M^y _{\mathcal D_y }$
is $ T ^{\ss} _{\mathcal D_y }$ (the missing $T_v  $ and $T_{v, v}$-operators for $v \in S$ are recovered from $\rho ^{\modular}_{\mathcal D_y } $). Let $\tilde T^y_{\mathcal D_y }   $ be the $T ^{\ss} _{\mathcal D_y }$-algebra generated by $  U ( p_v), U (p_v, p_v ) $ for $ v \in P ^{\bold u}_{\mathcal D_y} \setminus P ^{\exc}$ and $  \tilde U ( p_v), \tilde U (p_v, p_v ) $ for $ v \in P ^{\bold {n.o.}}$ in $  \End _{o _{\mathcal D_y }} \tilde M^y _{\mathcal D_y }$. If we define the maximal ideal $\tilde m ^y _{\mathcal D_y }   $ of $\tilde T_{\mathcal D_y }(y)$ similarly as in \S\ref{subsec-modular2}, by Lemma \ref{lem-modular41}, the localization of $\tilde T^y_{\mathcal D_y }  $ at $\tilde m ^y_{\mathcal D_y } $ is identified with $ \tilde T _{\mathcal D_y } $.
\begin{dfn}\label{dfn-modular42} $M _{\mathcal D}$ is a $  \tilde T _{\mathcal D_y } $-module given by
$$
 M _{\mathcal D } =  ( \tilde M^y _{\mathcal D_y } ) _{\tilde m ^y _{\mathcal D_y}  }. 
$$ 
\end{dfn}
Note that $M_{\mathcal D }$ is a $T_{\mathcal D} $-module by the identification $T_{\mathcal D}\simeq T_{\mathcal D_y }  $ given in Proposition \ref{prop-modular41}, and $M_{\mathcal D , E_{\mathcal D}  }= M _{\mathcal D } \otimes _{o _\mathcal D} E_{\mathcal D}  $ is free of rank $ 2 ^ {q_{\bar \rho} } $ over $T _{\mathcal D, E_{\mathcal D}} $ by the decomposition in Proposition \ref{prop-modular415}, Lemma \ref{lem-modular23}, and Lemma \ref{lem-modular41}.

\begin{rem}\label{rem-modular5}
It is possible to show that $M_{\mathcal D } $ is independent of a choice of auxiliary places $y$ and $S$ up to isomorphisms.
\end{rem}

The following theorem, which is proved in the next subsections, is the main result in this section. 
\begin{thm} \label{thm-modular41}
Assume that $ \bar \rho  $ is an irreducible modular representation of type $\mathcal D$, and the validity of Hypothesis \ref{hyp-nearlyordinary21} for the representations $\pi $ of $ G _D (\mathbb A _{\mathbb Q, f} ) $ which contributes to $T^{\ss} _{\mathcal D, \bar E_{\mathcal D}} $. Then we have the following:
\begin{enumerate}
\item The representation
$\rho ^{\modular} _{\mathcal D }  $ is a deformation of $\bar \rho$ of type $ \mathcal D $. 
\item $T^{\ss} _{\mathcal D}  = T _{\mathcal D} $.
\end{enumerate}
\end{thm}
Assume the validity of this theorem for the moment. 
Let $R_{\mathcal D } $ be the universal deformation ring of $\bar \rho$ of deformation type $\mathcal D$. 
By the universality of $R _{\mathcal D} $, we have the $o_{\mathcal D}$-algebra homomorphism 
$$
 \pi _{\mathcal D} : R _{\mathcal D} \longrightarrow  T_{\mathcal D} 
$$
which corresponds to $ \rho^{\modular} _{\mathcal D}$ by Theorem \ref{thm-modular41}, (1). $\pi _{\mathcal D} $ is surjective by Theorem \ref{thm-modular41}, (2).
So $M_{\mathcal D}$ is regarded as an $R_{\mathcal D }$-module, and we will construct a Taylor-Wiles system for the pair $ (R _{\mathcal D} , M_{\mathcal D} )$ in the next section.

\subsection {Local conditions outside $\ell$}\label{subsec-modular5}

In this subsection, we prove Theorem \ref{thm-modular41} at a finite place $ v \nmid \ell$.  If $ \deform_{\mathcal D}  (v) =
\bold u $, (1) is clear. For (2), we show that $U( p_v) $ and $ U( p_v, p_v)$ belong to $T^{\ss}_{\mathcal D}$.
The $U( p_v)$-operator is zero. $ U (p_v , p_v ) = q^{-1}_v \det \rho ^{\modular}_{\mathcal D} ( p_v)$ by the identification $G^{ \ab}_{F_v} \simeq \widehat {F^{\times}_v } $ by the local class field theory. \par
When $ \deform_{\mathcal D}  ( v) = \bold f$, (2) is clear, so we show that $\rho ^{\modular}_\mathcal D $ is
a finite deformation of $\bar \rho$ at $v$. \par
Take an element $\pi  $ of $\mathcal A_{(k, w)}  (G_D) $ which appears in the component of  $T ^{\ss}_{\mathcal D , \bar E_{\mathcal D} }$. By extending $E_{\mathcal D} $ if necessary, we may assume that $ \pi _f $ is defined over $ E_{\mathcal D} $, and the associated Galois representation $\rho _{\pi , E_{\mathcal D} } : G_F \to \GL_2 (o_{\mathcal D}  ) $ is a deformation of $\bar \rho$.\par
By Corollary \ref{cor-modular21}, $ e_{\bar \rho } I _{(k, w ) , \mathcal D } (\pi ) $ is isomorphic to $\otimes _{u \in \vert F \vert _f} I '_{\deform_{\mathcal D } (u)  }(\bar \rho \vert _{F_u} , \pi _u  )  $, and $I '_{\deform_{\mathcal D } }(\bar \rho \vert _{F_u} , \pi _u  ) = I _{\deform_{\mathcal D }(u)  }(\bar \rho \vert _{F_u} , \pi _u  )  $ if $\deform _{\mathcal D } (u ) = \bold f $. By the compatibility of the local and the global Langlands correspondences, $\pi _v $ corresponds to $\rho _{\pi , v} = \rho _{\pi , E_{\mathcal D} }\vert _{F_v} $, which is a deformation of $\bar \rho \vert _{F_v}$. Since we assume that $ e_{\bar \rho } I _{(k, w ) , \mathcal D } (\pi )$ is non-zero, $ I _{\bold f} (\bar \rho \vert _{F_v} , \pi _v ) $ is non-zero. By Proposition \ref{prop-galdef81}, $\rho _{\pi ,v}$ is a finite deformation of $\bar \rho \vert _{F_v} $. \par
If $\bar \rho \vert _{F_v} $ is of type $0_{NE}$, we need to show that $\det \rho ^{\modular} _{\mathcal D}  $ is unramified at $v$  (see Definition \ref{dfn-galdef21}). It suffices to check this over $E_{\mathcal D}$, and for any element $\pi $ of $\mathcal A_{(k, w)}  (G_D) $ that contributes to $e_{\bar \rho } I _{(k, w ) , \mathcal D } $, $\det \rho _{\pi , v} $ is unramified since it is a finite deformation of $\bar \rho \vert _{F_v} $ (Definition \ref{dfn-galdef21}). \par
If $\bar \rho \vert _{F_v} $ is of type $0_{E}$, it is of the form $ \Ind^{G_{F_v} } _{G_{\tilde F_v} }\bar \psi$ for the unramified quadratic extension $\tilde F_v $ of $F_v$. We must show that $ \rho ^{\modular} _{\mathcal D}  \vert _{I_{\tilde F_v}} $ is isomorphic to the sum of $ \psi \vert _{I_{\tilde F_v}}$ and the Frobenius twist (cf. Definition \ref{dfn-galdef21}). Here $\psi$ is the Teichm\"uller lift of $\bar \psi $. Again it suffices to see this over $E_{\mathcal D}$, and the claim follows since $\rho _{\pi ,v}$ is a finite deformation of $\bar \rho \vert _{F_v} $ for any element $\pi $ which appears in the component of $T ^{\ss}_{\mathcal D , \bar E_{\mathcal D} }$.\par
If $\bar \rho \vert _{F_v} $ is absolutely reducible with the twist type $\bar \kappa _v$, we need to show that the $I_{F_v}$-fixed part $\mathcal L _v  $ of $\rho ^{\modular} _{\mathcal D}  \vert _{I_{F_v} } \otimes \kappa ^{-1} _v $ is $T^{\ss}_{\mathcal D}$-free, is a $T^{\ss}_{\mathcal D}$-direct summand, and the $T^{\ss}_{\mathcal D}$-rank of $\mathcal L_v$ is equal to $d_v = \dim _{k} ( \bar \rho \vert _{I_{F_v}} \otimes \bar \kappa ^{-1} _{\mathcal D, v} ) ^{I_{F_v}} $. By the same argument as the previous cases, $\mathcal L _v \otimes _{o_{\mathcal D} } E_{\mathcal D} $ is $T^{\ss}_{\mathcal D , E_{\mathcal D} }$-free of rank $d_v$. If $d_v= 2$, that is, in the case of $2_{PR} $, it follows that $ \mathcal L_v = (T^{\ss}_{\mathcal D} ) ^{\oplus 2}$, and the claim is shown. \par
So the case of $d_v= 1$ is left. It is clear that $\mathcal L_v $ is an $o_{\mathcal D}$-direct summand of $ (T^{\ss}_{\mathcal D} ) ^{\oplus 2}$, and $G_{F_v}$-stable as a subrepresentation of $\rho^{\modular} _{\mathcal D} \vert _{F_v}  $. Thus we have an exact sequence of $G_{F_v}\times T ^{\ss}_{\mathcal D} $-modules
$$
0 \longrightarrow \mathcal L _v \longrightarrow \rho^{\modular} _{\mathcal D} \vert _{F_v} \longrightarrow \mathcal L ' _v\longrightarrow 0 , 
$$
where $\mathcal L ' _v $ is $o_{\mathcal D}$-free, and has the $T^{\ss}_{\mathcal D  } $-rank one. The $I_{F_v}$-action on $\mathcal L ' _v\ $ is $\det \rho ^{\modular} _{\mathcal D} \vert _{I_{F_v}} \otimes \kappa^{-1} _{\mathcal D, v}  $ by calculating it over $E_{\mathcal D} $. \par
$\mathcal L ' _v / m _{\mathcal D }\mathcal L ' _v$ is a quotient of $\bar \rho $, and the $G_{F_v}$-action on $\mathcal L ' _v / m _{\mathcal D }\mathcal L ' _v$ is $\det \bar \rho \vert _{F_v} \otimes \bar \kappa^{-1}_v $. Since $\bar \rho \vert_{F_v} $ has a unique one dimensional quotient where $G_{F_v}$ acts via $\det \bar \rho \vert _{F_v} \otimes \bar \kappa^{-1}_v $, $ \dim _{k }\mathcal L ' _v / m _{\mathcal D }\mathcal L ' _v = 1$. By
Nakayama's lemma, there is a surjection $T ^{\ss}_{\mathcal D} 
\twoheadrightarrow 
 \mathcal L' _v$. Since the rank of $\mathcal L' _v $ over $E_{\mathcal D}$ is one, $
\mathcal L' _v $ is $T^{\ss} _{\mathcal D}$-free of rank one, and hence $ \mathcal L_v$ is also $T^{\ss} _{\mathcal D}$-free of rank one. This shows that $ \rho ^{\modular}_{\mathcal D}$ is a finite deformation at $v$. 

\subsection{Local conditions at $\ell$} \label{subsec-modular6}
We prove Theorem \ref{thm-modular41} at $ v \vert \ell$ when $q_{\bar \rho} $ is one, under the Hypothesis \ref{hyp-nearlyordinary21} if $\deform_{\mathcal D} (v) = \bold{n.o.}$. \par
We denote the residual characteristic at $v$ by $p$.  To emphasize the difference with the previous cases, we mainly use $p$ instead of $\ell$. \par

Assume that $\deform _{\mathcal D}(v)  = \bold { n.o.} $. Let $ ( k  _v , w ) = ( ( k_\iota )_{\iota \in I_{F, v}} , w) $ be the $v$-type, and $k'_v = (k'_\iota )_{\iota \in I_{F, v}}$ is as in Definition \ref{dfn-nearlyordinary11}.\par
By Hypothesis \ref{hyp-nearlyordinary21},  for any representation $\pi \in \mathcal A_{(k , w ) }  ( G_{\bar \rho} ) $ which appear as a component of $ T^{\ss} _{\mathcal D , E_{\mathcal D} }$, $\rho _{\pi , E_{\mathcal D} }\vert _{F_v} $, contains a subrepresentation $ L _{k ' _v} ( \chi_v)$ with multiplicity one. Here
$\chi_v$ is the nearly ordinary character of $\pi$ at
$v$, and $\chi_v \vert _{o^{\times}_{F_v}}  $ is the automorphic nearly ordinary type of $\bar \rho $ at $v$ by Definition \ref{dfn-nearlyordinary34}. By choosing this subrepresentation for each $\pi$, we may assume that 
 $
\rho ^{\modular}_{\mathcal D  , E_{\mathcal D} } \vert_{F_v}$ has the form
$$
0 \longrightarrow \mathcal L_{E_\lambda}  \longrightarrow   \rho ^{\modular}_{\mathcal D  , E_{\mathcal D} } \vert_{F_v} \longrightarrow \mathcal L ' _{E_\lambda}
\longrightarrow 0 ,
$$
where $ \mathcal L_{E_\lambda},\   \mathcal L ' _{E_\lambda} $ are both free of rank one over $
T_{\mathcal D, E_{\mathcal D} }$. By the same argument as in the case of $d_v = 1 $ of \S\ref{subsec-modular5} using that $\bar  \rho \vert _{F_v}$ is $G_{F_v}$-distinguished, $ \mathcal L = (T ^{\ss}_{\mathcal D }
) ^{\oplus 2} \cap \mathcal L _{E_\mathcal D} 
$ and $ \mathcal L ' = (T^{\ss} _{\mathcal D } )^{\oplus 2}/ \mathcal L $ are free of rank one over $T ^{\ss}_{\mathcal D} $. This shows that
$ \rho_{\mathcal D}^{\modular}\vert _{F_v}$ is a nearly ordinary deformation of $\bar \rho\vert _{F_v} $, and the restriction of the nearly ordinary
character of $  \rho_{\mathcal D}^{\modular}\vert _{F_v}$ to $I_{F_v}$ is $\kappa_{\mathcal D, v} $ (\S\ref{subsec-modular2}, {\bf D2}). The elements $\tilde U(p_v )$ and $\tilde U (p_v, p_v ) $ of $T _{\mathcal D} $ are equal to the action of $
\rho^{\modular}  _{\mathcal D} (\sigma)  $ on $\mathcal L $ and $ \det \rho ^{\modular}_{\mathcal D }\cdot  \chi^w _{\cycle} (\sigma)$ respectively, where $\sigma$ is the element of $G^{\ab}_{F_v} $ corresponding to $p_v$. Thus $\tilde U(p_v )$ and $\tilde U(p_v , p_v)$ belong to $ T ^{\ss} _{\mathcal D} $. \par
\bigskip

Next we consider the case when $ \deform (v) =\bold {n.o.f.}  $, by assuming that $q_{\bar \rho}$ is one. 
In this case, we assume that the $v$-type $(k  _v, w) $ is $( ( 2, \ldots, 2  ) , w) 
$ (see Definition \ref{dfn-nearlyordinary32}), in particular $w$ is even.  By twisting by $\chi _{\cycle} ^{ \frac w 2} $, we may assume that 
$w= 0$. We prepare several lemmas.\par
\begin{lem} \label{lem-modular61} Let $D$ be a division quaternion algebra over $F$  with $q_D= 1$ which is split at the places dividing $p$, $K$ an $F$-factorizable small compact open subgroup of $G_D ( \mathbb A_{\mathbb Q, f} ) $ such that the $v$-component $K_v$ is a maximal hyperspecial subgroup of $D^\times  (F_v) $ for some place $v$ dividing $p$. 
\begin{enumerate}
\item The canonical model $S_{K, F_v} = S _K ( G_D (\mathbb A _{\mathbb Q , f } ) , X_D )_{F_v}$ over $F_v$ has a projective smooth model $ S_ {K, o_{F_v}}$ over $ o_{F_v}$, which is unique up to canonical isomorphisms. For an inclusion $K ' \hookrightarrow K $ of compact open subgroups with $ K ' _v = K _v $, the induced morphism $S_{K ' , o_{F_v} } \to S_{K  , o_{F_v} } $ is \'etale. 
\item For a discrete infinity type $(k, w)$, assume that the $v$-type $ (k_v , w ) $ is 
$ ( (2, \ldots, 2) , 0 )$. Then the $o_{E_\lambda}$-smooth sheaf $ \mathcal F ^K _{(k, w)} $ on $S_{K, F_v}$ admits a unique extension $ \tilde {\mathcal F} ^K _{(k, w)} $ as an $o_{E_\lambda}$-smooth
\'etale sheaf to $S_ {K, o_{F_v}} $. 
\end{enumerate}
\end{lem}
\begin{proof}[Proof of Lemma \ref{lem-modular61}]
(1) is a consequence of \cite{Car1}. To be precise, in \cite{Car1}, to discuss integral models over $o_{F_v} $, compact open subgroups are assumed to be ``sufficiently small'' outside $v$, which is stronger than our notion of smallness, though we can reduce to the case stated in the reference. See \cite{Fu1}, \cite{Ja1} for the argument. \par
We show (2). The covering $\pi^v_p : \tilde S \to S_{K , F_v} $ which corresponds to $\overline {(F ^\times \cap K )} \backslash
\prod_{ u \vert p, u \neq v } K _u $ extends to a (pro-) \'etale coverings of $\tilde S_K$ by (1). As in \S\ref{subsec-shimura2}, $ \mathcal F _{(k, w)} $ is defined by an $o_{E_\lambda} $-lattice of
 $ G_D (\mathbb Q_{p} )$-representation
$V _{(k, w ) , E_\lambda }= \otimes _{ \iota : F_u \hookrightarrow E_{\lambda}, u \vert p }( \iota \det ) ^ { -k ' _\iota} \cdot \Sym ^{ k _\iota -2 }   ( V _{o_u } \otimes _{\iota} E_\lambda  ) ^{\vee} 
$. By the assumption on the $v$-type $  (k_v , w ) $, $V _{( k, w ) , E_\lambda } $ is regarded as a representation of $\prod _{ u \vert \ell, u \neq v } K _u $. Thus $ \mathcal F_{(k, w) } $ is defined by the covering $\pi ^v _p $, and the claim follows. 
\begin{lem}
\label{lem-modular62} Assumptions are as in Lemma \ref{lem-modular61}. Let $\Sigma $ be a set of finite places which contain $ \Sigma _K $ and $\{ v : v \vert p  \} $, $T _\Sigma = H ( D ^{\times} (\mathbb A_{F, f} ) , K ^\Sigma ) _{o_{E_{\lambda}} }$ the convolution algebra over $o_{E_\lambda} $, and $m _{\Sigma }$ a maximal ideal of $T_{\Sigma } $, which is not of residual type.  
Then there is a unique $p$-divisible group $ E$ over $o_{F_v}$ such that $T_p (E_{F_v} ) (-1)  $ is isomorphic to 
$ H ^ 1_{\et}  (S_{K, \bar F _v }, \  \bar {\mathcal F}_{(k, w) , o _{E_\lambda} } ) _{m_\Sigma } $ as a $G_{F_v}
$-module.
\end{lem}
\begin{proof}[Proof of Lemma \ref{lem-modular62}]
Fix $n \geq 1 $. By taking a conjugation, we may assume that $ K _u $ for $ u \vert p $ is a subgroup of $ \GL_2 ( o_{F_u } ) $. Let $K _n $ be the compact open subgroup 
defined by $K _ n =K _v \cdot  \prod _{u \vert p , u \neq v } (K (p ^n ) \cap K _u ) \cdot K ^p$, and  $X _ n = S _{K _n , o _{F_v} } $ the integral model for $S_{K_n } $. For $ X=  S_{K , o_{F_v} }$, $\pi _n : X_n \to X $ is the covering of $X$ induced by the inclusion $K _n \hookrightarrow K$. $ G _n = \overline{F ^{\times} }\cap K _n \backslash K  $ acts on $X_n$, and $\pi _n $ is a $G_n$-torsor since $K$ is small.  \par
For $\mathcal G _n = \tilde { \mathcal F} ^K _{(k,w) , o_{E_\lambda } / p ^n o_{E_\lambda }  }$, by the Hochshild-Serre spectral sequence
$$
H ^1 ( G_n , H ^0 ( X_{n , \bar F_v }, \pi ^* _n  \mathcal G_ n (1)_{\bar F_v} ) ) \longrightarrow    H ^1 ( X_{ \bar F_v }, \mathcal G_n (1) _ { \bar F_v} )\longrightarrow H ^1 ( X_{n , \bar F_v }, \pi ^* _n  \mathcal G_ n (1)_{\bar F_v} ) ^{G_n}
 $$
$$\longrightarrow  H ^2( G_n , H ^0 ( X_{n , \bar F_v }, \pi ^* _n  \mathcal G_n (1)_ {\bar F_v} ) )
$$
is exact. The action of $T_\Sigma $ on $H ^0  ( X_{n , \bar F_v } , \mathcal G _n (1)_{\bar F_v } ) $ is of residual type by Lemma \ref{lem-coh41}. After localization at $ m_{\Sigma} $, we have
$$
  H ^1 ( X_{ \bar F_v }, \mathcal G_n (1) _ { \bar F_v} )_{m_\Sigma } \stackrel{\sim}{\longrightarrow} (H ^1 ( X_{n , \bar F_v }, \pi ^* _n  \mathcal G_n (1) _ { \bar F_v} ) _{m_\Sigma } )^{G_n} .
\leqno{(\ast)}
$$
By the definition of $\tilde {\mathcal F }^K _{(k,w)}$, $\pi ^* _n \mathcal G _n $ is trivialized on $X_n $ in $D ^{\times} ( \mathbb A^{\Sigma }_{F, f} ) $-equivariant way. Take a trivialization $  \pi ^* _n \mathcal G _n \simeq M_n  $, where $M_n $ is a free $o_{E_\lambda } / p ^n o_{E_\lambda }   $-module.  \par
Let $ P_n = Pic ^0 (X_n / o_{F_v}) $ be the connected component of the Picard scheme of $ X_n $ over $o_{F_v}$. We regard the $p^n $-division points $ P _n [ p ^n ]$ of $P_n$ as a finite flat group scheme over $o_{F_v} $. The convolution algebra $T_{\Sigma } $ and $G_n$ acts on $ P _n [ p ^n ]$ by the Picard functoriality. We regard $M_n  \otimes _{\mathbb Z_p / p ^n } P _n [ p ^n ] $ as a finite flat group scheme with an action of $G_n$ and $T_\Sigma$. $G_n$ acts diagonally on $ M_n$ and $ P _n [ p ^n ] $, and $ T_{\Sigma}$ acts trivially on $M_n$. 
These two $G_n $ and $T_{\Sigma }$-actions commute. 
It is easily checked that $H ^1 (X_{ \bar F_v }, \mathcal G_n (1) _ { \bar F_v} ) $ is $G_n \times T_{\Sigma} $-isomorphic to $M_n  \otimes _{\mathbb Z_p / p ^n } P _n [ p ^n ] (\bar F_v )   $ with the actions defined as above.  \par
The $G_n$-invariant $(M_n  \otimes _{\mathbb Z_p / p ^n } P _n [ p ^n ] ) ^{G_n }$ is represented by a closed subgroup scheme $\tilde E _n $ of $ M_n  \otimes _{\mathbb Z_p / p ^n } P _n [ p ^n ]$ ($\tilde E_n $ is equal to the Zariski-closure of $ (M_n  \otimes _{\mathbb Z_p / p ^n } P _n [ p ^n ] )_{F_v} ^{G_n }$ in $M_n  \otimes _{\mathbb Z_p / p ^n } P _n [ p ^n ] $). Then we define
$$
\mathcal E_n = (\tilde E _n ) _{m _\Sigma}.
$$
$\mathcal E_n$ is a finite flat group scheme over $o_{F_v}$. By ($\ast$), $\mathcal E_n (\bar F_v) $ is canonically isomorphic to $ H ^1 ( X_{ \bar F_v }, \mathcal G_n (1) _ { \bar F_v} )_{m_\Sigma } $ as a $G_{F_v}$-module. \par 
For an integer $ m \geq 0 $, the projection $X_{m+n }\overset{\pi _{m+n , n }} \to  X_n  $ induces $\mathcal E_n\to \mathcal E_{m+n } $ by the Picard functoriality. 
By Proposition \ref{prop-coh21}, $ \{ \mathcal E_{n , F_v} \} _{n \geq 1} $ forms a $p$-divisible group over $F_v$.
\begin{sublem}\label{sublem-modular61}
$\mathcal E_n $ is equal to $\mathcal E_{m+n } [ p ^ n ] =\ker (\mathcal E_{m+n } \overset{p ^m } \rightarrow \mathcal E_{m+n } ) $. 
\end{sublem}
\begin{proof}[Proof of Sublemma \ref{sublem-modular61}]
$\pi _{m+n , n }$ induces $\pi ^* _{m+n , n }: P _n \to P_{m + n} $ by the Picard functoriality. The kernel $K$ of $\pi ^* _{m+n , n } $ is finite flat over $o_{F_v}$. We show that 
$$
(K \cap P _n [ p^{\infty} ] ) _{m_\Sigma }= \{0 \}.
\leqno{(\dagger)}
$$
 It is sufficient to see $ (K \cap P _n [ p ] ) _{m_\Sigma }= \{0 \}$. $(K \cap  P _n [ p ]) (\bar F_v ) $ is the kernel of $ H ^1 ( X_{n , \bar F_v } , \mathbb Z/p (1) ) \to  H ^1 ( X_{m + n, \bar F_v } , \mathbb Z/p (1) )  $. By Lemma \ref{lem-coh41}, the localization at $m_{\Sigma }$ is zero, and ($\dagger$) is shown. \par
By ($\dagger$), $P _n [ p ^n ] _{m_\Sigma} $ is regarded as a closed subgroup scheme of $ P_{m + n}  $. This implies that $\mathcal E_n $ is a closed subgroup scheme of $ \mathcal E_{m+n} $. Since $\mathcal E_n $ and $\mathcal E _{m+n}[p^n]  $ are two finite flat closed subgroup schemes of $ \mathcal E _{m+n}$, it is enough to see that the generic fibers over $F_v$ are the same.  $ \{ \mathcal E_{n , F_v} \} _{n \geq 1} $ forms a $p$-divisible group over $F_v$, so the claim is shown.
\end{proof}
It is clear that the image of $\mathcal E_{m+n } \overset{p ^m } \to \mathcal E_{m+n } $ is contained in $\mathcal E_n $. 
Now we have an increasing sequence $ \{ \mathcal E_{n } \} _{n \geq 1}  $ of finite flat group schemes over $o_{F_v}$ which satisfies the following properties:
\begin{itemize}
\item For $m, n \geq 1 $, 
$$
0\longrightarrow \mathcal E _n \longrightarrow \mathcal E_{m +n } \overset{ p ^n } \longrightarrow \mathcal E _{ m}
$$
is exact. 
\item $ \{ \mathcal E_{n } \} _{n \geq 1}  $ forms a $p$-divisible group over $F_v$. 
\end{itemize}
By a celebrated argument of Tate, there is an integer $n_0 \geq 0 $ such that $\{ \mathcal E _{n + n_0 } / \mathcal E _{ n_0 } \} _{n \geq 1 }  $ forms a $p$-divisible group $E$ over $o_{F_v}$ (\cite{Ta}, p.181--182: the affine algebra of $ \mathcal E _{n+1} / \mathcal E_n $ forms an increasing sequence of orders as $n$ varies, and $n_0 $ is taken so that it becomes stationary for $ n \geq n_0 $) .  By the construction, the Tate module $T_p (E) $ is isomorphic to $ H ^ 1(S_{K, \bar F _v }, \  \bar {\mathcal F}_{(k, w) , o _{E_\lambda}} (1) ) _{m_\Sigma } $.
The uniqueness of $E$ follows from Tate's theorem \cite{Ta}, Theorem 4. 
\end{proof}
We go back to the proof of Theorem \ref{thm-modular41}. 
We reduce to the case when the nearly ordinary type $\bar \kappa_{\mathcal D, v} $ is trivial. 
 Take a character $ \tau : G_F \to o^{\times}_{E_\mathcal D}$ of finite order which satisfies the following conditions:
\begin{itemize}
\item The order of $\tau $ is prime to $\ell$. 
\item $\tau \vert _{ I_{F_v}} = \bar \kappa_{\mathcal D, v} $.
\item  $\tau$ is unramified at any finite place $u$ in $(\Sigma_{\mathcal D
}\setminus \{ v\} )$. 
\end{itemize}
We can always find such a character $\tau$ with these properties by the global class field theory, by replacing $E_{\mathcal D}$ by a finite extension if necessary. Let $ \Sigma _{\tau}$ be the ramification set of $\tau$. \par

Let $\bar \tau$ be the mod $m _{\mathcal D}$-reduction of $\tau$, $ \bar \rho ' = \bar \rho \otimes \bar \tau ^{-1} $ the twist of $\bar
\rho$ by $ \bar \tau ^{-1}$. Note that $\bar  \rho ' $ has type $2_{PR} $ at the places in $\Sigma _\tau$.\par
The twist type of $\bar \rho ' \vert _{F_u } $ at a finite place $u $ is $\nu _u \otimes  \bar \tau\vert _{I_{F_u} }  ^{-1}$, where $\nu _u  $ is the twist type of $\bar \rho$ at $u$. 
 \par
We define a deformation condition $\mathcal D '  $ for $\bar \rho '$. $\mathcal D '  $ is the same as $\mathcal D $ except for the twist types at finite places, and the twist types for $\bar \rho'$ is the one defined as above. In particular $\deform _{\mathcal D }= \deform _{\mathcal D' } $ holds (at the places in $ \Sigma _\tau$, finite deformations are considered). \par
It suffices to show the claim for $\bar \rho '$ and $\mathcal D'$, since the twist of $\rho ^{\modular}_{\mathcal D}  $ by
$\tau ^{-1} $ gives $\rho ^{\modular} _{\mathcal D' }$ lifting $\bar \rho '$. \par

Thus we may assume that the nearly ordinary type $\bar \kappa _{\mathcal D, v} $ of $\bar \rho$ is trivial. 
We choose an auxiliary place $y$ and a finite set $S$ for $\bar \rho$ as in \S\ref{subsec-modular4}. Recall that $\tilde K = \tilde  K _{\mathcal D _y  } ( y, S) $ defined in \S\ref{subsec-modular4} is small. \par
$\tilde  K _v =\GL_2 (o_{F_v})$, and the $K$-character $\nu_{\bold{n.o.f.}} (\bar \rho \vert _{F_v} ) $ at $v$ is trivial by our assumption. \par
For $ K' =\ker \nu _{\mathcal D_y } \vert _{\tilde K } $, let $S_{K', F_v}= S _{K'} ( G _{\bar \rho} , X_{D(\bar \rho)} )_{F_v}$ be the Shimura curve over $F_v$. By the definitions of $\tilde M ^y _{\mathcal D_y} $ and $M _{\mathcal D} $ (Definition \ref{dfn-modular405}, \ref{dfn-modular42}), $M_{\mathcal  D} $ is a $G_F$-stable $o_{\mathcal D}$-direct summand of $ H ^ 1 _{\stack} (S_{K'} , \bar{ \mathcal F } ^{K'} _{(k, w)} )_{\tilde m _{\Sigma , \bar \rho }}  $. 
By Lemma \ref{lem-modular62}, there is a $p$-divisible group $ E ' $ over $o_{F_v} $ such that $T _p ( E') (-1)  $ is isomorphic to $ M _{\mathcal D}$ as a $G_{F_v}$-module. Since $ M _{\mathcal D, E_{\mathcal D}}$ is isomorphic to $\rho ^{\modular}_{\mathcal D, E_\mathcal D}  $, by Proposition \ref{prop-galdef61}, it follows that $\rho ^{\modular} _{\mathcal D } \vert _{F_v} $ is nearly ordinary and finite. \par
Assume that $\deform _{\mathcal D } (v) = \bold {fl}$.  We may assume moreover that the $v$-type is $ ((2, \ldots, 2) , 0)$ and the twist type at $v$ is trivial. As in the case of $\bold {n.o.f.}$, there is a $p$-divisible group $ E ' $ over $o_{F_v} $ such that $T _p ( E') (-1)  $ is isomorphic to $ M _{\mathcal D}$ as a $G_{F_v}$-module. Since $F_v $ is absolutely unramified and $p \geq 3$, for a $G_{F_v}$-representation $\rho $ on a finite dimensional $\mathbb Q_p $-vector space $V$ and a $G_{F_v}$-stable lattice $L$ in $V$, $L$ is associated to a $p$-divisible group if and only if there is one such lattice $L_0$ in $V$. Thus $ \rho ^{\modular}_{\mathcal D} \vert _{F_v} $ is a flat deformation. \par
\medskip
We have shown that $\rho ^{\modular} _{\mathcal D} \vert _{F_v}  $ is a deformation of type $\deform_{\mathcal D} (v) $ at any finite place $v$, and $U ( p _v ) ,\  U (p_v, p_v )  $ for $v \in P ^{\bold u } _{\mathcal D} \setminus P ^{\exc}$, and $\tilde U ( p _v ) , \ \tilde U (p_v, p_v )  $ for 
$v \in P ^{\bold {n.o.}} $ belong to $ T ^{\ss}_{\mathcal D } $. This proves Theorem \ref{thm-modular41} when $q_{\bar \rho} $ is one. 

\subsection{Approximation method of Taylor}\label{subsec-modular7}
 
We prove Theorem \ref{thm-modular41} at $\ell$ when $q_{\bar \rho } $ is zero (under Hypothesis \ref{hyp-nearlyordinary21} if $\deform(v) = \bold{n.o.}$). We use the approximation method of Taylor
\cite{T1}, \cite{T2}. \par
For a finite place $u $ such that $ u \not \in \Sigma _{\mathcal D }$, we assume the following condition:
\begin{itemize}
\item[{\bf{SP} }] $\bar \rho (\Fr _u )$ is a regular semi-simple element of $\GL_2 (k) $, and the two distinct eigenvalues of $\bar \rho( \Fr _u )$ are of the form $ \bar \alpha _u, \  q_u \bar \alpha  _u $.
\end{itemize}
 In the following, we assume that $ \bar \alpha _u \in k= o_{\mathcal D}/ m_{\mathcal D}  $.\par

$$
I ^{u-\parabolic } _{(k, w) } = \Hom _{\tilde K_{\mathcal D   }\cap K_0 ( u )   } (\nu _{\mathcal D } \vert  _{\tilde K_{\mathcal D   }  \cap K_0 ( u ) }, \mathcal S_ {(k, w) } ( G _{\bar \rho}  )  ).
$$
As in \S\ref{sec-coh}, $I ^{K_0 (u) } _{(k, w) }   $ contains $I ^{u-\old}_{(k,w)} = I_{( k, w ) , \mathcal D}(G_{\bar \rho} )  ^{\oplus 2}  $ as a subspace (the space of the forms which are old at
$u$). The inclusion is induced by $ ( f_1,\ f_2 ) \mapsto \pr  _1 ^* f_1 +
\pr_2^* f_2 $. Here $\pr _i $ ($i = 1, 2$) are degeneracy maps as in \S\ref{subsec-coh3}:
$\pr_1$ corresponds to the inclusion $ \tilde K_{\mathcal D   }\cap K_0 ( u ) \hookrightarrow \tilde K_{\mathcal D   }$, $\pr_2$ is the projection twisted by the conjugation by 
$\begin{pmatrix}
 1& 0  \\ 0 & p_u 
 \end{pmatrix}
$ at $u$. The quotient space 
$$
I ^{u-\new } _{(k, w) } = I ^{u-\parabolic } _{(k, w) } / I ^{u-\old}_{(k,w)}$$
is the space of the forms which are new at $u$. \par

\begin{dfn}\label{dfn-modular705} Under the condition {\bf SP} on $u$, the Hecke rings $ T^{u-\parabolic } _{\mathcal D}  $ and $ T^{u-\sp } _{\mathcal D} $ are defined as follows. 
\begin{enumerate}
\item 
$ \tilde T^{u-\parabolic }_{ \mathcal D } $ is the $o_{\mathcal D}$-algebra generated by the following elements in $ \End_{o_\mathcal D} I ^{u-\parabolic}_{(k, w)} $:
$ T_v $ and $T_{v,v} $ for $v \not \in \Sigma_{\mathcal D}\cup \{ u\} $, $ U ( p_v )$ and 
$ U(p_v, p_v) $ for $ v \in P ^{\bold u} _\mathcal D \setminus P ^{\exc}$, 
$\tilde U ( p_v ) $, and $\tilde U(p_v, p_v) $ for $ v \in P ^{\bold {n.o.}} $, $U (p_u ) $ and 
$ U ( p_u , p_u ) $ at $u$ as in Definition \ref{dfn-modular22}.  
\item The maximal ideal $\tilde  m^{u-\parabolic} _{\mathcal D  } 
$ of $\tilde T^{u-\parabolic }_{ \mathcal D }  $ is generated by 
$m_{o_{\mathcal D} }$,
$T_v - f_{\pi _{\min}}  ( T_v)  $ for $v
\not \in  \Sigma_{\mathcal D}\cup \{ u\} $, $  \tilde  U ( p_v ) - \alpha _{\pi _{\min} , v}  $ for $v  \in P ^{\bold {n.o.}}$, $ U (p_v) $ for $ v \in  P ^{\bold u } _{\mathcal D} \setminus P ^{\exc} $, $ U ( p _u ) - \bar \alpha _u
$, and $U (p_v, p_v )- \chi _{\bar \rho} (p_v)$-operators for $v \not \in P ^{\bold {n.o.}} \cup P ^{\exc}$ in the notation of \ref{subsec-modular2}. 
\item $ T ^{u-\parabolic} _{\mathcal D}$ is the localization of $\tilde T^{u-\parabolic }_{ \mathcal D }  $ at $ \tilde  m^{u-\parabolic} _{\mathcal D  } $. \par
\item $\tilde T^{u-\sp}_{\mathcal D} $ (resp. $\tilde T^{u-\unr}_{\mathcal D} $) is the image of $ \tilde T ^{u-\parabolic} _{\mathcal D} $ in $ \End _{o_\mathcal D} I ^{u-\new}_{(k,w)} $ (resp. $ \End _{o_\mathcal D} I ^{u-\old}_{(k,w)} $).
$\tilde m ^{u-\sp} _{\mathcal D} $ (resp. $\tilde m ^{u-\unr} _{\mathcal D} $) is the image of $ \tilde m ^{u-\parabolic} _{\mathcal D} $ in $\tilde T^{u-\sp}_{\mathcal D} $ (resp. $\tilde T^{u-\unr}_{\mathcal D} $). Then define $T^{u-\sp}_{\mathcal D} $ (resp. $T^{u-\unr}_{\mathcal D}$ ) by $(\tilde T^{u-\sp}_{\mathcal D}) _{\tilde m ^{u-\sp} _{\mathcal D}}$ (resp. $ (\tilde T^{u-\unr}_{\mathcal D}) _{\tilde m ^{u-\unr} _{\mathcal D}}$), which are quotients of $  T ^{u-\parabolic} _{\mathcal D}$. 
\end{enumerate}
\end{dfn}
\begin{lem}\label{lem-modular71}
\begin{enumerate}
\item $T^{u-\unr}_{\mathcal D}$ is isomorphic to $ T_{\mathcal D}$.
\item We regard $M _{\mathcal D} $ as a $T^{u-\parabolic}_{\mathcal D}$-module by $T^{u-\parabolic}_{\mathcal D}\twoheadrightarrow T ^{u-\unr}_{\mathcal D} \simeq  T_{\mathcal D}$. Then $T^{u-\parabolic}_{\mathcal D}$-action on $M _{\mathcal D} /(  U (p_u) ^ 2 - T_{u, u }) M _{\mathcal D} $ factors through $T^{u-\sp} _{\mathcal D }$. In particular $T^{u-\sp} _{\mathcal D } $ is non-zero.
\end{enumerate}
\end{lem}
\begin{proof} [Proof of Lemma \ref{lem-modular71}] This is well-known in the theory of congruence modules, so we briefly sketch the proof by explaining how to construct the isomorphisms in Lemma \ref{lem-modular71}.  \par
We choose an auxiliary place $y$ and $S$ as in \S\ref{subsec-modular4} such that $ u \not \in \{ y \} \cup S$. By condition {\bf SP}, the polynomial
$$f (X) = X ^2 - T_{v} X + q_v \cdot T_{v,v} 
$$ 
has two distinct roots in $ T ^{\ss}_{\mathcal D}/ m^{\ss}_{\mathcal D} $, and hence in $T^{\ss}_{\mathcal D}$ by Hensel's lemma. We choose the root $\alpha$ which lifts $\bar \alpha _u $. By the same argument as in the case of $\bold {n.o.f. }$ discussed in \S\ref{subsec-congruent3}, we obtain an injective homomorphism
$$
\xi : M _{\mathcal D} \longrightarrow M ^{u-\parabolic}_{\mathcal D }.
$$
Here $K = \ker \nu _{\mathcal D _y } \vert _{\tilde K_{  \mathcal D_ y }}    $ in the notation of \S\ref{subsec-modular4}, and 
$$
M ^{u-\parabolic}_{\mathcal D } 
=  \Hom _{\tilde K_{\mathcal D_y } \cap K_0 (v)  } 
( \nu _{\mathcal D_y } , 
H ^ 0 _{\stack}( S_{K\cap K_0 (v) } ,\ \bar {\mathcal F} _{(k, w)} )
_{\tilde m _{ \Sigma _{\mathcal D _y } \cup S \cup  \{ v \}  } }
).
$$
Since $U (p_u ) $ satisfies $ f (U ( p_u ) ) = 0 $, $U(p_v)$ acts via the multiplication by $\alpha $ on $ M_{\mathcal D}$.
\par
 Since $q_{\bar \rho} = 0 $, by Theorem \ref{thm-coh31}, $\xi$ is universally injective. 
$$
\xi ^{\vee} : M ^{u-\parabolic}_{\mathcal D }\longrightarrow M _{\mathcal D} 
$$
is the map induced by $\xi $ by Poincar\'e duality.  
By Lemma \ref{lem-congruent21}, 
$$\xi ^{\vee} \circ
\xi ( M _{\mathcal D}  )= (  U (p_u) ^ 2 - T_{u, u} ) M_{\mathcal D} 
$$ 
as in the case of $\bold {n.o.f.}$ in \S\ref{subsec-congruent2}. By a standard argument in the theory of congruence modules
\cite{Ri1},  the congruence module $ M _{\mathcal D} / \xi ^{\vee}\circ \xi ( M_{\mathcal D}  )$ is a $T^{u-\parabolic}_{\mathcal D} $-module, and Lemma \ref{lem-modular71} is shown. 
\end{proof}

By $T _{\mathcal D} ^{\ss, u-\parabolic}$ we denote the $o_{\mathcal D }$-subalgebra of $T _{\mathcal D} ^{ u-\parabolic} $ generated by $ T_v,\ v
\not \in
\Sigma _{\mathcal D} \cup \{ u , y \}   \cup S$ as in \ref{subsec-modular2}, and denote by $T _{\mathcal D} ^{\ss, u-\sp} $ (resp. $T _{\mathcal D} ^{\ss, u-\unr} $) the image of $ T _{\mathcal D} ^{\ss, u-\parabolic}$ in $T _{\mathcal D} ^{ u-\sp} $ (resp. $T _{\mathcal D} ^{ u-\unr} $). By Lemma \ref{lem-modular31}, $T _{\mathcal D} ^{\ss, u-\unr} = T ^{\ss}_{\mathcal D } $.

\begin{lem} \label{lem-modular72} For any integer $N \geq 1 $, there is a finite place $u$ which satisfies the following properties:
\begin{enumerate}
\item $u \not \in \Sigma _{\mathcal D} $, $q_u \equiv -1 \mod \ell$, and the eigenvalues of $\bar \rho( \Fr _u )$ are of the form $ \bar \alpha, \ -\bar \alpha  $.
In particular {\bf SP} holds at $u$.
\item We define $ T ^{u-\parabolic}_{\mathcal D} $ and $ T ^{u-\sp}_{\mathcal D} $ by choosing an eigenvalue of $\bar \rho( \Fr _u )$. Then $T ^{\ss}_{\mathcal D}/ m^N _{\mathcal D } T ^{\ss}_{\mathcal D} =T _{\mathcal D} ^{\ss, u-\unr}/m^N  _{\mathcal D}  T _{\mathcal D} ^{\ss, u-\unr}$ is a quotient of $T _{\mathcal D} ^{\ss, u-\sp} $ as a $T _{\mathcal D} ^{\ss, u-\parabolic}$-algebra.
\end{enumerate}
\end{lem}
\begin{proof}[Proof of Lemma \ref{lem-modular72}] By taking a finite unramified extension of $o _{\mathcal D} $, we may assume that all eigenvalues of $\bar \rho ( g ) $, $g \in G_F$ belong to $k = o_{\mathcal D}/ m_{\mathcal D}$. \par
$M _{\mathcal D , E_\mathcal D} $ is free of rank one as a $
T_{\mathcal D , E _{\mathcal D} }^{\ss}  $-module by Proposition \ref{prop-modular415}, Lemma \ref{lem-modular23}, and Lemma \ref{lem-modular41}. We choose an injective $T _{\mathcal D }^{\ss} $-homomorphism with a finite cokernel
$$
 i  : T _{\mathcal D }^{\ss}   \longhookrightarrow M _{\mathcal D} , 
$$
and take $ c _0  \geq 0
$ such that $m_{\mathcal D } ^ {c_0} $ annihilates $ \coker i  $. \par

Let $\tilde T$ be the normalization of $T _{\mathcal D }^{\ss} $. Since $T _{\mathcal D }^{\ss} $ is reduced, the canonical $o_{\mathcal D}$-algebra homomorphism gives an embedding
$$
j : T _{\mathcal D }^{\ss}   \longhookrightarrow \tilde T.
$$
$ c_1 \geq 0 $ is an integer such that $ m ^{c_1} _{\mathcal D} $ annihilates $\coker j $. \par
Define $ \rho _{\tilde T } : G_F \to \GL _2 (\tilde T ) $ as the composition of $\ G_F \overset {\rho ^{\modular} _{\mathcal D} } \to \GL _2 (T ^{\ss} _{\mathcal D}  )\to  \GL_2 ( \tilde T)$.  \par

We fix an element $c$ of order $2$ in $G_F$ which corresponds to the complex conjugation for some embedding $ F \hookrightarrow \mathbb C$. 
We view $\mathcal L =\tilde T ^{\oplus 2}  $ as a $G_F$-module by $\rho _{\tilde T}$. For $\gamma = c_1 + c_2 $, let $ H_N $ be the image of $G_F $ in $\Aut ( (\mathbb Z/ \ell ^ {N + \gamma}\mathbb Z) ( 1) \oplus  \mathcal L/ m _{\mathcal D } ^ { N  +\gamma } \mathcal L  ) $, $C $ the conjugacy class of $ c$ in $H_N $. 
We take $u$ so that $ \Fr _{u } $ belongs to $C$ by the Chebotarev density theorem. We show that $u$ satisfies the desired properties. \par
By our choice of $u$, \ref{lem-modular72} (1) is satisfied, $1+ q_u \equiv 0 \mod \ell ^{N +\gamma } $ holds, and the image of $T _u $ is zero in $\tilde T  / m^{N + \gamma}  _{\mathcal D}\tilde T $. \par
Since
$\coker j$ is annihilated by
$m_{\mathcal D } ^{c_1}$, $ T_u$ is contained in 
$m^{N +\gamma - c_1 } _{\mathcal D } T_{\mathcal D } ^{\ss} = m ^{N + c_0} _{\mathcal D } T_{\mathcal D } ^{\ss}  $.\par
\medskip
 We choose an eigenvalue $\bar \alpha $ of $\bar \rho (\Fr _u ) $, and define Hecke algebras $ T ^{u-\parabolic}_{\mathcal D} $ and $ T ^{u-\sp}_{\mathcal D} $ using $\bar \alpha$. Let $\alpha $ be the roots of $f (X) = X ^2 - T_u X + q _u \cdot T _{u, u}  $ which lift $\bar \alpha$. The action of $ U ( p_u ) $ on $ M _{\mathcal D} $ is the multiplication by $\alpha  $ in $T^{\ss} _{\mathcal D } $.  \par
 As
$$
A_u = \alpha ^2 - T _{u, u } = T_u \cdot \alpha - (1+ q_u ) T _{u, u } \in m^{N +c_0 } _{\mathcal D } T_{\mathcal D } ^{\ss}  ,
$$
$ M _{\mathcal D} / (U ( p_u ) ^2 - T_{u, u } )M _{\mathcal D } = M _{\mathcal D} / A_u  M _{\mathcal D}  $ admits $M _{\mathcal D} / m^{N +c_0  } _{\mathcal D } M _{\mathcal D}$ as a quotient. By Lemma \ref{lem-modular71} (2), $M _{\mathcal D} / m^{N +c_0  } _{\mathcal D } M _{\mathcal D}$ is a $T ^{\ss, u-\sp}_\mathcal D $-module. 
\par
By our choice of $c_0 $, the image of $ i (T_{\mathcal D } ^{\ss} ) $ in $M _{\mathcal D} / m^{N +c_0  } _{\mathcal D } M _{\mathcal D}$ admits $T^{\ss}  _{\mathcal D } / m_{\mathcal D} ^NT^{\ss}  _{\mathcal D }$ as a quotient. Thus \ref{lem-modular72}, (2) is shown. 
 
\end{proof}
For an integer $N \geq 1 $, we choose a finite place $u_N $ which fulfills the conclusion of Lemma \ref{lem-modular72}. 

We have a surjective $o_{\mathcal D }$-algebra homomorphism 
$$
\alpha _N : T _{\mathcal D} ^{\ss, u_N-\sp } \twoheadrightarrow T ^{\ss} _{\mathcal D}/ 
m _{\mathcal D} ^ N T ^{\ss} _{\mathcal D}
$$
which maps $T_v $ to $T_v$ for $ v \not \in \Sigma _{\mathcal D} \cup \{ u_N \} $. 

For
$T _{\mathcal D} ^{\ss, u_N-\sp } $, there is a Galois representation
$\rho_{\mathcal D , u_N-\sp }  ^{\modular} : G_F \to \GL_2 (T _{\mathcal D} ^{\ss, u_N-\sp } )$ as in Proposition \ref{prop-modular32}. Since any representation
$\pi$ which appear as a component of $ T _{\mathcal D}   ^{\ss, u_n -\sp }$ has a special representation $\pi _{u_N }$ at $u_N$, we can apply Hypothesis \ref{hyp-nearlyordinary21} to $\pi$, and show that $\rho_{\mathcal D , u_N -\sp}  ^{\modular} $ is nearly
ordinary, nearly ordinary finite or flat according to $ \deform _{\mathcal D} ( v) =
\bold { n.o.}, \ \bold {n.o.f.},\ \bold {fl}$ by the argument in \S\ref{subsec-modular6}. \par
By the surjectivity of $\alpha _N $ and the uniqueness (Lemma \ref{lem-modular31}), $ \rho ^{\modular } _{\mathcal D}
\mod
m _{\mathcal D}  ^ n = \alpha _N \circ \rho_{\mathcal D , u_N -\sp }  ^{\modular}  $ is nearly ordinary or nearly
ordinary finite at a place $v$ such that 
$ \deform _{\mathcal D} ( v) =
\bold { n.o.}, \ \bold {n.o.f.}, \bold{fl}$. 
 Since this is true for any $N$, the claim for $ T_{\mathcal D}$ is shown.
\end{proof}

\begin{rem}\label{rem-modular3}
In the {\bf fl}-case, one may use Taylor's result (\cite{T2}, theorem 1.6 and lemma
2.1, 3) if the infinity type is $( (2, \ldots, 2), w)$.  The demonstration there can be modified so that it fits to our setting (we only assume that the $v$-type is $( ( 2, \ldots, 2), w)$). 
\end{rem}
\begin{prop}\label{prop-modular71}
The following operators in $\End_{o_{\mathcal D}} M_{\mathcal D}$ belong to $ T^{\ss}_\mathcal D $. $  U ( p_v ), U (p_v , p_v )$ for $v \in P ^{\bold f} _{\mathcal D } $, $\tilde U ( p_v ) , \tilde U (p_v, p_v ) $ for $v \in P ^{\bold {n.o.f.}}\cup P ^{\bold {fl}}$.
\end{prop}
\begin{proof}[Proof of Proposition \ref{prop-modular71}] By Lemma \ref{lem-modular22}, the operators in question belong to $ T^{\ss}  _{\mathcal D , E_\mathcal D }$. We find these operators by using $\rho ^{\modular}_{\mathcal D}$. For a finite place $v$, let $\sigma _v$ be an element in $ I_{F_v}$ which lifts the geometric Frobenius element, and corresponds to $p_v$ by $(I^{\ab} _{F_v}) _{G_{F_v} } \simeq o^{\times}_{F_v}$.
 \par

For $v \in P ^{\bold f} _{\mathcal D } $, we use the notation of \S\ref{subsec-modular5}. $U ( p_v, p_v ) =q_v ^{-1}  \det \rho ^{\modular } _{\mathcal D} (\sigma _v) $. 
When $d_v= 2$, $U (p_v ) = \tr _{T^{\ss}_\mathcal D} \rho ( \sigma _v ) $. When $d_v= 1 $, $U(p_v) $ is equal to the action of  $\rho^{\modular} _{\mathcal D} (\sigma _v)$ on
$\mathcal L $. Note that it is sufficient to verify the equalities after tensoring $E_{\mathcal D} $. \par

For $v \in  P ^{\bold {n.o.f.}} \cup P ^{\bold {fl}}$, $\tilde U(p_v, p_v  ) = (\det \rho _{\mathcal D} ^{\modular} 
\cdot
\chi _{\cycle}^ {w+1} )  ( p_v )
$
holds. For $\tilde U (p_v)$-operator, by twisting, we may assume that $w= 0 $. By the argument (twisting by a finite order character) in \S\ref{subsec-modular6}, we may assume that the twist type at $v$ is trivial.\par
For $\rho_v =   \rho^{\modular} _{\mathcal D }
   \vert _{G_{F_v} } $, $\rho_v (1)$ extends to a unique $
\ell$-divisible group $E$ over $o_{F_v}$, and it is a crystalline representation.\par
If $v \in P ^{\bold {n.o.f.}}$, let $\tilde U  $ be the action of $\rho _v (\sigma _v )  $ on $\mathcal L $ in the notation of \S\ref{subsec-modular6}. We define $ \tilde T_v $ by
$$
\tilde T_v = \tilde U + q_v \tilde U (p_v, p_v ) \cdot \tilde U ^{-1}.
$$
\par
If $v \in  P ^{\bold {fl}}$, let $D(\rho_v) $ be the filtered module associated to
$\rho_v $ of filtration type $[
0,\ 1 ] $ by \cite{FL}. $D(\rho_v)$ is a free $T_{\mathcal D} ^{\ss}\otimes _{\Z_\ell} o_{F_v}
$-module with an action of the absolute Frobenius $\varphi$. We define $\tilde T_v $ by 
$$
\tilde T_v =  \tr _{T_{\mathcal D} ^{\ss}\otimes _{\Z_\ell} o_{F_v}} (
\varphi ^ f,
\ D(\rho_v ) ) \in T_{\mathcal D} ^{\ss}\otimes _{\Z_\ell} o_{F_v} . 
$$
Here $ f = [ k(v) : \F_\ell]$. $\varphi ^f$ is $T_{\mathcal D} ^{\ss}\otimes _{\Z_\ell} o_{F_v}$-linear, thus $\tilde T_v$ is well-defined. \par
We show 
$$
\tilde U (p_v) = \tilde T_v 
$$
in $T_{\mathcal D} ^{\ss}\otimes _{\Z_\ell} o_{F_v}$ in the both cases. It suffices to see it in $ T_{\mathcal D} ^{\ss}\otimes _{\Z_\ell} F_v$. If $q_{\bar \rho } = 1 $, this is a consequence of the main theorem of
\cite{Sa2}, which shows a weaker version of the compatibility of the local and the global Langlands correspondence at
$v \vert \ell$. \par¡¡¡¡
When $q_{\bar \rho}$ is zero, we use the
approximation method using Lemma \ref{lem-modular72}, and the equality holds in $T_{\mathcal D} ^{\ss}\otimes _{\Z_\ell} o_{F_v}/ m_{F_v } ^ N$ for any $N$, and hence in $T_{\mathcal D} ^{\ss}\otimes _{\Z_\ell} o_{F_v}$. \par
Since $ o_{F_v} $ is faithfully flat over $\mathbb Z_\ell$, $ \tilde U (p_v )\in T^{\ss} _{\mathcal D} $. 
\end{proof}

\section{Construction of Taylor-Wiles systems}\label{sec-const}
In this section, we
construct the Taylor-Wiles system for a minimal Hecke algebra
$ T_{\mathcal D}$.

\begin{thm}\label{thm-const0}
Let $\ell \geq 3$ be a prime, $\bar \rho  $ an absolutely
irreducible representation. Assume that $\bar
\rho \vert _{F(\zeta_\ell)}$ is absolutely irreducible if $ [F(\zeta_\ell):F] = 2$.

For a minimal deformation type $\mathcal D $, define $o^{\associate}_{\mathcal D} $ as $ 
o_{\mathcal D }  [ C_{F, \ell } ] $, where $ C_{F, \ell}  $ is the $\ell$-part of the class group of $F$.\par
 Then there is a Taylor-Wiles system
$\{ R _ Q,\ M _Q 
\} _{ Q
\in \mathcal X_{\mathcal D} } $ for $( R_{\mathcal D} , M _{\mathcal D} )$ and torus $ \G _{m , F} $ over $o^{\associate}_{\mathcal D}$. The index set $\mathcal X _{\mathcal D }  $ is defined by 
$$
Q\in  \mathcal X _{\mathcal D } \Longleftrightarrow Q \text{ is a finite subset of }Y_{\mathcal D} . 
$$
Here, $Y _{\mathcal D } $ is the set of the finite places $v \not\in \Sigma_{\mathcal D }$, $q_v \equiv 1 \mod \ell$, and $\bar \rho (\Fr _v ) $ is a regular semi-simple element of $\GL _2 (k )$. 
\end{thm}

In \cite{TW}, TW5 was obtained by a study of group cohomology. We use an argument based on a property of perfect complexes. 

\subsection{Hecke algebras}\label{subsec-const1}
For an absolutely irreducible representation $\bar \rho $, we fix a deformation type $\mathcal D$ for $\bar \rho$.

For an element $Q \in \mathcal X_{\mathcal D} $, we construct a Hecke algebra
$T_Q$ and a
$T_Q$-module
$M_Q$, define a ring homomorphism $R_Q \to T_Q$ from some deformation ring $R_Q$, and
verify the conditions of Taylor-Wiles systems.\par

By taking an extension of $k$, we may assume that any eigenvalues of $\bar \rho (g) $ for $g \in G_F $ belong to $k$. For a finite place $ v \in Y _{\mathcal D}$, let
$\bar \alpha _ v $ and $\bar  \beta _v $ be the eigenvalues of $ \bar \rho ( \Fr _v) $ at
$ v$. We make a choice, and specify one of the eigenvalues of $\bar \rho (\Fr_v ) $, say
$\bar \alpha _v$. We denote by $\Delta _v $ (resp. $ \Delta^{\vee}  _v $) the $\ell$-Sylow subgroup of $ k ( v ) ^{\times} $ (resp. the maximal prime to $ \ell$-subgroup of $k (v) ^\times $).\par

For $Q\in \mathcal X_{\mathcal D} $, we choose an auxiliary place $ y $ and a set $S$ as in \S\ref{subsec-modular4} which are disjoint from $ \Sigma _{\mathcal D } \cup Q$. 
We define the deformation function $\deform _{\mathcal D_Q }$ of $\bar \rho $ by
$$
\deform _{\mathcal D_Q } (v ) = \deform_{\mathcal D } ( v) \text{ for }v \not \in \Sigma _{\mathcal D} , \quad \bold u \text{ for } v \in Q ,
$$
and let $\mathcal D_Q $ be the deformation type of $\bar \rho$ which has $\deform _{\mathcal D_Q} $ as the deformation function, and the same data as $\mathcal D$ except for the deformation function. \par

Let $ K = \tilde K _{\mathcal D _y }= K _{11} ( m_y ^2 ) \cdot K ^y _{\mathcal D }  $ be the compact open subgroup of $ G _{\bar \rho }( \mathbb A _{\mathbb Q , f} ) $ defined in \S \ref{subsec-modular4} for $\mathcal D_y $. 
We define two compact open subgroups $K_{0, Q} $, $K_Q$ of $K $ for $Q$: 

$$
K_{ 0, Q}=\prod _{v \in Q} K_{0 }(m_v) \cdot K ^Q , 
$$

$$
K _Q = \prod _{v \in Q} K_{Q, v} \cdot  K ^Q .
$$
Here 
$$
 K _{ Q, v} = \{ g \in K_0 ( m_v),\ \ g \equiv 
\begin{pmatrix}
\alpha    & * \\ 
0 & \alpha \cdot h    \\
\end{pmatrix}
\mod m_v , \ \alpha \in k(v) ^\times,\ h , \in \Delta ^{\vee}_v
\} .
$$

There are inclusions $K_Q \subset K_{0, Q} \subset K $. $K_Q$ is a normal subgroup in
$K_{0, Q}$, and the quotient $(\overline{Z ( \mathbb Q  ) }\cap K_Q )\cdot K _Q  \backslash K_{0, Q}$ is isomorphic to 
$$
\Delta_Q :=
\prod_{v\in Q}
\Delta_v . 
$$

We define the Hecke algebras for $K_{0, Q} $ and $K_Q$ which are slightly different from the one defined for $\mathcal D $ and $\mathcal D_Q$ in \S\ref{sec-modular} at the places in $Q$.  \par

For $ \delta \in \Delta _Q $, we take an element $ \delta ' \in \prod _{v \in Q} o_{F_v} ^\times $ lifting $\delta
$ and $ c (\delta ' ) \in G_{\bar \rho} (\A  _{\mathbb Q , f} )  $ is the element which has  $ c (\delta ' ) _Q =
\begin{pmatrix} 
{\delta '} ^{-1}  &0\\ 
 0&{\delta ' }
\end{pmatrix} $,
at the $Q$-component and the components other than $Q$ are
one. Then the double coset $ K_Q c(\delta ' ) K_Q$ depends only on $\delta $, and defines an element $ V (\delta ) $ in $ H ( G_D (\A _{\mathbb Q , f} ) ,\ K _Q) $. \par

By our assumption, there is a
minimal lift $\pi _{\min} $ of
$\bar \rho $ defined over $o _{\mathcal D}$. $f_{\pi _{\min}} : T^{\ss} _{K , \mathcal o_{\mathcal D} } \to o _{\mathcal D} $ is the corresponding  $o _\mathcal D 
$-algebra homomorphism. 
\begin{dfn}\label{dfn-const11}
\begin{enumerate}
\item  $ I^0  _{(k, w ),  \mathcal D _Q } $ (resp. $ I  _{(k, w ),  \mathcal D _Q } $ ) is the intertwining space 
$$
I^0  _{(k, w ),  \mathcal D _Q }  = \Hom _{K_{0, Q}  } ( \nu _{\mathcal D} \vert _{K_{0, Q}  } , \mathcal S _{(k, w) } ( G_{\bar \rho} ) ) \quad (\text{resp. } I  _{(k, w ),  \mathcal D _Q }  = \Hom _{K_{Q}  } ( \nu _{\mathcal D} \vert _{K_{ Q}  }  , \mathcal S _{(k, w) } ( G_{\bar \rho} )  ).
$$
\item For an element $\pi \in \mathcal A _{(k, w ) } ( G _{\bar \rho } ) $, 
 $ I^0  _{(k, w ),  \mathcal D _Q } (\pi )  $ (resp. $ I  _{(k, w ),  \mathcal D _Q } (\pi )  $ ) is the intertwining space 
$$
I^0  _{(k, w ),  \mathcal D _Q } ( \pi )   = \Hom _{K_{0, Q}  } ( \nu _{\mathcal D} \vert _{K_{0, Q}  } , \pi _f  ) \quad (\text{resp. } I  _{(k, w ),  \mathcal D _Q }  = \Hom _{K_{Q}  } ( \nu _{\mathcal D} \vert _{K_{ Q}  }, \pi _f  )).
$$
\end{enumerate}
\end{dfn}
\begin{dfn}\label{dfn-const12}
\begin{enumerate}
\item  The Hecke algebra $\tilde  T_{ Q }$ is the $o _{\mathcal D}$-algebra generated by the following elements in $\End _{o_{\mathcal D} } I  _{( k, w ), \mathcal D_Q } $: 
$ T_v, $ and $T_{v, v} $ for $ v \not  \in \{ y \} \cup \Sigma _{\mathcal D_ Q } $, $ U ( p_v )$ and $  U(p_v, p_v ) $ for $v \in \{  y \} \cup S \cup P ^{\bold u } _{\mathcal D_Q } \setminus P ^{\exc} $,  $ \tilde U ( p_v )$ and $\tilde U ( p_v, p_v ) $ for $ v \in P ^{\bold {n.o.} } $, $V (\delta ) $ for $ \delta \in \Delta _Q$. \par
 $ \tilde T_{0, Q} $ is the image of $\tilde  T_{ Q } $ in $ \End _{o_{\mathcal D} } I ^0 _{( k, w ), \mathcal D_Q }$.
\item By $\tilde m  _{Q} $, we denote the maximal ideal of
$ \tilde T _{Q} $ generated by $ m _{o _{\mathcal D}} $ and the following operators:
$T_v - f_{\pi _{\min} }  ( T_v)$, $ T_{v,v} - f _{\pi _{\min} } ( T_{v,v}) $ for $  v \not  \in \{ y \} \cup S \cup  \Sigma _{\mathcal D_ Q }$, $U (p_v )  $, $ U ( p _v , p_v ) - \chi _{\bar \rho }(p_v)  $ for $ v \in ( \{ y \} \cup P ^{\bold u } _{\mathcal D} ) \setminus P ^{\exc}$,  $  \tilde  U ( p_v ) - \alpha _{\pi _{\min } , v } $, $\tilde U ( p_v, p_v ) - \gamma _{\pi _{\min}} $ for $ v \in P ^{\bold {n.o.} } $, $U (p_v )- \bar \alpha _v  $, $ U ( p _v , p_v ) - \chi _{\bar \rho }(p_v)  $ for $ v \in Q$,  $ V ( \delta ) -1$ for $ \delta \in \Delta _Q $. Here $ \alpha _{\pi _{\min } , v } $ and $\gamma _{\pi _{\min}} $ are elements of $ o_{\mathcal D} $ defined in \S\ref{subsec-modular2}. \par
$ \tilde m _{0, Q} $ is the image of $ \tilde m _Q $ in $\tilde T _{0, Q} $.
 
\item  The Hecke algebras $ T_{0, Q   } $ and $ T_Q$ associated to $Q$ are defined by
$$
T_{0, Q   } = (\tilde T_{0, Q} ) _{\tilde m_{0, Q}} , \quad T_ Q = (\tilde T_Q) _{\tilde m_Q } .
$$ 
\end{enumerate}
\end{dfn}
Note that the natural $o_{\mathcal D}$-algebra homomorphism $ T_{Q} 
\twoheadrightarrow T_{0, Q  } $ maps $ V(\delta ) $ for $\delta \in \Delta _Q$ to $1$.  

As in \S\ref{subsec-modular2}, we define $T _{0, Q } ^{\ss} $ (resp.
$T _{ Q } ^{\ss}$) as the $o_{\mathcal D }$-subalgebra generated by $T _v  $ and $T_{v, v} $ for $ v \not \in \{ y \} \cup S \cup \Sigma _{\mathcal D _Q} $. 
\begin{prop}\label{prop-const11}
$T ^{\ss} _{0, Q } =  T_{0, Q}  $ and $T ^{\ss} _{Q } =  T_{Q} 
$ hold. In particular $T_{0, Q}$, $T_{Q} $ are reduced.  
\end{prop} 
\begin{proof}[Proof of Proposition \ref{prop-const11}] We proceed as in the proofs of Lemma \ref{lem-modular23} and Proposition \ref{prop-modular41}. First we show $T ^{\ss} _{Q , E_{\mathcal D} } 
= T_{Q , E_{\mathcal D} } $. \par
Let $e _{\bar \rho}$ be the idempotent of $ \tilde T_Q $ corresponding to $T_Q $. Take a representation $\pi \in \mathcal A _{(k, w ) } ( G_{\bar \rho } ) $ such that $e _{\bar \rho}I _{(k, w), \mathcal D_Q} (\pi )\neq \{ 0 \} $. We show the dimension over $\bar E _{\mathcal D} $ is one. Since $ I _{( k, w ) , \mathcal D _Q } ( \pi ) $ is the tensor product local spaces, it suffices to study the local intertwining space at any finite place. If $v \not \in \{ y \} \cup Q $, this is similarly treated as in Lemma \ref{lem-modular23}. At $y$, this is done in Proposition \ref{prop-modular41}.  It suffices to consider the case of $ v \in Q$. \par
If $ v \in Q$, the $v$-compenent $\pi _ v$ belongs to principal series or a special representation twisted by a character since it
has a non-zero fixed vector by $ K_{11} ( m_v)$. Since we assume that $ q_v
\equiv 1
\mod
\ell$ and two eigenvalues $\bar \alpha _v ,\ \bar \beta _v $ of $\bar \rho (\Fr _ v)$ are
distinct modulo $m _{\mathcal D} $, the latter case does not occur. Since
$\pi _v$ has a non-zero fixed vector by $ K_v$, $\pi_v$ is a tamely ramified principal series. 

\begin{lem}\label{lem-const11}
For a local field $F$ and a principal series representation $
\pi =
\pi (
\chi _1 ,
\chi _2 ) $, assume that $\chi_1 $ and $\chi_2$ are tamely ramified characters of $ F ^\times$. 
Then there is a unique non-zero vector $w$ up to scalar in the representation space such that 
$$ 
\begin{pmatrix}
 a & * \\  
p_F \cdot *  &d 
\end{pmatrix}
 \cdot w = \chi _1 ( a )
\chi _2 ( d) w \quad (a, \ d  \in o_F ^\times) .
$$
\end{lem}
\begin{proof}[Proof of Lemma \ref{lem-const11}] Consider $\pi ' = \pi \otimes \chi ^{-1} _2 = \pi ( \chi _1 / \chi_2 , 1 ) $.
A new vector $ w'$ fixed by $K_1(m_F)$ of $\pi ' $ corresponds to $w$ as in the claim by twisting by $\chi ^{-1} _2$. 
\end{proof}

We choose an ordering of $\chi_1 $ and $\chi _2 $ so that $ \chi _2 (p _v) $ is a lift of $\bar \alpha _v$. Since $ \bar \alpha _v \neq
\bar \beta _v$, this determines the ordering uniquely. By Sublemma \ref{lem-const11} and Lemma \ref{lem-galdef81}, $I _{\bold u } ( \bar \rho \vert _{F_v} , \pi _v ) =  \pi _v ^{K_{11} ( m_v) } $ is isomorphic to $ \bar E _{\mathcal D } [ U ] / ( (U - \chi _1 ( p_v ) ) (U - \chi _2 ( p_v ) )) $, where the action
of $U$ is $U(p_v) $.    \par
It becomes one dimensional after localization at $  (U -
\chi _2 ( p_v ) )$. Since the local intertwining spaces have dimension one, it follows that $e _{\bar \rho}I _{(k, w), \mathcal D_Q} (\pi )$ is one dimensional. \par

This proves that $ T ^{\ss} _{Q , E_\mathcal D }  = T_{Q,  E_\mathcal D  }$ as in Lemma \ref{lem-modular23}. In particular $T _{Q}$ is reduced. From this, we construct Galois representation 
$$
\rho_Q ^{\modular}: G _{\Sigma _{\mathcal D_Q} } \longrightarrow \GL_2 ( T^{\ss} _{Q})
$$
as in Proposition \ref{prop-modular32}.  We show $T_{Q }  ^{\ss} = T_{Q} $ by an argument using Galois representations as in Theorem \ref{thm-modular41}. The only point we must discuss is to recover $V (\delta )  $-operators for $\delta \in \Delta _Q$ from $\rho_Q ^{\modular} $. 

\begin{lem}\label{lem-const12}
\begin{enumerate}
\item For $ v \in Q$, $\rho ^{\modular} _Q \vert_{ G_{F_v}} $ is a sum of two
characters $ \chi _{1, v} ,\ \chi _{2 , v}: G _{F_v}  \to T ^{\times}_Q 
$.
\item  Fix the ordering of the characters such that $ \bar \chi _{2, v} (\Fr _v ) =\bar  \alpha _v$.
For an element $ \sigma _\delta $ of $ I ^{\tame} _{F_v }$ which lifts $\delta \in \Delta _v$ via $ I ^{\tame} _{F_v }\to  (I ^{\tame} _{F_v }) _{G_F } \simeq k (v ) ^{\times} $, $ \chi _{2, v} ( \sigma _{\delta}) = V(\delta )$. 
\end{enumerate}
\end{lem}
\begin{proof}[Proof of Lemma \ref{lem-const12}] 
(1) follows from Faltings' theorem \ref{claim-galdef2}.  \par

For (2), take a representation $\pi \in \mathcal A _{(k,  w) } ( G _{\bar \rho } ) $ such that $e _{\bar \rho}I _{(k, w), \mathcal D_Q} (\pi )\neq \{ 0 \}$. It suffices to check the identity in (2) on the component of $T_Q$ corresponding to $\pi$. For a finite extension $E_{\lambda} $ of $E_{\mathcal D } $ such that $\pi _f $ is defined over $E_{\lambda } $, let $f _{\pi } :  T^{\ss} _{Q} 
\to o _{E _\lambda }$ be the $o_{\mathcal D}$-homomorphism corresponding to $\pi$. $f_{\pi }  $ induces $\rho : G_F \to \GL _2 (T^{\ss} _{Q} ) \to \GL _2 (o_{E_\lambda}  )  $, which is isomorphic to the $\lambda $-adic representation associated to $\pi$. By the compatibility of the local and the global Langlands
correspondence, $\pi _v$ corresponds to 
$\rho \vert _{G_{F_v}} $ by the local Langlands correspondence. $\rho \vert _{G_{F_v}} \simeq \tilde  \chi _{1, v} \oplus \tilde \chi_{2, v}$, where $ \chi _{i, v} $ is the composition of $G_{F_v} \overset{\chi _{i, v} } \to T ^\times _{Q} \overset {f_\pi } \to o^{\times}_{E_\lambda} $. It follows that $ \pi _v$ is isomorphic to the principal series $ \pi ( \tilde \chi _1 , \tilde \chi_2) $. By Lemma \ref{lem-const11}, the action of $V(\delta ) $ on $\pi ^{K_{11} ( m_v) } _v  $ is $\tilde \chi _2 (\sigma _{\delta}) $, and the claim follows. 
\end{proof}
By Lemma \ref{lem-const12}, $V(\delta ) $ belongs to $ T_{Q }  ^{\ss} $. We conclude that $T_{Q }  ^{\ss} = T_{Q} $, and 
$T_{0, Q }  ^{\ss} = T_{0, Q} 
$ follows from this as a consequence. 
\end{proof}

Let $T ^{\ss} _{\mathcal D _y }$ and $T _{\mathcal D _y  }  $ be the $\ell$-adic Hecke algebras defined in Definition \ref{dfn-modular22}. $T ^{\ss} _{\mathcal D _y } = T _{\mathcal D _y  }  $, and $ T _{\mathcal D _y  } \simeq T_{\mathcal D}$ by Theorem \ref{thm-modular41} and Proposition \ref{prop-modular41}. \par

For $\ell$-adic reduced Hecke algebras we have defined, there are $o _{\mathcal D }$-algebra homomorphisms 
$ T^{\ss} _{Q}  \to T^{\ss}_{0, Q } \to T^{\ss} _{\mathcal D_y }$. By the definition and Lemma \ref{lem-modular31}, these homomorphisms are surjective, and induce surjective homomorphisms between $\ell$-adic Hecke algebras
$$
T _{Q}  \twoheadrightarrow T_ {0, Q } \twoheadrightarrow T_{\mathcal D_y }  .
$$
by Proposition \ref{prop-const11}. 
\begin{prop}
\label{prop-const12}
The above $o_{\mathcal D}$-algebra homomorphism gives an isomorphism
$$
T_{0, Q  } \simeq  T _{\mathcal D_y }  .
$$
\end{prop}
\begin{proof}[Proof of Proposition \ref{prop-const12}] First we show $T_{0, Q , E_{\mathcal D} } \simeq  T _{\mathcal D_y , E_{\mathcal D}  }  $. \par
Let $e _{0, \bar \rho}$ be the idempotent of $ \tilde T_{0, Q} $ corresponding to $T_{0, Q }$. Take a representation $\pi \in \mathcal A _{(k, w ) } ( G_{\bar \rho } ) $ such that $e^0 _{\bar \rho}I _{(k, w), \mathcal D_Q} (\pi )\neq \{ 0 \} $. We show $\pi _v$ is spherical for $v \in Q$.\par
By our choice of $K _{0, Q}$, the central character of $\pi_v$ is unramified.  If $\pi_v$ is not spherical, it is
an unramified special representation since it is fixed by $K_0 ( m_v)$. Let $\rho _{\pi , \bar E _{\mathcal D} }$ be the Galois representation associated to $\pi$. By the compatibility of the local and the global Langlands
correspondence, $\rho _v = \rho _{\pi, \bar E _{\mathcal D} }\vert _{G_{F_v} }\simeq
\sp (2) \otimes \chi_v
$, where $\sp (2) $ is the special representation, and $\chi_v$ is an unramified character. This implies that for any Frobenius lift $\sigma \in G_{F_v}$,  the two eigenvalues $\alpha_v$ and $\beta_v$ of $\rho _v (\sigma) $ must satisfy $\alpha _v /
\beta _v = q_v ^{\pm1}$. Since $ q_ v
\equiv 1 \mod \ell $ and
$\bar \alpha_v $,
$\bar \beta_v$ are distinct, this is a contradiction. Thus $\pi _v$ is spherical at $v \in Q$, and $\pi $ appear as a component of $T _{\mathcal D_y }  $. \par
As in the proof of Proposition \ref{prop-const11}, $e^0 _{\bar \rho}I ^0 _{(k, w), \mathcal D_Q} (\pi )$ is one dimensional over $\bar  E _{\mathcal D}$, and the action of $T_v \in T _{0, Q}$ for $ v \not \in \{ y \} \cup S \cup \Sigma _{\mathcal D _Q}$ is the same as the action of $T_v \in T _{\mathcal D_y } $ on $e _{\bar \rho}I ^y _{(k, w), \mathcal D_y } (\pi ) $ in the notation of \S\ref{subsec-modular4}. Thus we have an isomorphism $ T_{0, Q , E_{\mathcal D} } \simeq  T _{\mathcal D_y , E_{\mathcal D}  }  $, which induces $ T^{\ss}_{0, Q  } \simeq  T ^{\ss}_{\mathcal D_y  } $. \par
To show $ T_{0, Q , E_{\mathcal D} }  \simeq T _{\mathcal D_y , E_{\mathcal D}  }  $, it suffices to see $U ( p_v )  $ and $U(p_v, p_v ) $ belong to $ T^{\ss}_{0, Q  } $. For $v \in Q$, note that there are elements $ T_v $ and $T_{v, v}$ in $T^{\ss}_{0, Q  } \simeq T _{\mathcal D_y , E_{\mathcal D}  }  $.  \par
For $v \in Q$, $U ( p_v, p_v ) = T_{v, v} $, and $U (p_v ) $ satisfies the equation
$$
U (p_v)^2  - T_v \cdot U (p_v) +  q_v T_{v, v}  =0 
\leqno{(\ast)}
$$
by Lemma \ref{lem-galdef81}. By our assumption, it has two
distinct roots in the residue field of $T ^{\ss}_{\mathcal D_y }$, and hence there are two distinct roots in
$T ^{\ss}_{\mathcal D_y }$ by Hensel's lemma. By our choice of $\tilde m _{0, Q} $, 
$U ( p_v) $ is the element  $\alpha _v $ in $T ^{\ss} _{\mathcal D_y } $ which satisfies ($\ast$), and lifts $\bar \alpha _v$. The claim is shown. 
\end{proof}

\subsection{Twisted sheaves}\label{subsec-const2}

To control the central character when the infinity type $ ( k,w ) $ is not $ ( (2, \ldots, 2) , 0 )$, a passage to the
adjoint group $G_{\bar \rho}  ^{\ad}  $ is necessary.  We slightly modify the definitions of
$S_K$ and the sheaf $\bar {\mathcal F}_{ (k, w) } $. \par

We return to the general setting. For a division algebra $D$ over $F$ with $ q _D \leq 1 $, let $Z$ be the center of $G_D$. For an $F$-factorizable compact open subgroup $U $ of $G
(\A _{\mathbb Q, f}  )$, we define 
$$
C_U = Z (
\mathbb Q )
\backslash Z (\A _{\mathbb Q} )/  U\cap Z (\A  _{\mathbb Q , f} )  \cdot Z(\R )  .
$$

$ Z (\A _{\mathbb Q} ) $ acts on $S_U = G _D ( \mathbb Q) \backslash G_D ( \A _{\mathbb Q}  )/ U \times K_\infty 
$ from the right. This action of $Z(\A _{\mathbb Q} )$ induces a
$C_U
$-action on
$S_U$.

\begin{lem}\label{lem-const21} We assume that there is a finite place $y$ and a finite set $S$ of finite places which satisfy conditions (1)-(3) of Proposition \ref{prop-shimura51}, and $U= U ( y, S) $.  Then the $C_{U}  $-action on $S_U $ is free.
\end{lem}
\begin{proof} [Proof of Lemma \ref{lem-const21}]
 Take an element $ z\in Z (\A _{\mathbb Q _f })$ which represents a class $c$ in $ C_U$. We may assume that the finite part of $z$ belongs to $U  $ by an approximation by an element of $F^{\times} $. By \ref{lem-shimura50}, there is an element $\epsilon \in F^ \times $ such that $\epsilon ^{-1} z \in U  $. Then the class of $z$ in $C_U$ is trivial. 
\end{proof}

\begin{dfn}\label{dfn-const21} Notations are as in Lemma \ref{lem-const21}. For a prime $\ell$, let $C_{U, \ell} $ be the $\ell$-Sylow subgroup of $C_U$. We denote by
$
\tilde S_ U $ the quotient of $S_U $ by $  C_{ U, \ell } $.
\end{dfn}

We view the pro-$\ell$-part $\chi_{\cycle, \ell }$ of $\chi_{\cycle}$ as a
continuous character of $F^ \times \backslash ( \A_{F, f}  )^\times
$. We also view it as a continuous character of $G_D (\mathbb Q)  \backslash  G_D (\mathbb A _{\mathbb Q, f} ) $ via the reduced norm $ G_D (\mathbb A _{\mathbb Q, f}  ) \overset{ \Norm _D} \to ( \A_{F, f}  )^\times
$. In particular $\chi_{\cycle, \ell }$ can be viewed as a $ G_D (\mathbb A _{\mathbb Q, f}  ) $-equivariant sheaf on $S= \varprojlim _U S_U$. Since $\chi _{\cycle, \ell }$ takes values in $ 1 + \ell \mathbb Z_{\ell} $ and $\ell \geq 3 $, the square root $\chi ^{\frac 1 2} _{\cycle, \ell } $ is defined as an $\ell$-adic character. 
\par

Assume that $U$ satisfies the condition of Lemma \ref{lem-const21}, and the sheaf
$\bar {\mathcal F}^U _{(k , w) , E_\lambda  }
$ is defined on
$S_U$. 
\begin{dfn} \label{dfn-const22} For a prime $\ell \geq 3$, 
$$ 
\bar {\mathcal F } ^{U, \tw}  _{(k , w) } =  \bar {\mathcal F} ^{U}  _{(k , w) }\otimes \chi _{\cycle,
\ell } ^ { -
\frac w 2}.
$$
\end{dfn}
By the definition of $ \bar {\mathcal F } ^{U, \tw}  _{(k , w) } $, it is $ G_D ( \mathbb A^{\ell}_{\mathbb Q , f} )$-equivariant. We analyze the action of $ Z(\A _{\mathbb Q _f}  ) $ in more detail. \par

\begin{lem}\label{lem-const22} For a prime $\ell \geq 3$, let $ \mathcal Z _{\ell} $ be the inverse image of $ C _{U , \ell }$ by $Z ( \mathbb A _{\mathbb Q , f} ) \to C_U $. Then the $\mathcal Z_{\ell}$-action on $ \bar {\mathcal F } ^{U, \tw}  _{(k , w) }$ induces a structure of a $C_U$-equivariant sheaf  on $ \bar {\mathcal F } ^{U, \tw}  _{(k , w) }$, and $ \bar {\mathcal F } ^{U, \tw}  _{(k , w) }$ descends to an $o_{E_{\lambda }} $-smooth sheaf $\tilde {\mathcal F } ^{U, \tw}  _{(k , w) }$ on $\tilde S_U$. 
\end{lem}
\begin{proof}[Proof of Lemma \ref{lem-const22}]
We use the notations in \S\ref{subsec-shimura2}. Let $ \pi _\ell : \tilde S_\ell \to S_U $ be the $\overline {Z ( \mathbb Q) } \cap U_{\ell} \backslash  U _{\ell}$-torsor, $V_{(k, w), E_\lambda } $ the
representation of $G _D (\mathbb Q_{\ell})  $ defining $\bar {\mathcal F}^U_{(k, w), E_\lambda } $.  
We twist $ V_{(k, w ) , E_\lambda } $ by  $\chi _{\cycle, \ell } ^ { - \frac w 2}$, and define
$$
V ^{\tw}_{ (k, w), E_\lambda  } =V _{(k, w), E_\lambda  }\otimes \chi _{\cycle, \ell } ^ { - \frac w 2} 
$$
as a continuous representation of $G _D (\mathbb Q_{\ell})  $. Since $ \chi _{\cycle, \ell } ^ { - \frac w 2} $ is $o^{\times}_{E_\lambda }$-valued, $V _{(k, w), o _{E_\lambda}  } $ gives an $o_{E_\lambda }$-lattice $ V ^{\tw}_{ (k, w),  o_{E_\lambda } } $ of $ V ^{\tw}_{ (k, w), E_\lambda  } $. 
$\mathcal F^{U, \tw} _{(k,w) , E_\lambda} $ on $S_U $ is obtained from $\pi _\ell$ and $G _D (\mathbb Q_{\ell})$-representation $ V ^{\tw} _{(k, w) , E_\lambda } $. \par
For $v \vert \ell$, we denote $\ker (o_{F_v}  ^\times \to  k (v) ^{\times} )$ by $ U ^ 1 _v $. Since $\prod _{v \vert \ell } U ^1 _v   $ acts trivially on $ V ^{\tw}_{ (k, w), E_\lambda  }  $, any element $ z \in Z( \mathbb A _{\mathbb Q, f} ) $ such that the $\ell$-component $z_\ell$ belongs to $\prod _{v \vert \ell } U ^1 _v $ acts trivially on $ \mathcal F^{U, \tw} _{(k,w) , E_\lambda} $ . In particular an equivariant action of  $Z ( \mathbb Q) \backslash Z (\A  _{\mathbb Q , f} )/ H 
$ is induced for 
$$ 
H = ( \prod _{ v \vert \ell} U _v \cap U ^ 1 _v )\cap Z(\mathbb Q_\ell) 
\cdot ( U
\cap Z (\A _{\mathbb Q, f } ) ) ^\ell. 
$$
The $\ell$-part of $ Z ( \mathbb Q) \backslash Z (\A _{\mathbb Q, f}  )/ H $ is canonically isomorphic to $ C_{U, \ell}$, and hence we obtain a $ C_{U,
\ell }$-action on $\bar { \mathcal F} ^{U, \tw}  _{(k, w), E_\lambda} $ which lifts the $C_{U,
\ell}$-action $S_U$. By Lemma \ref{lem-const21} and by descent, we have a sheaf
$
\tilde
{\mathcal F}^{U, \tw}  _{(k, w), E_\lambda}
$ on
$\tilde S_U 
$. Note that all classes of $C_{U, \ell} $ are represented by an element in $Z (\A ^{\ell} _{\mathbb Q, f} )$, and hence the $\mathcal Z_{\ell}$-action on $ \bar { \mathcal F} ^{U, \tw}  _{(k, w), E_\lambda} $ preserves the $o_{E_\lambda} $-structure
$\bar {\mathcal F}  ^{U , \tw} _{(k, w) }  $. Thus we have the $o_{E_\lambda}$-structure $\tilde {\mathcal F}^{U , \tw}  _{(k, w)} $. 
\end{proof}

$$
H ^ * ( S_U, \ \bar {\mathcal F} ^{U , \tw} _{(k, w)}  ) =  H ^ * ( S_U,  \bar {\mathcal F} ^{U } _{(k, w)}  )\otimes
\chi _{\cycle, \ell}^{-\frac w 2} 
$$
holds by the definition of $ \bar {\mathcal F} ^{U , \tw} _{(k, w)} $, and this equality is compatible with the action of the convolution algebra $H ( G _D ( \mathbb A^{\ell}  _{\mathbb Q, f} ) , U ^{\ell }  )_{o_{E_\lambda} }  $, by regarding $
\chi _{\cycle, \ell}^{-\frac w 2}$ as an $\ell$-adic character of $ G_D ( \mathbb A^{\ell}_{\mathbb Q, f} ) $. \par

We also obtain an action of $H ( G _D ( \mathbb A^{\ell}  _{\mathbb Q, f} ) , U ^{\ell }  )_{o_{E_\lambda} }  $ on $H ^ * ( \tilde S_U, \ \tilde {\mathcal F}^{U , \tw}_{(k, w) } ) $ since the $ C_{U , \ell }$-action on $S_U$ comes from the center $Z(\mathbb A_{\mathbb Q, f} )$, and commutes with the
Hecke correspondences. 

\subsection{Twisted algebras and modules}\label{subsec-const3}
 We go back to the situation in \S\ref{subsec-const1}. 
We use the same notation as in \S\ref{subsec-modular4}. For our choice of the auxiliary place $y$ and set $S$, $\Sigma = \Sigma _{\mathcal D _ Q } \cup\{ y \} \cup  S $, $H _{K ^{\Sigma } } $ and $ \tilde m _{\Sigma , \bar \rho }$ are the convolution algebra and the maximal ideal corresponding to $\bar \rho$, respectively.  For $ K ' _{ Q } =\ker \nu _{\mathcal D_ y } \vert _{K _{Q}  }$ (resp. $K ' _{ 0, Q } =\ker \nu _{\mathcal D_ y } \vert _{K _{0, Q}  } $), we
define $\tilde M _{ Q} $ (resp. $\tilde M _{0, Q} $) by 
$$
\tilde  M _{Q   } = \Hom _{K _{ Q}  }  ( \nu _{\mathcal D_ y } \vert _{K_{ Q}  }, H ^{q _{\bar \rho} } _{\stack} ( S_{K ' _{Q}  } , \bar {\mathcal F} ^{K' _{Q}   } _{(k, w) } ) )_{\tilde m _{\Sigma ,  \bar \rho } }
$$
$$
(\text{resp. } \tilde  M _{0, Q   } = \Hom _{K _{ 0, Q}  }  ( \nu _{\mathcal D_ y } \vert _{K_{0, Q}  }, H ^{q _{\bar \rho} } _{\stack} ( S_{K ' _{0, Q}  } , \bar {\mathcal F} ^{K' _{0, Q}   } _{(k, w) } ) )_{\tilde m _{ \Sigma , \bar \rho } }) .
$$
Similarly, $\tilde  M^{\tw} _{Q   } $ and $\tilde M ^{\tw}_{0, Q} $ are defined as above by using $\bar {\mathcal F} ^{\tw} _{(k, w) } $ instead of $\bar {\mathcal F} _{(k, w) }  $.  \par
The isomorphism $\tilde M ^{\tw} _{Q} = \tilde M_{Q} \otimes \chi ^{- \frac w 2 }_{\cycle, \ell}  $ (resp. $\tilde M ^{\tw} _{0, Q} = \tilde M_{0, Q} \otimes \chi ^{- \frac w 2 }_{\cycle, \ell}  $) respects the action of $ H ( G _{\bar \rho} ( \mathbb A _{\mathbb Q, f} ) , K ^{\ell} _{0 , Q} ) $ (resp. $ H ( G _{\bar \rho} ( \mathbb A _{\mathbb Q, f} ) , K ^{\ell} _{ Q} ) $ ) by viewing $\chi ^{- \frac w 2 }_{\cycle, \ell} $ as a one dimensional representation of the convolution algebra. This is also true for $\tilde U (p_v ) $ and $\tilde U (p_v, p_v )$-operators. In particular for the $o_{\mathcal D} $-algebra $\tilde T^{\tw}  _{Q}  $ (resp. $\tilde T ^{\tw} _{0, Q}$) generated by $ T_v^{\tw} ,\  T_{v, v} ^{\tw} \ (v \not \in \Sigma_{\mathcal D_{Q}}  \cup \{ y \} \cup S)$, $ U ^{\tw} ( p_v )$, $ U ^{\tw} ( p_v , p_v )$ fo $v \in P ^{\bold u } _{\mathcal D_Q}\setminus P ^{\exc}$,  $ \tilde U ^{\tw} ( p_v )$, $ \tilde U ^{\tw} ( p_v , p_v )$ for $v \in P ^{\bold {n.o.}} $, 
$ V^{\tw}  (\delta ) $ for $\delta \in \Delta _Q$ ($ (- ) ^{\tw}$ denotes the Hecke operators defined by the twisted action on $ \bar {\mathcal F} ^{\tw} _{(k, w) } $), there is an $o_{\mathcal D}$-isomorphism
$$
\tilde T_Q  \simeq \tilde  T ^{\tw}_Q \quad (\text{resp. } \tilde T_{0, Q}  \simeq \tilde  T ^{\tw}_{0, Q })
$$
which maps the standard operators $(-)$ as above ($T_v$, for example) to $ ( - ) ^{\tw} \cdot \chi ^{\frac w 2 } _{\cycle, \ell } $. Let $ \tilde m^{\tw}_Q $ (resp. $ \tilde m^{\tw}_{0, Q}$) be the maximal ideal of $  \tilde T^{\tw}  _{Q} $ (resp. $ \tilde T ^{\tw} _{0, Q}$), $T^{\tw}  _{Q}$ (resp. $ T^{\tw}  _{0, Q}$) the localization $( \tilde T^{\tw}  _{Q} ) _{\tilde m^{\tw}_Q} $ (resp. $( \tilde T^{\tw}  _{0, Q} ) _{\tilde m^{\tw}_{0, Q}} $). \par
The corresponding $T_Q$ (resp. $T_{0, Q} $)-module $M_Q = ( \tilde M _Q )_{\tilde m_Q }  $ (resp. $M_{0, Q} = ( \tilde M _{0, Q} )_{\tilde m_{0, Q} }  $) and $T^{\tw}_Q$ (resp. $T^{\tw}_{0,Q} $)-module $ M^{\tw} _Q= ( \tilde M ^{\tw} _Q) _{\tilde m^{\tw} _Q} $ (resp. $ M^{\tw} _{0, Q}= ( \tilde M ^{\tw} _{0, Q}) _{\tilde m^{\tw} _{0,Q}}$) are $o_{\mathcal D}$-free as in \S\ref{subsec-modular4} by using Proposition \ref{prop-coh21}, and the $o_{\mathcal D} $-algebra isomorphism 
$$
T_Q \simeq T^{\tw} _Q  \quad (\text{resp. }T_{0 , Q} \simeq T^{\tw} _{0, Q}  )
$$ 
induces
$$
M _Q \simeq M ^{\tw} _Q  \ \quad (\text{resp. }M _{0, Q} \simeq M ^{\tw} _{0, Q } ).
$$

\par
The same construction applies also to $ T _{\mathcal D _y } $ and $M _{\mathcal D} $. $T ^{\tw}_{\mathcal D _y }  $ and $M ^{\tw}_{\mathcal D}  $ denotes the twisted algebra and the module which satisfy $T _{\mathcal D _y } \simeq T ^{\tw}_{\mathcal D _y }  $ and $M _{\mathcal D} \simeq M ^{\tw}_{\mathcal D}  $.\par

\begin{prop}\label{prop-const31}
There is a natural isomorphism 
$$
 M_{\mathcal D}  \simeq  M_{0, Q}, 
 $$ 
which is compatible with the isomorphism
$T_{\mathcal D} \simeq T_{\mathcal D_y } \simeq  T_{0, Q }$ obtained by Proposition \ref{prop-modular41} and \ref{prop-const12}. 
The same is true for $ M ^{\tw} _{\mathcal D}
$ and $ M_{0, Q} ^{\tw}$.
\end{prop}
\begin{proof}[Proof of Proposition \ref{prop-const31}] We construct such an isomorphism by induction on $ n = \sharp Q$. For $n= 0 $ we take it as the identity.
Assume that $ Q =Q' \cup  \{ v \}$. The degeneracy maps $\pr _{i} : S_{K ' _{0, Q} } \to S_{ K ' _{0, \emptyset} } $ for $i = 1, 2$ defined in \S\ref{subsec-coh3} induce 
$$
H ^{q _{\bar \rho} } _{\stack} ( S_{K ' _{0, Q' }  } , \bar {\mathcal F} ^{K' _{0, Q ' }   } _{(k, w) }  )^{ \oplus 2} \overset{\pr_1 + \pr _2   }\longrightarrow H ^{q _{\bar \rho} } _{\stack} ( S_{K ' _{0, Q}  } , \bar {\mathcal F} ^{K' _{0, Q}   } _{(k, w) } ) . 
$$
Thus we have a $\tilde T_Q$-homomorphism $\tilde \beta _Q:  \tilde M _{Q' } ^{\oplus 2} \to \tilde M _{0, Q}$. 

By the decomposition of $ H ^ {q_{\bar \rho }} $ in \S\ref{subsec-shimura4}, we have an isomorphism
$$
(I ^ 0_{(k, w ) , \mathcal D_Q }  )^{\oplus 2 ^{q _{\bar \rho}  }  } _{\tilde m _{\Sigma ,  \bar \rho } } \simeq  \tilde  M _{0, Q   } \otimes _{o _\mathcal D } \bar E_{\mathcal D} , 
$$
which induce $\beta _{Q, E_\mathcal D} : \tilde M _{Q' , E_{\mathcal D}} ^{\oplus 2} \simeq \tilde M _{0, Q , E_{\mathcal D} }$ since the $v$-component of the representations which contibutes to $I ^ 0_{(k, w ) , \mathcal D_Q }  $ is spherical as in the proof of Proposition \ref{prop-const12}. In particular $\tilde T _Q = \tilde T_{Q'} [ U(p_v )  ]$, and $\tilde \beta_Q$ is regarded as a $\tilde T_{Q'} [ U(p_v )   ]   $-homomorphism. \par
We need a special case of cohomological universal injectivity. 
\begin{sublem}
\label{sublem-const31} Let $D$ be a division algebra with $q_D\leq 1$ which is split at a finite place $v$. Assume that $K$ is an $F$-factorizable compact open small subgroup of $G_D ( \mathbb A_{\mathbb Q_f} ) $ with the $v$-component $K_v = \GL_2 (o_{F_v} ) $. Then the kernel of the homomorphism defined by the degeneracy maps at $v$
$$
 H ^{q_D }  ( S _K  ,\  \bar {\mathcal  F} ^ K  _{ ( k, w ), k_\lambda  } )  ^ {\oplus 2} 
\overset { \pr _1 ^* + \pr _2 ^* } \longrightarrow H  ^{q_D}  ( S _{ K\cap K_{0 } (v)  },\ \bar 
{\mathcal F}^{K \cap K_0 (v) }   _{  (k, w) , k_{\lambda} }  )
$$
is annihilated by $ T ^2 _v - (1+q_v)^2T_{v, v} $.
\end{sublem}
\begin{proof}[Proof of Sublemma \ref{sublem-const31}] 
Consider the homomorphism 
$$
H  ^{q_D} ( S _{ K\cap K_{0 } (v)  },\
\bar {\mathcal F} _{  (k, w) , k_\lambda  }  ) \overset \gamma \longrightarrow H ^{q_D}  ( S _K  ,\ \bar {\mathcal F} _{  (k, w) , k_{\lambda}}
)  ^ {\oplus 2}
$$
induced by the trace map. The composite
$$
\delta :  H  ^{q_D }  ( S _K  ,\  \bar {\mathcal F} _{  (k, w), k_\lambda   } )  ^ {\oplus 2} 
\overset { \pr _1 ^* + \pr _2 ^* } \longrightarrow H  ^{q_D }  ( S _{ K\cap K_{0 } (v)  },\
\bar {\mathcal F} _{  (k, w)  , k_{\lambda} }  )   \overset \gamma \longrightarrow  H  ^{q_D } ( S _K  ,\  \bar {\mathcal F} _{(  k, w), k_{\lambda} }
)  ^ {\oplus 2}
$$
written in the matrix form is 
$$
\delta = 
\begin{pmatrix}
 1 + q _v &  T _v\\
  T_{v,v} ^{-1}\cdot T_v  & 1+q_v \\
\end{pmatrix}
,
$$
which is shown in \S\ref{subsec-congruent2}.
If 
$
\begin{pmatrix}
 a\\
b
\end{pmatrix}
$
is in the kernel of $\delta$, $T_v ^ 2 -( 1+ q_v ) ^
2 T_{v, v}  $ annihilates
$
\begin{pmatrix}
 a\\
b
\end{pmatrix}
$
by a direct calculation. The claim follows. 
\end{proof}
By Sublemma \ref{sublem-const31}, the image of $ \tilde \beta _Q$ is an $o_{\mathcal D}$-direct summand. 
 We already know that $\tilde \beta_{Q, E_{\mathcal D} }$ is an isomorphism. Thus $ \tilde \beta _Q $ is an isomorphism. \par
The localization of $\tilde M _{ Q'}$ at $\tilde M _{Q'} $ is isomorphic to $M_{\mathcal D}$ by our assumption on the induction. $M_Q$ is isomorphic to the localization of $M ^{\oplus 2} _{Q'}  $ at the maximal ideal of $ T _{Q'} [U ( p_v)] $ generated by $m _{T_{Q' }}  $ and $U(p_v) -\bar \alpha _v$. As in the proof of Proposition \ref{prop-const12}, $ U ( p_v)$ safisfies $f(U (p_v )) = 0 $ for  $f (X) = X ^2 - T_v X + q_v T_{v, v} $ and $f$ has two distinct roots in $  T _{Q'}$. For the root $\alpha_v $ of $f(X)$ which lifts $\bar \alpha _v$, $M_Q$ is the submodule of $M ^{\oplus 2} _{Q'}   $ where $U(p_v ) $ acts as $ \alpha _v$, which is identified with $M _{Q'} $. 
\end{proof}

\subsection{Perfect complex argument} \label{subsec-const4}

In this subsection, we verify the axiom TW5 for $M_Q$. 

If we set $K' _Q= \ker \nu _{\mathcal D_y } \vert _{K_Q ( y, S) }$ and $K' _{0,Q}= \ker \nu _{\mathcal D_y } \vert _{K_{0,Q} (y, S)} $, then
$$
 \Delta _Q  \overset \alpha  \simeq  \overline {Z ( \mathbb Q ) } \cdot K '_{Q}  \backslash K '_{0, Q} 
$$
holds. Here, $\alpha (\delta ) $ for $\delta \in \Delta_Q$ is given by
the class of $ \begin{pmatrix}
 \delta  ^{-1} &0\\ 
 0 &{\delta }
\end{pmatrix}
$. In particular $ M_Q$ has a structure of an $o_\mathcal D [
\Delta_Q]$-module. $\delta \in \Delta _Q$ acts by $V(\delta)$-operator. 

\begin{lem}\label{lem-const41}
$\Delta _Q\times   C_{K'_Q,
\ell}$-action on $S_{K'_Q}$ is free. 
\end{lem}

\begin{proof}[Proof of Lemma \ref{lem-const41}] Let $ c$ be a class in $ C_{K'_Q, \ell}$ represented by $ z\in Z ( \mathbb A_{\mathbb Q, f} ) $. We may assume that the $v$-component of $z$ belongs to $ K_{11} (m_v)$ for $v \in Q$, and $ z \in K '_Q$ by an approximation by an element of $F^\times$. For $\delta \in \Delta _Q$, $ c (\delta ' ) \in G_{\bar \rho} (\A _{\mathbb Q, f} ) $
is the element with 
$ c (\delta ' ) _Q =
\begin{pmatrix}
 {\delta '} ^{-1} &0\\ 
 0 &{\delta '}
\end{pmatrix}
$ for some $\delta'  \in \prod _{ v \in Q} o^{\times}_{F_v}$ lifting $\delta $.  \par
Assume that $z \cdot c (\delta ') $ has a fixed point in $S_{K_Q} $. By Lemma \ref{lem-shimura50}, and there is some unit $\epsilon \in o^{\times } _F $ such that $ \epsilon ^{-1} z c (\delta ' ) \in K_Q$. This implies that $\delta = 1$, and the class $c$ of $z$ is trivial. 
\end{proof}

 $C_{K'_{0, Q} ,\ell} $, $C_{K'_Q , \ell} $ are
isomorphic to $C_{F, P ^{\bold u } _{\mathcal D} , \ell}$, the $\ell$-part of the class group of $F$ which is unramified outside $P ^{\bold u } _{\mathcal D} $. 
\begin{prop}\label{prop-const41}
$ M^{\tw} _Q $ is an
$
o_{\mathcal D}  [C_{F, P ^{\bold u } _{\mathcal D} ,  \ell}\times 
\Delta _Q  ]$-free module. The same is true for $ M_Q$.
\end{prop}
\begin{proof}[Proof of Proposition \ref{prop-const41}]
 For $F = \mathbb Q$ and $ D = M_2 (\mathbb Q) $, this is essentially proved in \cite{TW}, proposition
1. We may assume that
$ D$ is a division algebra. We prove \ref{prop-const41} by the perfect complex argument as in \cite{Fu1}, \S3. \par

For $\Sigma = \Sigma _{ \mathcal D_Q} \cup \{ y \} \cup S$, let $ H = H ( D ^\times ( \A ^ {\Sigma} _f ),\ K ^ {\Sigma } ) _{o_{\mathcal D} } $ be the convolution algebra, and $\tilde m_{\Sigma , \bar \rho} $ is the maximal  maximal ideal corresponding to $\bar \rho$. 

Since $M _Q^{\tw}$ is a direct summand of $H ^{q_{\bar \rho} }   ( S_{K' _Q} ,\  \bar {\mathcal F} ^{\tw} _{(k, w)})
_{\tilde m_{\Sigma , \bar \rho}} $, it suffices to show the following claim. 
\end{proof}
\begin{lem}\label{lem-const42}
$ H ^ {q_{\bar \rho} } ( S_{ K'  _ Q} ,\
\bar { \mathcal F} _{ (k, w)  } ^{\tw} ) _{\tilde m_{\Sigma , \bar \rho}}
$ is
$o_{\mathcal D} [C_{F, P ^{\bold u } _{\mathcal D} ,  \ell} \times  \Delta _Q] $-free.  
\end{lem}
\begin{proof}[Proof of Lemma \ref{lem-const42}]
We need two sublemmas. 

\begin{sublem} \label{sublem-const41}
Let $A$ be a noetherian local algebra with the maximal ideal
$m_A 
$ and the residue field $k_A$, $B$ a commutative $A$-algebra. Let 
$L^\cdot 
$ be a complex of
$B$-modules bounded below. Assume the following conditions:
\begin{enumerate}
\item $ L ^\cdot  $ defines an object of the derived category $ D ^b_{\coh} (A) $ of bounded $A$-complexes with finitely generated cohomology groups (in particular $ L ^\cdot \otimes^{\mathbb L} _A k _A$ is defined).
\item For a maximal ideal
$m$ of
$B$ above $m_A$, $ H ^i (L ^\cdot \otimes^{\mathbb L} _A k _A) \otimes _B B_m  $ is zero for
$i \neq 0$. Here $B_m $ is the localization of $B$ at $m$.
\end{enumerate}
Then
$ H ^i( L^\cdot  ) \otimes _B B_m $ is $A$-free for $i = 0 $, and zero for $i \neq 0$. 
\end{sublem}
\begin{proof}[Proof of Sublemma \ref{sublem-const41}]  This is lemma 3.2 of \cite{Fu1}. The proof is simple so we include it here.
By (1), by replacing $L ^\cdot $ by the canonical truncation $ \tau _{\leq N } L ^\cdot $ for a sufficiently large integer $N$, we may assume that $L ^\cdot $ is bounded above.  Since
$B_m
$ is
$B$-flat, we may assume that 
$B= B_m
$.  $H ^i (L ^\cdot \otimes^{\mathbb L} _Ak_A)$ is zero except $i = 0 $. By taking the minimal resolution of $L^\cdot$ as $A$-complexes, the claim follows. 
\end{proof}

\begin{sublem}
\label{sublem-const42}
Let $\pi : X \to Y $ be an \'etale morphism between schemes of finite type over an algebraically closed field, which is a $G$-torsor for a finite group $G$, $\mathcal F$ a smooth $\Lambda $-sheaf on $Y$, where $\Lambda$ is a finite $\mathbb Z_\ell $-algebra. Then
$ R
\Gamma (X,\
\pi ^*
\mathcal F )$ is a perfect complex of $\Lambda [G] $-modules, and
$$ 
 R \Gamma (X,\  \pi ^* \mathcal F ) \otimes^{\mathbb L} _{ \Lambda [ G] } \Lambda [G] / I _G
\rightarrow  R\Gamma ( Y,\ \mathcal F ) 
$$
holds. Here $I_G $ is the augmentation ideal, and
the map is induced by the trace map.
\end{sublem}
This is known (usually in the dual form) in any standard cohomology theory.\par
\bigskip
 
We go back to the proof of Lemma \ref{lem-const42}.  $X=  S_{ K' _ Q}$, $Y=\tilde  S_{  K '_{0, Q}
}$, $\mathcal F = \tilde {\mathcal F}_{(k,w) }$. We view $X$ and $Y$ as a scheme over $\mathbb C$. 
The natural projection $
\pi : X\to Y
$ is an
\'etale $G$-torsor by Lemma \ref{lem-const41} for $G= C_{F, P ^{\bold u} _{\mathcal D } ,  \ell} \times \Delta _Q  $.\par
By the same argument as in Proposition \ref{prop-shimura31}, $R \Gamma (X, \ \pi ^* \mathcal F ) $ is represented by a complex $L^\cdot$ of
$H \otimes_{o_{\mathcal D}}o_{\mathcal D}[ G]$-modules bounded below. \par

If $q_{\bar \rho } = 0 $, Lemma \ref{lem-const41} is clear since the $G$-action on $ X$ is free. We may assume that $q_{\bar \rho }= 1 $. For $A = 
o_{\mathcal D}[ G]$, and $B = A \otimes_{o_{\mathcal D}} H $, let $m$ be a maximal ideal of $B$ which is above
$m_A$ and $\tilde m_{\Sigma , \bar \rho }  $. As in Lemma \ref{lem-coh21}, the action of $ H$ on $ H ^ q ( Y, \
\mathcal F
\otimes_{o_{\mathcal D}} k ) $ is of residual type for $ q = 0, 2 $, and hence vanishes after the localization at $m$. 
By Sublemma \ref{sublem-const42}, $ H ^ q( L ^\cdot \otimes _A A/ m_A) \simeq H ^ q ( Y,
\  \mathcal F
\otimes_{o_{\mathcal D} } k )$ for $ q \in \mathbb Z$, and the localization $ H ^ q( L ^\cdot \otimes _A A/ m_A)_m $ at $m$ is zero except for 
$ q = 1$. By Sublemma \ref{sublem-const41}, $ H ^ 1 ( L ^\cdot )_m  = H ^ 1 (X, \ \pi ^* \mathcal F ) _{\tilde m_{\Sigma , \bar \rho } }$ is
$o_{\mathcal D}  [ G] $-free. 
\end{proof}
\bigskip 
As a byproduct of the proof of Lemma \ref{lem-const42}, we have an isomorphism
$$
H ^ {q_{\bar \rho} } ( S_{K'_Q} ,\  \bar {\mathcal F} ^{\tw} _{(k,w) } ) _{\tilde m_{\Sigma , \bar \rho }  }\otimes_{o _\mathcal D  [C_{F,
P ^{\bold u} _{\mathcal D} , \ell}
\times 
\Delta_Q] }
o _\mathcal D  [C_{F,
P ^{\bold u} _{\mathcal D} , \ell}
\times 
\Delta_Q]  / I_{C_{F,
P ^{\bold u} _{\mathcal D} , \ell}
\times 
\Delta_Q}
\leqno{(\ast)}$$
$$ \stackrel{\sim} {\longrightarrow} H ^ {q_{\bar \rho}} (
\tilde S_{K'_{0, Q}} ,\
 \bar {\mathcal F} ^{\tw} _{(k,w) } )_{\tilde m_{\Sigma , \bar \rho }  } . 
$$
Here $I _{C_{F,
P ^{\bold u} _{\mathcal D} , \ell}
\times 
\Delta_Q}$ is the augmentation ideal of $ o _\mathcal D  [C_{F, P ^{\bold u} _{\mathcal D} , \ell}
\times  \Delta_Q] $. Similarly, by using the covering $ S _{K'_Q} \to S_{K'_{0 , Q}}
$ instead of $ S _{K'_Q} \to \tilde S_{K'_{0 , Q}}$, 
$$
H ^ {q_{\bar \rho} } ( S_{K'_Q} ,\  \bar {\mathcal F} ^{\tw} _{(k,w) } ) _{\tilde m_{\Sigma , \bar \rho }  }\otimes_{o ^{\associate}_\mathcal D  [\Delta_Q] }
 o ^{\associate}_\mathcal D  [
\Delta_Q]  / I_{\Delta _Q} \stackrel {\sim} {\longrightarrow} H ^ {q_{\bar \rho}} (
 S_{K'_{0, Q}} ,\
 \bar {\mathcal F} ^{\tw} _{(k,w) } )_{\tilde m_{\Sigma , \bar \rho }  } 
\leqno{(\ast \ast)}
$$
holds for 
$$
o ^{\associate}_{\mathcal D} =o _\mathcal D  [C_{F, P ^{\bold u} _{\mathcal D} , \ell}], 
$$
and the augmentation ideal $I_Q$ of $o ^{\associate}_{\mathcal D} [\Delta_Q]$ over $ o ^{\associate}_{\mathcal D} $. It is easy to check that the
isomorphisms ($\ast$) and ($\ast \ast$) commute with the action of the standard Hecke operators. By using Proposition \ref{prop-const31}, we have
\begin{prop}\label{prop-const42} $M_Q ^{\tw} $ is $o ^{\associate}_{\mathcal D}  [ \Delta _Q] $-free, and 
$M_Q ^{\tw} \otimes _{ o ^{\associate}_{\mathcal D} } o ^{\associate}_{\mathcal D}  [ \Delta _Q] /
I_Q\simeq  M^{\tw}_{\mathcal D} $ as
a $T_Q$-module.  The same is true for $M_Q$ and $M_{\mathcal D}$.  
\end{prop}

\subsection {Construction of the system}\label{subsec-const5}
 Finally we prove Theorem \ref{thm-const0} by summarizing the result of the previous subsections. \par
Let $\chi^{\Hecke} _{\bar \rho } $ be the algebraic Hecke character of weight $2w$ attached to $\det \bar \rho $ (see Definition \ref{dfn-nearlyordinary31}). $\chi = (\det \bar \rho)_{\lift}\cdot  \chi _{\cycle , \ell }^{-w} $ is the $\ell$-adic representation attached to $\chi^{\Hecke} _{\bar \rho } $, where  $(\det \bar \rho)_{\lift} $ is the Teichm\"uller lift of $\bar \rho$.\par
For a minimal deformation condition $\mathcal D $, $ o^{\associate}_{\mathcal D}= o _{\mathcal D} [C_{ F, \ell} ]$ as in Theorem \ref{thm-const0}. 
$o^{\associate}_{\mathcal D} $ is seen as a universal deformation ring of $\det \bar \rho$ by the class field theory: for the character 
$$
 \mu ^{\univ}_{F, \ell}: G^{\ab} _F\stackrel{\sim}{\longrightarrow} \pi _0 ( F ^{\times } \backslash \mathbb A^{\times } _F ) \twoheadrightarrow C_{ F, \ell} \hookrightarrow (o _{\mathcal D} [C_{ F, \ell} ]) ^{\times}, 
$$
$ \chi \cdot  \mu ^{\univ}_{F, \ell}$ is the universal deformation of $\det \bar\rho $, where the local conditions are finite at $v \nmid \ell$, and the restriction to $I_{F_v} $ is $\chi \vert _{I_{F_v}}$ for $v \vert \ell$.  \par
For the universal deformation ring $R_{\mathcal D} $ and the universal Galois representation $\rho ^{\univ} _{\mathcal D}$, 
$ \det \rho _{\mathcal D } ^{\univ}$ satisfies the deformation conditions of $ \det \bar \rho  $ as above. Thus
$R_{\mathcal D  } $ has a structure of $ o^{\associate}_{\mathcal D}$-algebra, and $ \det \rho _{\mathcal D } ^{\univ}$ is identified with $ \chi \cdot \mu_F
^{\univ}
$. By the algebra homomorphism $\pi _{\mathcal D} : R_{\mathcal D } \to T_{\mathcal D} $ given by $\rho ^{\modular}_{\mathcal D}$, $T _{\mathcal D} $ is also regarded as an $o^{\associate}_{\mathcal D} $-algebra. \par
For $ Q \in \mathcal X_{\mathcal D} $, $(T_Q , M_Q)$ is the pair of the Hecke ring and the Hecke module defined in \S\ref{subsec-const1} and \S\ref{subsec-const3}, and $\rho^{\modular}_Q  $ is the universal modular representation. By our construction, $\det \rho ^{\modular}_Q   $ is unramified at $v \in Q$. Thus we have an $o^{\associate}_{\mathcal D}$-algebra structure on $T_Q$ as above. 
 $R_Q$ is the universal deformation ring of $ \bar \rho $ of deformation type $\mathcal D_Q $ over $o^{\associate}_{\mathcal D } $ with the {\it determinant fixed} to $\chi\cdot \mu_F ^{\univ}$. The universal deformation is denoted by $\rho^{\univ} _Q $. The canonical algebra homomorphism $\pi _Q : R_Q \to T_Q $ is an $o^{\associate}_{\mathcal D} $-homomorphism, which is surjective since $ \tr _{T_Q} ( \Fr _v ) $ for $ v \not \in \Sigma_{\mathcal D_Q} \cup \{ y \} \cup S$ generate $T_Q$ by Proposition \ref{prop-const11}. \par

Note that this $o^{\associate}_{\mathcal D } $-structure on $T_Q$ is compatible with
the $o^{\associate}_{\mathcal D } $-module structure on $ M_Q$ given by $M_Q \simeq M_Q ^{\tw} $ and the $C_{F,
\ell}$-action on $M_Q ^{\tw}  $ given in \S\ref{subsec-const2}: the $C_{F,
\ell}$-action is defined by the action of $Z(\A_{\mathbb Q, f} )$ on $ S _{K ' _Q}$, and $Z(\A_{\mathbb Q, f})$-action yields the central character for any representation $\pi$ which appear as a component of $T_Q$. The central character  corresponds to the determinant (twisted by $\chi_{\cycle}$) by the Langlands correspondence, so the compatibility is shown.  \par

For the natural projections $ \alpha _Q : R _Q
\twoheadrightarrow R _{\mathcal D} $, $\beta _Q : T_Q \twoheadrightarrow T _{\mathcal D} $, the diagram
$$
\CD R_Q  @> \pi _Q >> T _Q\\
 @V \alpha _Q VV @V \beta_Q VV\\ 
R _{\mathcal D} @> \pi _{\mathcal D }>> T_{\mathcal D}  \\
\endCD
$$
is commutative: $\beta_Q \circ \pi _Q \circ  \rho ^{\univ}  _Q=\beta_Q \circ \rho ^{\modular}  _Q $ and $\pi _{\mathcal D } \circ \alpha_Q  \circ \rho ^{\univ} _Q = \pi _{\mathcal D } \circ  \rho ^{\univ}_{\mathcal D }= \rho ^{\modular} _{\mathcal D } $ have the same trace at $\Fr _v$ for $ v \not \in \Sigma_{\mathcal D_Q}\cup \{ y \} \cup S $ (the standard Hecke operators at $v$). By the Chebotarev density theorem, the two representations have the same trace function on $G_F$, and hence are isomorphic by Lemma \ref{lem-modular32}. \par

TW1 is clear from the definition of $X_\mathcal D$. We check axioms TW2-TW5 of Definition \ref{dfn-TW11}. By Proposition \ref{claim-galdef2}, $o^{\associate}_{\mathcal D}  [ \Z_\ell \times \Delta _v ] $ is regarded as a versal hull
of the unrestricted deformation of $\bar \rho \vert _{G_{F_v}} $ with a fixed determinant (we use the identification given there). We have a natural $o ^{\associate}_{\mathcal D} $-algebra homomorphism 
$$
\tau _Q : o ^{\associate}_{\mathcal D} [ \Delta _Q] =\bigotimes _{ v \in Q }o ^{\associate}_{\mathcal D}  [\Delta _v]
\to R_Q ,
$$
which corresponds to the restriction $\prod _{v \in Q} \rho ^{\univ}_Q \vert _{I_{F_v}}$. Thus $R_Q$ has an $o ^{\associate}_{\mathcal D} [ \Delta _Q ] $-algebra structure, and we view $ M_Q$ as an $R_Q$-module by $f_Q$. We use these structures in TW2. \par
By our normalization in Proposition \ref{claim-galdef2} and Lemma \ref{lem-const12}, 
$$
\pi _Q \circ \tau _Q ( \delta
 ) = V(\delta)  \quad (\delta \in \Delta _Q). 
$$
So $ o ^{\associate}_{\mathcal D} [ \Delta _Q] \overset {\pi_Q \circ\tau _Q} \to
T_Q$ is equal to the homomorphism $  o ^{\associate}_{\mathcal D}[ \Delta _Q] \to T _Q$
obtained by the $ V(\delta) $-operators. It is clear that $R _Q / I _Q R_Q \simeq R $, since $R_Q/ I _QR_Q$ is the maximal quotient of $R_Q$ where $ \rho ^{\univ} _Q$ is unramified at $Q$ by Faltings' theorem \ref{claim-galdef2}. TW3 is verified. \par
We regard $M_Q$ as an $R_Q$-module by $\pi _Q$, and the action of $o ^{\associate}_{\mathcal D}[ \Delta _Q]$ through $R_Q$ is the same as the $C_{F, \ell } \times \Delta _Q$-action geometrically constructed on $M ^{\tw}_Q$ by the isomorphism $M _Q \simeq M ^{\tw}_Q$.
By Proposition \ref{prop-const42}, $ M _Q / I _Q M _Q\simeq M _{\mathcal D}$, and TW4 is satisfied.   \par

Finally we check TW5. By Proposition \ref{prop-const41}, $M_Q $ is a free
$ o ^{\associate}_{\mathcal D}[\Delta_Q]
$-module. The $ o ^{\associate}_{\mathcal D}[\Delta_Q]$-rank of $ M _Q
$ is the
$ o ^{\associate}_{\mathcal D}$-rank of $ M _Q / I_Q M _Q $, and which is non-zero and independent of $Q$ by Proposition \ref{prop-const42}.  Thus TW5 is satisfied for
$(R_Q, M _Q )_{ Q
\in
\mathcal X_{\mathcal D} } $.

\begin{rem}\label{rem-const41}
Similarly, one has a Taylor-Wiles system $ (R_{Q, \chi} , \ \tilde M_{Q } )_{Q \in \mathcal X_{\mathcal  D} } $ for the deformation ring $ R_{\mathcal D, \chi}$ with determinant fixed to $\chi$. Here $\tilde M_{Q }$ is the image of $M_Q $ in $H ^ {q_{\bar \rho} } ( \tilde S_{K'_Q} ,\  \bar {\mathcal F} ^{\tw} _{(k,w) } ) _{\tilde m_{\Sigma , \bar \rho }  } \otimes \chi^{\frac w 2}_{\cycle, \ell} $ in the notation of \S\ref{subsec-const4}.
\end{rem}

\section{The minimal case} \label{sec-minimal}
The main purpose of this section is to prove the following theorem.
\begin{thm}\label{thm-minimal1}
 Let $F$ be a totally real number field, $ \bar \rho: G_F \to \GL_2( k_\lambda) $ an absolutely irreducible mod $\ell$-representation. We take a minimal deformation condition $\mathcal  D$, and assume the
following conditions.
\begin{enumerate}
\item  $\ell \geq 3$, and $\bar \rho
\vert _{ F(\zeta_\ell) }$ is absolutely irreducible. When $\ell= 5$, the following case is excluded: the projective image $\bar G$ of $ \bar \rho $ is isomorphic to $\PGL_2 (\mathbb F_5)$, and the mod $\ell$-cyclotomic character $\bar \chi _{\cycle} $ factors through $G_F \to  \bar G ^{\ab} \simeq \mathbb Z/ 2$ (in particular $[F(\zeta_5): F ] = 2$). 

\item  For $ v\vert
\ell$, the deformation condition for 
$\bar
\rho
\vert _{G_{F_v} } $ is either nearly ordinary or flat (cf. \ref{subsec-galdef6}). When the condition is nearly ordinary (resp. flat) at $v$, we assume that $\bar \rho \vert _{G_{F_v}}$ is $G_{F_v}$-distinguished (resp. $F_v$ is absolutely unramified).
\item There is a minimal modular lifting $\pi$ of $ \bar \rho $ over $o_{\mathcal D } $ as in Definition \ref{dfn-nearlyordinary32}.
 \par 

\item Hypothesis \ref{hyp-nearlyordinary21} is satisfied.
\end{enumerate}
  Then the universal deformation ring $R_{\mathcal  D} $ of $ \bar
\rho$ with the deformation condition $\mathcal  D$ is a complete intersection of relative dimension zero over $o^{\associate}_{\mathcal D}=o_{\mathcal D}[ C_{F, \ell} ]$, and $M_{\mathcal D} $ is a free $R_{\mathcal  D}  $-module. In particular $R_\mathcal D$ is isomorphic to the Hecke algebra $T_{\mathcal  D} $ attached to
$\bar \rho$. Here $C_{F, \ell}  $ is the $\ell$-part of the class group of $F$. 
\end{thm}

 Note that there is an exceptional case which does not happen in \cite{W2}.

\subsection{Preliminary on finite subgroups of $\GL _2 (\overline{\mathbb F}_\ell )$}
 
Let $\ell$ be a prime. Consider an absolutely irreducible representation $\rho: G \to \GL _2 (\overline{\mathbb F} _\ell ) $ of a finite group $G$. By the classification of finite subgroups of $  \GL_2 (\overline{ \mathbb F }_{\ell}) $, the projective image $ \bar G $ of $\rho$ is in the following list:
\begin{enumerate}
\item  $\ell\geq 5$, $\bar G \simeq A_4$. $\bar G ^{\ab} $ is isomorphic to $\mathbb Z / 3  $.
\item  $\ell\geq 5$, $\bar G \simeq S_4$. $\bar G ^{\ab}  $ is isomorphic to $\mathbb Z / 2  $.
\item  $\ell \geq 7$, $\bar G \simeq A_5$. $\bar G ^{\ab} $ is trivial. 
\item (Dihedral case) $\ell \geq 2$, $\bar G \simeq D_{2n}$, where $n$ is prime to $\ell$, $n\geq 2 $, and $D_{2n}$ is a dihedral group of order $2n$. $ \bar G ^{\ab} $ is isomorphic to $\mathbb Z/2  $ (resp. $ \mathbb Z/2 \times \mathbb Z/2 $) if $n$ is odd (resp. even). 
\item ($\PSL_2$-case) $\ell^n \geq 3$, $\bar G \simeq  \PSL _2 (\mathbb F_{\ell ^n}   )   $. $\bar G ^{\ab} $ is trivial except for $\mathbb F_{\ell ^n} = \mathbb F_3$, and is isomorphic to $\mathbb Z / 3  $ in the exceptional case.
\item ($\PGL_2$-case) $\ell ^n \geq 3$, $\bar G \simeq \PGL _2 (\mathbb F_{\ell ^n}   )  $. $\bar G ^{\ab} \simeq \mathbb Z/2 $.
\end{enumerate}
The above 6 cases are exclusive to one another. In the case of (5) (resp. (6)),  $\bar G$ is $\GL_2 (\overline{\mathbb F}_\ell )$-conjugate to the standard embedding of $ \PSL_2 (\mathbb F _{\ell^n} ) $ (resp. $ \PGL_2 (\mathbb F _{\ell^n} ) $) via the embedding $\mathbb F _{\ell ^n } \ \hookrightarrow  \overline{\mathbb F}_\ell $. 

 When $\ell \geq 3$, $\bar G$ is classified into (5) or (6) if and only if $\ell$ divides the order of $\bar G$.

Assume that $\ell\geq 3$, and let $\rho : G \to \GL_2 (\overline {\mathbb F} _\ell) $ be an absolutely irreducible representation. Then the following conditions are equivalent:
\begin{itemize}
\item $\rho$ is monomial.
\item $ \ad^0 \rho $ is absolutely reducible. 
\item The projective image $\bar G$ of $\rho $ is in the dihedral case. 
\end{itemize}
If one of the above conditions holds, then for an index $2$-subgroup $H$ of $G$ such that $ \rho$ is induced from a character $ \chi $ of $H$, $\ad^0 \rho \simeq \nu  \oplus \Ind _H ^ G \chi / \chi ^c $. Here $\nu$ is the character defined as the composition of $G \to G/H \simeq \{ \pm 1\} $, and $\chi ^c $ is the $c$-twist of $\chi$ for the generator $c$ of $G/H$. $ \Ind _H ^ G \chi / \chi ^ c$ is absolutely reducible if and only if $ \bar G \simeq \mathbb Z/ 2 \times \mathbb Z/ 2$.\par
\bigskip
We prepare several group theoretical lemmas which are used in the Chebotarev density argument in \S\ref{subsec-minimal2}
. 
\begin{lem} \label{lem-group11}
Assume that $\ell \geq 3$. Let $
\rho : G \to \GL_2 (\overline {\mathbb F}_\ell)$ be an absolutely irreducible representation
of a finite group $G$, and $H$ a normal subgroup of $G$ such that $G/H$ is cyclic.

 Then there is an element $\sigma \in H$ of order prime to $\ell$ such that 
$ \rho (\sigma ) $ is regular and semi-simple. 
\end{lem}

\begin{proof}[Proof of Lemma \ref{lem-group11}]  Assume that $\rho (h)$ has the same eigenvalues for any $h\in H$. Consider the
$H$-representation $ \ad ^0  \rho \vert_H$. Then the $H$-action factors through an
$\ell$-group, since the all eigenvalues of any element of $H$ are one. So $\ad ^0 \rho \vert _H$ has only
trivial constituents.

 $G/H $ acts on the invariant subspace $(\ad ^0  \rho ) ^H \neq \{ 0\}$. Since $G/H$ is
abelian, the
$G$-action on $(\ad ^0  \rho ) ^H$ is absolutely reducible, and hence $\ad ^0  \rho $ is also absolutely reducible. 
Since $\ell \neq 2$, $ \rho$ is monomial by the absolute reducibility of $\ad ^0 \rho$. Since the projective image $\bar G$ of $G$ has order prime to $\ell$, the $H$-action on $ \ad ^0  \rho $ is semi-simple. $H$ acts trivially on $  \ad ^0  \rho $, and hence $\bar G$ is a quotient of $G/H$. This is impossible since $\bar G$ is dihedral, and $G/H$ is cyclic. \par
So there is an element $h\in H$ such that $ \rho (h) $ is regular and semi-simple. By replacing $h$ by a suitable $\ell ^m$-th power $\sigma$ for some $m\geq 0$, the order of $\sigma $ is prime to $\ell$, and $\sigma$ is regular and semi-simple. 
\end{proof}
\begin{rem} By the proof of \ref{lem-group11}, the consequence remains true under the assumption that $G/H $ is abelian unless $\bar G$ is dihedral and abelian.
\end{rem}
In the dihedral case, we need a slightly stronger result than Lemma \ref{lem-group11}.

\begin{lem}\label{lem-group12}
Assume that $\ell \geq 3$. Let
$
\rho : G \to \GL_2 ( \overline {\mathbb F}_\ell )$ be a monomial representation of a finite group $G$, $H$ a subgroup of $G$ such that $
\rho\vert _H$ is absolutely irreducible. 
Then for any non-zero absolutely irreducible subspace $N $ of
$
\ad ^0 \rho  $, there is an element $\sigma \in H$ of order prime to $\ell$ such that
$\rho (\sigma ) $ is regular and semi-simple, and $N^{<\sigma>}\neq \{ 0 \} $.
\end{lem}
\begin{proof}[Proof of Lemma \ref{lem-group12}]
This is essentially \cite{W2}, lemma 1.12.\par
By the assumption, the restriction of $\rho$ to $H$ is absolutely irreducible. Since the $H$-action on $  \ad ^0 \rho $ is absolutely reducible, it is also monomial. We may assume that $H=G$, since it suffices to find $\sigma \in H$ for a non-zero absolutely  irreducible $H$-subrepresentation $N'$ of $N$. \par
First assume that the projective image $\bar G$ of $\rho$ is not abelian. Then $ \ad ^0 \rho$
is isomorphic to $\nu \oplus V$, where $\nu$ is a character of order $2$ such that $\rho $ is induced from a character $\chi$ of $\ker \nu$, and $V$ is irreducible. When $N$ is the subspace which is isomorphic to $V$, we take any element $ \sigma \in G $ of order prime to $\ell$ such that $\nu (\sigma )\neq 1$.  When $N$ is the subspace where $G$ acts by $\nu$, one takes any element $\sigma \in \ker \nu $ of order prime to $\ell$ such that $ \chi / \chi ^{c} (\sigma) \neq 1$. 
 \par
Finally we treat the case when $\bar G  $ is dihedral and abelian, that is, $\bar  G \simeq \mathbb Z/ 2 \times \mathbb Z/2 $. 
$\ad ^0 \rho \simeq \alpha  \oplus \beta \oplus \alpha \cdot  \beta  $ is isomorphic to the direct sum of the non-trivial characters of $\bar  G $. By replacing the role of characters if necessary, we may assume that $N$ is the subspace where $G$ acts by $ \alpha$.  Any element $\sigma $ of order prime to $\ell$ such that $ \alpha(\sigma) = 1$ and $\beta(\sigma )\neq 1$ satisfies the desired condition. 
\end{proof}

\begin{lem} \label{lem-group13} Assume that $\ell \geq 3 $. Let $ \rho: G \to \GL_2 (\overline {\mathbb F}_\ell )$ be a representation of a finite group $G$, $H$ a normal subgroup of $G$ of $\ell$-power index. 
\begin{enumerate}
\item $\rho$ is absolutely irreducible if and only if $\bar \rho \vert_H $ is absolutely irreducible. 
\item Assume that $\rho$ is absolutely irreducible. Then the projective image $\bar H  $ of $\rho \vert_H $ is equal to the projective image $\bar G$ of $ \rho $ except when $\ell = 3 $, $\bar G \simeq \PSL _2 (\mathbb F_3 ) $, and $ \bar  H \simeq \mathbb Z /2 \times   \mathbb Z /2 $.
\end{enumerate}

\end{lem}
\begin{proof}[Proof of Lemma \ref{lem-group13}] For (1), it suffices to show that $\rho $ is absolutely reducible if $\rho \vert _H$ is absolutely reducible. Let $X$ be the set of all constituants of $\rho \vert _H$. The cardinality of $X$ is at most $2$. $g \in G$ acts on $X$ by sending $\chi$ to the $g$-twist $\chi ^g$. Since $G/H$ is an $\ell$-group and $\ell \geq 3$, the action is trivial, and any element $\chi$ in $X$ is the restriction of a character $\tilde \chi$ of $G$. Assume that $\chi$ appears as a subrepresentation of $ \rho \vert _H$. Then $(\rho \otimes \tilde \chi ^{-1})^{H}\neq \{ 0 \} $.  Since $G/H$ is an $\ell$-group, $(\rho \otimes \tilde \chi ^{-1})^G$ is non-zero, and hence $\rho$ is absolutely reducible. \par
For (2), note that $\bar H \neq \bar G  $ implies that the order of ${\bar G }^{\ab}$ is divisible by $\ell$. Then the result follows from the classification of subgroups of $\GL_2 (\overline{\mathbb F }_\ell ) $.
\end{proof}

\begin{prop}
\label{prop-group11}
Let $\rho : G \to \GL_2 (\overline{ \mathbb F }_{\ell}) $ be an absolutely irreducible
representation of a finite group, $\ell \neq 2$,  $\mu : G
\to\overline{ \mathbb F }_{\ell}  ^\times$ a character of even order. Then 
$$
H ^1 (G/(\ker (\ad ^0 \rho \otimes \mu)  \cap \ker \mu),\ \ad ^0 \rho  \otimes \mu )= 0  
$$
except when 
\begin{itemize}
\item $\ell = 5$, the projective image of $\rho $ is isomorphic to $\PGL _2 ( \mathbb F _5 )$, and $\mu$ factors through the character of order $2$ of $\PGL _2 ( \mathbb F _5 )$.
\end{itemize}
\end{prop}
\begin{proof}[Proof of Proposition \ref{prop-group11}]
 Note that $(\ker ( \ad ^0 \rho\otimes \mu  )\cap \ker \mu =( \ker \ad ^0 \rho) \cap \ker \mu  $. 
Define subgroups $H '$ and $H$ of $G$ by $\ker \ad ^ 0  \rho $ and $\ker \ad ^0  \rho \cap \ker \mu$ respectively. By the Hochshild-Serre spectral sequence, 
$$
0 \longrightarrow H ^ 1 (G / H ' , (\ad ^0 \rho  \otimes \mu  ) ^{H ' / H } ) \longrightarrow H ^1 (G / H  , \ad ^0 \rho \otimes \mu  ) \longrightarrow H ^1 ( H' / H  , \ad ^0 \rho \otimes \mu  ) 
$$
is exact. If $ \mu \vert _{H'} $ is non-trivial, then $H ^1 (G / H , \ad ^0 \rho \otimes \mu  ) $ is zero since the order of $H  '/ H \simeq \mu ( H ' )  $ is prime to $\ell$, and $(\ad ^0 \rho \otimes \mu  ) ^{H ' / H  } $ vanishes. \par
So we may assume that $\mu $ is a character of $\bar G =  G/H '$, that is, a character of the projective image of $\rho$.

If the order of $\bar G $ is prime to $\ell$, then the statement of \ref {prop-group11} is clear, and hence we may assume that $\bar G $ contains an element of
order
$\ell$. By the classification, $ \bar G $ is either $ \PSL _2 (\mathbb F_{\ell ^n}   )  $ or $  \PGL _2 (\mathbb F_{\ell ^n}   ) $ for some $n\geq 1$. Since the order of $\mu$ is even, it follows that $\bar G$ is conjugate to the standard embedding of $ \PGL _2 (\mathbb F_{\ell ^n}   )  $, the order of $\mu $ is $2$, $\mu $ is obtained from $( \det ) ^{\frac { \ell ^ n -1 }{2}}$ of $\GL _2 (\mathbb F_{\ell ^n} )$, and $\ad ^0 \rho \otimes \mu $ is isomorphic to $ \Sym ^2 \otimes ( \det ) ^{\frac { \ell ^ n -1 }{2} -1} $ as $\PGL_2 (\mathbb F_\ell )$-representations over $ \overline { \mathbb F } _\ell$.\par
Then we apply the following lemma.
\begin{lem}\label{lem-group14}
Assume that $\ell$ is an odd prime. Then
$$
H ^1 ( \SL _2 (\mathbb F_{\ell ^n}  ) ,   \Sym ^2 )= 0
$$
except for $\mathbb F_{\ell ^n}  = \mathbb F _5$. 
\end{lem}
This follows from \cite{CPS}, p. 185, Table (4.5) if $ \mathbb F_{\ell ^n} \neq \mathbb F_3$.
In the $\mathbb F _3$-case, it is shown in \cite{W2}, p. 478. \par
\bigskip

The vanishing of the cohomology now follows from two exact sequences
$$
0 \longrightarrow H ^1 (\mathbb F_{\ell  ^n } ^\times ,  (\Sym ^2 \otimes ( \det ) ^{\frac { \ell ^ n -1 }{2} -1} ) ^{\SL_2 ( \mathbb F_{\ell ^n} ) })   \longrightarrow H ^1 ( \GL_2 ( \mathbb F_{\ell ^n} ) , \Sym ^2 \otimes ( \det ) ^{\frac { \ell ^ n -1 }{2} -1}  )$$
$$
\longrightarrow H ^1 ( \SL_2 ( \mathbb F_{\ell ^n} )     , \Sym ^2 \otimes ( \det ) ^{\frac { \ell ^ n -1 }{2} -1}  ) ^{\mathbb F_\ell ^\times} $$
$$
0 \longrightarrow H ^1 (   \PGL_2 ( \mathbb F_{\ell ^n} )  ,     (\Sym ^2 \otimes ( \det ) ^{\frac { \ell ^ n -1 }{2} -1} ) ^{ \mathbb F ^\times _{\ell^n }   } )  \longrightarrow H ^1 (\GL_2 ( \mathbb F_{\ell ^n} ) ,  \Sym ^2 \otimes ( \det ) ^{\frac { \ell ^ n -1 }{2} -1}  )  $$
and Lemma \ref{lem-group14}, except for the $\mathbb F_5 $-case.

\end{proof}

As for the application to elliptic curves, the exceptional case in Proposition \ref{prop-group11} does not happen. 
\begin{prop}\label{prop-exceptional}
Assume that $G$ contains an element $c$ of order $2$ such that $ \mu (c) = -1$, and $\rho (G )$ is a subgroup of $\GL _2 (\mathbb F _5)  $. Then the exceptional case in Proposition \ref{prop-group11} does not occur for the triplet $(G, \rho , \mu)$.
\end{prop}
\begin{proof} $\mu$ is identified with $(\det   ) ^ 2$. The equation $\mu (c) = -1 = (\det\ \rho( c)  ) ^2 $ 
is impossible since $\det\  \rho (c )\in \{ \pm 1  \}$. 
\end{proof}
\begin{rem}
\label{rem-minimal1}
It can be shown that $H ^1 (\GL_2 (\mathbb F _{\ell ^n }) , \Sym ^2     \otimes ( \det ) ^{-1})   $ is zero if $\ell \neq 2$. By \cite{CPS}, $ H ^ 1 (\SL _2 (\mathbb F _5), \Sym ^2   )$ is one dimensional. Using this, one shows that $H^1 (\PGL _2 (\mathbb F _5) , \Sym ^2 \otimes   \det )  $ is one dimensional. 
\end{rem}

\subsection{Chebotarev density argument}\label{subsec-minimal2}

We formulate an application of the Chebotarev density theorem following \cite{W2}, proposition 1.11.
\begin{dfn} \label{dfn-chebotarev11}
Let $F$ be a global field, $ M $ a finite
$\mathbb Z_{\ell}[G_F]$-module given by $ \rho _M : G_F \to \Aut _{\mathbb Z_\ell } M $, and $H$ an open normal subgroup of $G_F$.
For a finite Galois extension $F'$ of $F$, let $F ' _{(M , H) }$ be the Galois extension which corresponds to $(\ker
\rho _M \vert _{F ' } ) \cap H $. 
A cohomology class $x$ in $ H ^1 (F , \ M ) $ is exceptional for $(M, H) $ over $F'$ if the restriction $x_{F'} $ of $x$ to $F'$ belongs to $H ^1 ( \Gal (F ' _{(M , H ) }/ F' ) , \ M) = \ker ( H ^1 ( F ' ,\  M ) \to H ^1 ( F ' _{(M , H ) }, \ M ) ) $. 
\end{dfn} 

\begin{prop}
\label{prop-chebotarev11}

Notations are as in Definition \ref{dfn-chebotarev11}. Take a non-zero cohomology class $x$ in $ H ^1 (F , \ M ) $. Assume the following conditions.
\begin{enumerate}
\item $F'$ is a Galois extension of $F$, and $x$ is not exceptional for $(M, H)$ over $F'$ (cf. \ref{dfn-chebotarev11}). 

\item For any irreducible $\mathbb Z_\ell [G_F]$-submodule $N$ of $M$, a non-trivial element $ g_N \in \Gal (F ' _{(M , H ) }/F'  ) $ is given. The order of $g_N $ is prime to $\ell$, $g_N $ belongs to the image of $G_{F'}\cap H$, and 
$ N ^{<g_N >}\neq \{ 0\}
$. 
\end{enumerate}
Then there is a finite place $v$ of $F$ that satisfies:
\begin{itemize} 
\item $v$ is of degree 1 and splits completely in $F'$. 
\item $M$ is unramified at $v$, and the
restriction of 
$x
$ to
$ H^1 ( F_v,\ M)$ is not zero.
\item $\Fr _v $ acts on $M$ by a $G_F$-conjugate of $g_N$ for some non-zero irreducible submodule $N$.  In particular the image of $\Fr _v $ in $ G_F $ belongs to $ G_{F'} \cap H $. 
\end{itemize}
Such a place $v$ exists with a positive density. 
\end{prop}
\begin{proof}[Proof of Proposition \ref{prop-chebotarev11}] $L=  F ' _{(M , H ) }$, $\mathcal G = \Gal (L/ F ' ) $. The image $y$ of $ x _{F'}$ under the restriction map $H ^ 1 (F ' ,\ M ) \to
H ^1 ( L,\ M   )  $ does not vanish by assumption (1) of \ref{prop-chebotarev11}.  We regard $y$ as a non-zero $G_F$-equivariant homomorphism
$c (y) : G_L \to M  $. Let $ L_y$ be the abelian extension corresponding to the kernel of $ c (y)$ over $L$. $\Gal ( L_y / L ) $ with the conjugate action of $G_F$ is identified with $c (y) (G_L )$ as $\mathbb Z_\ell [G_F]$-modules.

For a non-zero irreducible $\mathbb Z_\ell [G_F]$-submodule $N$ in $c (y) (G_L ) 
$, take the element $ g_N$ in $\mathcal G$ of order $\alpha $ that is given for $N$ in \ref{prop-chebotarev11} (2), and an element $\tilde \sigma $ in $\Gal ( L_y /F ') $ which lifts $g _N$ and has the same order as $g_N$ (this is possible since  $\alpha$ is prime to $[ L_y : L ] $). 
$g_N \in H $, and $N ^{<g_N>} \neq \{ 0 \} $. 
Choose a non-zero element $\tilde \tau $ in $\Gal ( L_y /L)$ that belongs to $N  ^{<g_N>} $ by the identification $\Gal ( L_y /L)\simeq c (y) (G_L )  $.  \par

Since $\tilde \tau $ belongs to $N ^{<g_N>}$, $\tilde \tau $ is fixed by $\tilde \sigma $. This implies that $\tilde \tau$ commutes
with $\tilde \sigma$. 
We define a conjugacy class $C$ in $ \Gal ( L_y /F )$ as the class of $  \tilde \sigma \tilde \tau$. $C$ is contained in $\Gal ( L_y /F ' )$ since $F'$ is Galois over $F$.\par
Take a finite place $v$ of $F$ of degree one such that $M$ is unramified at $v$, and $ \Fr _v \in C$ by the Chebotarev density theorem. $v$ splits completely in $F'$ since the image of $\Fr _v $ in $\Gal (F' / F ) $ is trivial. 
$c (y)( \Fr _v ^\alpha )$ is in the $\Gal ( L_y /F ' )$-orbit of $\tilde \tau ^\alpha $, which is non-zero by the construction, and hence the restriction of $x$ to $ H ^1 ( F_v , M ) $ is non-zero. 
\end{proof}
\begin{cor} \label{cor-chebotarev11} Notations are as in Definition \ref{dfn-chebotarev11}. Assume moreover that $M$ is a $k_\lambda[G_F]$-module for some finite field $k_\lambda$, and the following conditions.
\begin{enumerate}
\item $F'$ is a Galois extension of $F$, and $x$ is not exceptional for $(M, H)$ over $F'$ (cf. \ref{dfn-chebotarev11}). 

\item For any irreducible $k_\lambda [G_F]$-submodule $N$ of $M$, a non-trivial element $ g_N \in \Gal (F ' _{(M , H ) }/F'  ) $ is given. The order of $g_N $ is prime to $\ell$, $g_N$ belongs to the image of $G_{F'}\cap H $, and 
$ N ^{<g_N >}\neq \{ 0\}
$. 
\end{enumerate}
Then there is a finite place $v$ of $F$ that satisfies:
\begin{itemize} 
\item $v$ is of degree 1 and splits completely in $F'$. 
\item $M$ is unramified at $v$, and the
restriction of 
$x
$ to
$ H^1 ( F_v,\ M)$ is not zero.
\item $\Fr _v $ acts on $M$ by a $G_F$-conjugate of $g_N$ for some non-zero irreducible submodule $N$.  In particular the image of $\Fr _v $ in $ G_F $ belongs to $  G_{F'} \cap H$. 
\end{itemize}
\end{cor}

\begin{proof}[Proof of corollary \ref{cor-chebotarev11}] We construct $g_N$ as in Proposition \ref{prop-chebotarev11} for any non-zero irreducible $\mathbb F _\ell [G_F]  $-submodule $N$ of $M$. Let $W $ be the image of $ N \otimes _{\mathbb F_\ell} k _\lambda $ in $M$, which is non-zero since it contains $N$. We take a non-zero $ k _\lambda [G_F]  $-irreducible submodule $N' $ in $W$, and the element $g_{N ' } $ given for this $N'$ in \ref{cor-chebotarev11} (2). 
Since the order of $g_{N ' } $ is prime to $\ell$, $( N \otimes _{\mathbb F_\ell} k _\lambda )  ^{<g_{N'} >}\to{ W }  ^{<g_{N'} >} $ is surjective. $  {N ' }  ^{<g_{N'} > }\subset W   ^{<g_{N'} >}$, and hence $( N \otimes _{\mathbb F_\ell }k _\lambda )  ^{<g_{N'} >} = N ^{<g_{N'} > } \otimes _{\mathbb F_\ell} k _\lambda$ is non-zero. By taking $g_{N' } $ as $g_N$, Proposition \ref{prop-chebotarev11}, (2) is satisfied, and Proposition \ref{prop-chebotarev11} is applied. 
\end{proof}

\subsection{Proof of Theorem \ref{thm-minimal1}}\label{subsec-minimal1}
We prove Theorem \ref{thm-minimal1}. 
We use the same notation as in \S \ref{subsec-const5}. We have already constructed a Taylor-Wiles system $
\{ R _Q , \ M_Q \}_{Q \in \mathcal X_{\mathcal  D}}  $ for
$(R_{\mathcal  D  } , M_{\mathcal D} )$ over $
o  ^{\associate}_{\mathcal D}$. \par
Recall that 
$$
\Hom _{k_\lambda} 	(m_{ R_Q}/(m_{ R_Q}^2, m_{o^{\associate}_{\mathcal D}} ), \ k_\lambda) \ = \ 
	H_{\mathcal  D_Q}^1(F,\ \ad^0 \bar \rho),
$$
for the universal deformation ring $R_{ \mathcal  D _Q } $, and $R_{ \mathcal  D _Q } $ is
generated by $\dim_{k_\lambda} H_{\mathcal  D_Q}^1(F,\ \ad^0 \bar \rho)$ many elements over $
 o  _{E_\lambda }$. Here
$$
H_{\mathcal  D_Q}^1(F,\ \ad^0 \bar \rho) =  \ker ( H ^ 1 (F ,\
\ad^0 
\bar  \rho ) \longrightarrow \bigoplus _{v \in \vert F\vert _f }  H ^ 1 ( F_v ,\ \ad ^0 \bar \rho   ) /L'_v ) 
$$
by Proposition \ref{prop-galdef92}, with the subspace $L'_v = H ^1 _{\deform _{\mathcal  D_Q} (v) } (F_v,\ \ad ^0
\bar
\rho )
$ for any  finite place $v$. \par  
For the Tate dual $ \ad ^0 \bar \rho ^* =  \ad ^0 \bar \rho  (1) $ of
$ \ad ^0 \bar
\rho$, the dual Selmer groups are defined by
$$
H_{{ \mathcal  D   } ^*  } ^1(F, \ \ad^0 \bar
\rho  (1)  )= \ker ( H ^ 1 (F,\
\ad ^0 
\bar  \rho (1) ) \longrightarrow \bigoplus _{v \in \vert F\vert _f } H ^ 1 ( F_v ,\ \ad ^0 \bar \rho (1)   ) /L_v ^* )
$$
$$ 
H_{{ \mathcal  D  _Q } ^*  } ^1(F, \ \ad^0 \bar
\rho  (1)  )= \ker ( H ^ 1 ( F, \
\ad ^0 
\bar  \rho (1) ) \longrightarrow \bigoplus _{v \in \vert F\vert _f } H ^ 1 ( F_v ,\ \ad ^0 \bar \rho (1)   ) /{L'}_v ^* )
$$
Here $ L ^* _v $ (resp. $ {L'} _v^ * $) is the annihilator of $H ^ 1 _{\deform _{\mathcal  D '}(v) }
(F_v,\
\ad ^0 \bar
\rho)
$ (resp. $H ^1 _{\deform _{\mathcal  D_Q} (v) } (F_v,\ \ad ^0 \bar
\rho )$) by the local duality. \par

To prove Theorem \ref{thm-minimal1}, it is sufficient to find a non-empty subset
$
\mathcal X
\subset
\mathcal X  _{\mathcal  D }
$ which satisfies the assumptions of the complete intersection-freeness criterion (Theorem \ref{thm-TW21}) under the minimality of $\mathcal  D $. For this, it suffices to show the following proposition. 

\begin{prop} \label{prop-minimal10}
Assume that $\mathcal D$ is minimal. For any integer
$n \geq 1 $, there is an element $ Q \in  \mathcal X  _{ \mathcal  D  } $ such that 
\begin{enumerate}
\item $ q_v
\equiv 1
\mod
\ell ^n $ for $v \in Q$,
\item
$\dim _{k_\lambda}  H_{\mathcal  D _Q}^1 ( F,\ \ad^0 \bar \rho) \leq \sharp Q 
 $, 
\item $ \sharp Q =  \dim_{k_\lambda}  H_{{ \mathcal  D } ^*} ^1(F, \ \ad^0 \bar \rho (1) )$
\end{enumerate}
hold. 
\end{prop}

\begin{lem}\label{lem-minimal11} Assume that $\mathcal  D$ is minimal. 
\begin{enumerate}
 \item The inequality
$$
\dim _{k_\lambda}  H_{\mathcal  D  } ^1 ( F ,\  \ad^0 \bar \rho ) \leq \dim _{k_\lambda} 
H_{{\mathcal  D } ^*}^1 ( F,\ 
\ad^0
\bar
\rho  (1) ) 
$$ 
holds. 
\item For any element $Q $ of $\mathcal X _{\mathcal D}$, 
$$
\dim _{k_\lambda}  H_{\mathcal  D _Q   } ^1 ( F ,\  \ad^0 \bar \rho ) \leq \dim _{k_\lambda } H_{\mathcal  D _Q^*}^1 (F , \ \ad^0 \bar \rho  (1) )  + \sharp Q.
$$
\end{enumerate}
\end{lem}

\begin{proof}[Proof of Lemma \ref{lem-minimal11}] We prove (1). By the formula in \cite{W2}, proposition 1.6, 
$$ 
\dim _{k_\lambda}  H ^1_{\mathcal  D } (  F ,\ \ad ^0\bar \rho ) - \dim _{k_\lambda}  H
^1_{\mathcal  D ^*} ( F ,\ \ad^0 \bar \rho (1)) 
$$
$$
= \sum _{ v \in \Sigma_\mathcal  D  } h _v  +
\sum _{ v
\in I_F } h _v  +
\dim _{k_\lambda}  H ^0( F,\ \ad ^0 \bar \rho ) - \dim _{k_\lambda}  H ^0 ( F,\ \ad ^0
\bar
\rho (1)) 
$$
holds. Here 
$$
 h _v = \dim _{k_\lambda} H ^0 ( F_v,\ \ad ^0 \bar \rho (1) ) -
\dim _{k _\lambda} H ^1 ( F_v ,\ \ad ^0 \bar \rho) /  H ^1 _{\deform (v) } ( F_v , \
\ad ^0
\bar
\rho)
$$
for
$v \in \Sigma _{\mathcal  D }$, and 
$$
h _v =
\dim _{k _\lambda} H ^0 ( F_v , \ad ^0 \bar \rho (1) )
$$
for $ v \in I_F $. Since
$\bar \rho \vert _{ G _{ F(\zeta_ \ell
 )} }$ is absolutely irreducible, 
$$
\dim _{k_\lambda}  H ^0( F,\ \ad ^0 \bar \rho ) =\dim _{k_\lambda} H ^0 ( F,\
\ad ^0 \bar \rho(1))= 0 . 
$$
Since $\bar \rho$ is odd for any complex conjugation, 
$$
\sum _{ v \in I_F } h _v =2 [ F: \Q] . 
$$
If $ v \in \Sigma_{\mathcal  D },\ v \nmid \ell $, then $ h_v=0$.  In this case, the
local condition is finite, and this follows from Proposition \ref{prop-galdef41}, since $ \dim
_{k_\lambda} H_f^1 ( F_v ,\ \ad ^0 \bar \rho ) = \dim _{k_\lambda } H ^ 0 ( F_v, \ad ^0
\bar \rho  )$.  When
$ v
\vert
\ell$, we show $h_v \leq  -2 [F_v : \mathbb Q_\ell]$. In the nearly ordinary and nearly
ordinary finite case, this follows from Theorem \ref{thm-galdef71}, and Tate's local Euler
characteristic formula. In the flat case we have 
$$
\dim _{k_\lambda}  H ^1 _{\bold {fl} }  ( F _v ,\
\ad ^0 \bar \rho\vert _{F_v} ) = \dim _{k_\lambda}  H ^0 ( F_v ,\ \ad ^0 \bar \rho\vert_{F_v}  ) + [ F_v:
\mathbb Q _\ell]
$$
by Theorem \ref{thm-galdef72}, and the claim follows again by Tate's local Euler
characteristic formula. 

 So the total sum can not be strictly positive.\par

 For (2), by \cite{W2}, proposition 1.6 and Lemma \ref{lem-minimal11}
$$
\dim _{k_\lambda}  H_{\mathcal  D_Q } ^1 (F ,
\
\ad^0
\bar
\rho ) /H_{\mathcal  D _Q^*}^1 (F , \ \ad^0 \bar \rho  (1) )  = \dim _{k_\lambda}  H_{\mathcal  D} ^1
(F ,
\
\ad^0
\bar
\rho ) /H_{\mathcal  D ^*}^1 (F , \ \ad^0 \bar \rho  (1) ) 
$$
$$
+ \sum_{v\in Q}\dim _{k _\lambda} H ^0 ( F_v  ,\  \ad^0
\bar
\rho (1) ) \leq \sharp Q
$$
holds. Here we have used (1), and the fact that 
$ \dim _k H^0(F_v ,\ \ad^0 \bar \rho  (1) ) = 1$ for each $ v \in Q $, since $ q_v \equiv 1 \mod \ell$, and $
\bar
\rho ( \Fr _v ) $ is a regular semi-simple element in $\GL_2 (k _\lambda)  $. 
\end{proof}
We prove Proposition \ref{prop-minimal10}.
By Lemma \ref{lem-minimal11}, (2), an element $Q $ in $\mathcal X _{\mathcal D} $ satisfies the inequality 
$$
 \dim _k  H_{\mathcal  D_Q}^1 (F ,
\
\ad^0
\bar
\rho )  \leq  \# Q
\leqno{(\ast)}
$$

if 
$$
\dim _{k_\lambda}  H_{\mathcal  D _Q ^{*}}^1(F ,
\
\ad^0
\bar
\rho (1) )= 0 . 
\leqno{(\ast \ast)}
$$

For any integer $n\geq 1$, it is clear that a finite set $Q$ given by the following lemma satisfies conclusions (1)-(3) of Proposition \ref{prop-minimal10} since ($\ast$), and hence ($\ast \ast$), is satisfied for this $Q$ with equality $  \# Q = \dim _{k_\lambda }H ^1 _{ \mathcal  D ^* } ( F,\ \ad ^0 \bar
\rho (1) )$.

\begin{lem}
\label{lem-minimal12}(cf. \cite{W2}, 3.8)
For a prime $\ell \geq 3$, let $ F$ be a
number field such that
$[F(\zeta_\ell):F] $ is even. Assume that $
\bar
\rho
\vert _{F(\zeta_\ell)}$ is absolutely irreducible.  When $\ell= 5$, the following case is excluded: the projective image $\bar G$ of $ \bar \rho $ is isomorphic to $\PGL_2 (\mathbb F_5)$, and the mod $\ell$-cyclotomic character $\bar \chi _{\cycle} $ factors through $G_F \to  \bar G ^{\ab} \simeq \mathbb Z/2  $.

Then for any
$n
\geq 1
$, there is an element
$Q
\in  \mathcal X  _{ \mathcal  D }$ which satisfies the following conditions.
\begin{enumerate}
\item  For $ v \in Q $, $v$ is of degree one, and $ q_v 
\equiv  1 \mod \ell ^n    $ holds.

\item The restriction map 
$$
H^1_{{\mathcal  D} ^*}(F ,\ \ad^0 \bar \rho (1)) \longrightarrow \bigoplus_{v \in Q} H _f ^1(F_v ,\ 
\ad^0 \bar \rho (1))
$$ 
is bijective. In particular the kernel $H^1_{\mathcal  D_Q ^{*}}(F ,\ \ad^0
\bar \rho (1)) $ vanishes.
\end{enumerate}
\end{lem}
\begin{proof}[Proof of Lemma \ref{lem-minimal12}]
Fix an integer $n \geq 1$. We apply Proposition \ref{prop-group11}, and Corollary \ref{cor-chebotarev11}.
We take $ \bar \chi _{\cycle} $ as $\mu$.  Let
$F$ be the number field in consideration, $F ' $ the $n$-th layer of the cyclotomic $\Z_\ell$-extension of $F$, $M = \ad ^0 \bar \rho ( 1) $, $H =\ker \mu$. Since $F'/F$ is an $\ell$-extension, $\mu \vert_{F'}$ has the same order as $\mu$. \par
 Let $x$ be a non-zero cohomology class  in $H ^ 1_{{\mathcal  D  }^*} (F,\ M )$, which is not exceptional for $ (M, H )$ over $F$. 
First we show that $x$ is not exceptional for $(M, H) $ over $F'$.  Since $\Gal (F' / F)$ is a cyclic group of $\ell $-power order, $\rho \vert _{G_{F'} } $ remains absolutely irreducible by Lemma \ref{lem-group13}, (1), and hence $ H ^ 0 ( F',
\ M ) = \{ 0\} 
$, and the restriction $x_{F'} $ of $x$ is non-zero. If $x_F $ is exceptional, by Proposition \ref {prop-group11}, $\ell = 5$ and the projective image of $\bar \rho \vert _{G_{F'}} $ is isomorphic to $\PGL_2 (\mathbb F _5 )$ and $\mu \vert _H$ factors through $ \PGL_2 (\mathbb F _5 )^{\ab}$. By Lemma \ref{lem-group13}, (2), the projective image of $ \rho $ is isomorphic to $ \PGL_2 (\mathbb F _5 )$, and hence $x$ is exceptional over $F$. \par

So the condition (1) of Corollary \ref{cor-chebotarev11} is satisfied. For the condition (2)
of Corollary \ref{cor-chebotarev11}, for any irreducible subspace $N$, we take $g_N$ as the image of an element $\sigma$ given by Lemma \ref{lem-group11} in the non-dihedral case ($ F'/F$ is a cyclic extension), and by Lemma \ref{lem-group12} in
the dihedral case (the restriction of $\bar
\rho
$ to $F'(\zeta_\ell)$ is absolutely irreducible by the assumption on $\bar \rho$ and Lemma \ref{lem-group13}, (1)). 

 Thus Corollary \ref{cor-chebotarev11} is applied, and we have a degree 1 place
$v$ of $F$ such that  $v$ splits completely in $F'$, the image of $\Fr _v$ in $G_F$ belongs to  $ G_{F'} \cap H  $, the restriction of $x$ to $ H ^1 _f (F_v , M ) $ is non-zero, and $\ad ^0 \bar \rho (\Fr _v) $, and hence $\bar \rho(\Fr _v )$, are regular and semi-simple. Since $v$ splits completely in $F'$ and $\mu ( \Fr_v)= 1 $, $ q_v
\equiv 1 \mod \ell ^ n $.

$\dim _{k_\lambda
} H ^ 1 _f (F_v,
\ M ) = 1$ because $\bar \rho (\Fr_v )$ is a regular semi-simple element in $\GL_2 (k_\lambda) $. Then the image of $x$ under the restriction map 
$$
H ^ 1_{{\mathcal  D  }^*} (F,\ M ) \longrightarrow H ^1_f ( F_v ,\ M )  
$$
spans $H ^1_f ( F_v ,\ M )$. Continuing
successively by taking a non-zero element in the kernel which is not exceptional for $(M, H )$, we have a finite set $Q$ which belongs to $ \mathcal X _{\mathcal D}$, such that 
$$
H ^ 1 _{\mathcal D ^*} (F, M ) \longrightarrow \bigoplus  _{v \in Q} H ^1 _f ( F_v ,\ M )  
$$
is surjective, and the kernel consists of classes which are exceptional for $ (M, H ) $ over $F$. We have excluded the exceptional case, so we finish the proof. 

\end{proof}
\section{Congruence modules} \label{sec-congruent}
We calculate the cohomological congruence modules made from modular varieties. 
For elliptic modular curves and constant sheaves, the idea is due to Ribet
\cite{Ri1}. In
\cite{W2},
\S2 this is generally discussed, to reduce the estimate of the Selmer group to the
minimal case. We proceed in the same way as in \cite{W2}, assuming cohomological universal injectivity in the Shimura curve case.
\medskip

\subsection{Dual construction}\label{subsec-congruent1}
For an absolutely irreducible representation $\bar \rho $ and a deformation type $\mathcal D$, we define an $R_{\mathcal D} $-module $ \hat M _{\mathcal D} = M_{\mathcal D}\oplus M^{\op}_{\mathcal D}$ which contains $M_{\mathcal D} $ as a direct summand. \par

In general, for a Galois representation $ \rho : G_F \to \GL_2 (A  ) $, we define $\rho ^{\op} $ as
$$
\rho ^{\op} =\rho ^{\vee} (-1) =  ( \det \rho ) ^{-1} \cdot \rho  (-1) .
$$
It is easily checked that $( \rho ^{\op})^{\op}= \rho  $, and for a place $v\nmid \ell$ where $\rho $ is unramified, 
$$
\begin{cases}
\tr \rho^{\op} ( \Fr _v) =q_v  \cdot ( \det \rho  (\Fr _v ) ) ^{-1} \cdot  \tr \rho  (\Fr_v)    
\\
\det \rho ^{\op} (\Fr _v  ) =q^2_v \cdot  (\det \rho (\Fr_v) ) ^{-1}  
\end{cases}
$$ hold. \par
We fix a discrete infinity type $(k, w )$. For a deformation type $\mathcal D $ of $\bar \rho$, $ M_{\mathcal D}  $ is regarded as a submodule of $ H ^{q_{\bar \rho} } _{\stack}(S_K ,    \bar {\mathcal F} _{( k, w) }  ) $ as in \S\ref{subsec-modular4}.
$$
M ^{\op} _{\mathcal D}  =\Hom_{o_{\mathcal D} }  (M _{\mathcal D} , o_{\mathcal D} )   
$$ 
is regarded as a submodule of $ H ^{q_{\bar \rho} }_{\stack} (S_K ,    \bar {\mathcal F} _{( k, -w) }  )(q_{\bar \rho }) $ by Poincar\'e duality discussed in \S\ref{subsec-shimura3}. \par
We define the deformation type $\mathcal D ^{\op} $ for $\bar \rho ^{\op} $. The deformation functions and the coefficient rings are the same for $\mathcal D $ and $\mathcal D ^{\op}$. At $v \vert \ell$, the nearly ordinary type (resp. flat twist type) of $\mathcal D ^{\op} $ is $(\det \bar \rho )^{-1} _{\lift} \vert _{F_v} \cdot   \kappa _{\mathcal D , v}  (-1)   $, where $   \kappa _{\mathcal D , v}$ is the nearly ordinary type (resp. flat twist type) of $\mathcal D $. \par
The $K$-type for $\bar \rho ^{\op}$ and $\mathcal D ^{\op} $ at a finite place $v$ is the same as that of $\bar \rho$ except for the $K$-character. The $K$-character is defined as  
$$
\nu _{\deform _{\mathcal D ^{\op} } }(\bar \rho \vert _{F_v} )  = \nu _{\deform _{\mathcal D } }(\bar \rho \vert _{F_v} ) ^{-1} . 
$$
Thus we have the $\ell$-adic Hecke algebra $T_{\mathcal D ^{\op}}$ and $ T_{\mathcal D ^{\op} }$-module $ M _{\mathcal D ^{\op } } $ for $\bar \rho ^{\op} $. 
\begin{prop}\label{prop-congruent11}
\begin{enumerate} 
\item $M ^{\op } _{\mathcal D} $ is canonically isomorphic to the module associated to $\bar \rho ^{\op} $ of deformation type $\mathcal D ^{\op} $. 
\item There is an $o_{\mathcal D}$-algebra isomorphism $T_{\mathcal D } \simeq T_{\mathcal D ^{\op} } $ between the $\ell$-adic Hecke algebras attached to $\bar \rho $ and $ \bar \rho ^{\op}$, which maps $ T _v $ to $T_{v, v} ^{-1}\cdot  T_v $, $T _{v, v} $ to $T _{v, v} ^{-1} $ for $v \not \in \Sigma_{\mathcal D} $. 
\end{enumerate}
\end{prop}
\begin{proof}[Proof of Proposition \ref{prop-congruent11}]
The compatibility of standard Hecke operators $T_v$ and $T_{v,v}$ follows from Proposition \ref{prop-shimura32}. 
One also checks the compatibility of
for $\tilde U(p_v)$ and $\tilde U(p_v, p_v)$-operators under Poincar'e duality. 
\end{proof}

By Proposition \ref{prop-congruent11}, $\hat M_{\mathcal D} = M _{\mathcal D } \oplus M ^{\op} _{\mathcal D }$ is regarded as a $T_{\mathcal D }$-module, and 
self-dual as a $T_{\mathcal D} $-module. Thus $\hat M _{\mathcal D} $ has a pairing
$$
\langle \ , \
\rangle _{\mathcal  D} : \hat M_{\mathcal  D} \times \hat M_{\mathcal  D} \longrightarrow o_{\mathcal D} 
$$
which induces $ \hat M_{\mathcal  D} \simeq \Hom _{o_{\mathcal D} } ( \hat M _{\mathcal  D} ,o_{\mathcal D}   ) $
as $T_{\mathcal  D}$-modules. 

 \subsection{Congruence modules (I)}\label{subsec-congruent2}
For a division quaternion algebra $D$ as in \S\ref{subsec-shimura3}, we assume that $D$ is split at all $ v\vert \ell$. 
For an $F$-factorizable compact open subgroup 
$ K$ of $G_D ( \mathbb A _{\mathbb Q, f} ) $ and a finite set $\Sigma $ of finite places which contains $\Sigma _K $, take
$ v
\not
\in
\Sigma$ such that 
$$
K_v =\GL_2 (o_{F_v} ) ^{\ell} = \text{ the inverse image of }\Delta _v \text{ by } \GL_2 (o_{F_v} )
\overset{
\det }
\to
o_{F_v }^\times
\to k(v) ^\times 
.
$$
Consider two degeneracy maps 
$\pr _1 ,
\
\pr _2 : S_ { K
\cap K_0 (v)}
\to S_K
$ defined in \S\ref{subsec-coh3}, and 
$$
\varphi_{(k, w)}  : H ^ { q_D } ( S_K ,\ \bar {\mathcal  F} _{(k, w)} ) ^{\oplus 2} \overset { \pr _1 ^*
+ 
\pr _2 ^*}\longrightarrow
 H ^ { q_D  } ( S_{K\cap K_0 ( v)}  ,\ \bar { \mathcal  F} _{(k, w) } ) . 
$$ 
is the map considered there. 
When $q_D = 0$, $\varphi_{(k, w) }$ is universally injective, that is, injective and the image is an $o_{E_\lambda}$-direct summand up to the modules of residual type (with respect to the action of the convolution algebra $H _\Sigma $). If $q_D = 1$,
$\varphi_{(k,w) }$ satisfies the same properties under Hypothesis \ref{hyp-coh31}. 

By taking the dual of $\varphi_{( k, -w)} $ induced by Poincar\'e duality, with appropriate
Tate twists, we have
$$
\tilde \varphi_{(k, w) }  :  H ^ { q_D } ( S_{K\cap K_0 ( v)}  ,\  \bar {\mathcal  F} _{(k, w)} ) \overset
{ \varphi _{(k, -w)}^\vee ( -q_D ) }\longrightarrow
 H ^ { q_D } ( S_K ,\
\bar {\mathcal  F} _{(k, w)} ) ^{\oplus 2} .
$$
By the construction of duality pairing, $\tilde \varphi_{(k, w) }   $ is the map $( \pr _1 ) _! + (\pr _2 ) _!$ obtained from the trace map. 
The following lemma is originally due to Ribet \cite{Ri2}.
\begin{lem}\label{lem-congruent21}
Assume that the $v$-type is $( (2, \ldots, 2) , 0 )  $ if $v \vert \ell$. Then 
$$
 \tilde \varphi
_{( k, w) } \circ \varphi _{(k, w) } (H ^ { q_D  } ( S_K ,\
\bar {\mathcal  F} _{(k, w) } ) ^{\oplus 2}) = ( U (p_v) ^ 2 - T _{v,v} ) H ^ { q_D } ( S_K ,\
\bar {\mathcal  F} _{(k, w)} ) ^{\oplus 2} 
$$
holds up to the modules of residual type (Hypothesis \ref{hyp-coh31} is unnecessary for this statement). 
\end{lem}
\begin{proof}[Proof of Lemma \ref{lem-congruent21}] First we calculate $\tilde \varphi
_{ (k, w) } \circ \varphi _{(k, w)} $ explicitly. 
$$
\tilde \varphi
_{ (k, w)} \circ \varphi _{(k, w)}= 
\begin{pmatrix}
 1+ q_v  &  T_v \\ 
T_{v,v} ^{-1} \cdot T_v &1+ q_v \\ 
\end{pmatrix} 
\leqno{(\ast)}
$$
Here the multiplication is on the left. ($\ast$) is seen as follows. By Proposition \ref{prop-coh21}, the localization of
$H ^ { q_D  } ( S_K ,\
\bar {\mathcal  F} _{(k, w)} ) ^{\oplus 2}$ at some maximal ideal of $ H_K $ which is not of residual type is $o_{E_\lambda}$-free, so it suffices to check ($\ast$) after tensoring with $ E_\lambda$. \par

Take an element $ \begin{pmatrix} f_1\\ f_2 \end{pmatrix}
 \in  H ^ { q_D} ( S_K ,\
\bar {\mathcal  F} _{(k, w) , E_\lambda } ) ^{\oplus 2} $.
Then 
$$
  \tilde \varphi
_{ (k, w)} \circ \varphi _{(k, w)} \begin{pmatrix} f_1\\ f_2 \end{pmatrix}
$$
$$
=( \pr _1) _!   \pr _1 ^*  f_1 + ( \pr _2) _! \pr _1 ^*  f_1 + ( \pr _1) _!   \pr _2 ^*  f_2 + ( \pr _2) _! \pr _2 ^*  f_2
$$
$$
= ( q_v + 1) f_1 + T _{v, v} ^{-1} f_1 + T _v f_2 + ( q_v + 1 ) f_2 .
$$
Here we have used that $(\pr_1  ) _! \pr_2 ^{*}= T_v $, and $(\pr_2  ) _! \pr_1 ^{*}$ is 
the dual correspondence $T_{v, v} ^{-1}\cdot T_v $. Thus ($\ast$) is shown. \par

Similarly, we have
$$
U_v\circ  \varphi _{(k, w)}  =  \varphi _{(k, w)}   \circ 
\begin{pmatrix} 
 0 & -T_{v, v}   \\ 
q_v &T_v  \\ 
\end{pmatrix} , 
$$
$$
(U_v ^2 - T_{v,v} )\circ  \varphi _{(k, w)} =
\begin{pmatrix}
 - T_{v,v} & 0 \\
 T_v  & -T_{v,v} 
\end{pmatrix}
 \circ ( \tilde \varphi
_{( k, w )} \circ \varphi _{(k, w)})  .
$$
\end{proof}

Consider an absolutely irreducible representation $\bar \rho$ with a deformation type $\mathcal  D$. We assume that the deformation function satisfies $ \deform _{\mathcal  D}(v) = \bold{n.o.f.}$. Let $
\mathcal  D ' $ be the deformation type of $\bar \rho$ which has the same data as $\mathcal D$ except the deformation function $\deform_{\mathcal D'} $.  $ \deform _{\mathcal  D ' } ( v) = \bold
{n.o.} $, and $\deform _{\mathcal D ' } $ takes the same value as $\deform _{\mathcal  D}$ at the other places.

We have defined $ T_{\mathcal  D}$ and $T_{\mathcal  D' }$-modules $M_{\mathcal  D} $ and $M_{\mathcal  D'}$. 

We show that there is a canonical surjective homomorphism $\beta : T_{\mathcal  D ' } \to T_{\mathcal  D }
$ and
$T_{\mathcal  D ' }
$-homomorphism
$$
\xi : M_{\mathcal  D} \longrightarrow M _{\mathcal  D ' } ,
$$
with the calculation of the congruence modules. We may assume that the
$v$-type satisfies
$ (k_v , w ) = (( 2,
\ldots, 2), w)$ by our assumption \ref{dfn-nearlyordinary32}, and hence $w$ is even. We may twist $\bar {\mathcal  F}_{(k, w) } $ by $\Norm_{D/F} ^{\frac
w 2} 
$, and assume that $ w= 0 $.  \par
We use the same notation as in \S\ref{subsec-modular4}. In particular we make a choice of an auxiliary place $y$ and a set of finite places $S$ to make compact open subgroups in consideration small. \par
For $K = \ker \nu _{\mathcal D_y } \vert _{\tilde K _{\mathcal D_y }  }$, there is a canonical map  
$$
H ^ { q_{\bar \rho}  } _{\stack } ( S_{ K  } ,\ \bar { \mathcal  F} _{(k,w) } )^{\oplus 2}  \longrightarrow H ^ { q_{\bar \rho} }_{\stack} 
( S_{(K _v \cap K _0( m_v) ) \cdot K ^ v } ,\  \bar {\mathcal  F }_{(k,w)} ) 
\leqno{(\dagger_1)}
$$
induced by the degeneracy maps. For $\Sigma ' = \Sigma _{\mathcal D ' _y } \cup S = \Sigma _{\mathcal D _y } \cup S \cup \{ v\} $, let $\tilde m _{ \Sigma  ' , \bar \rho  } $ be the maximal ideal of $H _{K ^{\Sigma ' } } = H ( G _D ( \mathbb A _{\mathbb Q , f} ), K ^{\Sigma '}  )  _{o _{\mathcal D }} $ corresponding to $\bar \rho$. \par

We localize ($\dagger_1$) at $\tilde m_{\Sigma ' , \bar \rho  } $, and take the part where $ \tilde K _{\mathcal D ' _y } $ acts via $ \nu _{\mathcal D'_y  } $. Thus we have an injective $o_{\mathcal D } $-homomorphism
$$
( \tilde M ^ y  _{\mathcal D _ y }  ) ^{\oplus 2} \longhookrightarrow  \tilde M ^ y  _{\mathcal D ' _ y }.
\leqno{(\dagger_2)}
$$

Let $\tilde T ^{\ss}_{\mathcal D}  $ (resp. $\tilde T ^{\ss}_{\mathcal D '}$) be the image of $ H _{K ^{\Sigma  } } $ in $\End _{o _{\mathcal D} } \tilde M ^ y  _{\mathcal D _ y }  $ (resp. $\End _{o _{\mathcal D} } \tilde M ^ y  _{\mathcal D '_ y } $). In the case of {\bf n.o.f.},
$\tilde U (p_v) = U(p_v)$-operator on the right hand side of ($\dagger _2$) satisfies the relation
$$
U ^ 2 - T _v \cdot U  + q_v T_{v,v} = 0, 
$$ 
over $\tilde T ^{\ss}_{\mathcal D}$, since
$ U(p_v) $ is given by 
$$
\begin{pmatrix}
 T_v & - T_{v,v} \\ 
q_v  & 0 
\end{pmatrix} 
$$
as in Lemma \ref{lem-congruent21}. Let $A =\tilde T ^{\ss}_{\mathcal D}  [U] / (U ^ 2 - T_v \cdot
U  + q_v T_{v,v} )  $. $A$ acts faithfully on
$ ( \tilde M ^ y  _{\mathcal D _ y }  ) ^{\oplus 2}$, where $ U$ acts as $U (p_v) $. 
Then ($\dagger_2 $) is compatible with the action of Hecke algebras (including $U (p_v ) $-operator), we have a canonical surjective homomorphism
$$
\tilde T ^{\ss}_{\mathcal D '} [ U ( p_v ) ] \twoheadrightarrow A. 
$$
Since $\deform _{\mathcal D '} (v) = \bold {n.o.}$, there is only one non-zero root of ($\dagger_3$) mod $\tilde m^{\ss}  _{\mathcal D_y } $, which defines a maximal ideal $m_A $ of $A$. The localization $ A_{m_A}$ is isomorphic to $T ^{\ss}_{\mathcal D }  $, and $(  ( \tilde M ^ y  _{\mathcal D _ y }  ) ^{\oplus 2} ) _{m_A} \simeq   \tilde M ^ y  _{\mathcal D _ y }  $. \par
By Theorem \ref{thm-modular41}, there is a maximal ideal ideal $ \tilde m ' $ of $ \tilde T ^{\ss}_{\mathcal D '} [ U ( p_v ) ] $, and the localization $(\tilde T ^{\ss}_{\mathcal D '} [ U ( p_v ) ] ) _{\tilde m '}  $ is equal to $ T _{\mathcal D '}$. 
Evidently $\tilde m ' $ is above $m_A$ since the condtion at $v$ is $\bold {n.o.} $ or $\bold {n.o.f. } $, and $ A_{\tilde m ' A} \simeq T _{\mathcal D_y}  $. \par
Thus we have 
$$
\beta : T_{\mathcal  D ' } \stackrel{\sim} {\longrightarrow}T_{\mathcal  D ' _y } 
\longrightarrow  T _{\mathcal  D_y }  \stackrel{\sim} {\longrightarrow}  T_{\mathcal  D} , 
$$
using Proposition \ref{prop-modular41}. Moreover, 
$$
 \xi : M_{\mathcal  D} \longrightarrow M _{\mathcal  D ' }
$$
is induced by localization at $\tilde m '$. By applying this result to $\bar \rho ^{\op} $ and $ \mathcal D ^{\op}$, we have a $T _{\mathcal D '} $-homomorphism 
$$
\hat \xi : \hat M _{\mathcal D } \longrightarrow \hat M _{\mathcal D ' } ,
$$
and we view the dual $(\hat \xi ) ^{\vee} : \hat M _{\mathcal D ' } \to \hat M_{\mathcal D } $ also as a $T _{\mathcal D '}$-homomorphism. 
The restriction of $(\hat \xi ) ^{\vee} $ to $ M _{\mathcal D }$ is $\xi ^{\vee} $ by the construction.  

\begin{prop}\label{prop-congruent21} 
Let $\bar \rho$ be an absolutely irreducible representation with a deformation type $\mathcal  D$ which satisfies $ \deform_{\mathcal  D}(v) = \bold {n.o.f.} $ at $v \vert \ell$.  Let $
\mathcal  D ' $ be the deformation type of $\bar \rho$ which has the same data as $\mathcal D$, and the deformation function $\deform _{\mathcal D'} $ satisfies $ \deform _{\mathcal  D ' } ( v) =
\bold {n.o.}$, and takes the same value as $\deform _{\mathcal  D}$ other than $v$.\par
\begin{enumerate}
\item There is a surjective homomorphism $\beta: T_{\mathcal  D '} \to
T_{\mathcal  D } $, and for the canonical map 
$\xi : M _{\mathcal  D} \to M _{\mathcal  D'} $, 
$$
\xi  ^{\vee} \circ  \xi ( M_{\mathcal  D}) =  \Delta \cdot M _{\mathcal 
D} , \quad \Delta =  \tilde U (p_v ) ^ 2 - \tilde U (p_v , p_v) 
.$$
Here $  \xi ^{\vee}: M _{\mathcal  D'} \to M_{\mathcal  D}$ is the map defined by self-duality pairings $\langle \ , \
\rangle _{\mathcal  D} $ and $\langle \ , \
\rangle _{\mathcal  D'}$.
\item Similarly, we have 
$$
(\hat \xi )  ^{\vee} \circ  \hat \xi ( \hat M_{\mathcal  D}) =  \Delta \cdot \hat M _{\mathcal 
D}.$$
Here $\Delta $ is the same as in (1). 
\item Assume Hypothesis \ref{hyp-coh31} if $q_{\bar \rho } = 1$. Then 
$\hat  \xi ( \hat M_{\mathcal  D})$ is an $o_{\mathcal D }$-direct summand of $\hat M _{\mathcal D'} $. 
\end{enumerate}
\end{prop}
\begin{proof}[Proof of Proposition \ref{prop-congruent21}] (1) follows from Lemma \ref{lem-congruent21}, since $
\xi^{\vee} \vert _{M _{\mathcal D ' } }
$ defined by the self-duality of $\hat M_{\mathcal  D ' }$ is equal to the map deduced from $ \tilde
\varphi _{(k,w)}$. For (2), we should check that the element $\tilde U ( p_v ) ^2 - \tilde U ( p_v , p_v )  $ in $ T _{\mathcal D ^{\op} } $ generate the same ideal as generated by $\tilde U ( p_v ) ^2 - \tilde U ( p_v , p_v )  $ in $ T _{\mathcal D} $ by the isomorphism $ T _{\mathcal D} \simeq T _{\mathcal D ^{\op} }$. This follows from the equality 
$\tilde U ( p_v ) ^2 - \tilde U ( p_v , p_v ) = \tilde U ( p_v , p_v ) ^2 ( (\tilde U (p_v , p_v  ) ^{-1} \tilde U ( p_v ) ) ^2 - \tilde U (p_v , p_v ) ^{-1}        )$. \par
(3) is clear from the definition of $\xi $. 
\end{proof}

\subsection{Congruence modules (II)}\label{subsec-congruent3}
From now on, we assume that the deformation type $\mathcal D $ of $\bar \rho$ satisfies $\deform _{\mathcal D} (v) = \bold f $. Let $\mathcal D '$ be the deformation type of $\bar \rho$ which has the same data as $\mathcal D$ except for the deformation function $\deform _{\mathcal D '}$. $\deform _{\mathcal D '} ( v) =\bold u $, and $\deform _{\mathcal D '}$ takes the same value as $\deform _{\mathcal D} $ at the other places. \par
In this subsection, we assume moreover that $\bar \rho$ is either of type $1_{SP}$ or $1_{PR}$ at $v$.
\par
 Let $D$ be a division algebra with $q_D \leq 1 $, $K$ an $F$-factorizable compact open subgroup of $G_D ( \mathbb A _{\mathbb Q, f} ) $, and for a place $v \nmid \ell$ and an integer $c \geq 1$, we assume that $ K _v = K_1 (m^c_{F_v} ) \cap \GL_2 (o_{F_v} ) ^{\ell}$ in the notation
of \S\ref{subsec-congruent2}. The degeneracy maps define
$$
\xi _{1, (k, w)} : 
H ^ { q_D } ( S_{ K } ,\ \bar { \mathcal  F} _{(k, w) } ) ^{\oplus 2}\overset {\pr_1 ^* +
\pr_2 ^*}
\longrightarrow  H ^ { q_D } ( S_{ K\cap K_0 (v   ^ {c+1}) } ,\
\bar {\mathcal  F } _{(k, w)} ).
$$

Note that $\xi _{1, (k, w) } $ may not be injective, even on the part which is not of residual type. $ \xi ^{\vee}_{1, (k, w) } $ is
defined similarly as in \S\ref{subsec-congruent2}. Then, as in the proof of Lemma \ref{lem-congruent21}, we have
$$
( \xi_{1, (k, w) } ) ^{\vee}  \circ \xi _{1, (k, w) }  =
 \begin{pmatrix}
 q_v &  U(p_v)\\
U ( p_v, p_v ) ^{-1} \cdot U(p_v) & q_v
\end{pmatrix} .
\leqno{(\ast)}
$$
Let $\tau _{(k, w) } :  H ^ { q_D } ( S_{ K } ,\ \bar { \mathcal  F} _{(k, w) } )\to H ^ { q_D } ( S_{ K } ,\
\bar {\mathcal  F}_{(k, w) } ) ^{\oplus 2}$ be the map defined by 
$$ 
\tau _{(k, w) } :  c \longmapsto \begin{pmatrix}  -q_v  \\   U ( p_v, p_v ) ^{-1} \cdot U(p_v )\end{pmatrix} \cdot c  .
$$

By ($\ast$),
$$
\tau_{(k, -w)}  ^{\vee}  \circ \xi ^{\vee}_{1, (k,w)} 
\circ\xi _{1, (k, w)}  \circ \tau _{(k, w)} =
\begin{pmatrix}
-q_v &  U(p_v)
\end{pmatrix}
\cdot
\begin{pmatrix}
 U(p_v, p_v ) ^{-1} \cdot U(p_v ) ^2- q_v ^ 2 \\ 
0
\end{pmatrix}
 $$
$$
= -q_v (U(p_v, p_v ) ^{-1} \cdot U(p_v ) ^2 - q_v ^2)   .$$

Since we allow the unrestricted condition on the determinant when $ \deform_{\mathcal D  '} = \bold u $, we need to discuss the passage from $K _0 ( m^{c+1} _v  )\cap GL _2 ( o_{F_v } ) ^{\ell}$ to $K_1(m_v ^{c+1})$. \par
We define $ \Delta  _{K , v} $ as the quotient of $ \Delta _v $ by the image of $Z(F) \cap K $. Then 
$$
\pi : S_{ K_ 1 ( m _{F_v} ^ {c+1}) \cdot K ^ v } \longrightarrow S_{ K \cap K_0 ( v ^ {c+1} ) }
$$
is a $ \Delta  _{K , v} $-torsor, and the composition
$$
H ^ { q_D  } (
S_{ K \cap K_0 ( v ^ {c+1}) },   \bar {\mathcal  F}_{(k, w)})\overset{\pi ^*} \longrightarrow  H ^ { q_D  } (
S_{ K_ 1 ( m _{F_v} ^ {c+1}) \cdot K ^ v } ,   \bar {\mathcal  F}_{(k, w)})\overset{ \pi _! } \longrightarrow H ^ { q_D  } (
S_{ K \cap K_0 ( v ^ {c+1}) } ,   \bar {\mathcal  F}_{(k, w)})
$$
induced by the trace map is the multiplication by $\sharp  \Delta  _{K , v}  $.

Thus we have
\begin{lem}\label{lem-congruent31} Assume that $ v \nmid \ell$.
\begin{enumerate}
\item  We obtain 
$$
\tau _{(k, -w)} ^{\vee}  \circ \xi^{\vee}_{1, (k,w)} 
\circ \xi _{1, (k, w)}  \circ \tau _{(k, w)} ( H ^ { q_D } ( S_{ K } ,\ \bar { \mathcal  F}_{(k, w)} ))= ( U(p_v, p_v
) ^{-1} \cdot U(p_v ) ^2 - q_v ^2 )  H ^ { q_D } ( S_{ K } ,\  \bar {\mathcal  F}_{(k, w)} )
$$
up to the modules of residual type.
\item The equality
$$
\tau _{(k, -w)} ^{\vee}  \circ \xi^{\vee}_{1, (k,w)}  \circ \pi _! \circ \pi ^* 
\circ \xi _{1, (k, w)}  \circ \tau _{(k, w)} ( H ^ { q_D } ( S_{ K } ,\ \bar { \mathcal  F}_{(k, w)} ))
$$
$$
=\sharp  \Delta  _{K , v} \cdot  ( U(p_v, p_v
) ^{-1} \cdot U(p_v ) ^2 - q_v ^2 )  H ^ { q_D } ( S_{ K } ,\  \bar {\mathcal  F}_{(k, w)} )
$$
holds up to the modules of residual type.
\end{enumerate}
\end{lem}

We apply the result to $K = \ker \nu _{\mathcal D _ y } \vert _{\tilde K _{\mathcal  D_y} } $. Recall that we assume that $\bar \rho$ is of type $1_{SP}$ or $1_{PR} $. Let $c = c (\bar \rho ) $ be the integer $\geq 1 $ defined in \S\ref{subsec-galdef7}. To define $\beta: T_{\mathcal  D'} \to T_{\mathcal  D}$ in this case,
$$
A= \tilde T^{\ss}  _{\mathcal  D }[U(p_v) ] [U] / (U ( U  - U(p _v)  ) ) .
$$
 The $A$-action on $H ^ { q _{\bar \rho}  } _{\stack} ( S _{ K } ,\ \bar { \mathcal  F} _{(k, w)} )^{\oplus 2}$ is given by 
$$
U \longmapsto  \begin{pmatrix}
 U (p_v ) & q_v \\ 
0 & 0 
 \end{pmatrix}
 .
$$
A direct calculation in this case shows that this $U$-action is compatible with the $U(p_v)$-action on $H ^
{ q _{\bar \rho} } ( S_{  K_1 (m_{F_v} ^ {c+1} )\cdot K ^v  } ,\
\bar { \mathcal  F} _{(k, w)} )$, and $A$ is regarded as a $\tilde T ^{\ss} _{\mathcal D'} [U (p_v) ]$-algebra. \par
By the same argument as in the case of $\bold{n.o.f.}$, there is a maximal ideal $\tilde m '$ of $\tilde T ^{\ss} _{\mathcal D'} [U (p_v) ]$, and the localization $(\tilde T ^{\ss} _{\mathcal D'} [U (p_v) ]) _{\tilde m ' } $ is equal to $ T _{\mathcal D ' } $. 
Since $\deform _{\mathcal D '} (v) = \bold u  $, $m _A = \tilde m ' A$ contains $U$, and the localization $ A_{m_A} $ at $m_A$ is isomorphic to $T _{\mathcal  D_y }  $ (though we have omitted
$U(p_v)$-operator in the definition of $T_{\mathcal  D _y} $, by Proposition \ref{prop-modular71} it is recovered from $\rho ^{\modular}_{\mathcal D} $). Thus we have a surjective homomorphism $\beta : T_{\mathcal D '} \to T_{\mathcal D}$. \par
Then $\xi_{1, (k,w) } $ becomes injective on the localization by $m_S $ by the cohomological universal exactness (Proposition \ref{prop-coh41} with $n= c$), since there is no contribution from \newline $ H ^{q_{\bar \rho} } _{\stack} ( S_{(\GL_2 (o_{F_v }) ^\ell \cap K_1 ( m^{c-1}_{F_v} )  \cdot
K^v},
\
\bar {\mathcal  F}_{(k, w)} )$ as $c = c(\bar \rho ) $ is the conductor of $G_{F_v}$-representation $\bar \rho\vert _{F_v}\otimes \bar \kappa^{-1}  _v $ and the level can not be smaller. It also follows that the image of $\xi_{1, (k,w) } $ is an $o_\lambda$-direct summand after localization. After localization at $m_A$, $\tau_{(k,w)}$  is an isomorphism, and $\pi $ is universally injective by Proposition \ref{prop-coh31}. \par
Thus we have a $T_{\mathcal D '}$-homomorphism
$$ 
\xi = (\pi \circ \xi_{1, (k , w)} \circ \tau _{(k, w)} )_{m_A}: M _{\mathcal D } \longrightarrow M _{\mathcal D '}. 
$$
\bigskip
As for the calculation of the congruence module, we treat the $1_{SP}$-case first.
Since all representations which contributes to $T_{\mathcal D}$ are special representations at $v$, $ U(p_v )^2 = U(p_v, p_v ) $ holds in $A_{m_A} $.  By Lemma \ref{lem-congruent31}, we have
\begin{prop}\label{prop-congruent31}
Assume that $\bar \rho
$ is of type $1_{SP} $ at $v$, and $\deform _{\mathcal D} (v) = \bold f $.
\begin{enumerate}
\item There is a surjective homomorphism $\beta: T_{\mathcal  D '} \to
T_{\mathcal  D } $, and for the canonical map 
$\xi : M _{\mathcal  D} \to M _{\mathcal  D'} $, 
$$
 \xi^{\vee} \circ  \xi ( M_{\mathcal  D}) = \Delta \cdot M _{\mathcal 
D} , \quad \Delta = \sharp (\ker ( C_{F, P ^{\bold u } _{\mathcal D '} , \ell  } \to C_{F, P ^{\bold u } _{\mathcal D } , \ell  }  )    )     \cdot ( q_v ^2-1 ) .$$
 Here, $ \xi ^{\vee}: M _{\mathcal  D'} \to
M_{\mathcal  D}$ is the map defined by Poincar\'e duality. 
\item 
$$
(\hat  \xi ) ^{\vee} \circ \hat  \xi ( \hat M_{\mathcal  D}) = \Delta \cdot \hat  M _{\mathcal 
D}  .
$$
Here, $\Delta $ is the same as in (1). 
\item $\hat  \xi ( \hat M_{\mathcal  D})$ is an $o_{\mathcal D }$-direct summand of $\hat M _{\mathcal D'} $. 
\end{enumerate}
\end{prop}

Next we treat the $1_{PR}$-case. Since 
$$
U (p_v) ^2 = q_v \cdot U (p_v, p_v ) 
$$
holds in $A_{m_A} $, by Lemma \ref{lem-congruent31}, we obtain
\begin{prop}\label{prop-congruent32}
Assume that $\bar \rho
$ is of type $1_{PR} $ at $v$, and $\deform _{\mathcal D } (v) = \bold f $.
\begin{enumerate}
\item There is a surjective homomorphism $\beta: T_{\mathcal  D '} \to
T_{\mathcal  D } $, and for the canonical map 
$\xi : M _{\mathcal  D} \to M _{\mathcal  D'} $, 
$$
 \xi^{\vee} \circ  \xi ( M_{\mathcal  D}) = \Delta \cdot M _{\mathcal 
D} , \quad \Delta = \sharp (\ker ( C_{F, P ^{\bold u } _{\mathcal D '} , \ell  } \to C_{F, P ^{\bold u } _{\mathcal D } , \ell  }  )    )     \cdot ( q_v -1 ) .$$
 Here, $ \xi ^{\vee}: M _{\mathcal  D'} \to
M_{\mathcal  D}$ is the map defined by Poincar\'e duality. 
\item 
$$
(\hat  \xi ) ^{\vee} \circ \hat  \xi ( \hat M_{\mathcal  D}) = \Delta \cdot \hat  M _{\mathcal 
D}.
$$
Here, $\Delta$ is the same as in (1). 
\item $\hat  \xi ( \hat M_{\mathcal  D})$ is an $o_{\mathcal D }$-direct summand of $\hat M _{\mathcal D'} $. 
\end{enumerate}
\end{prop}

\subsection{Congruence modules (III)}\label{subsec-congruent4}
In this subsection, we treat the $2_{PR}$-case. \par
For a division quaternion algebra $D$ with $q_D \leq 1 $, assume that $K $ is an $F$-factorizable small compact open subgroup of $G_D (\mathbb A _{\mathbb Q, f} )  $, and $ K _v = \GL_2 (o_{F_v} ) ^{\ell} $ in
the notation of \S\ref{subsec-congruent2} at $v \nmid \ell$. We denote $\pr _1, \ \pr _2 : S_{ K \cap K_0 ( v) } \to S_K  $ the degeneracy maps
defined in \S\ref{subsec-coh1}. Here $\pr_1 $ corresponds to the standard inclusion of the groups. 
$ \pr _3, \ \pr _ 4 : S_{ K \cap K_0 ( v^2) } \to S_{K \cap K_0 (v^ 2) } $ are
two degeneracy maps defined similarly. 

By Proposition \ref{prop-coh41} with $ n = 1 $,
$$
0 \longrightarrow H ^ { q_D  } ( S _{ K } ,\  \bar {\mathcal  F } _{(k, w)} )\overset {( \pr _2 ^*,\ -\pr _1 ^* )} \longrightarrow
  H ^ { q_D } ( S _{ K
\cap K_0 ( v)  } ,\ \bar {\mathcal  F} _{(k, w)} ) ^{\oplus 2} \overset {( \pr _3 ^*,\ \pr _4 ^* )} \longrightarrow
  H ^ { q_D  } ( S _{ K
\cap K_0 ( v^2)  } ,\  \bar {\mathcal  F }_{(k, w)} )\
$$
is universally exact up to the modules of residual type (in Proposition \ref{prop-coh41}, it is stated for $K_{11}
$-structure, but the above cohomology groups with $K_0$-structure is obtained by taking the invariants of a subgroup of $ k(v)
^\times \times k (v)^\times $ by using Lemma \ref{lem-coh21}).
By using Theorem \ref{thm-coh31} when $q_{\bar \rho} = 0 $, and Hypothesis \ref{hyp-coh31} when $q_{\bar \rho } = 1$, the homomorphism
$$
 H ^ { q_D } ( S _{ K } ,\ \bar {\mathcal  F} _{(k, w)} )^{\oplus 3} \overset {\xi _{2, (k, w) }  } \longrightarrow     H ^ { q_D } ( S _{ K
\cap K_0 ( v^2)  } ,\ \bar { \mathcal  F} _{(k, w)} )
$$
is universally injective up to the modules of residual type. Here $\xi_{2, (k,w) }   = ( \pr _1 \circ \pr _3 ) ^* \times  ( \pr _2 \circ \pr _3 ) ^*
\times  (
\pr _2
\circ
\pr _4 ) ^* $. We define $ \xi^{\vee}_{2, (k, w)}  $ similarly as in \S\ref{subsec-congruent2}.\par
The morphism
$$
\pi : S_{ K_ 1 ( m _{F_v} ^ {2}) \cdot K ^ v } \longrightarrow S_{ K \cap K_0 ( v ^ {2} ) }
$$
defined in \S\ref{subsec-congruent2} is a $ \Delta  _{K , v} $-torsor, and the composition 
$$
H ^ { q_D  } (
S_{ K \cap K_0 ( v ^ {2}) },   \bar {\mathcal  F}_{(k, w)})\overset{\pi ^*} \to  H ^ { q_D  } (
S_{ K_ 1 ( m _{F_v} ^ {2}) \cdot K ^ v } ,   \bar {\mathcal  F}_{(k, w)})\overset{ \pi _! } \to H ^ { q_D  } (
S_{ K \cap K_0 ( v ^ {2}) } ,   \bar {\mathcal  F}_{(k, w)})
$$
induced by the trace map is the multiplication by $\sharp  \Delta  _{K , v}  $.
Thus we have 
\begin{lem}\label{lem-congruent41} Assume that $v \nmid \ell$. 
\begin{enumerate}
\item We obtain
$$
\xi^{\vee}_{2,  (k, w ) }  \circ \xi_{2, (k, w)}  \ ( H ^ { q_D } ( S _{ K } ,\ \bar { \mathcal  F} _{(k, w) } )^{\oplus 3}  )= (
q_v -1 ) ( T_v ^ 2 - T_{v, v} ( 1+ q_v ) ^ 2 ) H ^ { q_D  } ( S _{ K } ,\  \bar {\mathcal  F} _{(k, w) }
)^{\oplus 3}  
$$
up to the modules of residual type. 
\item The equality
$$
\xi^{\vee}_{2,  (k, w ) }  \circ \pi _! \circ \pi ^* \circ  \xi_{2, (k, w)}  \ ( H ^ { q_D } ( S _{ K } ,\ \bar { \mathcal  F} _{(k, w) } )^{\oplus 3}  )=\sharp \Delta _{K, v} \cdot  (
q_v -1 ) ( T_v ^ 2 - T_{v, v} ( 1+ q_v ) ^ 2 ) H ^ { q_D  } ( S _{ K } ,\  \bar {\mathcal  F} _{(k, w) }
)^{\oplus 3}  
$$
holds up to the modules of residual type. 
\end{enumerate}
\end{lem}
This is shown as in \cite{W2}, proposition 2.6 by $3\times 3$-matrix calculation. So
we omit the details. \par

To define $\beta : T_{\mathcal  D'} \to T_{\mathcal  D} $, we discuss as in \S\ref{subsec-congruent2}, by taking an auxiliary place $y$ and a set $S$. $K =\ker \nu _{\mathcal D_y } \vert _{\tilde K _{\mathcal D_y } } $.
For the homomorphism $\tilde \xi$ and $\tilde T ^{\ss} _{\mathcal  D  }$-algebra $A$ defined by 
$$
\tilde \xi = \pi ^* \circ \xi_{2, (k, w) }  :H ^ { q_{\bar \rho}  } _{\stack} ( S _{ K } ,\ \bar { \mathcal  F} _{(k, w) } )^{\oplus 3} \overset
{\xi _{2, (k, w) }  } \longrightarrow
    H ^ { q _{\bar \rho}  } _{\stack}( S _{ K
\cap K_0 ( v^2)  } ,\  \bar {\mathcal  F} _{(k, w) } )
$$
$$
\overset {\pi ^* } \longrightarrow H ^ { q_{\bar \rho}  } _{\stack}(
S_{ K
\cap K_1 ( v ^ 2) } ,\ 
\bar {\mathcal  F}_{(k, w) })
$$
and
$$
A= \tilde T ^{\ss} _{\mathcal  D  } [U] / (U ( U ^ 2 - T _v U + q_v T_{v, v} ) ) ,
$$
we give the $A$-action on $H ^ { q_{\bar \rho} } _{\stack}( S _{ K } ,\ \bar {\mathcal  F} _{(k, w)} )^{\oplus 3}  $ by 
$$
\begin{pmatrix}
  T_v & -T_{v,v}& 0 \\
q_v  & 0 &0\\
0& q_v &0\\
\end{pmatrix}
 .
$$
One checks that the action of $U$ is the restriction of $U (p_v) $-action on $ H ^ { q_{\bar \rho}  } _{\stack}( S
_{ K } ,\ \bar { \mathcal  F} _{(k, w)} )^{\oplus 3}  $. 
As in \S\ref{subsec-congruent2}, there is a maximal ideal $\tilde m '$ of $\tilde T ^{\ss} _{\mathcal D'} [U (p_v) ]$, and the localization $(\tilde T ^{\ss} _{\mathcal D'} [U (p_v) ]) _{\tilde m ' } $ is equal to $ T _{\mathcal D ' } $ by using Proposition \ref{prop-modular71} and Theorem \ref{thm-modular41}.

Since $\deform _{\mathcal D '} (v) = \bold u  $, $m _A = \tilde m ' A$ contains $U$, and the localization $ A_{m_A} $ at $m_A$ is isomorphic to $T _{\mathcal  D }  $. 
Thus we have $ \beta:T_{\mathcal  D ' } \to T_{\mathcal  D } $, and a $T_{\mathcal D '}$-homomorphism $ \xi : M _{\mathcal  D} \to M _{\mathcal  D' } $ by localizing $\tilde \xi$. 

\begin{prop}\label{prop-congruent41}
Assume that $\bar \rho
$ is of type $2_{PR} $ at $v$, and $\deform_{\mathcal D} (v ) = \bold f $.
\begin{enumerate}
\item There is a surjective homomorphism $\beta: T_{\mathcal  D '} \to
T_{\mathcal  D } $, and for the canonical map 
$\xi : M _{\mathcal  D} \to M _{\mathcal  D'} $, 
$$
 \xi^{\vee} \circ  \xi ( M_{\mathcal  D}) = \Delta \cdot M _{\mathcal 
D} ,$$
$$
\Delta = \sharp (\ker ( C_{F, P ^{\bold u } _{\mathcal D '} , \ell  } \to C_{F, P ^{\bold u } _{\mathcal D } , \ell  }  )    )     \cdot ( q_v -1 ) \cdot (T_v ^ 2 - T_{v, v} ( 1+ q_v ) ^ 2  ) .$$
 Here $ \xi ^{\vee}: M _{\mathcal  D'} \to
M_{\mathcal  D}$ is the map defined by Poincar\'e duality. 
\item The equality
$$
(\hat  \xi ) ^{\vee} \circ \hat  \xi ( \hat M_{\mathcal  D}) = \Delta \cdot \hat  M _{\mathcal 
D} $$
holds. Here, $\Delta $ is the same as in (1). 
\item $\hat  \xi ( \hat M_{\mathcal  D})$ is an $o_{\mathcal D }$-direct summand of $\hat M _{\mathcal D'} $ under Hypothesis  \ref{hyp-coh31} if $q _{\bar \rho} = 1 $. 
\end{enumerate}
\end{prop}

\subsection{Congruence modules (IV)}\label{subsec-congruent5}

 Finally, we treat the cases of $0_{NE} $ and $0_E$.

In the case of $0_{NE}$, for $K =\ker \nu _{\mathcal D_y } \vert _{\tilde K _{\mathcal D_y } } $ and $c = c (\bar \rho ) =\Art \bar \rho \vert _{F_v}$, we consider the morphism
$$
\pi : S_{ K_ 1 ( m _{F_v} ^ {c}) \cdot K ^ v } \longrightarrow S_{ K \cap K_0 ( v ^ {c} ) }
$$
defined in \S\ref{subsec-congruent2} is a $ \Delta  _{K , v} $-torsor, and the composition 
$$
H ^ { q_D  }_{\stack}  (
S_{ K \cap K_0 ( v ^ {c}) },   \bar {\mathcal  F}_{(k, w)})\overset{\pi ^*} \longrightarrow H ^ { q_D  } _{\stack}(
S_{ K_ 1 ( m _{F_v} ^ {c}) \cdot K ^ v } ,   \bar {\mathcal  F}_{(k, w)})\overset{ \pi _! } \longrightarrow H ^ { q_D  } _{\stack}(
S_{ K \cap K_0 ( v ^ {c}) } ,   \bar {\mathcal  F}_{(k, w)})
$$
induced by the trace map is the multiplication by $\sharp  \Delta  _{K , v}  $ in the notation of \S\ref{subsec-congruent3}.\par

In the case of $0_E$, for $K =\ker \nu _{\mathcal D_y } \vert _{\tilde K _{\mathcal D'_y } } $ and $K '  =\ker \nu _{\mathcal D_y } \vert _{\tilde K _{\mathcal D_y } } $, 
$$ 
\pi : S_{ K ' } \longrightarrow
S_{ K  }
$$
is a torsor under $\tilde \Delta_v $. Here
$\tilde \Delta _v$ is the $\ell$-Sylow subgroup of $ \F_{q_v ^ 2}\simeq o_{D_{F_v} }/
m_{D_{F_v}}$ (since $q_v \equiv -1 \mod \ell$ and $\ell\geq 3$, 
$\Delta _{K, v } $ is trivial). 
The composition 
$$
H ^ { q_D  }_{\stack}  (
S_K,   \bar {\mathcal  F}_{(k, w)})\overset{\pi ^*} \longrightarrow  H ^ { q_D  } _{\stack}(
S_{ K '  } ,   \bar {\mathcal  F}_{(k, w)})\overset{ \pi _! } \longrightarrow H ^ { q_D  } _{\stack}(
S_K,   \bar {\mathcal  F}_{(k, w)})
$$
induced by the trace map is the multiplication by $\sharp  \tilde \Delta_v  $.

Thus we have
\begin{prop}\label{prop-congruent51} Assume that $\bar \rho $ is either of type $0_E$ or type $0_{NE}$ at $v$, and $ \deform _{\mathcal D } (v ) = \bold f $. 
\begin{enumerate}
\item There is a surjective homomorphism $\beta: T_{\mathcal  D '} \to
T_{\mathcal  D } $, and 
$$
\xi ^{\vee} \circ \xi (M_{\mathcal  D} ) = \Delta \cdot M_{\mathcal  D} .$$
Here $\xi $ the canonical map $ M _{\mathcal  D} \to M _{\mathcal  D'} $, 
 $\Delta = \sharp (\ker ( C_{F, P ^{\bold u } _{\mathcal D '} , \ell  } \to C_{F, P ^{\bold u } _{\mathcal D } , \ell  }  )    )   $ (resp. $ \sharp \tilde \Delta _v$) if $\bar \rho $ is of type $0_{NE}$ (resp. $0_E$) at $v$.
\item The equality
$$
(\hat \xi) ^{\vee} \circ \hat \xi (\hat M_{\mathcal  D} ) = \Delta \cdot \hat M_{\mathcal  D} 
$$
holds. Here, $\Delta $ is the same as in (1). 
\end{enumerate}
\end{prop}
The analysis of $0_E$-case was treated in \cite{DT2}, and by Diamond
\cite{D1} with applications to deformation rings in the case where $F=\Q$.

\section{The main theorem}\label{sec-final}
In this section, we generalize the results of \cite{W2} and \cite{D1}.
\begin{thm} ($R=T$ theorem)\label{thm-final1}
Let $F$ be a totally real number field of degree $d$, $ \bar \rho: G_F \to \GL_2( k) $ an absolutely irreducible mod $
\ell$-representation. We fix a deformation type $\mathcal  D$, and assume the following conditions.
\begin{enumerate}
\item $\ell \geq 3$, and $\bar \rho
\vert _{{ F(\zeta_\ell) }}$ is absolutely irreducible. When $\ell= 5$, the following case is excluded: the projective image $\bar G$ of $ \bar \rho $ is isomorphic to $\PGL_2 (\mathbb F_5)$, and the mod $\ell$-cyclotomic character $\bar \chi _{\cycle} $ factors through $G_F \to  \bar G ^{\ab} \simeq \mathbb Z/ 2$ (in particular $[F(\zeta_5): F ] = 2$). 
\item For $ v\vert
\ell$, the deformation condition for 
$\bar
\rho
\vert _{G_{F_v} } $ is either nearly ordinary or flat. 
When the condition is nearly ordinary (resp. flat) at $v$, we assume that $\bar \rho \vert _{G_{F_v} } $ is $G_{F_v} $-distinguished (resp. $F_v$ is absolutely unramified).
\item There is a minimal modular lifting $\pi$ of $ \bar \rho $ in Definition  \ref{dfn-nearlyordinary32}.

\item Hypothesis \ref{hyp-nearlyordinary21} is satisfied.
\item If $\mathcal  D $ is not minimal, we assume Hypothesis \ref{hyp-coh31} when $q_{\bar \rho } = 1$. 
\end{enumerate}
Then the universal deformation ring $R_{\mathcal  D} $ of $ \bar
\rho$ of type $\mathcal  D$ is of relative complete intersection of dimension zero over
$o_{\mathcal D} $, and is isomorphic to the Hecke algebra $T_{\mathcal  D} $.
\end{thm}

We prove the following theorem at the same time.
\begin{thm}(Freeness theorem)\label{thm-final2}
Under the same assumptions as Theorem \ref{thm-final1}, $ M_{\mathcal  D }  $ and $\hat M _{\mathcal D} $ are free $R_{\mathcal  D}
$-modules. 
\end{thm}

\begin{proof}[Proof of Theorem \ref{thm-final1} and \ref{thm-final2}]
Two theorems are already shown in \S\ref{sec-minimal} when $\mathcal D$ is minimal by using the Taylor-Wiles system constructed in \S\ref{sec-const}. We make the reduction to the minimal case by the level-raising formalism developped in \S\ref{sec-TW}.\par

For a given deformation type $ \mathcal  D$, assume that $\pi _{\mathcal  D} : R_{\mathcal  D} \simeq 
T_{\mathcal  D}
$,
$T_{\mathcal  D}
$ is a local complete intersection, and $\hat M_{\mathcal  D} $ is $T_{\mathcal  D }$-free. In the
terminology of Definition \ref{dfn-raising1}, the admissible quintet $(R_{\mathcal  D} , \ T_{\mathcal  D} , \pi _{\mathcal  D} ,\
\hat M_{\mathcal  D} ,\langle \ , \
\rangle_{\mathcal  D} ) 
$ is distinguished. \par
Take a finite place $v$, and assume that 
$ \deform _{\mathcal  D} (v) \in \{ \bold {n.o.f.} ,\
\bold {f}\}
$. In case of $\bold {n.o.f.}$, we twist by $\chi_{\cycle} ^{-\frac w 2} $, and assume
that $w= 0
$.\par
Let $
\mathcal  D ' $ be the deformation type of $\bar \rho$ which has the same data as $\mathcal D$ except the deformation function $\deform_{\mathcal D'} $.  $ \deform _{\mathcal  D ' } ( v) > \deform _{\mathcal D} (v)  $, and $\deform _{\mathcal D ' } $ takes the same value as $\deform _{\mathcal  D}$ at the other places.\par

We take a cuspidal
representation
$\pi
$ corresponding to a component of $T_{\mathcal  D } $. By extending $E_{\mathcal D}$ by taking an extension of scalars if necessary, we
may assume that 
$\pi _f $ is defined over $E_{\mathcal D}$. $\pi$ defines an $o_{\mathcal D} $-homomorphism $f:
T_{\mathcal  D} \to o_{\mathcal D} $, and a $G_F$-representation $\rho : G_F \to \GL_2 ( o_{\mathcal D} ) $ of type $\mathcal D$. \par

In \S\ref{sec-congruent}, we have defined an admissible morphism 
$$
(R_{\mathcal  D ' } , \ T_{\mathcal  D ' } , \pi _{\mathcal  D ' } ,\
\hat M_{\mathcal  D ' } ,\langle \ , \
\rangle _{\mathcal  D'} ) \longrightarrow 
(R_{\mathcal  D} , \ T_{\mathcal  D} , \pi _{\mathcal  D} ,\
\hat M_{\mathcal  D} ,\langle \ , \
\rangle _{\mathcal  D} ) 
$$ 
in the sense of Definition \ref{dfn-raising1} under Hypothesis \ref{hyp-coh31}. The morphism commutes with scalar extension $ o_{\mathcal D}  \to o _{\lambda} $, and the property of being distinguished is preserved under the scalar extension.\par

We apply Theorem \ref{thm-TW2}. First we check that $ \Delta $ given in \S\ref{sec-congruent} is a
non-zero divisor. In the cases of
$1_{SP }$ 
$1_{PR}$,
$0_E$ and $0_{NE} $, this is clear.  We consider the $\bold { n.o.f.}$ and $2_{PR}$-cases. Let $\pi $ be a cuspidal representation of $\GL_{2, F} $ appearing in the component of $T_{\mathcal  D} $. Then $
\pi _v$ belongs to principal series by the definition of $T_{\mathcal  D} $. Let $\chi _{1, v} $
and
$\chi _{2, v} $ be two quasi-characters such that $ \pi_v$ is isomorphic to $\pi ( \chi _{1, v} , \chi_{2, v} )$. We define $ \alpha _v =  \chi _{1, v}
(p_v)$, $\beta _v =  \chi _{2, v}
(p_v )$ for a uniformizer $p_v$ of $F_v$. By the next lemma, $\alpha
_v 
\neq q_v ^{\pm 1} \beta _v $ for any $\pi$, and hence  
$\Delta $ is a non-zero divisor.
\begin{lem}\label{lem-final1}
For any cuspidal representation $\pi $ of $\GL_{2, F} $ of a totally
real field $F$ with the infinity type
$ (k, w)
$, if
$\pi _v$ belongs to principal series, $\alpha _v \neq q_v ^{\pm1} \beta_v$ for the pair $(\alpha_v, \
\beta_v )$ defined as above.  
\end{lem}
\begin{proof}[Proof of Lemma \ref{lem-final1}]
Let $ \pi _{\bold u}  = \pi \otimes
\vert
\cdot \vert ^{\frac {w+1} 2 } $. $\pi _{\bold u}$ is a unitary cuspidal representation of $\GL_{2, F}$ by our
normalization. By taking the base change lift $\pi _{\bold u , F'} = \BC( \pi _{\bold u } )$ of $\pi _{\bold u } $ with respect to some abelian
extension $ F ' / F$ \cite{L}, one may assume that $\pi_{\bold u , F'} $ has a spherical component $(\pi _{\bold u , F'} )_{v'} $ at a place 
$v'
\vert v$. So it suffices to see $ \alpha _{v'} \neq q_{v'} ^{\pm 1} \beta _{v'} $ for the
Satake parameter of $(\pi _{\bold u , F' } ) _{v'} $. Since $(\pi _{\bold u , F' } )_{v'} $ is unitary, by \cite{JS}, 2.5, Corollary, $ \vert
\alpha _{v'}
\vert ,\  \vert \beta _{v'} \vert < q_{v '} ^{\frac 1 2}$. So the equality $ \alpha _{v'}= 
q_{v'} ^{\pm 1} \beta _{v'}  $ does not hold. 
\end{proof}

By Proposition \ref{prop-galdef101}, 
$$
\Hom _{o_{\mathcal D} } ( \ker f_{R_{\mathcal  D} } / (\ker f_{R_{\mathcal  D} } )^2, \ E_{\mathcal D} /
o_{\mathcal D} )=
\Sel _{\mathcal  D} ( F,\
\ad
\rho ) 
$$
and
$$
\Hom _{o_{\mathcal D} } ( \ker f_{R_{\mathcal  D ' } } /( \ker f_{R_{\mathcal  D ' } } )^2, \
E_{\mathcal D} /o_{\mathcal D} )=
\Sel _{\mathcal  D '} ( F,\
\ad
\rho ) 
$$
hold.
Since we assume that $ R_{\mathcal  D} \simeq T_{\mathcal  D} $, $\ker
f_{R_{\mathcal  D} } / (\ker f_{R_{\mathcal  D}} ) ^2 = \ker f/( \ker f) ^2$, and the finiteness of $\Sel
_{\mathcal  D }(F, \  \ad \rho )$ follows since $T_{\mathcal  D } $ is reduced. So the assumptions of level raising formalism (Theorem \ref{thm-TW2}) is satisfied if we check 
$$
\length_{o_{\mathcal D}}  \Sel
_{\mathcal  D '} ( F,\
\ad
\rho ) \leq \length_{o_{\mathcal D}}   \Sel
_{\mathcal  D } ( F,\
\ad
\rho ) + \length_{o_{\mathcal D}}  o _{\mathcal D}  / f (\Delta)o _{\mathcal D} 
$$
for the element $\Delta $ of $ T _{\mathcal D}$ given in \S\ref{sec-congruent}.\par
If $\deform _{\mathcal  D} (v) \neq \bold {n.o.f.}$, by Lemma \ref{lem-galdef102}, 
$$
\length_{o_{\mathcal D}}  \Sel
_{\mathcal  D '} ( F,\
\ad
\rho ) \leq \length_{o_{\mathcal D}}   \Sel
_{\mathcal  D } ( F,\
\ad
\rho ) +  \length_{o_{\mathcal D}}    H ^ 0 ( F_v,\ \ad \rho (1)\otimes _{o_{\mathcal D} } E_{\mathcal D} / o_{\mathcal D} ) ,
$$
and it suffices to see 
$$
\length_{o_{\mathcal D}}  H ^ 0 ( F_v,\ \ad \rho (1)\otimes _{o_{\mathcal D} } E_{\mathcal D} / o_{\mathcal D} ) 
 =\length_{o_{\mathcal D}}    o_{\mathcal D} / f
(\Delta) o_{\mathcal D}. 
$$
This is a consequence of Proposition \ref{prop-congruent31}, \ref{prop-congruent32}, \ref{prop-congruent41}, and \ref{prop-congruent51}. \par

If $\deform _{\mathcal  D} (v)= \bold {n.o.f.}$, one uses Lemma \ref{lem-galdef101} and Proposition \ref{prop-congruent21}, and the inequality
follows. \par

Thus we have shown that $(R_{\mathcal  D ' } , \ T_{\mathcal  D ' } , \pi _{\mathcal  D ' } ,\
M_{\mathcal  D ' } ,\langle \ , \
\rangle _{\mathcal  D'} ) $ is distinguished under the assumption that $(R_{\mathcal  D } , \ T_{\mathcal  D  } , \pi _{\mathcal  D  } ,\
M_{\mathcal  D  } ,\langle \ , \
\rangle _{\mathcal  D } ) $ is distinguished by Theorem \ref{thm-TW2}. Starting from the minimal case $\mathcal D= \mathcal D_{\min} $ (Theorem
\ref{thm-minimal1}), we raise the level place by place by enlarging $\mathcal D$, and the general case is shown. 
\end{proof}
\begin{cor}\label{cor-final1}
Under the same assumption as in \ref{thm-final1}, the Selmer group $\Sel
_{\mathcal  D} ( F, \ \ad \rho ) $ of
$\rho  =
\rho _{\pi, E_\lambda} $ for $\pi $ appearing in $T_{\mathcal  D} $ is finite. 
\end{cor}
This is a consequence of the reducedness of $T_{\mathcal D} $. \par
\bigskip
Theorem \ref{thm-final1} and \ref{thm-final2} can be applied to the nearly ordinary Hecke algebra of Hida. \par
Let $\mathcal D$ be a deformation type of $\bar \rho$ as in Theorem \ref{thm-final1}. For  $P  = \deform_{\mathcal D} ^{-1} ( \{\bold {n.o.} \} ) $, let $(o^{\times }_{F_v } ) _{\ell} $ be the pro-$\ell$ completion of $ o ^{\times }_{F_v} $, $ \mathcal X ^{\lc} _{\bold {n.o.}} = \prod _{  v \in P } ( o^{\times }_{F_v  })_{\ell}  $. 
We define the deformation type $\mathcal D ^{\bold {n.o.}} $ in the following way. 
\begin{itemize}
\item The deformation function of $\mathcal D^{\bold {n.o.}}$ is the same as that of $\mathcal D$. 
\item The coefficient ring $ o_{{\mathcal D } ^{\bold n.o.}} $ of $\mathcal D^{\bold {n.o.}}$ is $o_{\mathcal D} [[\mathcal X ^{\lc} _{\bold {n.o.}}  ]] $.
\item At $v \in P ^{\bold {n.o.}}$, the nearly ordinary type is 
$$
\kappa_{\mathcal D, v} \cdot \mu _v : (I^{\ab}_{F_v})_{G_{F_v}} \simeq o^{\times} _{F_v} \longrightarrow o^{\times}_{{\mathcal D } ^{\bold n.o.}} .
$$
Here, $ \kappa_{\mathcal D, v} $ is the neary ordinary type of $\mathcal D$, and $\mu _v :  o^{\times} _{F_v}\to (o_{\mathcal D} [[ o ^{\times}_{F_v , \ell }  ]] ) ^{\times} \to  o^{\times}_{{\mathcal D } ^{\bold {n.o.}}}$ is the universal character.  The flat twist type is the same as $\mathcal D$, by regarding it as a character with values in $ o^{\times}_{{\mathcal D } ^{\bold {n.o.}}}$.
\end{itemize}
Let $ R _{\mathcal D ^{\bold{ n.o.}}} $ be the universal deformation ring of $\bar \rho$ of type $\mathcal D ^{\bold {n.o.}} $. $ R _{\mathcal D ^{\bold {n.o.}}}$ is an $ o_{{\mathcal D } ^{\bold {n.o.}}} $-algebra.
Since any deformation of type $\mathcal D$ is a deformation of type $\mathcal D ^{\bold {n.o.}}$, there is a natural surjective map
$$
f_{\mathcal D}: R _{\mathcal D ^{\bold{ n.o.}}}\longrightarrow R_{\mathcal D},
$$
and $f_{\mathcal D} $ induces an isomorphism
$$ 
R _{\mathcal D ^{\bold {n.o.}}}/I^{\bold {n.o.} } R _{\mathcal D ^{\bold {n.o.}}} \stackrel{\sim}{\longrightarrow} R_{\mathcal D}.
$$ 
Here, $I^{\bold {n.o.} }  $ is the augumentation ideal of $  o_{{\mathcal D } ^{\bold {n.o.}}}$.

\begin{cor}\label{cor-final2} Assumptions are as in Theorem \ref{thm-final1}. Assume moreover that there is an $R _{\mathcal D ^{\bold n.o.}}$-module $M  _{\mathcal D ^{\bold n.o.}} $ which has the following properties:
\begin{enumerate}
\item $M  _{\mathcal D ^{\bold { n.o.}}} $ is a free $ o_{{\mathcal D } ^{\bold {n.o.}}}$-module. 
\item$ M  _{\mathcal D ^{\bold n.o.}} / I^{\bold {n.o.} }M  _{\mathcal D ^{\bold {n.o.}}}  \stackrel {\sim } {\rightarrow} M _{\mathcal D}$ as an $R_{\mathcal D} $-module. 
\item Let $T _{\mathcal D ^{\bold {n.o.}}}$ be the image of $R _{\mathcal D ^{\bold {n.o.}}}$ in $\End _{o_{{\mathcal D } ^{\bold { n.o.}}}} ( M  _{\mathcal D ^{\bold {n.o.}}} )  $. Then $ T _{\mathcal D ^{\bold {n.o.}}}$ is reduced, and $M  _{\mathcal D ^{\bold {n.o.}}} $ is generically free of the same rank as the $R_{\mathcal D} $-rank of $M_{\mathcal D}$ ($M_{\mathcal D}$ is $R_{\mathcal D}$-free by Theorem \ref{thm-final2}).
\end{enumerate}
Under these assumptions, $R _{\mathcal D ^{\bold { n.o.}}} \stackrel{\sim} {\rightarrow} T _{\mathcal D ^{\bold {n.o.}}}$, $R _{\mathcal D ^{\bold {n.o.}}} $ is of relative complete intersection of dimension zero over $o_{{\mathcal D } ^{\bold{ n.o.}}} $, and $M _{\mathcal D ^{\bold { n.o.}}}$ is a free $R _{\mathcal D ^{\bold {n.o.}} } $-module.
\end{cor}
\begin{proof}[Proof of Corollary \ref{cor-final2}] Let $ \alpha $ be the $R_{\mathcal D}$-rank of $M_{\mathcal D} $. By (2) and Nakayama's Lemma, there is a surjection 
$\beta: T^{\oplus \alpha }  _{\mathcal D ^{\bold {n.o.}}} \twoheadrightarrow M _{\mathcal D ^{\bold {n.o.}}}$. By (3), $\beta $ is an isomorphism generically, and hence is injective by the reducedness of $ T _{\mathcal D ^{\bold {n.o.}}}$. By (1), $T _{\mathcal D ^{\bold {n.o.}}}$ is $ o_{{\mathcal D } ^{\bold {n.o.}}} $-free. \par
Let $J $ be the kernel of $ R_{\mathcal D ^{\bold {n.o.}}} \twoheadrightarrow T _{\mathcal D ^{\bold {n.o.}}}$. Since $T _{\mathcal D ^{\bold {n.o.}}}$ is $  o_{{\mathcal D } ^{\bold {n.o.}}} $-free, 
$$
0 \longrightarrow   J   /  I^{\bold {n.o.} }  J \longrightarrow 
   R_{\mathcal D ^{\bold {n.o.}}}   /  I^{\bold {n.o.} }  R_{\mathcal D ^{\bold {n.o.}}} 
\overset{\gamma}\longrightarrow  T_{\mathcal D ^{\bold {n.o.}}}   /  I^{\bold {n.o.} }  T_{\mathcal D ^{\bold {n.o.}}} 
\longrightarrow 0 
$$
is exact. $\gamma$ is an isomorphism by Theorem \ref{thm-final1}. Thus $J = 0 $ by Nakayama's lemma, and $R _{\mathcal D ^{\bold { n.o.}}} \stackrel{\sim} {\rightarrow} T _{\mathcal D ^{\bold {n.o.}}} $ holds. Since $R _{\mathcal D ^{\bold { n.o.}}} $ is $o_{{\mathcal D } ^{\bold {n.o.}}} $-finite flat and $ R _{\mathcal D ^{\bold { n.o.}}} / I^{\bold {n.o.} } R _{\mathcal D ^{\bold { n.o.}}} $ is of relative complete intersection of dimension zero over $o_{{\mathcal D } ^{\bold {n.o.}}} / I^{\bold {n.o.} } = o _{\mathcal D}    $, $R _{\mathcal D ^{\bold {n.o.}}} $ is of relative complete intersection of dimension zero over $o_{{\mathcal D } ^{\bold{ n.o.}}} $.
\end{proof}

Using the results in \S\ref{sec-coh} and the perfect complex argument, the existence of $ T _{\mathcal D ^{\bold {n.o.}}}$ and $ M_{\mathcal D ^{\bold {n.o.}}}$ as in Corollary \ref{cor-final2} is shown (exact control theorem). If $\deform_{\mathcal D} (v) = \bold {n.o.}$ for any $v \vert \ell$, $T_{\mathcal D ^{\bold {n.o.}}}$ is the nearly ordinary Hecke algebra of Hida.


\begin{thebibliography}{99}



\bibitem {BK}
         Bloch, S., Kato, K.:
  The Grothendieck Festschrift, Vol. I,  333--400, Progr. Math., 86, Birkh\"auser Boston, Boston, MA, 1990.

\bibitem{Bou} 
 Bourbaki, N.:
 Alg\`ebre Commutative, chap.3--4,  Hermann  1961.


\bibitem {BR}
         Blasius, D., Rogawski, J.D.:
 Motives for Hilbert modular forms.  Invent. Math.  114  (1993),  no. 1, 55--87.

\bibitem{C}
         Conrad, B.:
 The flat deformation functor.  Modular forms and Fermat's last theorem (Boston, MA, 1995),  373--420, Springer, New York, 1997. 

\bibitem {Car1}
         Carayol, H.:
 Sur la mauvaise r\'eduction des courbes de Shimura, Compositio Math.  59  (1986),  no. 2, 151--230. 

\bibitem{Car2}
         Carayol, H.:
 Sur les repr\'esentations $p$-adiques associ\'ees aux formes modulaires de
Hilbert,
 Ann. Sci. \'Ecole Norm. Sup. (4)  19  (1986),  no. 3, 409--468. 

\bibitem {CPS}
         Cline, E., Parshall, B., Scott, L.:
 Cohomology of finite groups of Lie type. I.  Inst. Hautes \'Etudes Sci. Publ. Math. No. 45 (1975), 169--191. 




\bibitem{dJ} 
        de Jong, A.J.:
 Smoothness, semi-stability and alterations.  Inst. Hautes \'Etudes Sci. Publ. Math.  No. 83 (1996), 51--93. 

\bibitem{D1} 
        Diamond, F.:
 On deformation rings and Hecke rings.  Ann. of Math. (2)  144  (1996),  no. 1, 137--166.

\bibitem{D2}
      Diamond, F.:
 The Taylor-Wiles construction and multiplicity one.  Invent. Math.  128  (1997),  no. 2, 379--391. 

\bibitem{De}
        Deligne, P.:
 Travaux de Shimura
 S\'eminaire Bourbaki, 23\`eme ann\'ee (1970/71), Exp. No. 389,  pp. 123--165. Lecture Notes in Math., Vol. 244, Springer, Berlin, 1971.



\bibitem {DT1}
         Diamond, F., Taylor, R.:
Lifting modular mod $\ell$ representations.  Duke Math. J.  74  (1994),  no. 2, 253--269. 


\bibitem {DT2}
         Diamond, F., Taylor, R.:
 Nonoptimal levels of mod $\ell$ modular representations.  Invent. Math. 115  (1994),  no. 3, 435--462. 

\bibitem{FL}
           Fontaine, J.M., Lafaille, G.:
 Constructions de repr\'esentations $p$-adiques
  Ann. Sci. \'Ecole Norm. Sup. (4)  15  (1982),  no. 4, 547--608. 
 
 
\bibitem{Fu1}
           Fujiwara, K.:
 Level optimization in the totally real case,  preprint, arXiv.math.NT/0602586.

\bibitem{Fu2}
           Fujiwara, K.:
Galois deformations and arithmetic geometry of Shimura varieties, Proceedings of the International Congress of Mathematicians, Madrid, Spain, 2006, European Mathematical Society, to appear.



\bibitem{GD1}
Grothendieck, A., Dieudonn\'e, J.:
\'El\'ements de g\'eom\'etrie alg\'ebrique IV, Premi\`ere Partie, Publications Math\'ematiques, 20, 1964.

\bibitem{GD2}
Grothendieck, A., Dieudonn\'e, J.:
\'El\'ements de g\'eom\'etrie alg\'ebrique IV, Seconde Partie, Publications Math\'ematiques, 24, 1965.


\bibitem {Ger}
         G\'eradin, P.:
 Facteurs locaux des algebres simples de rang 4, I ,
 Groupes r\'eductif et formes automorphes, I, pp. 37--77, Publ. Math. Univ. Paris VII, 1, Univ. Paris VII, Paris, 1978. 

\bibitem {Gr1} Grothendieck, A.:
Le groupe de Brauer. III. Exemples et compl\'ements,
1968  Dix Expos\'es sur la Cohomologie des Sch\'emas  pp. 88--188 North-Holland, Amsterdam; Masson, Paris.


\bibitem {Gr2} Grothendieck, A.:
 Mod\`eles de N\'eron et monodromie
 Groupes de Monodromie en G\'eom\'etrie Alg\'ebrique,
S\'eminaire de G\'eom\'etrie Alg\'ebrique du Bois-Marie 1967--1969 (SGA 7 I). Dirig\'e par A. Grothendieck. Avec la collaboration de M. Raynaud et D. S. Rim. Lecture Notes in Mathematics, Vol. 288. Springer-Verlag, Berlin-New York, 1972. 

\bibitem{H0}
   Hida, H.:
 Galois representations into ${\rm GL}\sb 2(\mathbb Z \sb p[[X]])$ attached to ordinary cusp forms.  Invent. Math.  85  (1986),  no. 3, 545--613. 

\bibitem {H1}
       Hida, H.:
On $p$-adic Hecke algebras for ${\rm GL}\sb 2$ over totally real fields.  Ann. of Math. (2)  128  (1988),  no. 2, 295--384. 
 

\bibitem{H2}
     Hida, H.:
 On nearly ordinary Hecke algebras for ${\rm GL}(2)$ over totally real fields.  Algebraic number theory,  139--169, Adv. Stud. Pure Math., 17, Academic Press, Boston, MA, 1989.

\bibitem{Ja1}
       Jarvis, F.:
  Mazur's principle for totally real fields of odd degree, 
  Compositio Math.  116  (1999),  no. 1, 39--79. 

\bibitem{Ja2}
    Jarvis, F.:
Level lowering for modular mod $l$ representations over totally real fields.  Math. Ann.  313  (1999),  no. 1, 141--160.

\bibitem{JL}
      Jacquet, H., Langlands, R.P.:
 Automorphic forms on $\GL (2)$. Lecture Notes in Mathematics, Vol. 114. Springer-Verlag, Berlin-New York, 1970. 



\bibitem{JS}
      Jacquet, H., Shalika, J.A.:
On Euler products and the classification of automorphic forms. II.  Amer. J. Math.  103  (1981), no. 4, 777--815. 

\bibitem{K}
       Kutzko, P.:
 The Langlands conjecture for $\GL(2)$ of a local field.  Ann. of Math. (2)  112  (1980), no. 2, 381--412.


\bibitem{L}
       Langlands, R.P.:
 Base change for $\GL(2)$. Annals of Mathematics Studies, 96. Princeton University Press, Princeton, N.J.; University of Tokyo Press, Tokyo, 1980.

\bibitem{Len}
          Lenstra, H. W.:
 Complete intersections and Gorenstein rings.  Elliptic curves, modular forms, $\&$ Fermat's last theorem (Hong Kong, 1993),  99--109, Ser. Number Theory, I, Internat. Press, Cambridge, MA, 1995. 

\bibitem{M}
       Mazur, B.:
 Deforming Galois representations.  Galois groups over $\mathbb Q$ (Berkeley, CA, 1987),  385--437, Math. Sci. Res. Inst. Publ., 16, Springer, New York, 1989. 

\bibitem{O1}
       Ohta, M.:
 On $\ell$-adic representations attached to automorphic forms.  Japan. J. Math. (N.S.)  8  (1982),  no. 1, 1--47. 

\bibitem{O2}
 Ohta, M.:
On the zeta function of an abelian scheme over the Shimura curve.  Japan. J. Math. (N.S.)  9  (1983),  no. 1, 1--25. 

\bibitem{Raj}
        Rajaei, A.:
 On the levels of mod $\ell$ Hilbert modular forms.
J. Reine Angew. Math. 537 (2001), 33--65.

\bibitem{Ray}
        Raynaud, M.: 
 Sch\'emas en groupes de type $(p,\ldots, p)$.  Bull. Soc. Math. France  102  (1974), 241--280. 
 
\bibitem{Ram}
        Ramakrishna, R.:
 On a variation of Mazur's deformation functor.  Compositio Math.  87  (1993),  no. 3, 269-286. 

\bibitem{Ri1}
        Ribet, K.A.:
 Congruence relations between modular forms.  Proceedings of the International Congress of Mathematicians, Vol. 1, 2 (Warsaw, 1983),  503--514, PWN, Warsaw, 1984.

\bibitem{Ri2}
   Ribet, K.A.:
 On modular representations of $\Gal (\overline{\mathbb Q}/\mathbb Q)$ arising from modular forms.  Invent. Math.  100  (1990),  no. 2, 431--476. 


\bibitem{S}
      Shimura, G.:
 On canonical models of arithmetic quotients of bounded symmetric domains.  Ann. of Math. (2)  91  1970 144--222. 

\bibitem{Sa1}
        Saito, T.:
 Modular forms and $p$-adic Hodge theory.  Invent. Math.  129  (1997),  no. 3, 607--620. 

\bibitem{Sa2}
       Saito, T.:
 Hilbert modular forms and $p$-adic Hodge theory. 
Preprint (1999).


\bibitem{Sch}
        Schlessinger, M.:
Functors of Artin rings.  
Trans. Amer. Math. Soc.  130,  1968, 208--222. 

\bibitem{Sh}
          Shimizu, H.:
 Theta series and automorphic forms on $\GL (2)$.  J. Math. Soc. Japan  24  (1972), 638--683.

\bibitem{Ta}  
         Tate, J.:
$p$-divisible groups. 1967  Proc. Conf. Local Fields (Driebergen, 1966)  pp. 158--183 Springer, Berlin.
 
\bibitem{T1}
            Taylor, R.:
 On Galois representations associated to Hilbert modular forms.  Invent. Math.  98  (1989),  no. 2, 265--280. 



\bibitem{T2}
Taylor, R.:
  On Galois representations associated to Hilbert modular forms. II.  Elliptic curves, modular forms, $\&$ Fermat's last theorem (Hong Kong, 1993),  185--191, Ser. Number Theory, I, Internat. Press, Cambridge, MA, 1995. 

\bibitem{Tsu}
            Tsuji, T.:
$p$-adic \'etale cohomology and crystalline cohomology in the semi-stable reduction case.  Invent. Math.  137  (1999),  no. 2, 233-411. 



\bibitem{TW}
          Taylor, R,   Wiles, A.:
 Ring-theoretic properties of certain Hecke algebras.  Ann. of Math. (2)  141  (1995),  no. 3, 553--572. 

\bibitem{W1} Wiles, A.:
On ordinary $\Lambda$-adic representations associated to modular forms.  Invent. Math.  94  (1988),  no. 3, 529--573. 

\bibitem{W2} Wiles, A.:
Modular elliptic curves and Fermat's last theorem.  Ann. of Math. (2)  141  (1995),  no. 3, 443--551.
\end{thebibliography}
\end{document}